\input epsf.tex
\input amssym.def
\input amssym
\magnification=1000
\baselineskip = 0.17truein
\lineskiplimit = 0.01truein
\lineskip = 0.01truein
\vsize = 9.0truein
\voffset = 0.0truein
\parskip = 0.10truein
\parindent = 0.3truein
\settabs 12 \columns
\hsize = 6.4truein
\hoffset = -0.1truein

\setbox\strutbox=\hbox{%
\vrule height .708\baselineskip
depth .292\baselineskip
width 0pt}
\font\caps=cmcsc10

\font\bigtenrm=cmr10 at 14pt

\def\sqr#1#2{{\vcenter{\vbox{\hrule height.#2pt
\hbox{\vrule width.#2pt height#1pt \kern#1pt
\vrule width.#2pt}
\hrule height.#2pt}}}}
\def\square{\mathchoice\sqr46\sqr46\sqr{3.1}6\sqr{2.3}4}

\centerline{\bigtenrm THE EFFICIENT CERTIFICATION OF}
\centerline{\bigtenrm KNOTTEDNESS AND THURSTON NORM}
\tenrm
\vskip 14pt
\centerline{MARC LACKENBY}
\vskip 18pt
\centerline{\caps 1. Introduction}
\vskip 6pt

How difficult is it to determine whether a given knot is the unknot? The answer is not known.
There might be a polynomial-time algorithm, but so far, this has remained elusive. The complexity
of the unknot recognition problem was shown to be in NP by Hass, Lagarias and Pippenger [10].
The main aim of this article is to establish that it is in co-NP. This can be stated equivalently in
terms of the {\caps Knottedness} decision problem, which asks whether a given knot diagram represents
a non-trivial knot.

\noindent {\bf Theorem 1.1.} {\sl {\caps Knottedness} is in NP.}

In some sense, this result is not new. It was first announced by Agol [1] in 2002, but he has not provided
a full published proof. In 2011, Kuperberg gave an alternative proof of Theorem 1.1, but under the extra
assumption that the Generalised Riemann Hypothesis is true [19]. In this paper, we provide the first full proof
of the unconditional result.

Combined with the theorem of Hass, Lagarias and Pippenger [10], Theorem 1.1 gives the following corollary.

\noindent {\bf Corollary 1.2.} {\sl If either of the decision problems {\caps Unknot recognition} or {\caps Knottedness}
is NP-complete, then NP $=$ co-NP.}

This is because if any decision problem in co-NP is NP-complete, then the complexity classes NP and co-NP
must be equal. Since this is widely viewed not to be the case (see Section 2.4.3 in [7] for example), then it
seems very unlikely that either of these decision problems is NP-complete.

Decision problems that lie in both NP and co-NP are viewed as good potential candidates
for being solvable in polynomial time. However, it is very probable that some decision
problems in NP $\cap$ co-NP do not lie in P. A notable example is the problem of factorising an integer which,
when suitably rephrased as a decision problem, is in NP $\cap$ co-NP. It is a very interesting
consequence of Theorem 1.1 that {\caps Unknot recognition} now lies in that list of problems
lying in NP $\cap$ co-NP but that is not known to lie in P.

Our proof of Theorem 1.1 follows the outline given by Agol and, as in his argument, we establish more information
about the genus of the knot. Recall that the {\sl genus} of a knot $K$ in a 3-manifold is the minimal possible genus
of a {\sl Seifert surface} for $K$, which is a compact orientable embedded surface,
with no closed components and with boundary equal to $K$. When no such surface exists,
the genus of $K$ is defined to be infinite. We provide an algorithm to determine
the genus of a knot in the 3-sphere. More specifically, we consider the decision problem {\caps Classical knot genus}.
This takes, as its input a knot diagram $D$ and a positive integer $g$ (given in binary), and it asks whether the knot
specified by $D$ has genus $g$. We establish the following result.

\noindent {\bf Theorem 1.3.} {\sl {\caps Classical knot genus} is in NP.}

Theorem 1.1 is a consequence of Theorem 1.3 as follows. Given a diagram $D$ of a non-trivial
knot $K$, we need a certificate for its knottedness. Specifically, this certificate needs to
be verifiable in polynomial time, as a function of the crossing number $c(D)$ of $D$. Now, the genus $g$
of $K$ lies between 1 and $(c(D)-1)/2$. The certificate provided by Theorem 1.3 for this diagram $D$
and this genus $g$ is the required certificate for knottedness.

The term `classical' is meant to refer to knots in the 3-sphere. In fact, Theorem 1.3
is unlikely to generalise to knots in arbitrary closed 3-manifolds. Indeed, there is good reason
to believe that the problem of determining whether a knot in a closed 3-manifold has genus $g$ is 
not in NP (when both the knot and the 3-manifold are permitted to vary), because of the following
theorem, which is a well known consequence of work of Agol, Hass and Thurston [2]. 

\noindent {\bf Theorem 1.4.} {\sl If {\caps Knot genus in compact orientable 3-manifolds} is in NP, then NP $=$ co-NP.}

Here, the decision problem {\caps Knot genus in compact orientable 3-manifolds} takes as its input
a triangulation of a compact orientable 3-manifold $X$, a knot $K$ in $X$ given as a subcomplex of the 1-skeleton
and a positive integer $g$ (in binary), and it asks whether the minimal genus of a Seifert surface for $K$ is $g$.
We will give a short proof of this result in Section 16. As a result of Theorem 1.4, it is highly unlikely that
{\caps Knot genus in compact orientable 3-manifolds} is in NP.

The reason why knots in arbitrary compact orientable 3-manifolds appear to be significantly more
complicated than those in the 3-sphere is that when $M$ is the exterior of a
knot in the 3-sphere, $H_2(M, \partial M)$ is cyclic, and so there is only one possible homology
class (up to sign) that can be represented by a Seifert surface for $K$. However, when $M$ is the
exterior of a knot in a general compact 3-manifold, $H_2(M, \partial M)$ may have rank bigger than 1,
and so one must consider many different homology classes that could support a genus-minimising
Seifert surface. However, we will see shortly that if we restrict attention to a single homology
class, then an NP algorithm is available.

Although Theorem 1.3 is phrased in terms of the genus of a knot, it is not really the genus of
a surface $S$ that is its main measure of complexity in this paper. Instead, {\sl Thurston complexity}
$\chi_-(S)$ plays this role.  Recall that, for a compact connected orientable surface $S$, $\chi_-(S)$
equals $\max \{ - \chi(S), 0 \}$. When $S$ is a compact orientable surface with components $S_1, \dots, S_n$,
then $\chi_-(S) = \sum_{i=1}^n \chi_-(S_i)$. So, when $S$ is a Seifert surface for a knot $K$
and $S_1 = S - {\rm int}(N(K))$, then $S_1$ has minimal Thurston complexity in its
class in $H_2(S^3 - {\rm int}(N(K)), \partial N(K))$ if and only if $S$ has minimal possible genus.
However, the same is not true when $K$ is a link with more than one component.
In this case, $K$ may have disconnected Seifert surfaces as well as connected ones.

When $M$ is a compact orientable 3-manifold, and $z$ is a class in $H_2(M, \partial M)$, then
the {\sl Thurston norm} of $z$ is the minimal Thurston complexity of a compact oriented surface
representing $z$. Theorems 1.1 and 1.3 are special cases of a more general result which allows one to efficiently
determine the Thurston norm of a homology class.
We define the {\caps Thurston norm of a homology class} decision problem as follows.
The input is a triangulation of a compact orientable 3-manifold $M$,
a simplicial 1-cocycle $\phi$ representing an element $[\phi]$ in $H^1(M)$ and a non-negative integer $n$.
The problem asks whether the Thurston norm of the Poincar\'e dual to $[\phi]$ is $n$. The measure of complexity is
the size of the input, which is, up to a bounded linear factor, equal to the sum of the number of the tetrahedra in the
triangulation of $M$, the number of digits of the integer $n$ in binary and the sum of the
number of digits of $\phi(e)$, as $e$ ranges over all edges of the triangulation.

The following is the main theorem of this paper.

\noindent {\bf Theorem 1.5.} {\sl {\caps Thurston norm of a homology class} is in NP.}

%Note that the 3-manifold $M$ is required to have boundary a (possibly empty) collection of tori.
%The algorithm that we supply does not work in the case where $M$ has boundary components of higher genus.

Note that Theorem 1.3, and hence Theorem 1.1, follow quickly from Theorem 1.5. Given a diagram $D$
for a knot $K$, one can easily build a triangulation for its exterior $M$ in polynomial time, and where the number of
tetrahedra is at most a linear function of the crossing number of $D$. One can also readily build a simplicial 1-cocycle $\phi$
representing a generator $[\phi]$ of $H^1(M)$, where the maximal value of $|\phi(e)|$ for each edge $e$ is at most a linear function of the crossing number.
Setting $n = \max \{ 2g - 1, 0 \}$, we can apply Theorem 1.5 and thereby obtain an NP algorithm to
determine whether the genus of $K$ is $g$.

Our methods can also be used to certify the irreducibility of a 3-manifold. The decision problem {\caps Irreducibility of a compact orientable
3-manifold with toroidal boundary and } $b_1>0$ takes,
as its input, a triangulation of such a 3-manifold (which is permitted to have empty boundary) and asks whether it is irreducible.

\noindent {\bf Theorem 1.6.} {\sl The decision problem {\caps Irreducibility of a compact orientable
3-manifold with toroidal boundary and } $b_1>0$ is in NP.}

Since a compact orientable 3-manifold $M$ is irreducible and has incompressible boundary if and only if its double
is irreducible, we can use Theorem 1.6 to detect the incompressibility of $\partial M$. The decision problem
{\caps Incompressible boundary} takes as its input a triangulation of a compact orientable 3-manifold $M$ and it
asks whether $\partial M$ is incompressible. This was shown to be in co-NP by Ivanov (Theorem 4 in [12]).
We show that it is NP, and hence we have the following.

\noindent {\bf Theorem 1.7.} {\sl {\caps Incompressible boundary} is in NP $\cap$ co-NP.}

Theorem 1.7 has a consequence for the generalisation of unknot recognition to 3-manifolds other than the 3-sphere.
The decision problem {\caps Knottedness in 3-manifolds} takes, as its input, a triangulation of a compact
orientable 3-manifold $M$, and a knot $K$ given as a subcomplex
of the 1-skeleton, and it asks whether $K$ is knotted. This just means that $K$ does not bound an embedded disc.
Since the knottedness of $K$ is closely related to whether $M - {\rm int}(N(K))$ has compressible boundary, Theorem 1.7
can be used to establish the following result.

\noindent {\bf Theorem 1.8.} {\sl {\caps Knottedness in 3-manifolds} is in NP $\cap$ co-NP.}

\vskip 6pt
\noindent {\caps 1.1. Overview of the proof}
\vskip 6pt

In order to prove our main result, Theorem 1.5, one needs a method for certifying efficiently the Thurston norm of a class in $H_2(M, \partial M)$. 
By use of a doubling argument, we show in Section 14 that it suffices to consider the case where $M$ is closed and irreducible.
In fact, our methods work just as well when $\partial M$ is toroidal. So we now suppose that $M$ is a compact orientable irreducible 3-manifold with boundary a
(possibly empty) collection of tori. 

The method that both we and Agol use to certify the Thurston norm of a class in $H_2(M, \partial M)$ is based on work of
Thurston [27] and Gabai [3]. Thurston showed that if $S$ is a compact leaf of some taut foliation of $M$,
then $S$ minimises Thurston complexity in its class in $H_2(M, \partial M)$. Gabai showed that,
conversely, if $S$ is a compact orientable incompressible surface that minimises Thurston
complexity in its class and that intersects each component of $\partial M$ in a (possibly empty) collection of coherently
oriented essential curves, then there is some taut foliation of $M$ in which $S$ appears as a compact leaf.
Gabai's construction used a type of hierarchy for $M$, known as a sutured manifold hierarchy.
A sutured manifold structure on $M$ is a decomposition of $\partial M$ into two subsurfaces
$R_-(M)$ and $R_+(M)$, which meet along some simple closed curves $\gamma$. It is denoted $(M, \gamma)$.
The surface $R_-(M)$ is given a transverse
orientation pointing into $M$ and $R_+(M)$ is transversely oriented outwards. When $M$ is cut along a
transversely oriented properly embedded surface $S$ in general position with respect to $\gamma$, the new manifold inherits a sutured manifold structure.
A sutured manifold hierarchy is a sequence of decompositions
$$(M, \gamma) = (M_1, \gamma_1) \buildrel S_1 \over \longrightarrow (M_2, \gamma_2)
\buildrel S_2 \over \longrightarrow \dots \buildrel S_{n} \over \longrightarrow (M_{n+1}, \gamma_{n+1}),$$
where $(M_{n+1}, \gamma_{n+1})$ is a collection of 3-balls, each of which intersects $\gamma_{n+1}$ in
a single simple closed curve, and which satisfies some mild extra conditions. Gabai used these hierarchies
to construct taut foliations on the exteriors of many knots [4, 5], and was thereby able to determine
their genus. It was Scharlemann [25] who realised that much of Gabai's theory
could work without any reference to taut foliations; just the sutured manifold hierarchies
are enough to determine the Thurston norm of a homology class. For example, it is straightforward to verify that
the sequence of sutured manifold decompositions given in Figure 1 forms a taut sutured manifold
hierarchy, and hence the first decomposing surface minimises Thurston complexity in its
homology class. In particular, the existence of this hierarchy proves the knot $5_1$ is not the unknot.

It is sutured manifold hierarchies, such as the one in Figure 1, that we (and Agol) use to certify the
Thurston norm of a homology class. The existence of such a hierarchy was proved by Gabai,
but crucially, our certificate needs to be verifiable in polynomial time, as a function of the size of the input.
Essentially,  the sutured manifold hierarchy needs to be efficiently describable. The key to 
Agol's proof was to achieve this by placing some such hierarchy into `normal' form with respect
to a given triangulation for $M$. This normalisation procedure was based on work on Gabai [6],
who considered the related problem of normalising essential laminations.
In the present paper, we follow a similar approach, but instead
of using triangulations, we focus on handle structures. The machinery for placing
a sutured manifold hierarchy into `normal' form with respect to a handle structure
was developed by the author in [21]. This predates Agol's announcement of Theorems 1.1, 1.3 and 1.5,
and plays a central role in this paper.

\vskip 12pt
\centerline{
\epsfxsize=4in
\epsfbox{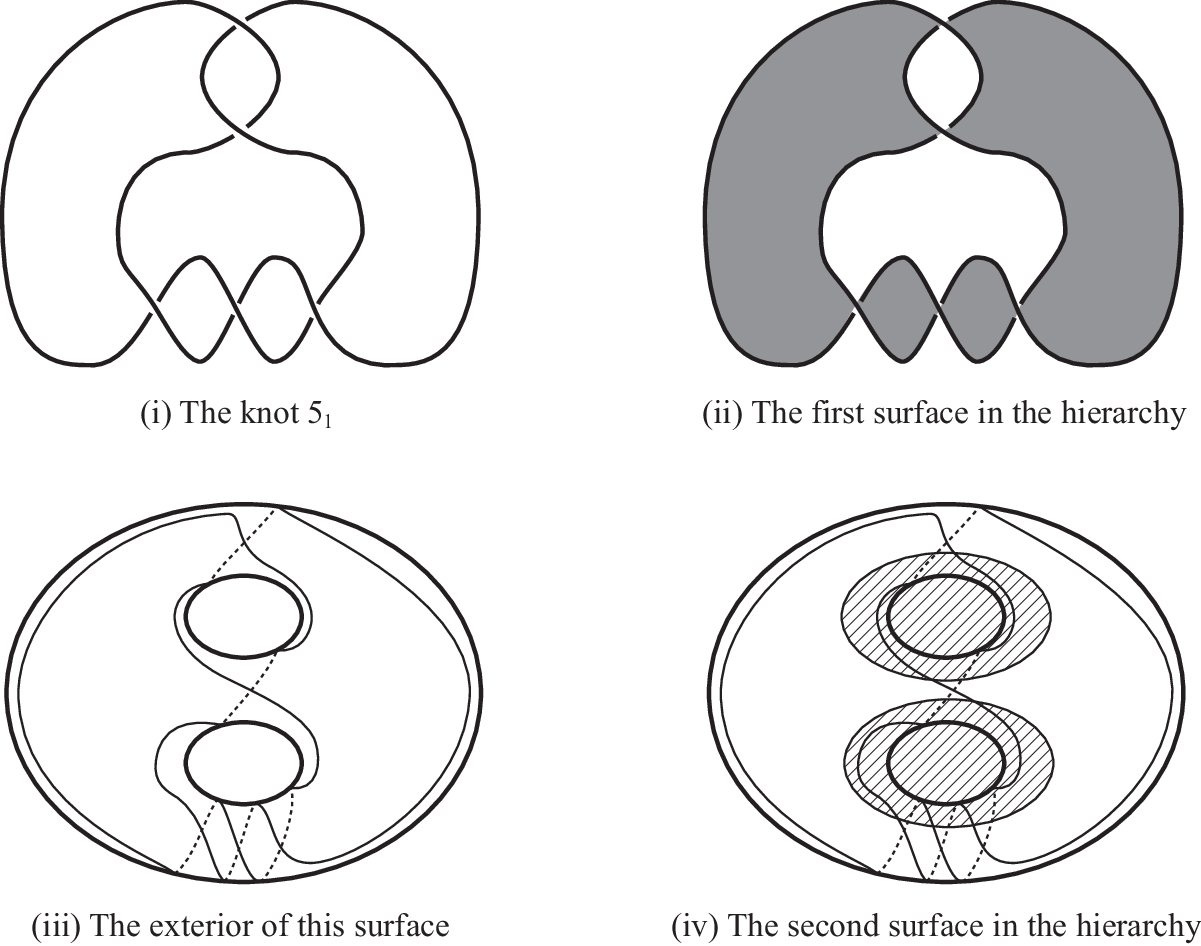}
}
\vskip 6pt
\centerline{Figure 1: A sutured manifold hierarchy}

\vskip 6pt
\noindent {\caps 1.2. Normalisation of sutured manifold hierarchies}
\vskip 6pt

Normal surfaces that minimise Thurston complexity in their homology class were studied by Tollefson and Wang [28].
So the first thing that we do is use their theory to realise the first surface $S_1$ in the hierarchy as a normal surface.
However, one difficulty with normal surfaces in triangulated 3-manifolds is that, when one cuts along them,
the resulting 3-manifold $M_2$ does not naturally inherit a triangulation. So we dualise the given triangulation
of the initial manifold $M$ to form a handle structure ${\cal H}$. There is a well-established theory of normal surfaces 
in handle structures [8, 13]. The next 3-manifold $M_2$ 
in the hierarchy then inherits a handle structure ${\cal H}_2$. 
Unfortunately, this may be much more complicated than ${\cal H}$, in two different senses.
For a start, ${\cal H}_2$ may simply have many more handles than ${\cal H}$. This happens
if $S_1$ intersects some handle of ${\cal H}$ in many discs, and then this handle of ${\cal H}$
is divided into many handles of ${\cal H}_2$. Fortunately, all but a bounded number of these
handles will be very simple copies of $D^2 \times I$, lying between two parallel normal discs.
These product regions patch together to form an $I$-bundle in $M_2$, known as its `parallelity bundle'.
It was shown in [21] how this parallelity bundle may be removed, primarily by decomposing 
along the annuli that form its intersection with the remainder of $M_2$. So, a key part of our argument
is the analysis of this $I$-bundle and an algorithmic method of removing it from $M_2$.

However, there is another reason why ${\cal H}_2$ may be more complicated than ${\cal H}$.
It is not at all clear that the local structure of the 0-handles of ${\cal H}_2$ is simpler than
that of ${\cal H}$. This issue was faced right at the very first use of hierarchies by Haken [8]
and Waldhausen [29]. They defined a notion of complexity for a handle structure of a 3-manifold,
and showed that, when it is decomposed along a normal surface, then the complexity does not go up. 
Unfortunately, their notion of complexity does not fit well with
sutured manifolds, primarily because it does not take account of the sutures. Fortunately, this issue 
was resolved in [21]. A variation of normal form more suited to sutured manifolds was introduced.
In this paper, we call such surfaces `regulated'. Also, a notion of complexity for a handle structure
of a sutured manifold was defined in [21]. It was also shown that this does not go up when the manifold
is decomposed along a regulated surface. Moreover, it was shown that sutured manifold
hierarchies can always be found where each decomposing surface is regulated.

Therefore, it is regulated surfaces that are used in this paper. Unfortunately, they come with their
own complications. Although it is the case that we may arrange for the decomposing surfaces to
be regulated, they may fail to satisfy one of the key technical requirements for a sutured manifold hierarchy.
Some curves of intersection with the surface $R_\pm$ may be simple closed curves bounding
discs in $R_\pm$. Such curves are called `trivial'. (See Remark 5.6 for an explanation for why
it seems hard to arrange that regulated decomposing surfaces have non-trivial boundary curves.)
Because surfaces with trivial boundary components are not permitted to be
part of a sutured manifold hierarchy, it is not clear that they can be used as part of a certificate
for Thurston norm. However, we develop a theory of `allowable hierarchies', where the decomposing
surfaces may have trivial boundary curves, but which can nonetheless be used to certify
Thurston norm.

So, allowable hierarchies of regulated surfaces will be used as part of our certificate. But these surfaces must be describable
in an efficient way. In particular, it is important that our surfaces intersect each handle
in at most $c^h$ discs, where $c$ is a universal constant and $h$ is the number of handles
in the initial handle structure. To be able to establish such a bound, we use methods from
linear algebra, that go back to Haken [8]. We show how regulated surfaces $S$ can be described
by means of solution $(S)$ to a system of linear equations, much in the same way that normal surfaces can be.
When $S$, $S_1$ and $S_2$ are regulated surfaces and $(S) = (S_1) + (S_2)$, then 
we say that $S$ is a `sum' of $S_1$ and $S_2$.
Just as in the normal surface case, one can place $S_1$ and $S_2$ into general position,
and then obtain $S$ by resolving the arcs and curves of intersection. It is crucial for our
purposes that we may find a decomposing surface that is `fundamental', which means
that it cannot then be expressed as a sum of other non-trivial surfaces. This is because 
fundamental surfaces have a bounded number of discs of intersection with each handle,
by methods that go back to Hass and Lagarias [9]. 

Therefore, we must analyse the case when $S$ is a sum of surfaces $S_1$ and $S_2$.
Now, $S_1$ and $S_2$ need not inherit orientations from $S$, and indeed they need
not even be orientable. But when they do inherit orientations, then the situation is
fairly easy to analyse. It turns out that we can generally show that decomposition
along $S_1$ or $S_2$ is taut, and so decompose along one of these instead.
These are `simpler' surfaces, and so in this way, one may arrange for decomposing
surfaces to be fundamental. When $S_1$ and $S_2$
do not inherit transverse orientations, then the aim is to show that $S$ was  not
as simple as it could have been. One can perform an `irregular switch' along
one of the arcs or curves of $S_1 \cap S_2$, creating a new transversely oriented
surface $S'$. We show that decomposition along $S'$ is also taut. This is possible
when the irregular switch takes place along a curve of $S_1 \cap S_2$. However, the
argument does not work in the case of an arc of $S_1 \cap S_2$, because the orientations
of $R_\pm$ at its endpoints may be problematic. Fortunately, regulated surfaces rescue us here,
because their boundaries are very tightly controlled, and in fact, it is possible to show
that the regular switch along arcs of $S_1 \cap S_2$ always respect the orientation of the surface.

Thus, we may arrange that the decomposing surfaces are fundamental, and hence
have an exponential bound on complexity. But decomposing along these surfaces is not
straightforward, because we never want to deal with handle structures having exponentially
many handles. They are too unwieldy to be efficiently describable within our certificate.
Fortunately, there is technology due to Agol, Hass and Thurston [2] which is applicable.
Using their methods, we show that, given a surface $S_i$ with an exponential bound on its complexity,
it is possible to determine the topological types of the components of the parallelity bundle for the manifold $M_{i+1}$
obtained by decomposing along $S_i$. So, we never need to construct the handle structure for $M_{i+1}$.
Instead, we can go directly to the handle structure for the manifold $M_{i+1}'$ that is obtained
by removing the parallelity bundle.

Thus, our certificate for Thurston norm is comprised (primarily) of the following pieces of information:
handle structures for a sequence of 3-manifolds, and regulated surfaces within these manifolds,
expressed as solutions to a system of equations. We verify this certificate by checking that the
next manifold is indeed obtained from the previous one by cutting along the regulated surface
and then purging the resulting manifold of its parallelity bundle.

\vfill\eject
\noindent {\caps 1.3. Comparison with Agol's strategy}
\vskip 6pt

The strategy of our proof, as explained in Section 1.1, is very similar to the one outlined by Agol [1].
However, the details, as given in Section 1.2, are very different. Instead of using handle structures,
Agol used triangulations. Since the result of decomposing a triangulation along a normal surface
is not in general a triangulation, Agol had to work hard to build a triangulation for each manifold
in the hierarchy. His technique was based on placing the surface into `spun' normal form.
The method of doing this was based on Gabai's method for normalising taut foliations,
as explained in [6], which used `sutured manifold evacuation'. In Agol's argument, the surface still
needed to be made `fundamental' in some way. Moreover, an analogue of the removal of
parallelity bundles would still need to have been achieved, again by the use of the algorithm
of Agol, Hass and Thurston [2]. Agol's strategy certainly has some advantages, but it seems
to us that the extensive use of the established techniques from [21], as followed in this paper, 
is also very convenient.

\vskip 6pt
\noindent {\caps 1.4. Structure of the paper}
\vskip 6pt

In Section 2, we recall some of the basic theory of sutured manifolds. In Section 3, we introduce
`decorated' sutured manifolds and `allowable' hierarchies. The key result here is Theorem 3.1 which implies that
allowable hierarchies can be used certify Thurston norm. In Section 4, we show that if a surface extends
to an allowable hierarchy, then certain modifications can be made to it, and the resulting surface
still extends to an allowable hierarchy. In Section 5, we introduce handle structures for sutured manifolds,
and regulated surfaces. We also introduce the complexity of a handle structure, and explain
how it behaves when a decomposition along a regulated surface is performed. In Section 6,
we give some simplifications that can be made to a handle structure. We also explain how
we need only to work with a finite universal list of 0-handle types. In Section 7, we develop an
algebraic theory of regulated surfaces. We introduce a minor variant, known as `boundary-regulated' surfaces,
and show how they may expressed as solutions to certain equations. We also explain how
summation of boundary-regulated surfaces can be interpreted topologically. Finally,
we show that decompositions can always be made along fundamental regulated surfaces.
This central result is Theorem 7.9.
In Section 8, we bound the complexity of fundamental surfaces, using tools from linear algebra.
In Section 9, we recall an algorithm of Agol, Hass and Thurston [2], and show how it can be used
to determine the parallelity bundle for the manifold obtained by decomposing along a surface.
In Section 10, we show to remove this parallelity bundle algorithmically. In Section 11, we show
how it suffices to focus our attention on 3-manifolds that are irreducible and atoroidal. This convenient hypothesis
occurs at several points in the preceding argument. In Section 12, we show how to certify efficiently
that a sutured manifold is a product. This is useful because the hierarchies that we use terminate
in products rather than 3-balls, for mostly technical reasons. We then go on to prove the main
theorem in the special case where the 3-manifold is Seifert fibred. In Section 13, we complete the
proof of Theorem 1.5 in the case of compact orientable irreducible 3-manifolds with (possibly empty)
toroidal boundary. We describe in detail the certificate for Thurston norm, we show why it always exists,
and how it may be verified in polynomial time. In Section 14, we show how to deal with the general case of 3-manifolds
that might be reducible or have non-toroidal boundary components. In Section 15, we give
the proofs of Theorems 1.7 and 1.8. In Section 16, we provide a proof of Theorem 1.4.

\vskip 6pt
\noindent {\caps 1.5. Acknowledgements}
\vskip 6pt

The author would like to thank the referee for many helpful suggestions which substantially improved the
paper. In particular, the referee suggested that one should double manifolds with non-toroidal boundary
components. This led to the versions of Theorems 1.5, 1.7 and 1.8 that appear in this paper, and 
that are improvements over the original versions. The author would also like to thank Mehdi Yazdi
for helpful conversations regarding the some of the material in Section 11.

\vfill\eject
\centerline{\caps 2. Sutured manifolds}
\vskip 6pt

A {\sl sutured manifold} is a compact orientable 3-manifold $M$, with its boundary decomposed into two compact
subsurfaces $R_-(M)$ and $R_+(M)$, in such a way that $R_-(M) \cap R_+(M)$ is a collection of simple closed curves
$\gamma$. These curves are called {\sl sutures}. The surfaces $R_-(M)$ and $R_+(M)$ are assigned
transverse orientations, with $R_-(M)$ pointing into $M$ and $R_+(M)$ pointing outwards. The sutured
manifold is usually denoted $(M, \gamma)$.

A compact oriented surface $S$ embedded in a 3-manifold $M$, with $\partial S$ in $\partial M$, is called {\sl taut}
if $S$ is incompressible and it has minimal Thurston complexity in its homology class
in $H_2(M, N(\partial S))$, where $N(\partial S)$ is a regular neighbourhood of $\partial S$ in $\partial M$.

Note that, in this definition, we consider the homology group $H_2(M, N(\partial S))$, not $H_2(M, \partial M)$. Hence,
it is possible for $S$ to be taut and yet not have minimal Thurston complexity in its class in $H_2(M, \partial M)$. However, the following
simple lemma (which appears as Lemma A.7 in [20])
implies that this phenomenon cannot occur in an important case, when $\partial M$ is a (possibly empty) union of tori.
%It is for this reason that we restrict, in Theorem 1.5, to manifolds of this form.

\noindent {\bf Lemma 2.1.} {\sl Let $M$ be a compact orientable 3-manifold with boundary a (possibly empty)
union of tori. Let $S$ be a compact oriented properly embedded surface, such that no component of $\partial S$
bounds a disc in $\partial M$. Then $S$ has minimal Thurston
complexity in its class in $H_2(M, N(\partial S))$ if and only if it has minimal Thurston complexity in its class
in $H_2(M, \partial M)$.}

A sutured manifold $(M, \gamma)$ is {\sl taut} if $M$ is irreducible and $R_-(M)$ and $R_+(M)$ are both taut.

Suppose that $(M, \gamma)$ is a sutured manifold and that $S$ is a compact, transversely oriented surface
properly embedded in $M$. Suppose that $\partial S$ intersects $\gamma$ transversely. Then 
$M' = {\rm cl}(M - N(S))$ inherits a sutured manifold structure, since both the
parts of $\partial M'$ arising from $\partial M$ and the parts of $\partial M'$ arising from $N(S)$
have a well-defined transverse orientation. This is called a {\sl sutured manifold decomposition} and is denoted
$$(M, \gamma) \buildrel S \over \longrightarrow (M', \gamma').$$
It is called {\sl taut} if $(M, \gamma)$ and $(M', \gamma')$ are taut. 
The surface $S$ is often called a {\sl decomposing surface}. Note that in a
taut decomposition, we do not require $S$ to be taut, although this will often be the case.

\vskip 18pt
\centerline{
\epsfxsize=5.5in
\epsfbox{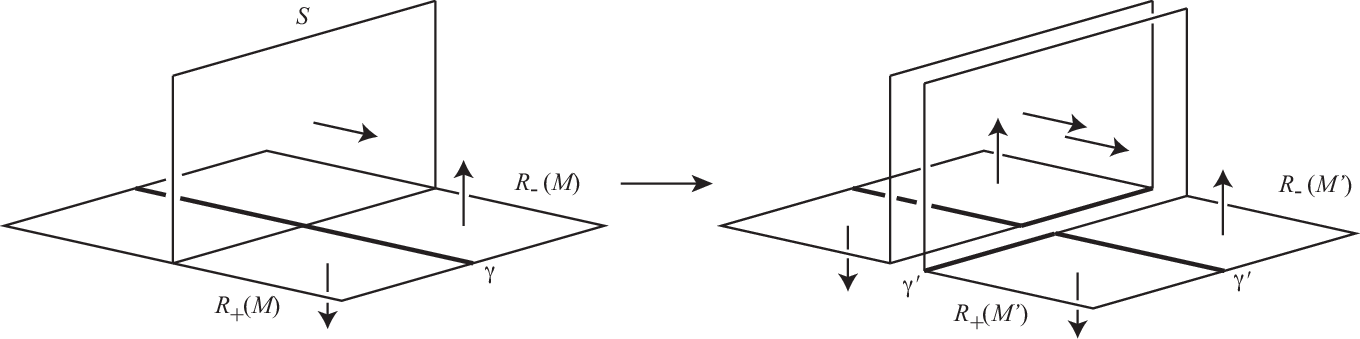}
}
\vskip 6pt
\centerline{Figure 2: Sutured manifold decomposition}

In general, decomposing a taut sutured manifold $(M, \gamma)$ along a taut surface $S$ does not result
in a taut decomposition. This is because the resulting sutured manifold $(M', \gamma')$ need not be taut.
However, there is one important situation where this is the case, which is summarised in the following
lemma.

\noindent {\bf Lemma 2.2.} {\sl Let $M$ be a compact orientable irreducible 3-manifold with boundary
a (possibly empty) collection of incompressible tori. Give $M$ a sutured manifepsold structure $(M, \gamma)$
where $R_+(M) = \partial M$ and $R_-(M) = \emptyset$. Let $S$ be a compact oriented properly embedded surface, such that
$S \cap \partial M$ is a (possibly empty) collection of essential simple closed curves. Suppose also that
no collection of annular components of $S$ is trivial in $H_2(M, \partial M)$, and that no component of $S$ is a 2-sphere.
Then, if $S$ is taut, so is the decomposition 
$$(M, \gamma) \buildrel S \over \longrightarrow (M', \gamma').$$
}

\noindent {\sl Proof.} The hypotheses on $M$ imply that $(M, \gamma)$ is taut. So, we need only show that $(M', \gamma')$ is taut.

Note first no component of $S$ is a disc, because $S \cap \partial M$ is essential and $\partial M$ is incompressible.
Hence, no component of $R_+(M')$ or $R_-(M')$ is a sphere or disc. Therefore,
$$\chi_-(R_+(M')) = -\chi(R_+(M')) = -\chi(S) - \chi(\partial M) = -\chi(S) = \chi_-(S) = -\chi(R_-(M')) = \chi_-(R_-(M')).$$
Hence, $R_+(M')$ has minimal Thurston complexity in its class in $H_2(M', N(\gamma'))$ if and only if
$R_-(M')$ does. But $R_-(M')$ is a copy of $S$. So if $R_-(M')$ was homologous in $H_2(M', N(\gamma'))$ to a surface with
smaller Thurston complexity, then $S$ would be homologous in $H_2(M, N(\partial S))$ to a surface
with smaller Thurston complexity, contradicting its tautness. 

We now show that $M'$ is irreducible. Consider a reducing sphere for $M'$. This bounds a ball in $M$.
The intersection between this ball and $S$ must be non-empty. Since $S$ is incompressible,
it must therefore contain a 2-sphere component, contrary to our assumption.

Now suppose that $R_-(M')$ is compressible in $M'$. But $R_-(M')$ is parallel in $M$ to $S$, and so this would imply
that $S$ is compressible. Finally suppose that $R_+(M')$ is compressible in $M'$. Then, since $R_+(M')$ minimises
Thurston complexity, this implies that $R_+(M')$ contains a component that is a compressible annulus or a 
compressible torus. Now, $R_+(M')$ consists of a copy of $S$, with annuli from $\partial M$ attached, plus possibly some toral
components of $\partial M$. So,
if $R_+(M')$ contains a compressible annulus, then $\partial M$ or $S$ compresses in $M$, which is contrary
to assumption. So, $R_+(M')$ contains a compressible torus, which therefore bounds a solid torus, by
the irreducibility of $M'$. This solid torus gives a homology in $H_2(M, \partial M)$ between
a collection of annular components of $S$ and a collection of annuli in $\partial M$. This again is
contrary to assumption. $\square$

One of the most important results in sutured manifold theory is summarised in the
phrase `tautness usually pulls back'. The precise theorem is as follows. (See Theorem 3.6 in [25]).

\noindent {\bf Theorem 2.3.} {\sl Let
$$(M, \gamma) \buildrel S \over \longrightarrow (M', \gamma')$$
be a sutured manifold decomposition. Suppose that no component of $\partial S$
bounds a disc in $\partial M$ disjoint from $\gamma$, and that no component of $S$
is a disc disjoint from $\gamma$ that forms a compression disc for a solid toral component
of $(M, \gamma)$ with no sutures.
Then if $(M', \gamma')$ is taut, so is $(M, \gamma)$ and so is $S$.}

As a consequence, if we have a sequence of sutured manifold decompositions, each
satisfying the requirements of Theorem 2.3, and we can show that the final sutured manifold
is taut, then the entire sequence of decompositions is taut, and each
decomposing surface is taut. It is a reasonably straightforward matter to verify
that a sutured manifold structure on a collection of 3-balls is taut. Indeed, it is
clear that it is taut if and only if each ball contains at most one suture. Thus,
this forms the basis for our certificate for Thurston norm. The existence of such a sequence
of decompositions is guaranteed by the following central result of Gabai [3]
(see also Theorems 2.6 and 4.19 in [25]).

\noindent {\bf Theorem 2.4.} {\sl Let $(M, \gamma)$ be a taut sutured manifold,
and let $z$ be a non-trivial element of $H_2(M, \partial M)$. Then there exists a
sequence of taut decompositions
$$(M, \gamma) = (M_1, \gamma_1) \buildrel S_1 \over \longrightarrow \dots
\buildrel S_n \over \longrightarrow (M_{n+1}, \gamma_{n+1})$$
such that
\item{(i)} no component of any $\partial S_i$ bounds a disc in $\partial M_i$ disjoint from $\gamma_i$;
\item{(ii)} no component of any $S_i$ is a disc disjoint from $\gamma_i$;
\item{(iii)} $[S_1, \partial S_1] = z \in H_2(M, \partial M)$;
\item{(iv)} $(M_{n+1}, \gamma_{n+1})$ is a collection of taut 3-balls.

}

The {\sl length} of this sequence of decompositions is $n$.

There are two types of decomposing surface that arise frequently.
A {\sl product disc} in a sutured manifold $(M, \gamma)$ is a properly embedded
disc that intersects $\gamma$ transversely at two points. A {\sl product annulus}
is a properly embedded annulus $A$ that is disjoint from $\gamma$, and that has
one component of $\partial A$ in $R_-(M)$ and one component of $\partial A$ in $R_+(M)$.
A product annulus $A$ is {\sl trivial} if some component of $\partial A$ bounds a
disc in $R_-(M)$ or $R_+(M)$. (Note that this latter definition is a minor variation of the
one given in Definition 4.1 of [25].)

Here, we have the stronger form of Theorem 2.3 (see Lemma 4.2 in [25]).

\noindent {\bf Proposition 2.5.} {\sl Let 
$$(M,\gamma) \buildrel S \over \longrightarrow (M', \gamma')$$
be a sutured manifold decomposition, where $S$ is either a product disc
or a non-trivial product annulus. Then $(M,\gamma)$ is taut if and only if $(M', \gamma')$ is taut.}

Finally, a {\sl product sutured manifold} is of the form $F\times [-1,1]$, where $F$ is some compact orientable surface and
where $\gamma = \partial F \times \{ 0 \}$. We will frequently refer to the product sutured manifold simply by $F \times [-1,1]$,
with the understanding that $\gamma = \partial F \times \{ 0 \}$.
We note that when a decomposition
$$(M,\gamma) \buildrel S \over \longrightarrow (M', \gamma')$$
is performed along a product disc or non-trivial product annulus, then $(M, \gamma)$ is a product sutured manifold
if and only if $(M' ,\gamma')$ is.

\vskip 18pt
\centerline{\caps 3. Decorated sutured manifolds}
\vskip 6pt

Unfortunately, we must consider decompositions
$$(M, \gamma) \buildrel S \over \longrightarrow (M', \gamma')$$
where some curves of $\partial S$ bound discs in $\partial M$ disjoint from $\gamma$. This leads
to some complications, because Theorem 2.3 does not apply, and so it is not obvious that $S$ can
be used as part of a certificate for Thurston norm. To get around this problem, we keep track
of the curves of $\partial S$ that bound discs in $\partial M$ disjoint from $\gamma$. These give rise to `special' sutures of $\gamma'$. 
We want to be able to distinguish these sutures, and so we now introduce a structure where this is possible, called a decorated sutured manifold.

\vskip 6pt
\noindent {\caps 3.1. Definition}
\vskip 6pt

We define a {\sl decorated sutured manifold} to be a sutured manifold $(M, \gamma)$ where some of the sutures are distinguished. These distinguished sutures
are called {\sl u-sutures}, where `u' stands for untouchable.

Given a decorated sutured manifold $(M,\gamma)$, its {\sl canonical enlargement}, denoted $E(M, \gamma)$ is an undecorated sutured manifold that is obtained by removing each u-suture and attaching a 2-handle $D^2 \times [0,1]$
along this curve. One component of $D^2 \times \{ 0,1\}$ lies in $R_-(E(M,\gamma))$ and the other lies in 
$R_+(E(M,\gamma))$. We denote the underlying 3-manifold of this enlargement by $E(M)$.

\vskip 6pt
\noindent {\caps 3.2. Pre-balls and pre-spherical products}
\vskip 6pt

A {\sl pre-ball} is a decorated sutured manifold of the form $P \times [-1,1]$, where $P $ is a compact connected planar surface.
The sutures are required to be $\partial P \times \{ 0 \}$ and exactly one is not a u-suture.
Note that the canonical enlargement of a pre-ball is a taut sutured ball.

A {\sl pre-spherical product manifold} is a decorated sutured manifold of the form $P \times [-1,1]$, 
where $P $ is a compact connected planar surface, where $\partial P \times \{ 0 \}$ forms the sutures, and every
suture is a u-suture. In this case, its canonical enlargement is a product sutured manifold
of the form $S^2 \times [-1,1]$.

\vskip 6pt
\noindent {\caps 3.3. Trivial curves}
\vskip 6pt

Let $(M, \gamma)$ be a decorated sutured manifold, and let $C \subset \partial M$ be a simple closed curve disjoint from the sutures.
Then $C$ is {\sl trivial} if there is a planar surface $P$ embedded in $R_\pm(M)$ such that
\item{(i)} $C$ is a component of $\partial P$; and
\item{(ii)} $P \cap \gamma$ is the remaining components of $\partial P$, and each is a u-suture.

\noindent The planar surface $P$ is called the {\sl trivialising planar surface}.

An alternative way of defining a trivial curve is that it is a simple closed curve in $R_\pm(M)$ that bounds a disc
in the boundary of $E(M)$ disjoint from the sutures of $E(M, \gamma)$. In this way, we see that a trivial curve
cannot bound trivialising planar surfaces on both sides, unless some component of $\partial E(M, \gamma)$ is a 
sphere with no sutures.

Let $S$ be a surface properly embedded in $M$.
Then a {\sl trivial boundary curve} of $S$ is a boundary component of $S$ that is trivial. Note that, by definition,
trivial boundary curves are disjoint from the sutures.

\vskip 6pt
\noindent {\caps 3.4. Allowable decompositions and hierarchies}
\vskip 6pt

If $(M, \gamma)$ is a decorated sutured manifold, then a sutured manifold decomposition 
$$(M, \gamma) \buildrel S \over \longrightarrow (M', \gamma')$$
is {\sl allowable} provided that
\item{(i)} $S$ is disjoint from the u-sutures of $\gamma$, and
\item{(ii)} if $C$ is any trivial boundary curve for $S$, then its trivialising planar surface inherits a transverse orientation from $R_\pm(M)$ that agrees with $S$ near $C$.

\noindent In this case, $(M', \gamma')$ inherits a decoration, by declaring that a component of 
$\gamma'$ is a u-suture if one of the following holds:
\item{(i)} it came from a u-suture of $\gamma$, or
\item{(ii)} it came from a trivial boundary curve of $S$.

\vskip 12pt
\centerline{
\epsfxsize=3.5in
\epsfbox{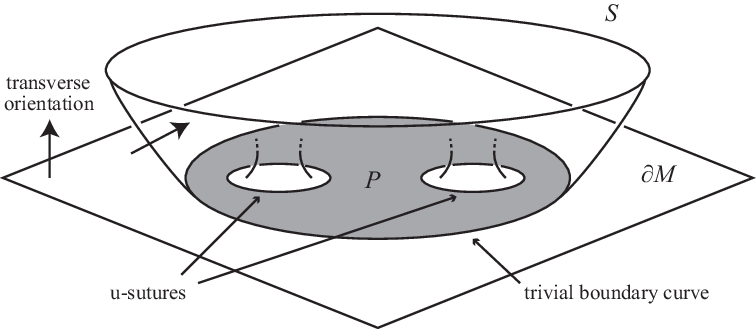}
}
\vskip 6pt
\centerline{Figure 3: Trivial boundary curve in an allowable decomposition}
\vfill\eject

A sequence of allowable decompositions is called an {\sl allowable hierarchy} if
\item{(i)} the final sutured manifold is a product manifold, no component of which pre-spherical; and
\item{(ii)} no component of any decomposing surface is planar, disjoint from the sutures and where all but at most one
of its boundary components is trivial.

We say that the allowable hierarchy is {\sl complete} if the final sutured manifold is a collection of pre-balls.

\vskip 6pt
\noindent {\caps 3.5. The utility of allowable hierarchies}
\vskip 6pt

Before we go any further with the development of the theory of allowable hierarchies, we state a result
which explains why they are useful.

\noindent {\bf Theorem 3.1.} {\sl Let $M$ be a compact orientable irreducible 3-manifold with boundary a
(possibly empty) union of tori. Give $M$ the sutured manifold structure with $R_+(M) = \partial M$ and $R_-(M) = \emptyset$.
Let $S$ be a compact oriented properly embedded surface with no 2-sphere or disc components, and give $\partial S$
the orientation that it inherits from $S$. Suppose that the intersection
between $S$ and each component of $\partial M$ is a (possibly empty) collection of essential curves that
are parallel as oriented curves in $\partial M$.
Then $S$ is incompressible and has minimal Thurston complexity in its class in $H_2(M, \partial M)$ if and only if
it is the first surface in an allowable hierarchy.}

Thus, allowable hierarchies will be the method that we employ to certify Thurston norm.
We will prove Theorem 3.1 in Section 3.10.

\vskip 6pt
\noindent {\caps 3.6. The canonical enlargement of a sutured manifold decomposition}
\vskip 6pt

Given an allowable decomposition
$$(M,\gamma) \buildrel S \over \longrightarrow (M',\gamma')$$
between decorated sutured manifolds, there is an associated decomposition
$$E(M,\gamma) \buildrel E(S) \over \longrightarrow E(M', \gamma'),$$
called its {\sl canonical enlargement} and
which is defined as follows. Consider any trivial simple closed curve of $\partial S$.
Attach its trivialising planar surface $P$ to $S$, and push this new part of the surface a little into the interior of $M$. 
Note that it is possible for two different trivialising planar surfaces $P$ and $P'$ to intersect. 
But their boundary components that are not u-sutures are disjoint. Hence, $P$ and $P'$ are nested. If $P'$ is contained in $P$, say,
then we push $P'$ a little further into the interior of $M$. 
We now extend this surface into the 2-handles that are attached to the u-sutures of $M$. For each such 2-handle, we insert a 
collection of discs, so that the resulting surface is properly embedded. The resulting surface is $E(S)$. Then  it is clear that decomposing $E(M, \gamma)$ along
$E(S)$ gives the same sutured manifold as $E(M', \gamma')$. For a pictorial proof of this fact, see Figure 4.

It is partly because of this canonical enlargement that u-sutures are `untouchable', in the sense that decomposing surfaces must avoid them.
If, on the contrary, a decomposing surface $S$ were allowed to run over a u-suture, then there would be no good way of defining $E(S)$, because
the boundary of $S$ would run over the attaching annulus of the 2-handle.

We now note that $E(S)$ satisfies one of the key hypotheses of Theorem 2.3.

\noindent {\bf Lemma 3.2.} {\sl No boundary curve of $E(S)$ is trivial in $E(M, \gamma)$.}

\vskip 18pt
\centerline{
\epsfxsize=6in
\epsfbox{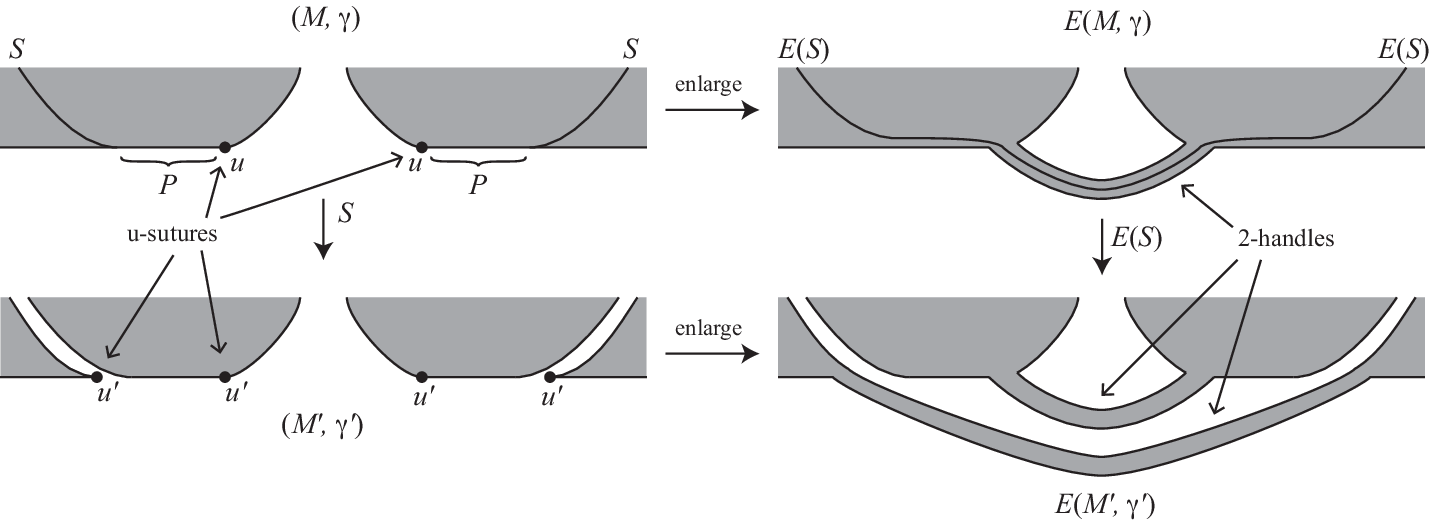}
}
\vskip 6pt
\centerline{Figure 4: Interchanging the order of enlargement and decomposition}

\noindent {\sl Proof.} Suppose that a boundary curve of $E(S)$ is trivial. It therefore bounds a disc $D$ in $\partial E(M)$
disjoint from the sutures. By choosing the curve of $\partial E(S)$ appropriately, we may ensure that the interior of $D$
is disjoint from $E(S)$. The intersection between this disc $D$ and $M$ is a planar surface $P$. All but one of its boundary curves are u-sutures of $\gamma$.
The remaining boundary component is $P \cap \partial S$. The orientations of $P$ and $S$ must agree locally near $P \cap \partial S$, by the definition of an
allowable decomposition. Hence, this boundary curve of $\partial S$ is removed in the construction of $E(S)$, which is a contradiction. $\square$

\vskip 6pt
\noindent {\caps 3.7. The canonical enlargement of an allowable hierarchy}
\vskip 6pt

\noindent {\bf Lemma 3.3.} {\sl The canonical enlargement of an allowable hierarchy of decorated sutured manifolds is taut
and each decomposing surface in this canonical enlargement is taut.}

\noindent {\sl Proof.} Let 
$$E(M_1, \gamma_1) \buildrel E(S_1) \over \longrightarrow E(M_2, \gamma_2) \buildrel E(S_2) \over \longrightarrow \dots
\buildrel E(S_n) \over \longrightarrow E(M_{n+1}, \gamma_{n+1})$$ 
be this enlargement. By assumption, $(M_{n+1}, \gamma_{n+1})$ is a product sutured manifold, no component
of which is pre-spherical. Hence, $E(M_{n+1}, \gamma_{n+1})$ is a taut product manifold.
By Lemma 3.2, no boundary curve of any $E(S_i)$ bounds a disc in $\partial E(M_i)$ disjoint from the sutures. 
Also, by (ii) in the definition of an allowable hierarchy, no component of any $E(S_i)$ is a disc disjoint from the sutures.
Hence, by Theorem 2.3, each decomposition in the sequence is taut and each decomposing surface is taut. $\square$

\noindent {\bf Corollary 3.4.} {\sl Let
$$(M_1, \gamma_1) \buildrel S_1 \over \longrightarrow (M_2, \gamma_2) \buildrel S_2 \over \longrightarrow \dots
\buildrel S_n \over \longrightarrow (M_{n+1}, \gamma_{n+1})$$ 
be an allowable hierarchy. Suppose that $(M_1, \gamma_1)$ has no u-sutures and that no boundary
component of $S_1$ is trivial in $(M_1, \gamma_1)$. Then $S_1$ is taut in $M_1$.}

\noindent {\sl Proof.} By Lemma 3.3, the canonical enlargement of this allowable hierarchy is taut and each decomposing surface is taut. 
But $(M_1, \gamma_1)$ has no u-sutures and therefore $E(M_1,\gamma_1)$ is just $(M_1,  \gamma_1)$. 
Also, since $S_1$ has no trivial boundary curves, $E(S_1)$ is just $S_1$.
Hence, $S_1$ is taut in $M_1$. $\square$

\vskip 6pt
\noindent {\caps 3.8. An alternative interpretation of complete allowable hierarchies}
\vskip 6pt

The following lemma provides a useful alternative way of understanding complete allowable hierarchies.
Its proof is immediate.

\noindent {\bf Lemma 3.5.} {\sl Let
$$(M_1,\gamma_1) \buildrel S_1 \over \longrightarrow \dots \buildrel S_n \over \longrightarrow (M_{n+1}, \gamma_{n+1})$$
be a sequence of allowable decompositions between decorated sutured manifolds. Then no boundary curve
of $E(S_i)$ is trivial in $E(M_i, \gamma_i)$. Furthermore, this forms a complete allowable hierarchy
if and only if its canonical enlargement satisfies the following conditions:
\item{(i)} the final manifold $E(M_{n+1}, \gamma_{n+1})$ is a collection of taut 3-balls, each of the form $D^2 \times [-1,1]$,
where the sutures are $\partial D^2 \times \{ 0 \}$ and where the co-cores of the attached 2-handles are vertical
in this product structure;
\item{(ii)} no component of $E(S_i)$ is a sphere or disc disjoint from the sutures.

}

\vskip 6pt
\noindent {\caps 3.9. Enlargement and tautness}
\vskip 6pt

\noindent {\bf Proposition 3.6.} {\sl Let $(M, \gamma)$ be a decorated sutured manifold. If $E(M, \gamma)$ is taut,
then so is $(M,\gamma)$.}

\noindent {\sl Proof.} We will prove the following stronger statement. If $(M, \gamma)$ is a decorated sutured manifold,
a single suture is removed, a 2-handle is attached along it and the resulting sutured manifold $(M', \gamma')$ is
taut, then $(M, \gamma)$ is also taut.

Note first that we may assume that $M$ is connected. For we may restrict attention
to the component of $M$ to which the 2-handle is attached.

Consider a 2-sphere $S$ properly embedded in $M$. Since $M'$ is irreducible,
$S$ bounds a ball in $M'$. This ball is disjoint from $\partial M'$, and so is disjoint from the attached
2-handle. Therefore, $S$ bounds a ball in $M$. Hence, $M$ is irreducible.

Suppose now that $R_\pm(M)$ is compressible, via a compression disc $D$. Since $R_\pm(M')$
is incompressible, $\partial D$ bounds a disc $D'$ in $\partial M'$ disjoint from the sutures. 
Note that at least one of the two discs of intersection between the 2-handle and
$R_\pm(M')$ misses $D'$, because one of these discs lies in $R_-(M')$ and the other lies in
$R_+(M')$. Since $M'$ is irreducible, $D \cup D'$ bounds a ball in $M'$. This ball cannot intersect the
2-handle, because this would imply that $D$ intersected this 2-handle, whereas $D$ lies in $M$.
Hence, we deduce that the ball lies in $M$, and in particular, $D'$ lies in $M$. Therefore,
this was not a compression disc for $R_\pm(M)$, which is a contradiction.

We now show that $R_\pm(M)$ has minimal Thurston complexity in its class in $H_2(M, N(\gamma))$.
Consider another surface $S$ in the same class. We may assume that this runs over the attaching
annulus of the 2-handle in a single essential curve. Hence, we may extend $S$ to a surface $\tilde S$
in $M'$, by attaching a disc in this 2-handle. Note that $\tilde S$ is in the same class as
$R_\pm(M')$ in $H_2(M', N(\gamma'))$. Since this has minimal Thurston complexity,
we deduce that $\chi_-(\tilde S) \geq \chi_-(R_\pm(M'))$. 

Now, we may assume that $R_\pm(M')$
has the same number of sphere and disc components as $R_\pm(M)$. For if a disc component
of $R_\pm(M')$ is created, then $(M',\gamma')$ is a taut 3-ball and so $R_\pm(M)$ consists of incompressible annuli,
and therefore has minimal Thurston complexity. If a sphere component of $R_\pm(M')$ is
created, this bounds a ball, by the irreducibility of $M'$, and hence one of $R_-(M')$
or $R_+(M')$ is empty, which is impossible, because they both have non-empty intersection with the
2-handle. Similarly, $\tilde S$ has the same number of sphere and disc components as $S$. 
Hence,
$$\chi_-(S) = \chi_-(\tilde S) + 1 \geq  \chi_-(R_\pm(M')) +1 = \chi_-(R_\pm(M))$$
which implies that $R_\pm(M)$ does indeed minimise Thurston complexity in its homology class.
$\square$

\noindent {\bf Corollary 3.7.} {\sl An allowable hierarchy of decorated sutured manifolds is taut.}

\noindent {\sl Proof.} By Lemma 3.3, the canonical enlargement of an allowable hierarchy is taut.
Hence, by Proposition 3.6, each of the original decorated sutured manifolds is taut. $\square$

Note that the converse of Proposition 3.6 is false. Consider the case where $E(M,\gamma)$ is a
3-ball with three parallel sutures in its boundary. Then $R_-(E(M,\gamma))$ is an annulus
and a disc, as is $R_+(E(M,\gamma))$. In particular, $E(M,\gamma)$ is not taut. However,
we may pick a properly embedded arc in $E(M)$, running between the disc component of $R_-(E(M,\gamma))$
and the disc component of $R_+(E(M,\gamma))$. Remove a regular neighbourhood
of this arc to form $(M,\gamma)$. Then $M$ is irreducible, and $R_\pm(M)$ consists of
incompressible annuli. Therefore, $(M, \gamma)$ is taut.
This example shows that admitting an allowable hierarchy is a stronger condition than
simply being taut.

\vskip 6pt
\noindent {\caps 3.10. Allowable hierarchies and tautness}
\vskip 6pt

We are now in a position to be able to prove Theorem 3.1.

\noindent {\sl Proof.} Let $M$ and $S$ be as in the statement of Theorem 3.1. Suppose that $S$ is incompressible
and has minimal Thurston complexity in its class in $H_2(M, \partial M)$. Then by Lemma 2.1,
$S$ is taut. So, by Lemma 2.2, the decomposition
$$(M, \gamma) \buildrel S \over \longrightarrow (M', \gamma')$$
is taut. Note that no collection of annular components of $S$ is trivial in $H_2(M, \partial M)$ by our
assumption about the orientations on $\partial S$, and so Lemma 2.2 does apply. Note also that
$(M', \gamma')$ has no u-sutures. By Theorem 2.4, $(M', \gamma')$ admits a taut sutured manifold
hierarchy, where no decomposing surface has a disc component disjoint from the sutures.
Also, at no stage does a boundary curve bound a disc in the boundary of the manifold disjoint from the sutures. So, no u-sutures
are created. Thus, this forms an allowable hierarchy.

Conversely, suppose that $S$ is the first surface in an allowable hierarchy. This hierarchy is taut by Corollary 3.7,
and in particular the first decomposition is taut. So, by Theorem 2.3, $S$ is taut, and so by Lemma 2.1,
$S$ is incompressible and has minimal Thurston complexity in its class in $H_2(M, \partial M)$. $\square$

\vskip 6pt
\noindent {\caps 3.11. Modifying the decoration}
\vskip 6pt

\noindent {\bf Lemma 3.8.} {\sl Let $(M, \gamma)$ be a decorated sutured manifold, and let $(M' , \gamma')$ be
the same sutured manifold, but where some of the u-sutures have been replaced by ordinary sutures.
Suppose that $(M, \gamma)$ has an allowable hierarchy. Then so does $(M' ,\gamma')$, with the
same reduced length.}

\noindent {\sl Proof.} Let
$$(M, \gamma) = (M_1, \gamma_1) \buildrel S_1 \over \longrightarrow \dots \buildrel S_n \over \longrightarrow (M_{n+1}, \gamma_{n+1})$$
be the allowable hierarchy for $(M, \gamma)$. We can view this as a hierarchy
$$(M', \gamma') = (M'_1, \gamma'_1) \buildrel S'_1 \over \longrightarrow \dots \buildrel S'_n \over \longrightarrow (M'_{n+1}, \gamma'_{n+1})$$
for $(M', \gamma')$. We show by induction the following:
\item{(i)} each $(M'_i, \gamma'_i)$ inherits a decoration;
\item{(ii)} each u-suture of $\gamma'_i$ is a u-suture of $\gamma_i$;
\item{(iii)} if a curve in $\partial M'_i - \gamma_i$ is trivial, then so is the corresponding curve in $\partial M_i - \gamma_i$, with the same trivialising planar surface;
\item{(iv)} decomposition along $S'_i$ is allowable;
\item{(v)} no component of $S'_i$ is planar, disjoint from the sutures and where all but at most one of its boundary components is trivial.

The induction starts with the case $i = 1$, where (i) and (ii) are part of the hypotheses of the lemma. The three remaining claims (iii), (iv) and (v) 
in the case $i = 1$ are proved in the same way as in the inductive step. So, suppose that $(M'_i, \gamma'_i)$ satisfies (i) and (ii).
Then any trivialising planar surface for a curve in $\partial M'_i - \gamma_i'$ is then a trivialising planar surface in $\partial M_i - \gamma_i$.
Hence, we have (iii). Then (ii) and (iii) imply that $S'_i$ is allowable, giving (iv). Also, by (ii), a component of $S'_i$ violating (v) would 
give a component of $S_i$ that could not be part of an allowable hierarchy. So, we obtain (v). Since decomposition along $S'_i$
is allowable, $(M'_{i+1},  \gamma'_{i+1})$ inherits a decoration, giving (i). The u-sutures of $\gamma'_{i+1}$ come either
from u-sutures of $\gamma'_i$ or from trivial boundary curves of $S'_i$. In the former case, (ii) for  $(M'_i, \gamma'_i)$
implies that it is also a u-suture of $\gamma_i$. In the latter case, this is also a trivial boundary curve of $S_i$, by (iii).
Hence, we obtain (ii) for $(M'_{i+1},  \gamma'_{i+1})$. This completes the induction.

The final manifold $(M_{n+1}, \gamma_{n+1})$ is a product sutured manifold, no component of which is pre-spherical. With the new decoration, some of the u-sutures
have become ordinary sutures. So, $(M'_{n+1}, \gamma'_{n+1})$ is still a product sutured manifold, no component of which is pre-spherical. 
$\square$

\vskip 6pt
\noindent {\caps 3.12. Canonical extensions that are balls}
\vskip 6pt

\noindent {\bf Lemma 3.9.} {\sl Let $(M, \gamma)$ be a decorated sutured manifold that
admits an allowable hierarchy. Suppose that $E(M)$ is a 3-ball. Then $(M, \gamma)$ is
a pre-ball.}

\noindent {\sl Proof.} Note first that, by Lemma 3.3, $E(M, \gamma)$ is taut, and hence
is a product sutured 3-ball $B$. The manifold $(M, \gamma)$ is obtained from $E(M, \gamma)$
by removing the attached 2-handles. The co-cores of these 2-handles form a tangle $t$ in $B$.
Thus, the lemma is equivalent to the assertion is that there is an ambient isotopy,
keeping the suture in $B$ fixed, after which the tangle respects the product structure on $B$.
Equivalently, there is a collection of disjoint embedded discs $D$ embedded in $B$, such that 
\item{(i)} the intersection between each component $D'$ of $D$ and $t$ is a single component of $t$
in $\partial D'$;
\item{(ii)} the remainder of the boundary of $D'$ lies in $\partial B$;
\item{(iii)} the intersection between each component of $D$ and the suture of $B$ is a single point;
\item{(iv)} each component of $t$ lies in a component of $D$.

Note that there is considerable flexibility in the choice of $D$. In particular, it may be chosen
so that its intersection with $R_+(B)$ is any given collection of disjoint embedded arcs, with the property
that each arc starts at a point of $t \cap R_+(B)$ and ends on the suture, and each point of $t \cap R_+(B)$ lies
at the endpoint of such an arc.

Now, $(M, \gamma)$ admits an allowable hierarchy where each surface is connected.
We prove the lemma by induction on the length of such a hierarchy.
The induction starts trivially, because a sutured manifold that has an allowable hierarchy with length
zero is a product, and the only product sutured manifold with canonical extension that is a ball is a pre-ball.
So, we consider the inductive step. Let 
$$(M, \gamma) \buildrel S \over \longrightarrow (M_S, \gamma_S)$$
be the first decomposition in the hierarchy. By assumption $E(M, \gamma)$ is a ball.
It is taut by Lemma 3.3. The surface $E(S)$ is taut in $E(M, \gamma)$ by Lemma 3.3.
In particular, it is incompressible, and is therefore a disc.
Therefore, $E(M_S, \gamma_S)$ is two 3-balls $B_1$ and $B_2$. By induction therefore,
$(M_S, \gamma_S)$ is two pre-balls. Let $t_1$ and $t_2$ be the tangles
in $B_1$ and $B_2$ forming the co-cores of the attached 2-handles. The union of $t_1$ and
$t_2$ forms a tangle $t$ in $B$ that is the co-cores of the 2-handles there. Suppose that
the copies of $E(S)$ in $B_1$ and $B_2$ lie in $R_-(B_1)$ and $R_+(B_2)$, say.
For each point of $t \cap E(S)$, pick an arc in $E(S)$ running from that point to $\partial E(S)$.
We may arrange that these arcs $\alpha$ are disjoint, and that they all end on a suture of $B_1$, say.
These arcs $\alpha$ are arcs in $R_-(B_1)$. By adding in extra arcs disjoint from $E(S)$,
extend them to a collection of arcs in $R_-(B_1)$ running from every point of $t \cap R_-(B_1)$
to the suture of $B_1$. The arcs $\alpha$ also lie in $R_+(B_2)$, but they do not end on the suture
there. Without changing their intersection with $E(S)$, extend them so that they do end on
the suture of $B_2$. Furthermore, add in extra arcs if necessary so that all point of
$t \cap R_+(B_2)$ are at the start of such an arc. We may find a collection of discs
$D_1$ in $B_1$ satisfying (i)-(iv) above such that $D_1$ intersects $R_-(B_1)$ in the given collection of arcs.
We may find a similar collection of discs $D_2$ in $B_2$. Then $D_1 \cup D_2$ is the required
collection of discs in $B$ satisfying (i)-(iv). Hence, $(M, \gamma)$ is a pre-ball,
as required. $\square$

\vfill\eject
\noindent {\caps 3.13. Some consequences of atoroidality}
\vskip 6pt

Recall that a compact orientable 3-manifold $M$ is {\sl atoroidal} if any properly embedded incompressible
torus is boundary parallel. This will be a useful hypothesis at various points in this paper.
In this subsection, we collate a few consequences of atoroidality. In Section 11, we will show
how the proof of the main theorems of this paper may be reduced to this case.

\noindent {\bf Lemma 3.10.} {\sl Let $(M, \gamma)$ be a connected sutured manifold with no u-sutures. Let
$$(M, \gamma) = (M_1, \gamma_1) \buildrel S_1 \over \longrightarrow \dots \buildrel S_n \over \longrightarrow (M_n, \gamma_n)$$
be an allowable hierarchy. Suppose that $M$ is atoroidal and that its boundary is not a single torus.
Then, for any $i > 1$, no component of $M_i$ has boundary a single incompressible torus with no u-sutures.}

\noindent {\sl Proof.} Suppose that $Y$ is such a component of $M_i$. 
Because $Y$ contains no u-sutures, it is its own canonical enlargement. The decomposing surfaces
$E(S_j)$ are all incompressible by Lemma 3.3. Therefore, the inclusion of $Y$ into
$E(M_1) = M$ is $\pi_1$-injective. Therefore $\partial Y$ is an incompressible torus in $M$. It is not boundary
parallel, because $\partial M$ would then be a single torus. This contradicts the
atoroidality of $M$. $\square$

\noindent {\bf Lemma 3.11.} {\sl Let $(M, \gamma)$ be a taut decorated sutured manifold.
Let
$$(M, \gamma) = (M_1, \gamma_1) \buildrel S_1 \over \longrightarrow \dots \buildrel S_n \over \longrightarrow (M_{n+1}, \gamma_{n+1})$$
be an allowable hierarchy.
Suppose that $E(M, \gamma)$ is atoroidal. Then $E(M_i, \gamma_i)$ is atoroidal for each $i$.
Furthermore, if the only Seifert fibred components of $E(M, \gamma)$
are solid tori and copies of $T^2 \times I$, then the same is true for $E(M_i, \gamma_i)$.}

\noindent {\sl Proof.} Consider an incompressible torus $T$ properly embedded in $E(M_i, \gamma_i)$.
Since the inclusion of $E(M_i, \gamma_i)$ into $E(M, \gamma)$ is $\pi_1$-injective, by
Lemma 3.3, we deduce that $T$ is incompressible in $E(M, \gamma)$. As $E(M,\gamma)$ is atoroidal, $T$
is therefore boundary-parallel in $E(M, \gamma)$. Let $T \times I$ be the product region between $T$
and a boundary component of $E(M, \gamma)$. Now the surfaces $E(S_j)$ are incompressible by Lemma 3.3,
and so their intersection with $T \times I$ is boundary-parallel. Hence, the result of
decomposing $T \times I$ along the surfaces $E(S_j)$ is to retain a copy of $T \times I$.
We deduce that $T$ is boundary-parallel in $E(M_i, \gamma_i)$. So, $E(M_i, \gamma_i)$ is
atoroidal. Suppose now that some component of $E(M_i, \gamma_i)$ is Seifert fibred
and neither a solid torus nor a copy of $T^2 \times I$. Its boundary tori must be boundary parallel in $E(M, \gamma)$,
which implies that this component of $E(M, \gamma)$ is a Seifert fibred space which is neither a solid
torus nor a copy of $T^2 \times I$. $\square$

%\noindent {\bf Lemma 3.11.} {\sl Let Let $(M, \gamma)$ be a taut decorated sutured manifold.
%Let
%$$(M, \gamma) = (M_1, \gamma_1) \buildrel S_1 \over \longrightarrow \dots \buildrel S_n \over \longrightarrow (M_n, \gamma_n)$$
%be an allowable hierarchy.
%Suppose that $M$ is atoroidal and has no u-sutures. Then any toral component of any decomposing surface $S_i$ is boundary-parallel in $M$.}

%\noindent {\sl Proof.} Consider a toral component $T$ of some decomposing surface $S_i$. This is incompressible in $E(M_i, \gamma_i)$
%by Lemma 3.3. Since the inclusion of $E(M_i, \gamma_i)$ into $E(M, \gamma)$ is $\pi_1$-injective, $T$ is
%boundary-parallel in $E(M, \gamma) = M$. $\square$

\noindent {\bf Lemma 3.12.} {\sl Let $(M, \gamma)$ be a connected taut decorated sutured manifold that admits an allowable hierarchy. Suppose that
$E(M, \gamma)$ is atoroidal and not a Seifert fibre space other than a solid torus or a copy of $T^2 \times I$. Let $A_1$ and $A_2$ be properly embedded
incompressible annuli in $M$ that intersect in a collection of essential simple closed curves and with non-trivial boundary disjoint from $\gamma$.
Suppose that there is no homeomorphism $h \colon M \rightarrow M$ fixed on $\partial M$ such that
$h(A_1) \cap A_2$ consists of fewer simple closed curves. Then $A_1 \cap A_2$ consists of at most two curves.}

\noindent {\sl Proof.} Suppose that $A_1 \cap A_2$ consists of at least two curves. Consider an annular component $A'_2$ of $A_2 - {\rm int}(N(A_1))$
with neither boundary component in $\partial M$. Let $A'_1$ be the annulus in $A_1$ bounded by the two curves
of $\partial A'_2$.

We first show that the two components of $A'_2$ emanate from the same
side of $A_1$. Suppose not. Then consider the torus $T' = A'_1 \cup A_2'$. 
A Dehn twist $h$ around $T'$ may be chosen so that $|h(A_1) \cap A_2| < |A_1 \cap A_2|$,
contradicting our minimality assumption.

Again consider the torus $T' = A'_1 \cup A_2'$. This forms a boundary component of the manifold $Y$ obtained
from $E(M, \gamma)$ by cutting along $A_1$, then along $A_2'$. Since $A_1$ and $A'_2$ are incompressible, the inclusion
of each component of $Y$ into $E(M, \gamma)$ is $\pi_1$-injective. Let $Y'$ be the component of $Y$ containing $T'$.
So, if $T'$ is incompressible in $Y'$, then $T'$ is
boundary-parallel in $E(M, \gamma)$. If the product region between $T'$ and a component of $\partial E(M)$ contains $A_1 - A'_1$,
then we deduce that $A_1$ is boundary parallel. In this case, the lemma is clear. So, we may assume that
this product region $Y'$ does not contain $A_1 - A_1'$. On the other hand, if $T'$ is compressible in $Y'$,
then $Y'$ is a solid torus, because $Y'$ is irreducible. Give this solid torus a Seifert fibration, in such a way that $A_1'$ and $A'_2$
are each a union of fibres. Then $Y'$ contains an exceptional fibre. For otherwise, we could isotope
$A'_2$ across $Y'$ (together with any other components of $A_2 \cap Y'$) and thereby reduce
the number of curves of $A_1 \cap A_2$.

Now consider an annulus $A''_2$ of $A_2 - {\rm int}(N(A_1))$ adjacent to $A_2'$, and also disjoint from $\partial M$.
Such an annulus exists if $A_1 \cap A_2$ consists of at least three curves. Note that $A''_2$ emanates from
the other side of $A_1$ from $A'_2$. Let $A''_1$ be the annulus in $A_1$ bounded by $\partial A''_2$.
It is possible that $A'_1$ and $A''_1$ intersect. Let $T''$ denote the torus $A''_1 \cup A''_2$. As argued above,
we deduce that either $T''$ is boundary parallel in $E(M, \gamma)$ and that the product region $Y''$ has interior disjoint from $A_1$,
or $T''$ bounds a Seifert fibred solid torus $Y''$ with an exceptional fibre and with interior disjoint from $A_1$.

Let $Y'''$ be a regular neighbourhood of $Y' \cup Y''$. This is Seifert fibred, with planar base space, and where the
sum of the number of exceptional fibres and boundary components is at least three. So, $\partial Y'''$ is incompressible
in $Y'''$. All but one of its boundary components is boundary-parallel in $E(M, \gamma)$. The remaining component of
$\partial Y'''$ separates off a subset of $M$ with non-empty boundary. So, by the atoroidality of $E(M, \gamma)$, 
this component of $\partial Y'''$ is boundary-parallel in $E(M, \gamma)$. We deduce that $E(M, \gamma)$ 
is Seifert fibred, and not a solid torus or a copy of $T^2 \times I$, which is contrary to assumption. $\square$

A very similar argument, which we omit, also gives the following result.

\noindent {\bf Lemma 3.13.} {\sl Let $(M, \gamma)$ be a connected taut decorated sutured manifold that admits an allowable hierarchy. Suppose that
$E(M, \gamma)$ is atoroidal and not a Seifert fibre space other than a solid torus or a copy of $T^2 \times I$. 
Let $A$ be an annulus properly embedded in $M$ with non-trivial boundary disjoint
from $\gamma$, and let $T$ be a properly embedded incompressible torus.
Suppose that there is no homeomorphism $h \colon M \rightarrow M$ fixed on $\partial M$ such that
$h(T) \cap A$ consists of fewer simple closed curves than $T \cap A$. Then $T \cap A$ consists of at most one curve.}

\vskip 18pt
\centerline{\caps 4. Surfaces that extend to an allowable hierarchy}
\vskip 6pt

In this section, we will consider surfaces $S$ that form the first surface in an allowable hierarchy.
We will show that certain modifications can be made to $S$ that preserve this property.

\vskip 6pt
\noindent {\caps 4.1. The reduced length of a hierarchy}
\vskip 6pt

Given a sequence of sutured manifold decompositions
$$(M_1, \gamma_1) \buildrel S_1 \over \longrightarrow \dots \buildrel S_{n} \over \longrightarrow (M_{n+1}, \gamma_{n+1})$$
its {\sl reduced length} is the number of surfaces $S_i$ that are not a union of product discs, annuli disjoint from the sutures and tori.

\noindent {\bf Lemma 4.1.} {\sl Let $(M,\gamma)$ be a decorated sutured manifold that admits an allowable
hierarchy with reduced length zero. Then each component of $(M,\gamma)$ is either a product sutured manifold, which is not pre-spherical,
or a copy of $T^2 \times I$ with no sutures, an orientable $I$-bundle over the Klein bottle with no sutures, the union of
two such $I$-bundles glued along their boundary or a torus bundle over the circle.}

\noindent {\sl Proof.} We may assume that $M$ is connected. We may also assume that each surface in the allowable hierarchy for $(M, \gamma)$ is
connected. We will prove the lemma by induction on the length of this hierarchy. Let
$$(M, \gamma) \buildrel S \over \longrightarrow (M', \gamma')$$
be the first decomposition in the hierarchy. By induction, each component of $(M', \gamma')$ is either a product sutured manifold
which is not pre-spherical or a copy of
$T^2 \times I$ with no sutures or an orientable $I$-bundle over the Klein bottle with no sutures.

Suppose first that  $S$ is a product disc or non-trivial product annulus. Then the components
of $(M', \gamma')$ must be products, because they have non-empty intersection with
both $R_-(M')$ and $R_+(M')$. So, $(M, \gamma)$ is also a product, as required. 

Suppose that $S$ is an
annulus, disjoint from $\gamma$, with both boundary components in $R_-(M)$, say.
Then a copy of $S$ becomes a component of $R_+(M')$. Let $X$ be the
component of $M'$ containing this copy of $S$. It has non-empty intersection
with both $R_-(M')$ and $R_+(M')$, and so must be a product.
It is therefore homeomorphic to $S \times I$. There are now two cases to consider.

Suppose that $X$ does not contain the other copy of $S$.
Then $S$ is parallel to an annulus in $R_-(M)$. Hence, decomposing along $S$ simply
peels off a copy of $S \times I$. Therefore, each component of $(M, \gamma)$ has
the required form.
On the other hand, if $X$ contains both copies of $S$, then we deduce that
$S$ was non-separating in a component of $M$. This component of $M$ is, up to homeomorphism, 
obtained from $S \times I$ by identifying $S \times \{ 0 \}$ and $S \times \{ 1 \}$. So,
it is either a copy of $T^2 \times I$ or the orientable $I$-bundle over the Klein bottle.
In both cases, it contains no sutures.

Finally suppose that $S$ is a torus. Then each component of $M'$ contains at least one toral boundary component disjoint
from the sutures. Hence, inductively, each component of $M'$ is an $I$-bundle over the Klein bottle or a copy of $T^2 \times I$.
Therefore, $M$ is a copy of $T^2 \times I$ with no sutures, an orientable $I$-bundle over the Klein bottle with no sutures, the union of
two such $I$-bundles glued along their boundary or a torus bundle over the circle.
$\square$

\vskip 6pt
\noindent {\caps 4.2. Controlling trivial boundary curves}
\vskip 6pt

Let $S$ be an allowable decomposing surface for a decorated sutured manifold $(M, \gamma)$.
We say that a curve $C$ of $\partial S$ is {\sl parallel towards a u-suture} if it is disjoint from $\gamma$ and there is an annulus $A$ in
$R_\pm(M)$ such that one component of $\partial A$ is $C$, the other component of $\partial A$ is a u-suture
and the orientations on $A$ and $S$ agree near $C$. Note that this annulus $A$ is a trivialising planar surface 
for $C$.

\noindent {\bf Lemma 4.2.} {\sl Let $(M, \gamma)$ be a taut decorated sutured manifold that admits an allowable hierarchy.
Then it admits an allowable hierarchy in which each trivial boundary curve of each decomposing surface is
parallel towards a u-suture. Moreover, we may ensure that each decomposing surface is incompressible.
The new hierarchy has the same length as the original hierarchy, the same canonical extension
and no greater reduced length. Moreover, if some surface in the original hierarchy is non-separating,
then this remains true of the corresponding surface in the new hierarchy.}

\noindent {\sl Proof.} Let
$$(M_i, \gamma_i) \buildrel S_i \over \longrightarrow (M_{i+1}, \gamma_{i+1})$$
be a decomposition in the allowable hierarchy. Consider a trivial boundary curve $C$ of $S_i$,
with trivialising planar surface $P$. Suppose that $C$ is not parallel towards  a u-suture,
and hence that $P$ is not an annulus. We may assume
that if there are any curves of ${\rm int}(P) \cap S_i$, then these are parallel in $P$ to u-sutures.
Let $S'_i$ be the result of attaching $P$ to $S_i$,
and isotoping a little so that it becomes properly embedded and so that each component of $\partial P - C$
is parallel to a u-suture. Decompose $(M_i, \gamma_i)$ along $S'_i$ instead,
giving a new sutured manifold $(M'_{i+1}, \gamma'_{i+1})$. 
Then $(M'_{i+1}, \gamma_{i+1}')$ is obtained from $(M_{i+1}, \gamma_{i+1})$ by attaching a copy of
$P \times [0,1]$ to the u-suture coming from $C$. All later decomposing surfaces beyond $(M_{i+1}, \gamma_{i+1})$ avoid this
u-suture. So we can view them as forming a hierarchy, starting at $(M_{i+1}', \gamma_{i+1}')$
and ending at $(M'_{n+1}, \gamma'_{n+1})$, say. We note
$(M'_{n+1}, \gamma'_{n+1})$ is also obtained from $(M_{n+1}, \gamma_{n+1})$
by attaching a copy of $P \times [0,1]$ to the suture that is the copy of $C$. Hence, $(M'_{n+1}, \gamma'_{n+1})$
is also a product sutured manifold, no component of which is pre-spherical.  Note also that, for each $j$,
a curve in $\partial M_j$ disjoint from the sutures is trivial if and only if the corresponding curve in
$\partial M'_j$ is trivial. Furthermore, the canonical extensions $E(S_i)$ and $E(S'_i)$ are homeomorphic,
via a homeomorphism that respects intersection with the sutures. So no component of $E(S'_i)$
is a disc disjoint from the sutures. Using these observations, it is easy to see that we have found
an allowable hierarchy extending $S'_i$. This surface has fewer trivial boundary curves that are not parallel
towards u-sutures. So, repeating this, we end with the required hierarchy. Note that these modifications to the
hierarchy do not change its length. Moreover, if a component of a decomposing surface $S_i$ was a product disc, torus or annulus disjoint from
$\gamma_i$, then this is true of the corresponding component of the new surface $S'_i$. In the annular
case, this follows from the fact that $S_i$ has no component that is planar, disjoint from $\gamma_i$ and with
all but at most one boundary component trivial.
So, the reduced length of the new hierarchy is at most that of the original one. Note also there is
a one-one correspondence between components of $M_i$ and $M'_i$, and so if a surface $S_i$
is non-separating, then this remains true of $S'_i$.

In this way, we may arrange that each trivial boundary curve of each decomposing surface is parallel towards a u-suture.
This implies that no component of $\partial S_i$ bounds a disc in $\partial M_i$ disjoint from the sutures.
Therefore the hypotheses of Theorem 2.3 hold, and so each decomposing surface is taut and, in particular, incompressible. $\square$

\vskip 6pt
\noindent {\caps 4.3. Slicing under a disc of contact}
\vskip 6pt

Suppose that there is a disc $D$ in $R_\pm(M)$ with $D \cap S = \partial D$, and with the orientation of $D$
matching that of $S$ near $\partial D$. Then $D$ is a {\sl disc of contact}. If we attach $D$ to $S$, and then
push it a little into the interior of $M$ so that the resulting surface $S'$ is properly embedded, this
is known as {\sl slicing under the disc of contact}.

\noindent {\bf Lemma 4.3.} {\sl Let 
$$(M, \gamma) \buildrel S \over \longrightarrow (M', \gamma')$$ be a taut allowable decomposition between decorated
sutured manifolds that extends to an allowable hierarchy. Let $S'$ be obtained from $S$ by slicing under a disc of contact. Then $S'$ also extends to a
allowable hierarchy, with reduced length at most that of the one starting with $S$, and with the same length.}

\noindent {\sl Proof.} Let $D$ be the disc of contact. Then $\partial D$ is a trivial boundary curve of $S$, and so
it gives rise to a u-suture $u$ of $\gamma'$. Therefore, all future decompositions in the given allowable hierarchy avoid
this u-suture. The manifold $(M'', \gamma'')$ obtained by decomposing along $S'$ is obtained from $(M' ,\gamma')$
by removing the u-suture $u$ and attaching a 2-handle. All the decompositions in the hierarchy for $(M', \gamma')$
avoid $u$, and so this sequence of decompositions may also be viewed as a hierarchy for $(M'' ,\gamma'')$.
Boundary curves of surfaces in this new hierarchy are trivial if and only if the corresponding boundary curves
are trivial in the original hierarchy. Hence, there is a one-one correspondence between the u-sutures in the
new and original hierarchy. Therefore, all the decompositions in the new hierarchy are allowable.
Let $(M_{n+1}, \gamma_{n+1})$ be the final manifold in the original hierarchy. The component containing $u$
is a product that is not pre-spherical. The corresponding component in the new hierarchy is therefore
also a product that is not pre-spherical. Hence, this forms an allowable hierarchy. The surfaces in the new
hierarchy have the same topological type and the same intersection with the sutures as in the original 
hierarchy, except the first surface $S'$. Note that $S$ could not have been a union of product discs, annuli disjoint from the
sutures and tori. So, the reduced length of the new hierarchy is at most that of the original one. $\square$

\vskip 6pt
\noindent {\caps 4.4. Decomposition along a product disc}
\vskip 6pt

Decomposition along an allowable product disc is a very common operation. In this section, we show that this
preserves the existence of an allowable hierarchy.

\noindent {\bf Lemma 4.4.} {\sl Suppose that $(M, \gamma)$ and $(M', \gamma')$ are decorated sutured manifolds that
differ by decomposition along an allowable product disc. Then $(M, \gamma)$ admits an allowable hierarchy if and only if $(M', \gamma')$ does.
Moreover, we may arrange that these hierarchies have the same reduced length.}

\noindent {\sl Proof.} Suppose that there is an allowable decomposition
$$(M, \gamma) \buildrel P \over \longrightarrow (M', \gamma')$$
where $P$ is a product disc. If $(M', \gamma')$ admits
an allowable hierarchy, then we may place $P$ at the beginning of this hierarchy,
and obtain an allowable hierarchy for $(M, \gamma)$ with the same reduced length. 

We now need to show that if $(M, \gamma)$ admits an allowable hierarchy, then so does $(M', \gamma')$, with the same reduced length.
We will also show that this hierarchy for $(M', \gamma')$ can be taken to have the same length as the one for $(M, \gamma)$.
We will prove this by induction. Our primary measure of complexity for this induction is the minimal reduced length of an allowable hierarchy for $(M, \gamma)$.
Our secondary measure of complexity is the length of such a hierarchy.

The induction starts as follows. Suppose that $(M, \gamma)$ has an allowable hierarchy with zero reduced length. Note that
$R_+(M)$ and $R_-(M)$ are both non-empty because $(M, \gamma)$ contains the product disc $P$. So, by Lemma 4.1, the
component of  $(M, \gamma)$ containing $P$ is a product. Moreover, no component is pre-spherical, because $E(M, \gamma)$ is taut by Lemma 3.3.
Hence, the new components of $(M', \gamma')$ are also products, and no component is pre-spherical, because each component created
by decomposition along an allowable product disc contains at least one suture that is not a u-suture.
So, $(M', \gamma')$ admits an allowable hierarchy with zero reduced length, and the same length as the hierarchy
for $(M, \gamma)$.

We now consider the inductive step. So, suppose that $(M, \gamma)$ admits an allowable hierarchy
$$(M,\gamma) = (M_1, \gamma_1) \buildrel S_1 \over \longrightarrow \dots \buildrel S_n \over \longrightarrow (M_{n+1}, \gamma_{n+1}).$$
We take this to have minimal reduced length, and subject to this condition, minimal length. By Lemma 4.2, we may assume that if any $S_i$ has a trivial boundary curve, 
then this is parallel in $R_\pm(M_i)$ towards a u-suture.
Let $S_i$ be the first surface to intersect $P$. We may assume that no trivial boundary curve intersects $P$.
This is because the sutures that $P$ runs over are not u-sutures, and so when we isotope
trivial boundary curves of $S_i$ towards u-sutures, we pull these curves away from $P$.

Suppose first that $P \cap S_i$ contains a simple closed curve that bounds a disc in $P$.
Then, we may find such a curve that is innermost in $P$, and so bounds a disc $D$ in $P$ with interior
disjoint from $S_i$. Now $(M_{i+1}, \gamma_{i+1})$ is taut, by Corollary 3.7, and so $\partial D$ bounds a disc $D'$ in $R_\pm(M_{i+1})$.
Now, $D'$ consists of parts of $R_\pm(M_i)$ and parts of $S_i$. Let $C$ be the intersection between these
two subsurfaces of $D'$, which is a collection of simple closed curves, because $\partial D'$ is disjoint from $R_\pm(M_i)$.
We claim that $C$ is empty. If there is such a curve, then an innermost one in
$D'$ bounds a disc, which is either a disc component of $S_i$ or a disc in $R_\pm(M_i)$. In the former case,
this violates our assumption that $S_i$ is part of an allowable hierarchy. In the latter case, the
boundary of this disc is a trivial boundary curve of $S_i$ not parallel towards a u-suture, which is contrary to assumption.
Hence, $D'$ is a disc in $S_i$, and we may isotope this disc onto $D$, and thereby reduce the
number of components of $P \cap S_i$. So, we may assume that $P \cap S_i$ contains no simple closed curves.

Suppose now $P \cap S_i$ contains an arc that has endpoints in the same component of $\partial P - \gamma$. Then there is one that is outermost in $P$,
and that separates off a disc $D$ disjoint from $\gamma$. Suppose that this disc $D$ is disjoint from $\gamma_{i+1}$. Then $\partial D$
bounds a disc $D'$ in $R_\pm(M_{i+1})$. Let $C$ be the intersection between $D' \cap R_\pm(M_i)$ and $D' \cap S_i$.
As argued above, $C$ contains no simple closed curves. Hence, it is a single arc. We may isotope
the disc $D' \cap S_i$ onto $D$, and thereby remove this component of $P \cap S_i$. Suppose now
that $D$ is not disjoint from $\gamma_{i+1}$. This is therefore a product disc in $(M_{i+1}, \gamma_{i+1})$.
We may boundary compress $S_i$ along this disc, giving a new surface $\hat S_i$. 
Let $(\hat M_{i+1}, \hat \gamma_{i+1})$ be the result of decomposing $(M_i, \gamma_i)$ along $\hat S_i$.
Then $(\hat M_{i+1}, \hat \gamma_{i+1})$ and $(M_{i+1}, \gamma_{i+1})$ differ by an allowable decomposition along a product disc.
Now, the part of the hierarchy after $(M_{i+1}, \gamma_{i+1})$ has no greater reduced length and shorter length
than the one for $(M, \gamma)$. So, by induction, $(\hat M_{i+1}, \hat \gamma_{i+1})$ extends to an allowable hierarchy with the same reduced length
and the same length.

In this way, we may assume that no arc of $P \cap S_i$  has endpoints in the same component of $\partial P - \gamma$.
Hence, $P \cap S_i$ consists of a collection of parallel arcs that divide $P$ into a collection of product discs $P'$. 
Let $(M_i', \gamma_i')$ be the result of decomposing $(M', \gamma')$ along $S_1, \dots, S_{i-1}$.
Note that there is a commutative diagram of sutured manifold decompositions
$$

\matrix{
(M_i, \gamma_i) & \buildrel S_i \over \longrightarrow & (M_{i+1}, \gamma_{i+1}) \cr
\Big\downarrow {\scriptstyle P} && \Big\downarrow {\scriptstyle P'} \cr
(M_i', \gamma_i') & \buildrel S_i - {\rm int}(N(P)) \over \longrightarrow & (M_{i+1}', \gamma'_{i+1}). \cr}
$$
The diagram also commutes as allowable decompositions of decorated sutured manifolds, for the following reasons.
Decompositions along $S_i$ and $P$ are allowable by assumption. Neither creates any u-sutures that intersect a
regular neighbourhood of $P$. In the case of $S_i$, this is because the trivial curves of $\partial S_i$ miss $P$. 
Hence, the decompositions along $P'$ and $S_i - {\rm int}(N(P))$ are allowable. Neither decomposition creates
any u-sutures that intersect a regular neighbourhood of $P$, because each boundary curve of $P'$
and $S_i - {\rm int}(N(P))$ that intersects this regular neighbourhood runs over a suture. Away from $P$,
$S_i$ and $S_i - {\rm int}(N(P))$ have the same trivial curves. So, both ways around the commutative diagram
create the same decorated sutured manifolds.

We need to verify that $S_i - {\rm int}(N(P))$ can be part of an allowable hierarchy. Note that no component of
$S_i - {\rm int}(N(P))$ is planar, disjoint from the sutures and where all but one of its boundary components is trivial.
This is because this was true of $S_i$, and when any component of $S_i$ is cut along $P$, the resulting
components of $S_i - {\rm int}(N(P))$ each runs over a suture. 

By induction, $(M_{i+1}', \gamma'_{i+1})$ admits an allowable hierarchy with the same reduced length and the same length as the one
starting with $(M_{i+1}, \gamma_{i+1})$. Hence, $(M', \gamma')$
admits an allowable hierarchy, by starting
with $S_1, \dots, S_{i-1}$, then $S_i - {\rm int}(N(P))$, and then the allowable hierarchy for $(M_{i+1}', \gamma'_{i+1})$. 
This clearly has the same length as the given hierarchy for $(M, \gamma)$.
We claim that this has the same reduced length as this hierarchy for $(M, \gamma)$.
We must check that each component of $S_i$ is a product disc, torus or annulus disjoint from the sutures if and only if the same is
true for each component of $S_i - {\rm int}(N(P))$. Note that each arc of $S_i \cap P$ intersects both $R_-(M_i)$
and $R_+(M_i)$. Hence if it lies in an annulus of $S_i$ disjoint from the sutures, then this is a product
annulus. So, if $S_i$ is a union of product discs, annuli disjoint from the sutures and tori, then the same is true of $S_i - {\rm int}(N(P))$. 
Conversely, no component of $S_i - {\rm int}(N(P))$ that intersects $N(P)$ can be an annulus
disjoint from the sutures or a torus, and if each such component of $S_i - {\rm int}(N(P))$ is a product disc, then these
patch together to form product discs and product annuli of $S_i$. This proves the claim. Hence,
the given hierarchy for $(M', \gamma')$ does indeed have the same reduced length as the original
one for $(M, \gamma)$. $\square$

\vskip 6pt
\noindent {\caps 4.5. Tubing along an arc}
\vskip 6pt

Suppose that $\alpha$ is an arc in $R_\pm(M)$ with $\alpha \cap S = \partial \alpha$. Then there is an embedding
of $\alpha \times [-1,1]$ in $R_\pm(M)$ with $\alpha \times \{ 0 \} = \alpha$ and 
$(\alpha \times [-1,1]) \cap S = \partial \alpha \times [-1,1]$. Suppose that the orientation
that $\alpha \times [-1,1]$ inherits from $R_\pm(M)$ agrees with the orientation of $S$
near $\partial \alpha \times [-1,1]$. Then $\alpha$ is called a {\sl tubing arc}.

We can construct a new surface $S'$ as follows. Attach $\alpha \times [-1,1]$ to $S$, and push the surface
a little into the interior of $M$, so that it becomes properly embedded. Then $S'$ is obtained from $S$ by
{\sl tubing along the arc $\alpha$}. (See Figure 5.)

Note that $S$ is obtained from $S'$ by boundary-compressing along a product disc $D$. This product disc
has boundary consisting of an arc in $\alpha \times [-1,1]$ running from $\alpha \times \{-1\}$ to $\alpha \times \{1\}$,
and an arc in $S'$. 

This type of modification played a key role in [21]. So it will be important to understand
how it behaves in the context of allowable hierarchies.

\noindent {\bf Lemma 4.5.} {\sl Let $S$ be an allowable decomposing surface for a decorated sutured
manifold $(M, \gamma)$ that extends to an allowable hierarchy. Let $S'$ be obtained from $S$ by 
tubing along an arc $\alpha$. Then $S'$ also extends to an allowable hierarchy. Moreover,
the reduced length of the hierarchy beyond $S'$ is equal to the reduced length of the hierarchy beyond $S$.}

\vskip 18pt
\centerline{
\epsfxsize=5in
\epsfbox{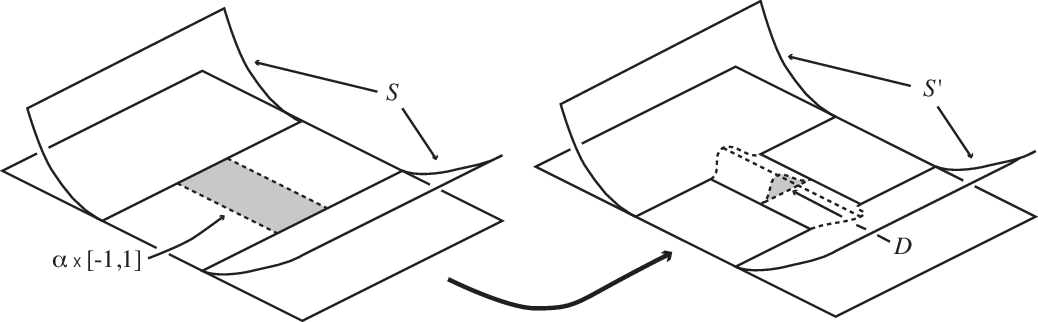}
}
\vskip 6pt
\centerline{Figure 5: Tubing along an arc}

\noindent {\sl Proof.} We first show that decomposition along $S'$ is allowable. Since $S$ was disjoint
from any u-sutures and so was $\alpha$, the same is true of $S'$. Suppose a component $C'$ of $\partial S'$
bounds a trivialising planar surface $P'$, but where the transverse orientations of $P'$ and $S'$ disagree near $C'$.
Since decomposition along $S$ was allowable, this component $C'$ must run along the new tube. Consider 
a component $C$ of $\partial S$ that is incident to an endpoint of $\alpha$. Then $C$ has non-empty intersection with
$P'$. It therefore bounds a planar subsurface $P$ of $P'$. This forms a trivialising planar surface for $C$
and its transverse orientation disagrees with that of $S$ near $C$. This contradicts the assumption that
$S$ was allowable.

We now check that $S'$ can be part of an allowable hierarchy, by verifying that no component of
$S'$ can be a planar surface disjoint from $\gamma$ with all but at most one boundary curve being trivial.
Suppose that $F'$ is such a component. This must contain the new tube, because otherwise it forms
a component of $S$. Thus $F'$ boundary-compresses to form either one or two components
of $S$. Suppose first that we get two components $F_1$ and $F_2$ of $S$. Then the arc
of intersection between the boundary-compression disc and $F'$ must have had endpoints on the
same component of $\partial F'$. 
If $\partial F'$ contains a non-trivial component, it is either this curve or
a boundary component of $F_1$ or $F_2$. In each case, at least one of $F_1$ or $F_2$ has 
at most one trivial boundary curve. This contradicts the fact that $S$ was part of an allowable hierarchy.
Suppose now that we obtain a single component $F$ of $S$ by boundary-compressing the tube.
Then the arc of intersection between the boundary-compression disc and $F'$ had endpoints
on distinct components of $\partial F'$. If they were both trivial, then their trivialising
planar surfaces combine to form a trivialising planar surface for the boundary component
of $F$ at the endpoints of $\alpha$. Hence, all but at most one boundary component of $F$
is trivial, and again this is a contradiction. On the other hand, if one of the components
of $\partial F'$ incident to the boundary-compression disc is non-trivial, then the remaining
components of $\partial F'$ are trivial. These end up forming all but one boundary component
of $F$. Again, this is a contradiction.

Thus, we have shown that $S'$ may be part of an allowable hierarchy. We now show that
such an allowable hierarchy exists.

Let $(M_S, \gamma_S)$ and $(M_{S'}, \gamma_{S'})$ be the result of
decomposing $(M, \gamma)$ along $S$ and $S'$ respectively.
Let $C_1$ and $C_2$ be the curves of $\partial S$ at the endpoints of
$\alpha$. It is possible that $C_1 = C_2$. For each curve $C_i$, let $A_i$ be the product annulus
properly embedded in $(M_{S'}, \gamma_{S'})$ that runs parallel to it. Orient $A_i$ so that its transverse
orientation near $\partial M$ is the same as that of $\partial S$. 

We claim that $A_i$ is trivial in $(M_{S'}, \gamma_{S'})$ if and only if $\alpha$ has endpoints on the same
trivial curve of $\partial S$. Suppose first that $\alpha$ has endpoints on the same trivial curve of $\partial S$.
Then $A_1 = A_2$. Removing a regular neighbourhood of $\alpha$ from the trivialising planar surface gives trivialising
planar surfaces for the two new components of $\partial S'$. These therefore give rise to u-sutures
of $\gamma_{S'}$. The two u-sutures and $A_1 \cap \partial M$ (which equals $A_2  \cap \partial M$) together bound a pair of pants in $\partial M_{S'}$,
which forms a trivialising planar surface, and hence $A_1$ (which equals $A_2)$ is trivial.

Suppose now that some $A_i$ is trivial. Then one of its boundary curves is trivial. It cannot be $A_i \cap S$, because
it would then lie in a component of $S$ that is planar, disjoint from the sutures and has all but at most one boundary curve trivial.
This contradicts the assumption that $S$ is part of an allowable hierarchy. Hence, $A_i \cap \partial M$ must be the trivial curve.
Its trivialising planar surface $P$ must have an orientation that is consistent with that of $A_i$, as otherwise $S$ is not allowable.
Therefore $P$ contains the new components of $\partial S'$ in its boundary. These are therefore u-sutures of $\gamma_{S'}$
and hence trivial boundary curves of $\partial S'$. They must be distinct curves of $\partial S'$, and hence $\alpha$
has endpoints on the same component of $\partial S$, as otherwise the other component of $\partial S$ would bound 
a trivialising planar surface with the wrong transverse orientation. The component of $\partial S$ at the endpoints
of $\alpha$ is trivial, with trivialising planar surface formed from the union of $P$ and the trivialising planar surfaces
for the two components of $\partial S'$. This proves the claim.

The lemma now divides into the cases of the claim. Suppose first that some $A_i$ is trivial. Then 
$\alpha$ has endpoints on the same trivial curve of $\partial S$. This curve becomes a u-suture
in $(M_S, \gamma_S)$, and later surfaces in the hierarchy avoid it. Call this hierarchy
$$(M_S, \gamma_S) = (M_2, \gamma_2) \buildrel S_2 \over \longrightarrow \dots \buildrel S_n \over \longrightarrow (M_{n+1}, \gamma_{n+1}).$$
We have a (non-allowable) decomposition
$$(M_{S'}, \gamma_{S'}) \buildrel D \over \longrightarrow (M_S, \gamma_S).$$
Hence, we can view $M_S$ as lying in $M_{S'}$. The hierarchy for $(M_S, \gamma_S)$ gives a sequence of
decompositions
$$(M_{S'}, \gamma_{S'}) = (M'_2, \gamma'_2) \buildrel S'_2 \over \longrightarrow \dots \buildrel S'_n \over \longrightarrow (M'_{n+1}, \gamma'_{n+1}).$$
A simple induction gives that 
\item{(i)} each $(M'_i, \gamma'_i)$ is decorated;
\item{(ii)} apart from the three sutures that intersect a regular neighbourhood of $\alpha$, the u-sutures of $\gamma_i$ and $\gamma'_i$ are equal;
\item{(iii)} a curve in $\partial M_i - \gamma_i$ is trivial if and only if the corresponding curve in $\partial M_i' - \gamma_i'$ is trivial;
\item{(iv)} decomposition along $S'_i$ is allowable;
\item{(v)} no component of $S'_i$ is planar, disjoint from the sutures and where all but at most one boundary component is trivial.

\noindent The argument is very similar to the proof of Lemma 3.8 and is omitted. There is a (non-allowable) decomposition
$$(M'_{n+1}, \gamma'_{n+1}) \buildrel D \over \longrightarrow (M_{n+1}, \gamma_{n+1}).$$
Hence, $(M'_{n+1}, \gamma'_{n+1})$ is also a product sutured manifold, no component
of which is pre-spherical. Thus, we have found the required allowable hierarchy for $(M_{S'}, \gamma_{S'})$.

The remaining case is where no $A_i$ is trivial. Then we may decompose $(M_{S'}, \gamma_{S'})$ along
these annuli. The resulting decorated sutured manifold is then a copy of $(M_S, \gamma_S)$ plus a product
sutured manifold that is not pre-spherical. Therefore, we may follow this decomposition with the given
hierarchy for $(M_S, \gamma_S)$. $\square$

We close this section with a technical lemma that will be useful later.

\noindent {\bf Lemma 4.6.} {\sl Let 
$$(M,\gamma) \buildrel S \over \longrightarrow (M_S, \gamma_S)$$
be a taut allowable decomposition between decorated sutured manifolds.
Let $S'$ be obtained from $S$ by tubing along an arc $\alpha$.
Let $C$ be a trivial curve in $\partial M_S$ disjoint from $\gamma_S$
and from the discs in $S$ to which the tube is attached.
Then $C$ corresponds to a curve trivial $C'$ in $\partial M_{S'}$ disjoint from $\gamma_{S'}$.}

\noindent {\sl Proof.} Suppose first that one of the curves of $\partial S$ to which the
tube is attached is trivial. It then bounds a trivialising planar surface $P$. This must contain $\alpha$.
The curve of $\partial S$ at the other endpoint of $\alpha$ cannot lie in the interior of $P$, because
it would then bound a trivialising planar surface with the wrong transverse orientation.
We therefore deduce that both endpoints of $\alpha$ lie in the same trivial curve of
$\partial S$. Therefore, this gives rise to two trivial curves of $\partial S'$.
Hence, in this case, $E(M_S, \gamma_S)$ is homeomorphic to $E(M_{S'}, \gamma_{S'})$,
and the homeomorphism takes $C$ to $C'$. Thus, $C$ is trivial if and only if $C'$ is trivial.

Suppose now that neither of the curves of $\partial S$ to which the tube is attached
is trivial. Then the trivialising planar surface $P$ for $C$ is disjoint from sutures of $\gamma_S$ corresponding to these curves.
Hence, it corresponds to a trivialising planar surface $P'$ for $C'$. $\square$

\vskip 6pt
\noindent {\caps 4.6. Boundary compressing along a product disc}
\vskip 6pt

Let $S$ be a surface properly embedded in $(M, \gamma)$. Suppose that there is a disc $D$ embedded in $M$,
such that $D \cap S$ is a single arc in $\partial D$, and $D \cap \partial M$ is a single arc in $\partial D$ and where
these two arcs intersect at their endpoints. Suppose $D \cap \gamma$ is empty and the
orientations of $S$ and $\partial M$ disagree at $\partial D$. We will consider
the surface $S'$ obtained from $S$ by boundary compressing along $D$. Thus,
$S$ is obtained from $S'$ by tubing along an arc.
Let $(M_S, \gamma_S)$ and $(M_{S'}, \gamma_{S'})$ be the manifolds obtained by decomposing $(M, \gamma)$
along $S$ and $S'$ respectively.
Note that $D$ is a product disc in $(M_S, \gamma_S)$, and that decomposing $(M_S, \gamma_S)$ along $D$
gives a sutured manifold homeomorphic to $(M_{S'}, \gamma_{S'})$. 

\noindent {\bf Proposition 4.7.} {\sl Let $(M, \gamma)$ be a decorated sutured manifold, and let $S$ be a decomposing surface that extends to an allowable
hierarchy. Let $S'$ be obtained from $S$ by boundary compressing along a product disc disjoint from $\gamma$.
Let $S''$ be obtained from $S'$ by removing
any component that is a planar surface disjoint from $\gamma$ with all but at most one boundary curve that is trivial. Then $S''$ also extends
to an allowable hierarchy. Moreover, the reduced length of this hierarchy beyond $S''$ is equal to the reduced
length beyond $S$.}

\noindent {\sl Proof.} Let $(M_S, \gamma_S)$ and $(M_{S'}, \gamma_{S'})$ be the result of decomposing $(M, \gamma)$ along $S$ and $S'$
respectively. There is a decomposition
$$(M_S, \gamma_S) \buildrel D \over \longrightarrow (M_{S'}, \gamma_{S'})$$
where $D$ is the product disc. Thus, one might hope to use Lemma 4.4 to prove the proposition.
However, there are a number of potential complications. Firstly, $S'$ might not be an allowable decomposing surface,
because it might have a trivial boundary curve, where the trivialising planar surface and $S'$
are incompatibly oriented. In this case, we have not even specified which sutures of $\gamma_{S'}$
are to be viewed as u-sutures. Secondly, $S'$ might have a planar component disjoint from 
$\gamma$ and with all but at most one boundary component trivial, and hence it might not be possible 
for $S'$ to be part of an allowable hierarchy. Thirdly, the product disc $D$ may run over a u-suture
of $\gamma_S$, and so decomposition along $D$ might not be allowable, and therefore
Lemma 4.4 might not apply.

\noindent {\sl Case 1.} $D$ does not run over a u-suture of $\gamma_S$.

In this case, Lemma 4.4 is applicable. Therefore, we declare that the sutures of $\gamma_{S'}$
incident to $D$ are not u-sutures. Then by Lemma 4.4, $(M_{S'}, \gamma_{S'})$ with this
decoration admits an allowable hierarchy with the same reduced length as the one starting from
$(M_S, \gamma_S)$. However, this does not complete the proof of the lemma, because the
first or second problems mentioned above may still hold.

Suppose that $S'$ has a trivial boundary curve $C'$, where the trivialising planar surface $P$ has transverse orientation
differing from that of $S'$ near $C'$. The planar surface $P$ may contain other boundary curves of $S'$. These
are all trivial, and by choosing $C'$ appropriately, we may assume that the trivialising subsurfaces for the curves
$S' \cap {\rm int}(P)$ are all correctly oriented. Let $P'$ be the subsurface of $P$ obtained by cutting along the
curves $S' \cap {\rm int}(P)$ and then taking the component containing $C'$. The curve $C'$ becomes a suture of $\gamma_{S'}$.
The surface $P'$ becomes a planar component of $R_\pm(M_{S'})$, and all its boundary curves apart from $C'$
are u-sutures. Hence, it extends to a disc in $E(M_{S'}, \gamma_{S'})$. Now, $E(M_{S'}, \gamma_{S'})$ is 
taut, by Lemma 3.3, and so the component of $E(M_{S'}, \gamma_{S'})$ is a taut 3-ball. Therefore,
the component of $S'$ containing $C'$ is planar, disjoint from the sutures and all but at most one of its boundary curves are
trivial. We therefore focus on this case.

Let $S'_1$ be a planar component of $S'$, disjoint from $\gamma$ and with all but at most one boundary component trivial.
Then $E(S'_1)$ is a disc disjoint from the sutures of $E(M, \gamma)$. By Lemma 3.3, this is taut. Hence,
$E(S'_1)$ is boundary parallel in $E(M, \gamma)$. It therefore separates off a 3-ball component of $E(M_{S'}, \gamma_{S'})$.
By Lemma 3.9, the corresponding component of $(M_{S'}, \gamma_{S'})$ is a product sutured manifold.
Therefore, $(M_{S'}, \gamma_{S'})$ is the disjoint union of this product sutured manifold
and the manifold obtained by decomposing $(M, \gamma)$ along $S - S'_1$. So, $S - S'_1$ extends to an allowable hierarchy,
with the same reduced length. Thus, the proposition is proved in this case.

\noindent {\sl Case 2.} $D$ runs over a u-suture of $\gamma_S$.

\noindent {\sl Case 2A.} $D$ runs over two distinct sutures of $\gamma_S$. 

Call these $C_1$ and $C_2$, where $C_1$ is trivial.
It is also possible that $C_2$ is trivial. Let $C'$ be the curve of $\partial S'$ incident to $D$. Then $C'$ is
trivial if and only if $C_2$ is trivial. Furthermore, the trivialising planar surfaces lie on the same side.
Thus, $C'$ cannot bound a trivialising planar surface with the wrong orientation, because this
would imply that $C_2$ did also, contradicting the assumption that $S$ was allowable.
Hence, decomposition along $S'$ is allowable. Moreover, boundary-compressing $S$ along $D$
does not disconnect this component of $S$. Hence, if the resulting component of $S'$ is a planar
surface disjoint from the sutures and with all but at most one boundary component trivial, then
the same is true of the component of $S$, which is a contradiction. So, $S'$ can be part of an
allowable hierarchy.

Now $S$ is obtained from $S'$ by tubing along an arc that starts and ends in $C'$. We now create a new surface $S''$ that is
obtained from $S'$ by tubing along a different arc. This has the same endpoints but runs along a regular neighbourhood of $C_1$.
Decomposing along $S''$ gives a decorated sutured manifold homeomorphic to $(M_S, \gamma_S)$. Hence, $S''$
extends to an allowable hierarchy with the same reduced length as the one starting with $S$. But $S''$ has a disc of
contact which lies within a regular neighbourhood of $C'_1$. Slicing under this disc of contact gives
a surface isotopic (relative to $\gamma$) to $S'$. By Lemma 4.3, $S'$ extends to an allowable hierarchy with the same reduced length 
as the one starting with $S$. 

\noindent {\sl Case 2B.} $D$ runs over the same u-suture of $\gamma_S$ twice. 

This comes from a boundary component
$C$ of $\partial S$. In this case, the boundary compression gives
rise to two distinct boundary components $C'_1$ and $C'_2$ of $S'$. Suppose first
that neither $C'_1$ nor $C'_2$ is trivial.

Let us change the decoration of $(M_S,\gamma_S)$ by declaring that the suture of $\gamma_S$ incident to $D$
is {\sl not} a u-suture.  Then, by Lemma 3.8, $(M_S, \gamma_S)$ with this new decoration still admits an allowable hierarchy.
Note that the decomposition 
$$(M_S, \gamma_S) \buildrel D \over \longrightarrow (M_{S'}, \gamma_{S'})$$
is now allowable, because $D$ now does not run over any u-sutures. Moreover, $(M_{S'}, \gamma_{S'})$ has
the correct decoration because the curves $C'_1$ and $C'_2$ are non-trivial. Hence,
Lemma 4.4 establishes the proposition in this case.

Now consider the case where $C'_1$ or $C'_2$ is trivial. In fact, if one of them is
trivial, then so is the other. Moreover, $C'_1$ (say) bounds a trivialising planar surface $P'_1$
with the incorrect orientation, whereas $C'_2$ bounds a trivialising planar surface $P'_2$ with
the correct orientation. We again change the decoration of 
$(M_S,\gamma_S)$ by declaring that the sutures of $\gamma_S$ incident to $D$
are not u-sutures. Then $(M_{S'}, \gamma_{S'})$ inherits a decoration where neither of
the sutures coming from $C'_1$ or $C'_2$ is a u-suture. Again, Lemma 4.4 implies that $(M_{S'}, \gamma_{S'})$
admits an allowable hierarchy. Therefore, by Lemma 3.3, $E(M_{S'}, \gamma_{S'})$ is taut.
Now, $P'_1$ extends to a disc in the boundary of $E(M_{S'}, \gamma_{S'})$,
and so by the tautness of $E(M_{S'}, \gamma_{S'})$, the component of $E(M_{S'}, \gamma_{S'})$
is a product sutured 3-ball. Hence, as argued in Case 1, the component $S'_1$ of $S'$ incident to $P_1'$ is a planar surface, disjoint
from $\gamma$ and that is boundary-parallel via a product region that is a product sutured manifold.
Therefore, we may remove this component $S'_1$ from $S'$. Let $S''$ be the union of the remaining components.
Let $(M_{S''}, \gamma_{S''})$ be obtained by decomposing
$(M, \gamma)$ along $S''$. Note that this decomposition is allowable. Note also 
that this sutured manifold differs from $(M_S, \gamma_S)$ by attaching the sutured manifold
$S'_1 \times I$ along a product disc. Hence, the allowable sutured manifold hierarchy for $(M_S, \gamma_S)$
gives an allowable sutured manifold hierarchy for $(M_{S''}, \gamma_{S''})$ with the
same reduced length. This is the required hierarchy. $\square$

\vskip 6pt
\noindent {\caps 4.7. Annular swaps and decompositions along annuli}
\vskip 6pt

Let $S$ be a surface properly embedded in $M$. Suppose that there is an annulus or M\"obius band $A$ 
embedded in the interior of $M$ such that $A \cap S = \partial A$. In the case where $A$ is an annulus,
suppose that the orientations on $S$ near $S \cap A$ either both point towards $A$ or both
point away from $A$. Let $N(A)$ be the orientable $I$-bundle over $A$, where $A$ is viewed as 
a section of $N(A)$. Then $N(A)$ embeds in $M$ in such a way that $N(A) \cap S$ is the $I$-bundle
over $\partial A$. Let $S'$ be the surface that results from removing $N(A) \cap S$ from $S$ and 
attaching ${\rm cl}(\partial N(A) - S)$. Then our assumption about the orientation on $S$ near $S \cap A$ implies
that $S'$ can be oriented in such a way that the orientations on $S$ and $S'$ agree on their intersection.
We say that $S'$ is obtained from $S$ by {\sl an annular swap} along $A$.

\vskip 12pt
\centerline{
\epsfxsize=4in
\epsfbox{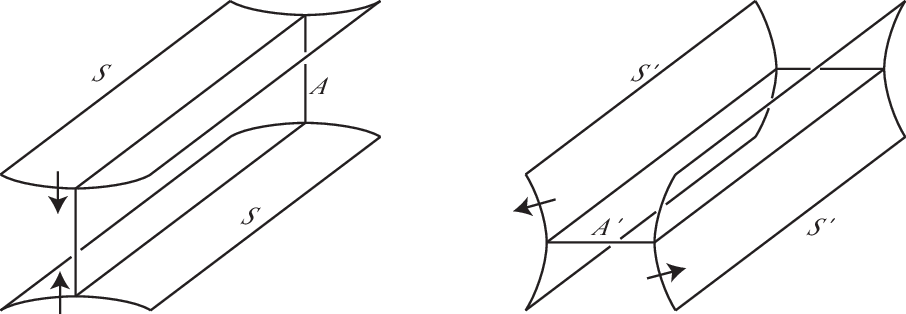}
}
\vskip 6pt
\centerline{Figure 6: An annular swap}

Note that there is an annulus or M\"obius band $A'$ embedded in $N(A)$, which is the $I$-bundle over a core curve
of $A$. So $A' \cap S'$ equals $\partial A'$. Then $S$ is obtained from $S'$ by an annular
swap along $A'$. (See Figure 6.) This modification to a decomposing surface will be important in this paper. In particular, it arises in the context of
`irregular switches' in normal surface theory. More details can be found in Sections 7.5 and 7.3.

The following result implies that this modification may be made, while preserving the existence of an allowable hierarchy,
under some fairly mild assumptions.

\noindent {\bf Proposition 4.8.} {\sl Let 
$$(M, \gamma) \buildrel S \over \longrightarrow (M_S, \gamma_S)$$
be a taut allowable decomposition between decorated sutured manifolds that extends to an allowable hierarchy.
Suppose that $E(M,\gamma)$ is atoroidal and
that no component of $E(M,\gamma)$ is a Seifert fibre space other than a solid torus or a copy of $T^2 \times I$.
Suppose also no component of $M_S$ has boundary a single torus with no sutures.
Assume that if $E(S)$ contains a component that is a torus or annulus disjoint from $\gamma$ and from the
attached 2-handles of $E(M, \gamma)$, then no other component of $E(S)$ is of this form.
Suppose that no component of $A \cap S$ bounds a disc in the canonical extension $E(S)$.
Let $S'$ be a surface obtained from $S$ by a swap along an annulus or M\"obius band $A$ in $M$.
Orient $S'$ so that the orientations of $S'$ and $S$ agree on their intersection.
Then at least one of the following holds:
\item{(i)} $S'$ extends to an allowable hierarchy, or
\item{(ii)} $A$ is an annulus and the two copies of $A$ in $S'$ lie in distinct components of $S'$,
one of these components $S'_1$ separates off a solid torus with no sutures,
and $S' - S'_1$ extends to an allowable hierarchy.

\noindent In both cases, the new hierarchy has the same reduced length as the original one.
Moreover, if some component of $S$ was not an annulus or torus disjoint from $\gamma$, then this remains true of at least one
component of the first surface of the new hierarchy.
}

At the same time, we will prove the following.

\vfill\eject
\noindent {\bf Proposition 4.9.} {\sl Let $(M, \gamma)$ be a taut decorated sutured manifold
that admits an allowable hierarchy. Suppose that $E(M, \gamma)$ is atoroidal and 
that no component of $E(M,\gamma)$ is a Seifert fibre space other than a solid torus or a copy of $T^2 \times I$.
Let $A$ be an annulus properly embedded in $M$ with boundary disjoint from $\gamma$. Suppose that
neither curve of $\partial A$ is trivial. Suppose that decomposition of $(M, \gamma)$ along $A$ does not
create a component that is a solid torus disjoint from the sutures. Then $A$ extends
to an allowable hierarchy with the same reduced length as the original one.}

\noindent {\sl Proof.}
We will prove these propositions simultaneously, by induction on the reduced length of the hierarchy.

We first prove that Proposition 4.9 in the case where the allowable hierarchy for $(M, \gamma)$ has
some reduced length $n$ implies Proposition 4.8 in the case where $(M_S, \gamma_S)$ has an allowable hierarchy
with reduced length $n$. 

So consider the situation of Proposition 4.9. Let $2A$ denote ${\rm cl}(\partial N(A) - S)$,
the $(\partial I)$-bundle over $A$. View $2A$ as one or two annuli embedded in $(M_S, \gamma_S)$.
Orient $2A$ so that the region it bounds containing $A$ inherits all  the new sutures. 
We start to verify that the hypotheses of Proposition 4.9 apply to the
components of $2A$ in $(M_S, \gamma_S)$. Note that no boundary curve of $2A$ is trivial in $(M_S, \gamma_S)$.
For if a boundary curve of $2A$ did bound a disc in the boundary of $E(M_S, \gamma_S)$ disjoint from the sutures, then this
would imply that either a component of $A \cap S$ bounded a disc in $E(S)$, contrary to assumption,
or some component of $E(S)$ is a disc disjoint from $\gamma$, contradicting the fact that $S$ is part
of an allowable hierarchy, or some boundary curve of $E(S)$ is trivial in $E(M,\gamma)$, contradicting Lemma 3.2. 
The same argument gives that when $2A$ is disconnected and 
$(M_S, \gamma_S)$ is decomposed along one component of $2A$,
then the other component of $2A$ remains non-trivial in the resulting manifold.
Note also that $E(M_S, \gamma_S)$ is atoroidal, by Lemma 3.11.  Also, no component of $E(M_S, \gamma_S)$
is a Seifert fibre space, other than a solid torus or a copy of $T^2 \times I$.
Thus, we may assume that all of the hypotheses of Proposition 4.9 apply to the
components of $2A$ in $(M_S, \gamma_S)$, with the possible exception that 
this decomposition may create some solid torus components with no sutures.
The proof divides into two cases, according to whether this hypothesis also holds.

\noindent {\sl Case 1.} Decomposing $(M_S, \gamma_S)$ along $2A$ does not create any solid torus components with no sutures.

Then by Proposition 4.9 applied once or twice, decomposing $(M_S,  \gamma_S)$ along $2A$ extends to an allowable hierarchy
with reduced length $n$. But decomposing $(M, \gamma)$ along $S$ then $2A$ gives the same sutured manifold as decomposing
$(M, \gamma)$ along $S'$ then $2A'$, where $2A' = S - {\rm int}(N(S'))$. Both the decomposition along $S'$
and the decomposition along $2A'$ are allowable, and the resulting sutured manifold has the same
u-sutures as the manifold obtained by decomposing $(M_S, \gamma_S)$ along $2A$. This is because the
boundary curves of $2A'$ are non-trivial, for the following reason. Suppose that a boundary curve of
$2A'$ had a trivialising planar surface $P$. The intersection $P \cap \partial S$ would give rise to
a (possibly empty) collection of u-sutures in $(M_S, \gamma_S)$. The component of $P - {\rm int}(N(\partial S))$
incident to $\partial A$ would then form a trivialising planar surface for this component of $\partial A$,
which is contrary to hypothesis. Also, no component of $E(S')$ is a disc, because this would imply
that $E(S)$ had a disc component disjoint from the sutures or some curve of $\partial A$ bounded a disc in $E(S)$,
neither of which is possible. Hence, $S'$ extends to an allowable hierarchy the same reduced
length as the one starting with $S$, as required.

\noindent {\sl Case 2.} Decomposing $(M_S, \gamma_S)$ along $2A$ does create at least one solid torus component with no sutures.

\noindent {\sl Case 2A.} The same solid torus with no sutures lies on both sides of $2A$.

Let $V$ be the component of $(M_S, \gamma_S)$ containing this solid torus. Then $\partial V$ consists
of tori with no sutures. It cannot be a single torus, by assumption. So it contains two tori.
Consider the intersection between $E(S)$ and the boundary of $E(V) = V$.
Since $E(S)$ has no trivial boundary curves, this intersection consists of some tori and some essential annuli.
Since $V$ has two boundary tori and each has non-empty intersection with $E(S)$, we deduce
that $E(S)$ contains at least two components that are tori or annuli disjoint from
$\gamma$ and from the attached 2-handles. This contradicts one of our assumptions.

\noindent {\sl Case 2B.} On both sides of $2A$, there are distinct solid tori disjoint from the sutures.

These solid tori then patch together to give a component of $M_S$. This has boundary a single torus,
that is disjoint from $\gamma_S$, contradicting one of our assumptions.

\noindent {\sl Case 2C.} On exactly one side of $2A$, there is a solid torus disjoint from the sutures.

Let $Y$ be this solid torus. If $S \cap Y$ contains any trivial curves, then one of these bounds a disc of contact in $M$.
We may slice under this, and maintain the existence of an allowable hierarchy, with the same reduced length, by Lemma 4.3.
In this process, the solid torus $Y$ with no sutures is preserved. So, in this way, we ensure that $S \cap Y$ is essential in $\partial Y$.
We note that the $(\partial I)$-bundle over $A$ must consist of two copies of $A$ in $S'$ lying in distinct components
of $S'$, since one lies in $Y$ and the other does not. Let $S'_1$ denote the components $Y \cap S'$.
We note that these are tori and annuli disjoint from the sutures, because they are essential subsurfaces of $\partial Y$.
By our assumption about annular and toral components of $E(S)$, $S'_1$ must be connected.
We will show that decomposing along $S' - S'_1$ extends an allowable hierarchy with the
same reduced length as the one starting with $S$. Orient $A$ so that when $(M_S, \gamma_S)$ is decomposed
along it, $Y$ inherits some sutures. Then by Proposition 4.9, this decomposition extends
to an allowable hierarchy with the same reduced length as the one for $(M_S, \gamma_S)$. But
decomposing $(M, \gamma)$ along $S$ then $A$ gives the same manifold as decomposing along
$S' - S'_1$ then along the annuli $S'_1 \cap S$. Hence, $S' - S'_1$ does extend to the required
allowable hierarchy. This establishes that Proposition 4.8 in the case where $(M_S, \gamma_S)$ has an allowable
hierarchy with reduced length $n$ follows from Proposition 4.9 in the case where 
$(M, \gamma)$ has an allowable hierarchy with reduced length $n$.

So, we will focus on Proposition 4.9. Therefore consider a taut decorated sutured manifold $(M, \gamma)$, 
where $E(M, \gamma)$ is atoroidal and where no component of $E(M, \gamma)$ is Seifert fibred other than a solid torus or a copy of $T^2 \times I$. 
By hypothesis, there is an allowable hierarchy
$$(M, \gamma)  = (M_1, \gamma_1) \buildrel S_1 \over \longrightarrow \dots \buildrel S_n \over \longrightarrow (M_{n+1}, \gamma_{n+1}).$$

We will make several hypotheses about this hierarchy, the first of which is:
\item{(1)} The reduced length of this allowable hierarchy is as short as possible. 

By Lemma 4.2, we may assume our second hypothesis:
\item{(2)} If any $S_j$ has a trivial boundary curve, then this is parallel in $R_\pm(M_j)$ towards a u-suture.

We note that, in general, we cannot assume that the surfaces in the hierarchy are connected, because
to break a surface down into its components may increase the reduced length of the hierarchy. However, decomposing
along tori and annuli disjoint from the sutures separately is possible:
\item{(3)} If any surface $S_j$ contains a component that is a torus or an annulus disjoint from the sutures, then
this is all of $S_j$.

Our final hypothesis is:
\item{(4)} Assuming (1) - (3), the length of the hierarchy is as short as possible.

We will prove the proposition by induction. We will consider only allowable hierarchies satisfying (1) - (4).
Our primary measure of complexity for such a hierarchy will be its reduced length. Our secondary
measure of complexity will be the length of the hierarchy.

The induction starts with the case where the hierarchy has reduced length zero. Let 
$$(M, \gamma) \buildrel A \over \longrightarrow (M', \gamma')$$
be the given decomposition.  Then by Lemma 4.1, the component of
$(M, \gamma)$ containing $A$ is a product sutured manifold, a copy of $T^2 \times I$ with no sutures or an orientable $I$-bundle over the
Klein bottle  with no sutures. Note that the other possibilities are ruled out because $A$ has non-empty boundary.
Since $A$ is non-trivial, it is incompressible in $M$. So, in all cases, $A$ is either
a union of $I$-fibres or is boundary parallel. In both cases, $(M', \gamma')$ has an allowable hierarchy with
reduced length zero, as required.

Before we embark on the main flow of the argument, we deal with the case where the component of $M$ containing $A$ is Seifert fibred.
We first show that this component of $M$ has no u-sutures. For if it has a u-suture, then the annular components of $R_\pm(M)$ 
that are adjacent to it extend to discs or spheres in the boundary of $E(M, \gamma)$. But $E(M, \gamma)$ is taut,
by Lemma 3.3, and so this component of $E(M, \gamma)$ is a 3-ball with a single suture. We deduce that this component of
$M$ had a single torus boundary component, containing two sutures, exactly one of which is a u-suture. But in this case, each component of
$\partial A$ is parallel to the u-suture and hence is trivial, which is contrary to assumption. So, this component of $M$ does indeed
have no u-sutures, and so this component of $(M, \gamma)$ is its own canonical extension. By our assumption
about Seifert fibred components of $E(M, \gamma)$, we deduce that it is a solid torus or a copy of $T^2 \times I$.
It is easy to check that in this case, the manifold obtained by decomposing this component
of $(M, \gamma)$ along $A$ has an admissible hierarchy with reduced length at most one.
Moreover, it admits an admissible hierarchy with reduced length zero if and only if the same if true of $(M, \gamma)$.
Thus, Proposition 4.9 is proved in this case.

We now consider the general inductive step. We will consider various possibilities for $S_1 \cap A$.

\noindent {\sl Case 1.} $S_1$ is disjoint from $A$.

Then we have a commutative diagram
$$

\matrix{
(M_1, \gamma_1) & \buildrel S_1 \over \longrightarrow & (M_2, \gamma_2) \cr
\Big\downarrow {\scriptstyle A} && \Big\downarrow {\scriptstyle A} \cr
(M_1', \gamma_1') & \buildrel S_1 \over \longrightarrow & (M_2', \gamma'_2). \cr}
$$
Note that $(M_2, \gamma_2)$ admits an allowable hierarchy satisfying (1)-(4),
with smaller complexity than the given hierarchy for $(M_1, \gamma_1)$.
Hence, we may apply induction to the decomposition of $(M_2, \gamma_2)$
along $A$. However, we need to check that the hypotheses of the proposition hold.
Note that $E(M_2, \gamma_2)$ is atoroidal by Lemma 3.11, and no component of $E(M_2, \gamma_2)$ is Seifert fibred,
unless it is a solid torus or a copy of $T^2 \times I$. Note also that $A$ 
has non-trivial boundary in $M_2$. If $(M_2', \gamma'_2)$ has no solid toral components
disjoint from the sutures, then inductively, $(M_2', \gamma'_2)$ admits an allowable
hierarchy with the same reduced length as the one for $(M_2, \gamma_2)$ starting with $S_2$.
Hence, in this case, the proposition is proved.

So, suppose that there is a solid torus component $Y$ of $M_2'$ disjoint from $\gamma_2'$.
Then $Y$ has non-empty intersection with $S_1$ because otherwise decomposing $(M_1, \gamma_1)$ along $A$
would also create $Y$, contrary to a hypothesis of the proposition.
So, $S_1$ is an annulus disjoint from the sutures, by hypotheses (2) and (3). 
Give $Y$ a Seifert fibration so that $A$ and $S_1$ are a union of fibres.
This must have an exceptional fibre. For otherwise, we can isotope $S_1$
across this solid torus. It is then easy to see that decomposing $(M_2, \gamma_2)$ along $A$ does not
then create a solid torus with no sutures. For if it did, this solid torus would have to meet
both $A$ and $S_1$, and we could then deduce that $M$ is Seifert fibred, 
and we have dealt with this case already. Hence, in this case, the proposition is
proved by induction, as above. So, $Y$ has an exceptional fibre.

Note that $S_1$ is not parallel to an annulus in $(M_1, \gamma_1)$ disjoint from the sutures.
For if it were, then decomposition along $S_1$ just creates a copy of $(M_1, \gamma_1)$ and a product
sutured manifold. In this case, we could remove $S_1$ from the start of the hierarchy,
contradicting our assumption (4).  Similarly, we may assume that $A$ is not parallel to an annulus in $(M_1, \gamma_1)$ disjoint from the sutures,
because in this case, Proposition 4.9 is trivial.

Let $\overline{S_1}$ be a parallel copy of $S_1$ but with opposite orientation, and arranged
so that it does not lie in $Y$. Define $\overline{A}$ similarly. Then decomposing 
$(M_2, \gamma_2)$ along $\overline{S_1} \cup \overline{A}$ does not create any solid
tori disjoint from the sutures, because this would imply either that $S_1$ or $A$ is parallel to an annulus
in $R_\pm(M_1)$, or that $E(M, \gamma)$ is toroidal, or that $E(M, \gamma)$ has a Seifert fibred component
other than a solid torus or a copy of $T^2 \times I$.
No boundary curve of $\overline{S_1} \cup \overline{A}$ is trivial in $(M_2,  \gamma_2)$.
So, inductively, decomposing $(M_2, \gamma_2)$ along
$\overline{S_1} \cup \overline{A}$ extends to an allowable hierarchy with the required reduced
length. In other words, decomposing $(M_1, \gamma_1)$ along $\overline{S_1} \cup S_1 \cup \overline{A}$
extends to such an allowable hierarchy. Let $(M_2'', \gamma_2'')$ be the result of making this
decomposition. This is homeomorphic to the manifold $(M_2''', \gamma_2''')$ 
obtained from $(M_1, \gamma_1)$
by decomposing along $\overline{S_1} \cup A \cup \overline{A}$. So, $(M_2''', \gamma_2''')$ also admits an allowable
hierarchy with the required reduced length. But then we see that $A$ does indeed extend
to the hierarchy that we are looking for.

Thus, we may assume that $S_1 \cap A$ is non-empty.
Exactly as argued in the proof of Lemma 4.4, we may ensure that $A \cap S_1$ contains no simple closed curves
or arcs that are inessential in $A$.

\noindent {\sl Case 2.} $A \cap S_1$ is a collection of essential simple closed curves in $A$, and $S_1$ is neither an annulus
disjoint from $\gamma_1$ nor a torus.

Note first that no component of $E(S_1)$ is a torus or annulus disjoint from $\gamma_1$ and from the attached
2-handles. This is because, in this situation, $S_1$ would either contain a trivial curve bounding a disc of contact,
contradicting hypothesis (2) or would itself be a torus or annulus disjoint from $\gamma_1$ and from the attached
2-handles, which we are assuming is not the case here.

We give the curves $A \cap S_1$ the transverse orientation coming from $S_1$. We also transversely
orient $\partial A$ by using the transverse orientation of $R_\pm(M_1)$.

Suppose that two curves of $A \cap S_1$ are adjacent in $A$ and incoherently oriented. They then bound an annulus $A_1$
such that $A_1 \cap S_1 = \partial A_1$. Let $S'_1$ be the result of performing an annular swap along $A_1$.
We now check that the conditions of Proposition 4.8 hold in this case. Note that $E(M, \gamma)$ is atoroidal,
and that no component of $E(M,\gamma)$ is a Seifert fibre space other than a solid torus or a copy of $T^2 \times I$.
Suppose that a component of $M_2$ has boundary a single torus with no sutures. 
At least one component of $S_1$ would be an essential subsurface of $\partial M_2$, and therefore 
$S_1$ would have a component that is a disc, annulus or torus. The former case is impossible since $S_1$ intersects $A$ in essential curves.
So $S_1$ contains a component that is an annulus or torus disjoint from the sutures. But we assuming that
this is not the case here. We verified above that no component of $E(S_1)$ is a torus or annulus disjoint
from the sutures and the attached 2-handles. Neither curve of $A_1 \cap S_1$ bounds a disc in $E(S_1)$. This is because we could extend this
to form a compression disc for a curve of $\partial A$. Since $E(M_1, \gamma_1)$ is taut, this curve
would bound a disc in $\partial E(M_1)$ disjoint from the sutures, and so $A$ would be trivial, which is contrary to hypothesis. 
Note that, because $S_1$ is neither an annulus disjoint from $\gamma_1$ nor a torus, the
reduced length of the hierarchy beyond $(M_2, \gamma_2)$ is strictly less than the one for $(M_1, \gamma_1)$.
Therefore, by induction, we deduce that, possibly after removing one of its components, $S'_1$ extends to an allowable
hierarchy with the same reduced length as the one starting with $S_1$. Note that $S'_1$ has fewer curves
of intersection with $A$. 

Suppose now all the curves of $A \cap S_1$ are coherently oriented, but that this orientation
disagrees with at least one component of $\partial A$. Then we may find 
such a curve that is outermost in $A$, which therefore separates off an annulus $A'$,
and so that the transverse orientations on $\partial A'$ either both point towards
$A'$ or both point away from it. We now add an annular component $A_+$ to $S_1$. This annulus $A_+$
is boundary parallel in $M_1$, and runs along a regular neighbourhood of $\partial A' \cap \partial M_1$. 
We orient $A_+$ so that the product region between it and $\partial M_1$ inherits two sutures.
Clearly, $S_1 \cup A_+$
extends to an allowably hierarchy with the same reduced length as the original one. This is because by cutting along $S_1 \cup A_+$ rather than $S_1$,
the resulting sutured 3-manifold inherits a new component, which is a taut product solid torus.
By adding $A_+$, we have introduced a new curve of intersection with $A$. But we may then perform
an annular swap to remove this curve and also the other curve of $\partial A'$. Let $S_1'$ be the
result of performing this annular swap. It is easy to check that the conditions of Proposition 4.8 are satisfied.
Note in particular that $E(S_1)$ contains no components that are annuli or tori disjoint from the sutures
and the attached 2-handles, and so only one component of $E(S_1 \cup A_+)$ has this form.
Hence, possibly after removing a component, $S'_1$ extends
to an allowable hierarchy with the same reduced length. 

In this way, we may replace the first surface $S_1$ of the allowable hierarchy with a surface $S'_1$,
such that $A \cap S_1'$ is a (possibly empty) collection of coherently oriented essential curves,
and so that these transverse orientations agree with both components of $\partial A$.
In doing so, we have not increased the reduced length of the hierarchy. We note that at least one component of
$S'_1$ is not an annulus disjoint from $\gamma_1$, a torus or product disc. Let $S''_1$ be the result of removing from $S'_1$
the components are tori and annuli disjoint from the sutures. Then $S''_1$ also extends to an allowable hierarchy with the
same reduced length. The surface $S''_1$ satisfies (2) and (3). 
We may further ensure that the later surfaces also satisfy (2) and (3), by Lemma 4.2.
The reduced length of this new hierarchy is at most that of the original one. We do not make any claim about the
length of the new hierarchy, however.  Note that $A' = A - {\rm int}(N(S_1''))$
is a collection of product annuli. We have a commutative diagram
$$

\matrix{
(M_1, \gamma_1) & \buildrel S_1'' \over \longrightarrow & (M_2'', \gamma_2'') \cr
\Big\downarrow {\scriptstyle A} && \Big\downarrow {\scriptstyle A'} \cr
(M_1', \gamma_1') & \buildrel S''_1 - {\rm int}(N(A)) \over \longrightarrow & (M_2', \gamma'_2). \cr}
$$
The decomposition along $A'$ does not create any solid tori with no sutures, because $A'$
consists of product annuli. No boundary curve of $A'$ is trivial. Note also that $E(M''_2,\gamma''_2)$ is
atoroidal by Lemma 3.11, and no component is a Seifert fibre space, other than a solid torus or a copy of $T^2 \times I$. 
The reduced length of the hierarchy for $(M''_2, \gamma''_2)$ is strictly less than that for $(M_1, \gamma_1)$.
Hence, inductively, $(M_2', \gamma'_2)$ extends to an allowable
hierarchy with the same reduced length as that for $(M_2'', \gamma_2'')$. So,
$(M_1' , \gamma_1')$ has the required allowable hierarchy. It has the
same reduced length as the original one.

\noindent {\sl Case 3.} $A \cap S_1$ is a collection of essential simple closed curves in $A$, and $S_1$ is a torus.

By Lemma 3.13, we may assume that $A \cap S_1$ is a single curve. Let $A' = A - {\rm int}(N(S_1))$. 
We have a commutative diagram of allowable decompositions
$$

\matrix{
(M_1, \gamma_1) & \buildrel S_1 \over \longrightarrow & (M_2, \gamma_2) \cr
\Big\downarrow {\scriptstyle A} && \Big\downarrow {\scriptstyle A'} \cr
(M_1', \gamma_1') & \buildrel S_1 - {\rm int}(N(A)) \over \longrightarrow & (M_2', \gamma'_2). \cr}
$$
It is easy to check that the hypotheses of Proposition 4.9 apply to $A'$. In particular, decomposition along $A'$
does not create a solid torus component disjoint from the sutures. So, inductively, 
$(M_2', \gamma'_2)$ extends to an allowable hierarchy with the same reduced length as the
one starting with $(M_2, \gamma_2)$. Therefore, $(M'_1, \gamma'_1)$ admits an allowable hierarchy
with the required reduced length.

\noindent {\sl Case 4.} $A \cap S_1$ is a collection of essential simple closed curves in $A$, and $S_1$ is an annulus
disjoint from $\gamma_1$.

We may assume that $A \cap S_1$ consists of as few curves as possible. 
By Lemma 3.12, this number is therefore one or two. Now,
$(M_2, \gamma_2)$ admits a hierarchy satisfying all the requirements (1) - (4). Its reduced length
is at most that of the original hierarchy, and it has shorter length. So, we may attempt to apply induction
to the decomposition of $(M_2, \gamma_2)$ along $A - {\rm int}(N(S_1))$. However, this decomposition may create
a solid torus disjoint from the sutures. Suppose that it does not. Then we have a commutative diagram
$$

\matrix{
(M_1, \gamma_1) & \buildrel S_1 \over \longrightarrow & (M_2, \gamma_2) \cr
\Big\downarrow {\scriptstyle A} && \Big\downarrow {\scriptstyle A - {\rm int}(N(S_1))} \cr
(M_1', \gamma_1') & \buildrel S_1 - {\rm int}(N(A)) \over \longrightarrow & (M_2', \gamma'_2). \cr}
$$
By induction,
$(M_2', \gamma'_2)$ extends to a hierarchy with reduced length equal to the
one starting with $(M_2, \gamma_2)$. Hence, we obtain a hierarchy for
$(M'_1, \gamma'_1)$, starting with $S_1 - {\rm int}(N(A))$,
then the hierarchy for $(M_2', \gamma'_2)$. Its reduced length is equal to
the one for $(M, \gamma)$. 

So, suppose that decomposing $(M_2, \gamma_2)$ along $A - {\rm int}(N(S_1))$ does create a solid torus $Y$ disjoint from the sutures.
Consider first the case where $A \cap S_1$ is a single curve.
The solid torus $Y$ has boundary consisting of a component of $A - {\rm int}(N(S_1))$,
one component of $S_1 - {\rm int}(N(A))$, and an annulus in $\partial M$. We give $Y$ a Seifert fibration
with at most one exceptional fibre, so that $A \cap S_1$ is a fibre. In fact, $Y$ must contain an exceptional fibre, because otherwise,
we may isotope $S_1$ off $A$, and we have dealt with this case already. Hence, by the atoroidality of $E(M, \gamma)$
and the assumption that no component of $E(M, \gamma)$ is a Seifert fibre space other than a solid torus or a copy of $T^2 \times I$, 
there is just one component of $M - {\rm int}(N(A \cup S_1))$ that is a solid torus
disjoint from the sutures. It is incident to just one of the four sides of $A \cap S_1$, because
of the orientations on $S_1$ and $A$.

Consider the annuli $A'$ shown in Figure 7, and oriented as shown. Note that decomposing $(M_2, \gamma_2)$ along $A'$
does not create a solid torus with no sutures. For if it did, this solid torus could not support an
exceptional fibre in its Seifert fibration, and we could therefore isotope $S_1$ off $A$. Hence, by induction, this decomposition extends
to an allowable hierarchy with the same reduced length as the one starting with $(M_2, \gamma_2)$. 
But instead of decomposing along $S_1$ then $A'$, we may start with $(M_{1}, \gamma_{1})$,
decompose along $A$, and then $A'$. This gives the same decorated sutured manifold,
up to homeomorphism. Hence, $A$ does extend to an allowable hierarchy, with the required reduced length, as claimed.

\vskip 18pt
\centerline{
\epsfxsize=5in
\epsfbox{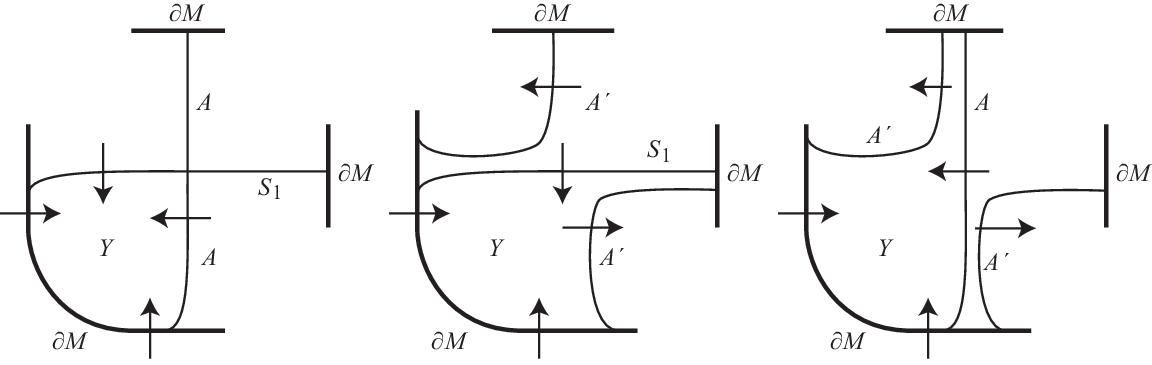}
}
\vskip 6pt
\centerline{Figure 7: When $A \cap S_1$ is a single curve}

Now suppose that $A \cap S_1$ consists of two simple closed curves.
Consider the annulus of $A - {\rm int}(N(S_1))$ disjoint from $\partial M$ and the annulus of $S_1 - {\rm int}(N(A))$
disjoint from $\partial M$. Their union is a torus $T$. We claim that $T$ bounds a solid torus with interior disjoint
from $S_1 \cup A$. Let $Y$ be the component of $E(M, \gamma) - {\rm int}(N(S_1 \cup A))$ with $T$ in its boundary.
The inclusion of $Y$ into $E(M, \gamma)$ is $\pi_1$-injective. So, if $Y$ is not a solid torus, then $T$ is boundary
parallel in $E(M, \gamma)$. But it is then clear that the component of $E(M, \gamma) - {\rm int}(N(S_1 \cup A))$
that is a solid torus with no sutures cannot be Seifert fibred with an exceptional fibre. Hence,
in this case, there is an isotopy reducing $|S_1 \cap A|$, and we have dealt with this already.
So, $Y$ is indeed a solid torus. It must be Seifert fibred with an exceptional fibre, for otherwise,
we can isotope $S_1$ off $A$. No other component of $E(M, \gamma) - {\rm int}(N(S_1 \cup A))$ can be 
an exceptionally fibred solid torus, because this would imply that $E(M, \gamma)$ itself was Seifert fibred and not a solid torus or a copy of $T^2 \times I$,
by atoroidality. Therefore, $Y$ is the solid torus component of $M - {\rm int}(N(S_1\cup A))$
with no sutures, which are assuming to exist. Therefore, the orientations on $S_1$ and $A$
near $Y$ are as shown in Figure 8 (up to reversing all of the transverse orientations). 
The three different rows in that figure illustrate the possible transverse orientations on
$\partial M$ near $\partial A$. 

In the figure are also shown various annuli $A'_1$, $A'_2$, $A'_3$ and $A'_4$.
Note that none of $A_1'$, $A_2'$ or $A_4'$ separate off a solid torus with no sutures, because in this case,
either we could isotope $S_1$ to reduce $|S_1 \cap A|$ or decomposition along $S_1$ would create a solid torus with no sutures. 
However, $A_3'$ might separate off a solid torus with no sutures, in which case it is boundary parallel. Let $A'$ denote the union of $A'_1$, $A'_2$, $A_4'$ and possibly
also $A'_3$. Include $A'_3$ if and only if it does not separate off a solid torus with no sutures. Then
decomposing $(M_2, \gamma_2)$ along $A'$ does not create any solid tori with no sutures.
Hence, inductively, this extends to an allowable hierarchy with same reduced length as the one
starting with $S_2$. Let $(M_2', \gamma_2')$ be the manifold that results from decomposing
$(M_2, \gamma_2)$ along $A'$. Similarly, let $(M_2'', \gamma_2'')$ be the manifold
that results from decomposing $(M_1, \gamma_1)$ along $A' \cup A$. Then, aside from two components
in each, $(M_2', \gamma_2')$ and $(M_2'', \gamma_2'')$ are equal. Furthermore, it is
easy to check that they have admissible hierarchies with the same reduced length. Hence,
$A'$ extends to the required admissible hierarchy, with the required reduced length.

\vskip 18pt
\centerline{
\epsfxsize=5in
\epsfbox{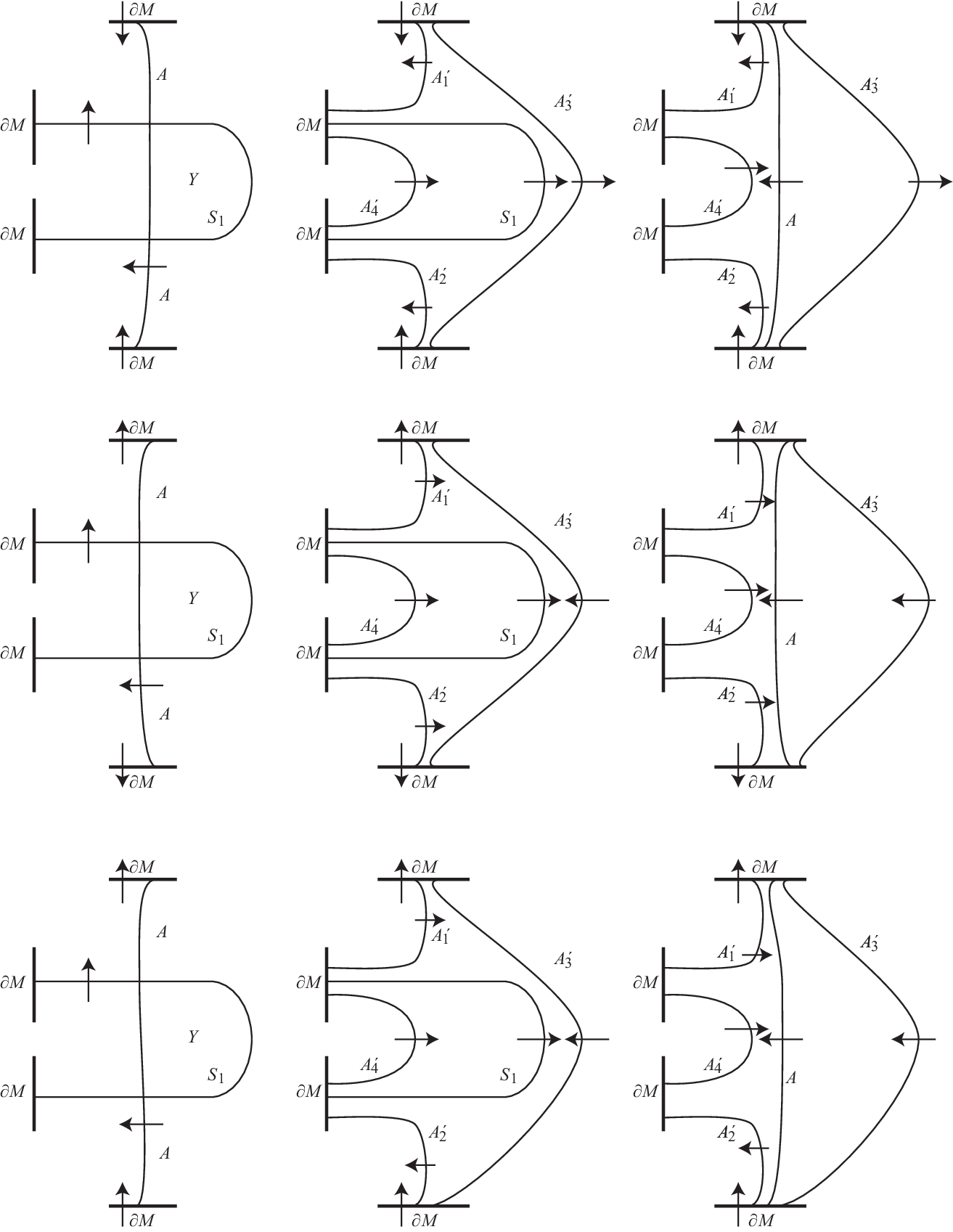}
}
\vskip 6pt
\centerline{Figure 8: When $A \cap S_1$ is two curves}

\noindent {\sl Case 5.} $A$ is not a product annulus, $A \cap S_1$ consists of essential arcs in $A$, and these arcs are not all coherently oriented. 

Then we may find a component of $A - {\rm int}(N(S_1))$
that is a disc $D$ disjoint from $\gamma_2$. Using the fact that $(M_2, \gamma_2)$ is taut,
$\partial D$ bounds a disc $D'$ in $R_\pm(M_2)$. Let $C$ be the intersection between $S_1 \cap D'$
and $R_\pm(M_1) \cap D'$. As argued in the proof of Lemma 4.4, $C$ is a collection of arcs. In this case there are two
of them. There are two cases. If the region in $D'$ between these two arcs lies in $S_1$, then
we may isotope this part of $S_1$ onto $D$, and so remove two curves of $A \cap S_1$.
On the other hand, if the region in $D'$ lies in $R_\pm(M_1)$, then we deduce that $A$ is
boundary-compressible in $M_1$ along a boundary compression disc disjoint from $\gamma_1$.
Therefore, $A$ is parallel to an annulus in $R_\pm(M_1)$. Hence, decomposing
$(M_1, \gamma_1)$ creates a solid torus with no sutures, or a solid torus with a product sutured
manifold structure. The former case is excluded by hypothesis, and in the latter case, the
proposition holds trivially.

\noindent {\sl Case 6.} $A$ is a product annulus, or $A \cap S_1$ consists of arcs that are essential in $A$ and coherently oriented. 

Hence, $A - {\rm int}(N(S_1))$ consists of a collection of product discs $P$ embedded in $(M_2, \gamma_2)$.
Note that we have a commutative diagram
$$

\matrix{
(M_1, \gamma_1) & \buildrel S_1 \over \longrightarrow & (M_2, \gamma_2) \cr
\Big\downarrow {\scriptstyle A} && \Big\downarrow {\scriptstyle P} \cr
(M_1', \gamma_1') & \buildrel S_1 - {\rm int}(N(A)) \over \longrightarrow & (M_2', \gamma'_2). \cr}
$$
By Lemma 4.4, $(M_2', \gamma'_2)$ admits an allowable hierarchy. 
But we may also obtain $(M_2', \gamma'_2)$ by first decomposing
$(M, \gamma)$ along $A$, then along
$S_1 - {\rm int}(N(A))$. Each of these surfaces meets the conditions of being
part of an allowable hierarchy, except possibly $S_1 - {\rm int}(N(A))$.
Its canonical extension may contain some discs disjoint from the sutures.
Suppose that there is such a disc. Then it is not boundary parallel in $M_1'$, because 
this would imply that $A$ is not a product annulus and that not all the arcs of $A \cap S_1$
are coherently oriented. Therefore we deduce that
$R_\pm(E(M'_1))$ is compressible. We deduce that at least one of the hypotheses
of Theorem 2.3 does not apply to the decomposition 
$$E(M_1', \gamma_1') \buildrel E(S_1 - {\rm int}(N(A))) \over \longrightarrow E(M_{2}', \gamma'_{2})$$
because $E(M_{2}', \gamma'_{2})$ is taut, whereas $E(M_1', \gamma'_1)$ is not.
We see that the only possibility is that $E(M_1', \gamma_1')$ contains a component $Y$
that is a solid torus disjoint from the sutures, and which is compressed by some disc components
of $E(S_1 - {\rm int}(N(A)))$. But then $(M_1', \gamma_1')$ contains a solid torus disjoint from the
sutures, and we are assuming that this is not the case.

We have therefore found the required hierarchy extending $A$.  $\square$

\vskip 6pt
\noindent {\caps 4.8. Product-separating surfaces}
\vskip 6pt

Let $(M, \gamma) \buildrel S \over \longrightarrow (M_S, \gamma_S)$ be an allowable decomposition between
decorated sutured manifolds. 
A connected surface $S$ is known as {\sl product-separating} if it separates the component of $M$ that
contains it, and one of the components of $E(M_S, \gamma_S)$ obtained by decomposing the corresponding component
of $E(M, \gamma)$ along $E(S)$ is a product sutured manifold $F \times I$, such that
\item{(i)} $F$ is homeomorphic to $E(S)$;
\item{(ii)} any attached 2-handles in $F \times I$ are of the form $D^2 \times I$, for a disc
$D^2$ in the interior of $F$;
\item{(iii)} each such 2-handle has non-empty intersection with the copy of $E(S)$ in the boundary of $F \times I$.

\noindent We say that a 2-handle satisfying (ii) and (iii) is {\sl strongly vertical} in $F \times I$. A disconnected surface $S$ is {\sl product-separating}
if each of its components is.

Note that this does not imply that $S$ is boundary parallel, even when $S$ is connected. For example, let $F$ be a connected surface properly embedded in $M$,
that is parallel to an essential subsurface of $R_\pm(M)$ and that has no trivial boundary curves. Orient $F$ so that the region between it and the subsurface of $R_\pm(M)$
is a product sutured manifold. Now pick an embedded arc $\alpha$ such that $\alpha \cap F$ is an endpoint of $\alpha$
lying in the interior of $F$, and such that $\alpha \cap \partial M$ is the other endpoint of $\alpha$. Choose $\alpha$
so that it is disjoint from the product region. Also arrange that at one endpoint, $F$ points towards $\alpha$ and
at the other endpoint, $\partial M$ points away from $\alpha$. Now remove a regular neighbourhood of $\alpha \cap F$
from $F$, and attach the annulus that encircles $\alpha$. Let $S$ be the new surface.
Then $E(S)$ is homeomorphic to $F$, and one component of $E(M_S, \gamma_S)$ is a copy of $F \times I$.
This product manifold contains none of the attached
2-handles of $E(M_S, \gamma_S)$. So, $S$ is product-separating. However, $S$ need not be boundary-parallel. 
This example shows that $|\chi(S)|$ may be greater than $|\chi(F)|$, where $F \times I$ is the product sutured manifold.

A related example is where $\alpha$ is chosen to lie within the product sutured manifold between $F$
and the subsurface of $R_\pm(M)$. Again, the resulting surface $S$ is product-separating. However,
if $\alpha$ is knotted, then $S$ need not be boundary-parallel.

A third example again starts with $F$ as above. We now isotope $F$, keeping its boundary in $\partial M$,
so that it approaches a component of $\gamma$ that is not a u-suture. We then slide $\partial F$ over this component of $\gamma$ so that
afterwards, a sub-arc of $\partial F$ and a sub-arc of $\gamma$ together bound a disc in $R_\pm(M)$. The resulting surface $S$
continues to have non-trivial boundary. Hence $E(S) = S$. Moreover, $E(S)$ separates off a product sutured manifold $F \times I$, where
$F$ is homeomorphic to $E(S)$ and that is disjoint from the 2-handles
of $E(M_S, \gamma_S)$. So $S$ is product-separating. In this case, $S$ is boundary-parallel. Note, however,
that in this case $S$ is not the entirety of a component of $R_\pm(F \times I)$.

The following lemma asserts that the examples given above essentially exhaust the possibilities
for a product-separating surface.

\noindent {\bf Lemma 4.10.} {\sl Let $S$ be a connected product-separating surface in the sutured manifold $(M, \gamma)$.
Then, after slicing under discs of contact, $S$ becomes a surface $\overline{S}$ that is parallel to a subsurface $G$
in $\partial M$. Moreover, one of $G \cap R_-(M)$ or $G \cap R_+(M)$ is a connected subsurface of $G$ that
is homeomorphic to $G$. The remainder of $G$ is a collection of discs, each of which intersects $\partial G$
in a single arc, and annuli, each of which intersects $\partial G$ in a single closed curve. Moreover,
decomposing $(M, \gamma)$ along $\overline{S}$ creates a copy of $(M, \gamma)$ and a component that is a product sutured manifold
homeomorphic to $\overline{S} \times I$.}

\noindent {\sl Proof.} By assumption, $S$ is product-separating. Modify $S$ by slicing under discs of contact
as many times as possible. This does not change $E(S)$ or $E(M_S, \gamma_S)$. Thus, the new surface $\overline{S}$
is still product-separating. So, one component of $E(M_S, \gamma_S)$
is a product sutured manifold $F \times I$. By (i) in the definition of product-separating,
$F$ is homeomorphic to $E(S)$. Now, a copy of $E(S)$ is a subsurface of $R_\pm(F \times I)$, 
$R_+(F \times I)$ say. But by (i), $R_+(F \times I)$ is homeomorphic to $E(S)$. Furthermore,
${\rm cl}(R_+(F \times I) - E(S))$ can contain no discs disjoint from $\gamma_S$ because
such a disc would form a disc of contact for $\overline{S}$, which we have removed.
Therefore, ${\rm cl}(R_+(F \times I) - E(S))$ consists of discs that intersect $\gamma_S$
in a single arc and annuli that intersect $\gamma_S$ in a simple closed curve. So, $E(S)$ is boundary parallel in $E(M, \gamma)$. Let $G'$ be the subsurface
of $\partial E(M, \gamma)$ to which it is parallel. The surface ${\rm cl}(R_+(F \times I) - E(S))$
is disjoint from the attached 2-handles, because if a 2-handle
did intersect ${\rm cl}(R_+(F \times I) - E(S))$, it would be disjoint from $E(S)$, contradicting (iii) in
the definition of product-separating. Hence, each component of ${\rm cl}(R_+(F \times I) - E(S))$
is a copy of a component of ${\rm cl}(G' - R_-(E(M,\gamma)))$. 
In the product region between $E(S)$ and $G'$, any attached
2-handles are strongly vertical. Thus, when we form $\overline{S}$ by removing from $E(S)$ any trivialising planar
surfaces for trivial boundary curves of $\overline{S}$, we see that $\overline{S}$ is parallel to a subsurface $G$ of $\partial M$
with the required properties. $\square$

The following result asserts that we may avoid product-separating surfaces in
an allowable hierarchy.

\noindent {\bf Proposition 4.11.} {\sl Let $(M, \gamma)$ be a decorated sutured manifold that
admits an allowable hierarchy. Then $(M, \gamma)$ admits an allowable hierarchy where each
decomposing surface $S$ is connected and not product-separating.}

\noindent {\sl Proof.}
We may arrange that each decomposing surface is connected, because we may replace a decomposition
along a disconnected surface by a sequence of decompositions along its components.
Consider a shortest allowable hierarchy for $(M, \gamma)$ where each decomposing
surface is connected. Then we claim that the first surface $S$ cannot be product-separating.
This will imply that every surface in the hierarchy is not product-separating, by induction.
Suppose $S$ is product-separating. Then we modify the hierarchy to a shorter one, as follows. First arrange,
using Lemma 4.2, for each trivial boundary curve of each decomposing surface to be
parallel towards to a u-suture. This does change the fact that the first surface (which we will still call $S$)
is product-separating. Since $S$ now has no discs of contact, Lemma 4.10 implies that $S$ is now boundary-parallel.
Moreover, decomposing $(M, \gamma)$ along $S$ gives a copy of $(M, \gamma)$ and a product sutured manifold.
So, we may discard $S$, and discard the product sutured manifold. The later surfaces in the hierarchy therefore form
an admissible hierarchy for $(M, \gamma)$ with shorter length. $\square$

\vfill\eject
\centerline{\caps 5. Sutured manifolds and handle structures}
\vskip 6pt

\noindent {\caps 5.1. Handles structures and their complexity}
\vskip 6pt

In this paper, handle structures on 3-manifolds will play a central role. If ${\cal H}$
is a handle structure, then ${\cal H}^i$ denotes the union of the $i$-handles.
We will always insist that our handle structures satisfy the following requirements:
\item{(i)} The intersection between each $i$-handle $D^i \times D^{3-i}$ and $\bigcup_{j < i} {\cal H}^j$
is $\partial D^i \times D^{3-i}$.
\item{(ii)} Any two $i$-handles are disjoint.
\item{(iii)} The intersection between any 2-handle $D^2 \times D^1$ and any 1-handle $D^1 \times D^2$
is $\alpha \times D^1$ in the 2-handle, where $\alpha$ is a collection of disjoint arcs in $\partial D^2$, and of the
form $D^1 \times \beta$ in the 1-handle, where $\beta$ is a collection of disjoint arcs in $\partial D^2$.
\item{(iv)} Every 2-handle runs over at least one 1-handle.

When $(M, \gamma)$ is a sutured manifold, we will also insist that a handle structure
${\cal H}$ for $M$ satisfies the following conditions:
\item{(i)} The intersection between $\gamma$ and any 0-handle is a (possibly empty) collection of arcs and simple closed curves.
\item{(ii)} The intersection between $\gamma$ and any 1-handle $D^1 \times D^2$ is
a (possibly empty) collection of arcs, each of the form $D^1 \times \{ \ast \}$, where $\ast$ is a point in $\partial D^2$.
\item{(iii)} The sutures $\gamma$ are disjoint from the 2-handles and 3-handles.

In a handle structure ${\cal H}$, the surface ${\cal F} = {\cal H}^0 \cap ({\cal H}^1 \cup {\cal H}^2)$
will be important. We view it as a subsurface of $\partial {\cal H}^0$. It inherits a
handle structure, where each 0-handle is a component of ${\cal H}^0 \cap {\cal H}^1$,
and each 1-handle is a component of ${\cal H}^0 \cap {\cal H}^2$. We denote the union
of these 0-handles of ${\cal F}$ by ${\cal F}^0$, and the union of the 1-handles of ${\cal F}$
by ${\cal F}^1$.

\vskip 6pt
\centerline{
\epsfxsize=4in
\epsfbox{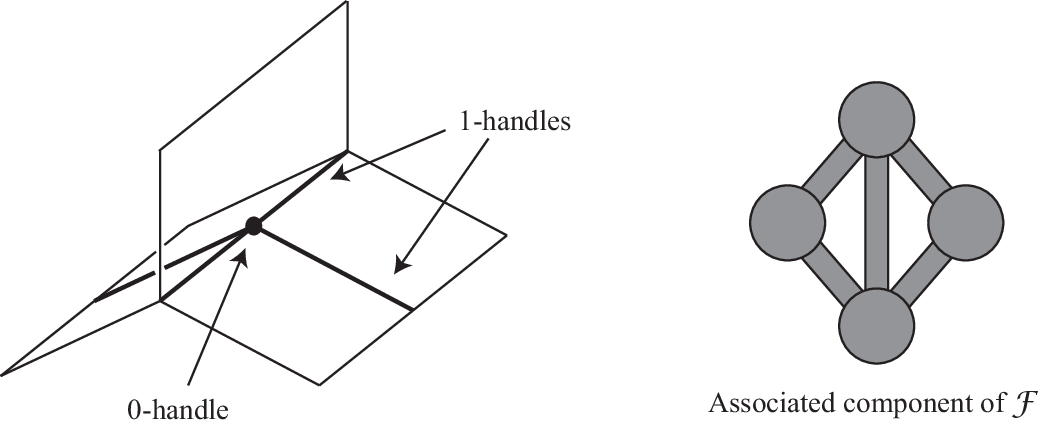}
}
\vskip 6pt
\centerline{Figure 9: Handles near a 0-handle}

In [21], the author defined a measure of complexity for a handle structure of a sutured manifold.
This will form a foundation for this paper. We recall it now.

For a component $F$ of ${\cal F}$, its {\sl index} $I(F)$ is
$$I(F) = - 2\chi(F) + |F \cap \gamma|.$$
Three integers were defined:
$$\eqalign{
C_1(F) &= |F \cap {\cal F}^1| + 1, \cr
C_2(F) &= I(F), \cr
C_3(F) &= |\partial F|.}$$
The {\sl ${\cal F}$-complexity set} $C_{\cal F}({\cal H})$ is the collection of ordered triples
$$\{ (C_1(F), C_2(F), C_3(F)): F \hbox{ is a component of } {\cal F}  \hbox{ with } I(F)>0\},$$
where repetitions are retained.

We order the triples within an ${\cal F}$-complexity set, using lexicographical ordering, as follows. We declare that
$$(C_1(F), C_2(F), C_3(F)) \leq (C_1(F'), C_2(F'), C_3(F'))$$
if and only if 
\item{(i)} $C_1(F) \leq C_1(F')$, or
\item{(ii)} $C_1(F) = C_1(F')$ and $C_2(F) \leq C_2(F')$, or
\item{(iii)} $C_1(F) = C_1(F')$ and $C_2(F) = C_2(F')$ and $C_3(F) \leq C_3(F')$.

We use this to place a total ordering on ${\cal F}$-complexity sets, again using lexicographical ordering, as follows. Given two
${\cal F}$-complexity sets $C_{\cal F}({\cal H})$ and $C_{\cal F}({\cal H}')$, for handle structures ${\cal H}$
and ${\cal H}'$, we first consider
the largest triple in each set. If one is greater than the other, then we declare that the
${\cal F}$-complexity set containing it is greater. If the two triples are equal, then
we pass to second-largest triples in each set, and compare these, and so on.
If at some stage, we run out of triples in one of the sets, then the other ${\cal F}$-complexity
set is declared to be the greater one.
This allows us to compare any two ${\cal F}$-complexity sets. It is clear that this is in
fact a well-ordering on ${\cal F}$-complexity sets (see Lemma 5.3 in [21]).

We can use this ordering on ${\cal F}$-complexity sets to define a notion of complexity
for a handle structure ${\cal H}$, as follows. The {\sl complexity} $C({\cal H})$ of a handle structure ${\cal H}$
is the ordered pair $(C_{\cal F}({\cal H}), n({\cal H}))$, where $n({\cal H})$ is the number
of 0-handles of ${\cal H}$ containing a component of ${\cal F}$ with positive index. We compare
the complexity of handle structures ${\cal H}$ and ${\cal H}'$ by asserting that $C({\cal H}) > C({\cal H}')$
if and only if
\item{(i)} $C_{\cal F}({\cal H}) > C_{\cal F}({\cal H}')$, or
\item{(ii)} $C_{\cal F}({\cal H}) = C_{\cal F}({\cal H}')$ and $n({\cal H}) < n({\cal H}')$.

\noindent It is shown in Lemma 5.3 in [21] that this ordering on complexity of handle structures
is a well-ordering.

The index $I(F)$ of a component $F$ of ${\cal F}$ may be computed as follows. If $F_0$ is a 0-handle of ${\cal F}$, then its {\sl index} is defined to be
$$I(F_0)  = -2 + |F_0 \cap {\cal F}^1| + |F_0 \cap \gamma|.$$
Then it is easy to check that
$$I(F) = \sum_{F_0} I(F_0),$$
where $F_0$ runs over every 0-handle of $F$.

In many circumstances, it will be useful to focus on handle structures ${\cal H}$ with some constraints.
We say that ${\cal H}$ is {\sl positive} if each 0-handle of ${\cal F}$ has positive index, and, for each 0-handle $H_0$ of ${\cal H}$,
$H_0 \cap ({\cal F} \cup \gamma)$ is connected.

\vskip 6pt
\noindent {\caps 5.2. Standard and regulated surfaces}
\vskip 6pt

A surface $S$ properly embedded in a 3-manifold $M$ with a handle structure ${\cal H}$ is {\sl standard} if
\item{(i)} it intersects each 0-handle in a collection of properly embedded disjoint discs;
\item{(ii)} its intersection with any 1-handle $D^1 \times D^2$ is $D^1 \times \alpha$, where
$\alpha$ is a collection of arcs properly embedded in $D^2$;
\item{(iii)} it intersects each 2-handle $D^2 \times D^1$ in discs of the form $D^2 \times \{ \ast \}$,
where $\ast$ is a point in the interior of $D^1$;
\item{(iv)} it is disjoint from the 3-handles.

\vskip12pt
\centerline{
\epsfxsize=4in
\epsfbox{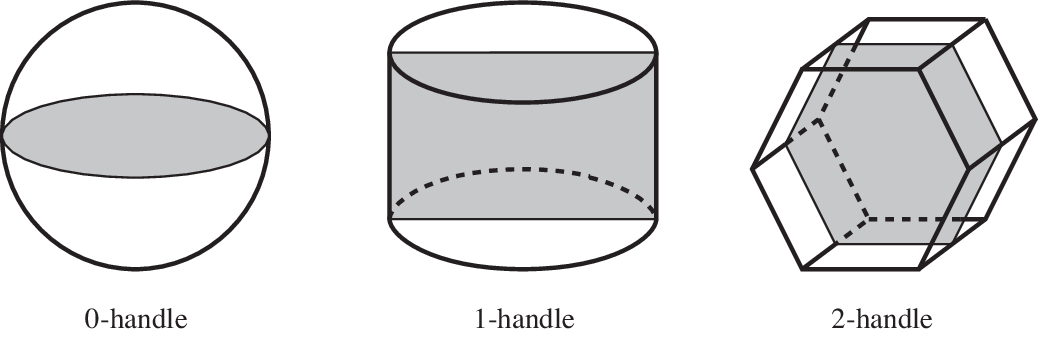}
}
\vskip 6pt
\centerline{Figure 10: Standard surface}

Let $(M, \gamma)$ be a sutured manifold with a handle structure ${\cal H}$. Let $S$ be a transversely oriented, standard surface properly embedded in $M$, 
with $\partial S$ transverse to $\gamma$. Then,
the sutured manifold $(M', \gamma')$ obtained by decomposing along $S$ inherits a handle
structure. Each component of intersection between $M'$ and an $i$-handle of ${\cal H}$ becomes
an $i$-handle for $M'$. Let ${\cal H}'$ denote this handle structure.

We now explain how, in certain circumstances, the complexity of ${\cal H}'$ must be at most that of ${\cal H}$.
In order to guarantee this, the standard surface $S$ must satisfy some conditions. These were explained in detail in
Section 9 of [21]. We recall them here.

\noindent {\sl Condition 1.} Each curve of $S \cap \partial {\cal H}^0$ meets any 1-handle of ${\cal F}$
in at most one arc.

This crucial condition is also one that appears in the theory of normal surfaces in handle structures.

\vskip 12pt
\centerline{
\epsfxsize=1.7in
\epsfbox{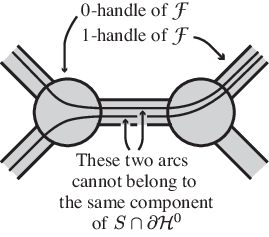}
}
\centerline{Figure 11: Condition 1}

\noindent {\sl Condition 2.} If $\alpha$ is an arc of $({\cal F}^0 \cap \partial {\cal F}) - (S \cup \gamma)$,
and both its endpoints lie in $S$, then the transverse orientations on $\alpha$
and on $S \cap {\cal F}$ near $\partial \alpha$ cannot all agree. Here, the transverse orientation on 
$\alpha$ is the one that it inherits from lying in $\partial {\cal F} \cap R_\pm(M)$. We mean that it
`agrees' with the transverse orientation on $S \cap {\cal F}$ near $\partial \alpha$ if these
three transverse orientations either all point out of the component of ${\cal F} - S$ containing $\alpha$,
or they all point in. (See Figure 12.)

This condition does not arise in normal surface theory. It was one of the distinctive features of
the surfaces considered in [21]. It imposes considerable constraints on the possible transverse orientations
on the discs of $S \cap {\cal H}^0$ that intersect $\partial M$. 

In [21], the arc $\alpha$ is known as a {\sl tubing arc}, and we use the same terminology here.

\vskip 6pt
\centerline{
\epsfxsize=2.2in
\epsfbox{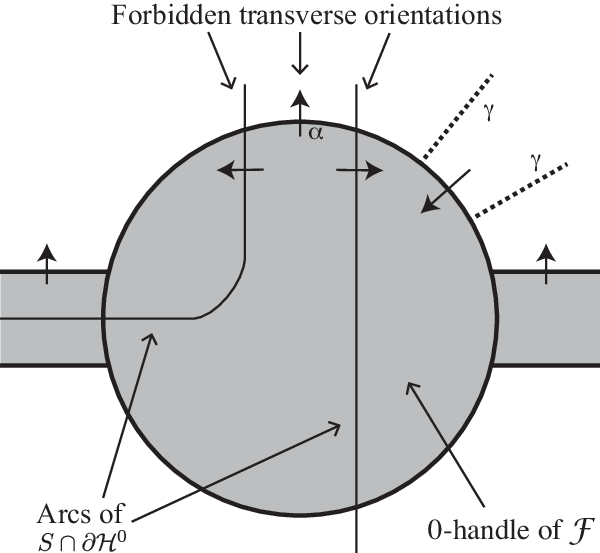}
}
\vskip 6pt
\centerline{Figure 12: Condition 2}

\noindent {\sl Condition $3'$.} Suppose that $D$ is a disc in ${\cal F}^0$ with $\partial D$
the union of two arcs $\alpha$ and $\beta$, where $\alpha = S \cap \partial D$
and $\beta = D \cap \partial {\cal F}^0$. Suppose that one endpoint of $\alpha$
lies in $R_\pm(M)$ and one endpoint lies in ${\cal F}^1$. Then at least one
of the following must hold:
\item{(i)} the interior of $\beta$ contains a component of ${\cal F}^0 \cap {\cal F}^1$;
\item{(ii)} $\beta$ has non-empty intersection with $\gamma$.

This is very close to Condition 3 of [21], except it is somewhat weaker. In Condition 3 of [21],
another configuration was ruled out, where $\beta$ had a single intersection with $\gamma$,
and $S$ was oriented in a certain way.

\vskip 6pt
\centerline{
\epsfxsize=1.9in
\epsfbox{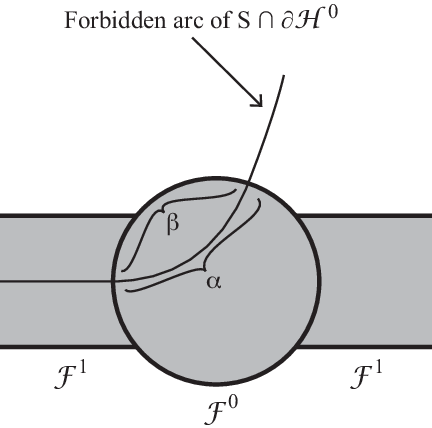}
}
\vskip 6pt
\centerline{Figure 13: Condition $3'$}

\noindent {\sl Condition 4.} Each component of $S \cap \partial {\cal H}^0$
meets any component of $R_\pm(M) \cap \partial {\cal H}^0$ in at most one arc.

\vskip 12pt
\centerline{
\epsfxsize=2.2in
\epsfbox{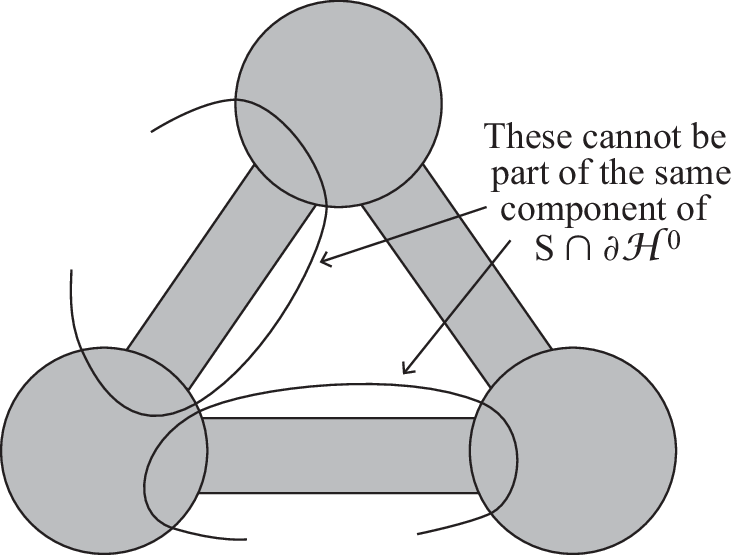}
}
\vskip 6pt
\centerline{Figure 14: Condition 4}

\noindent {\sl Condition 5$'$.} If $\beta$ is an arc of $S \cap {\cal F}^0$ with both endpoints
in $R_\pm(M)$, then each of the two arcs in $\partial {\cal F}^0$ joining $\partial \beta$
must either contain a component of ${\cal F}^0 \cap {\cal F}^1$ or hit $\gamma$.

\centerline{
\epsfxsize=1.6in
\epsfbox{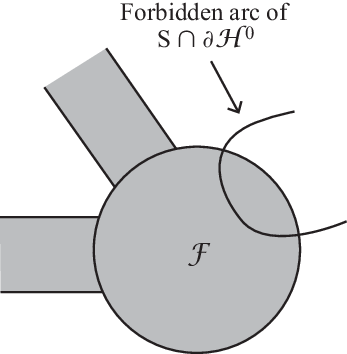}
}
\vskip 6pt
\centerline{Figure 15: Condition 5$'$}

Condition 5$'$ is very similar to Condition 5 in [21], except it is weaker. In Condition 5 of [21],
other configurations were also ruled out, where one of the arcs in $\partial {\cal F}^0$ joining $\partial \beta$
had at most two intersections with $\gamma'$.

We say that a standard surface satisfying Conditions 1, 2, 3$'$, 4 and 5$'$ is {\sl regulated}.

\vskip 6pt
\noindent{\caps 5.3. Comparison with normal surfaces}
\vskip 6pt

There is obviously a close connection between standard surfaces, regulated surfaces and the more well-known
normal surfaces. In this subsection, we briefly recall the definition of normal surfaces and compare them with
standard and regulated surfaces.

Normal surfaces also come in various flavours. In the setting of triangulated 3-manifolds, a properly embedded surface
is {\sl normal} if it intersects each tetrahedron in a collection of triangles and squares. However, although triangulations
will be used briefly in this paper, it is handle structures that will be our main concern.

So, let $M$ be a compact 3-manifold with a handle structure ${\cal H}$.
For simplicity, we will focus on the case where $\partial M$ is not decorated in any way.
In particular, we will disregard any sutured manifold structure. Then 
a {\sl normal surface} in $M$ is a standard surface satisfying Conditions 1 and 4 from the previous subsection
and in addition, the following constraints:
\item{(i)} No component of $S \cap {\cal F}^0$ has endpoints lying in the same component of
$\partial {\cal F}^0 \cap {\cal F}^1$, or in the same component of $\partial {\cal F}^0 - {\cal F}^1$, or in
adjacent components of $\partial {\cal F}^0 \cap {\cal F}^1$ and $\partial {\cal F}^0 - {\cal F}^1$.
\item{(ii)} No arc of intersection between $S$ and $\partial {\cal H}^0 - {\cal F}$ has endpoints in 
the same component of $\partial {\cal F}^0 - {\cal F}^1$.

%\vskip 6pt
%\noindent {\caps 5.4. Weakly normal surfaces}
%\vskip 6pt
%
%Unfortunately, the above definition of a normal surface is often a little too restrictive for our purposes.
%This is because when $S$ is a properly embedded normal surface in a triangulated 3-manifold $M$,
%and ${\cal H}$ is the handle structure on $M$ dual to the triangulation, then $S$ does not necessarily
%correspond to a normal surface in ${\cal H}$. 
%
%When $M$ is closed, such a correspondence does hold. This is because, in this case, for each
%0-handle $H_0$ of ${\cal H}$, the intersection $H_0 \cap {\cal F}$ is as shown in the left of Figure 14.
%The only elementary normal discs in $H_0$ with boundary in ${\cal F}$ correspond to triangles
%and squares in the triangulation.
%
%However, when $M$ has non-empty boundary, the situation is a little more complicated.
%In this case, the handles of ${\cal H}$ are in one-one correspondence with the simplices
%of the triangulation {\sl that do not lie wholly in $\partial M$}. So, in general, 0-handle $H_0$ 
%of ${\cal H}$, the intersection $H_0 \cap {\cal F}$ is  subset of the handle structure shown in
%the left of Figure 14. An elementary normal disc in the triangulation corresponds to a standard
%surface in ${\cal H}$, but it may violate conditions (i) and (ii) in the previous subsection.
%However, it does satisfy Conditions 1 and 4 from Section 5.2.
%
%So, we define a surface $S$ properly embedded in a 3-manifold $M$ with a handle structure ${\cal H}$ to be
%{\sl weakly normal} if it is standard and satisfies Conditions 1 and 4 from Section 5.2.

\vskip 6pt
\noindent {\caps 5.4. Elementary disc types}
\vskip 6pt

Let ${\cal H}$ be a handle structure of a 3-manifold, and let $H$ be one of its handles. Let $D$
be a disc properly embedded in $H$. Then a {\sl normal isotopy} of $D$ is an ambient isotopy
that preserves each of the handles of ${\cal H}$ and also $\gamma$ (when $M$ is a sutured manifold). 
When two discs properly embedded in $H$ are normally isotopic,
they are said to be of the same {\sl type}.

Typically, we will be interested in discs $D$ that are a component of $S \cap H$, for some handle $H$, 
where $S$ is a standard, normal or regulated surface. Such discs are called {\sl elementary discs}.

One of the limitations of standard surfaces is that there is, in general, no upper bound on the number of
elementary disc types in a handle. However, it is clear that there is usually such an upper
bound for regulated and for normal surfaces. In fact, for this to be true, we need to make the
hypothesis that, for each 0-handle $H$ of ${\cal H}$, $\partial H \cap (\gamma \cup {\cal F})$
is connected. This is equivalent to the statement that $\partial H \cap R_\pm(M)$ consists
of discs or a sphere. Under this hypothesis, we can then show that for each handle $H$ of ${\cal H}$,
$H$ contains only finitely many normal or regulated disc types, as follows.

It suffices to prove this for elementary discs in a 0-handle, because once these are controlled,
then the elementary discs in 1-handles are determined. Furthermore, in each 2-handle,
there is always only one elementary disc type. An elementary disc $D$ properly embedded in a
0-handle $H$ is determined, up to normal isotopy, by its boundary curve. This curve is determined
by its intersection with ${\cal F}^1$, with ${\cal F}^0$ and with $\partial H \cap R_\pm(M)$.
It can run over each 1-handle of ${\cal F}$ at most once by Condition 1 and it can run
over each component of $\partial H \cap R_\pm(M)$ at most once by Condition 4. Thus,
there are only finitely many possibilities for its intersection with ${\cal F}^1$
and $\partial H \cap R_\pm(M)$, because $\partial H \cap R_\pm(M)$ consists of discs. 
Hence, there are only finitely many possibilities for its intersection with $\partial {\cal F}^0$, and so there are only finitely many
possibilities for its intersection with ${\cal F}^0$. So, up to normal isotopy,
there are only finitely many possibilities for $D$. Moreover, it is clear that
once one is given the handle structure on ${\cal F}\cap \partial H$ and the arcs $\gamma \cap \partial H$,
the possible elementary disc types are readily computable.

\vskip 6pt
\noindent {\caps 5.5. Trivial and weakly trivial modifications}
\vskip 6pt

When a sutured manifold $(M, \gamma)$ with a handle structure ${\cal H}$ is decomposed
along certain surfaces, then sometimes very little happens to $(M, \gamma)$ and ${\cal H}$.
In [21], decompositions of this sort were considered, and the resulting changes to ${\cal H}$
were termed `trivial modifications'. In this subsection, we briefly recall this term in this context,
and then introduce a variant, known as a `weakly trivial modification'.

Let $S$ be a standard, transversely oriented surface properly embedded in $(M, \gamma)$. Let $H_0$ be a 0-handle
of ${\cal H}$. Then decomposition along $S$ results in a {\sl trivial modification} to $H_0$ if
\item{(i)} each component of $S \cap H_0$ is parallel to a disc in $\partial H_0$ 
that is disjoint from $\gamma$;
\item{(ii)} this component of $S \cap H_0$ has boundary entirely in ${\cal F}$;
\item{(iii)} the disc in $\partial H_0$ intersects ${\cal F}$ in an annulus $A$,
and, in the induced handle structure on $A$, each 0-handle of $A$ has index zero.

We say that decomposition along $S$ results in a {\sl weakly trivial modification} to $H_0$ if 
\item{(i)} each component of $S \cap {\cal F}$ is parallel
in ${\cal F}$ to an arc or circle in $\partial {\cal F}$;
\item{(ii)} in the resulting product region, each 0-handle of ${\cal F} - {\rm int}(N(S))$ has index zero;
\item{(iii)} no component of $S \cap \partial H_0$ separates components of ${\cal F} - {\rm int}(N(S))$
with positive index.

\vskip 18pt
\centerline{
\epsfxsize=5.8in
\epsfbox{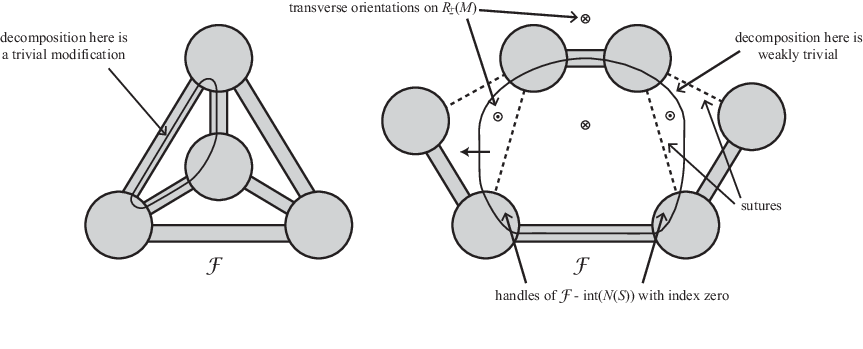}
}
\vskip 6pt
\centerline{Figure 16: Trivial and weakly trivial modifications}

\vfill\eject
\noindent {\bf Proposition 5.1.} {\sl Let ${\cal H}$ be a handle structure of a sutured manifold
$(M, \gamma)$. Suppose that ${\cal H}$ is positive. Let 
$$(M, \gamma) \buildrel S \over \longrightarrow (M_S, \gamma_S)$$
be a taut decomposition, where $S$ is a standard surface properly in $M$.
Suppose that this results in a weakly trivial modification to each 0-handle of ${\cal H}$.
Then $S$ is product-separating.}

\noindent {\sl Proof.} We will show that each component $S''$ of $S$ satisfies the following:
\item{(i)} $S''$ is separating in the component of $M$ that contains it;
\item{(ii)} one of the components of the sutured manifold obtained by decomposing this component of $(M, \gamma)$
along $S''$ is a product sutured manifold $F'' \times I$;
\item{(iii)} the inclusion of $S''$ into $R_+(F'' \times I)$ or $R_-(F'' \times I)$ is a homotopy equivalence.

This will imply that $S''$ is product-separating, for the following reasons. Note that (iii) implies
that any u-suture in $F'' \times I$ must be parallel to a trivial component of $\partial S''$. So when
we form $E(S'')$, the effect is to add a disc along this trivial curve, and when we form $E(F'' \times I)$,
the effect is to add a 2-handle along the corresponding suture. Hence, $E(F'' \times I)$ is still a product sutured manifold $F' \times I$,
and moreover $F'$ is homeomorphic to $E(S'')$. This verifies condition (i) in the definition of
product-separating. We also verified that any attached 2-handles in $E(F'' \times I)$ respect its
product structure and intersect $E(S'')$. Therefore, (ii) and (iii) in the definition of product-separating
are also verified.

Our aim is now to verify (i), (ii) and (iii) above.
We first modify $S$, forming a new surface $S'$. Consider any 0-handle $H_0$ of ${\cal H}$,
and let $F$ be ${\cal F} \cap \partial H_0$.
Consider any arc component of $S \cap F$. By definition of a weakly trivial modification, each such arc is
boundary parallel in $F$. Hence, we may find an arc that is outermost in $F$.
It separates off an index zero subdisc of $F$ which intersects $\partial F$ in a single arc $\alpha$.
Push $\alpha$ a little away from $F$. Then $\alpha$ forms part of the
boundary of a product disc in $M_S$. The rest of the boundary of this product disc
lies in the interior of a disc of $S \cap H_0$. Boundary compress $S$ along this disc.
Now repeat this process for all arcs of $S \cap F$, and for all 0-handles of ${\cal H}$.
Note that this new surface $S'$ will, in general, run over u-sutures, and so decomposition
along $S'$ need not be allowable, but it is taut by Proposition 2.5. Let $(M_{S'}, \gamma_{S'})$ be the
sutured manifold that results from this decomposition.

\vskip 18pt
\centerline{
\epsfxsize=3.5in
\epsfbox{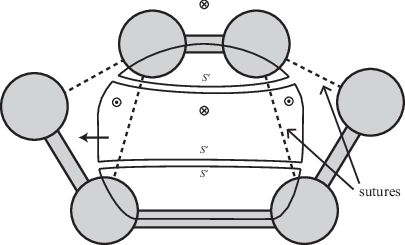}
}
\vskip 6pt
\centerline{Figure 17: The surface $S'$ that results from boundary-compressing $S$}

We will now show that every component of $S''$ of $S'$ satisfies (i), (ii) and (iii) above. 
Our proof follows Proposition 10.7 in [21] closely. Let $C$ be the collection
of arcs and simple closed curves of $S' \cap {\cal F}$ that are incident to components
of ${\cal F} - {\rm int}(N(S'))$ with positive index.
Then $C$ separates off from ${\cal F}$ a collection of annuli $A$ disjoint from $\gamma$ and product discs $P$
(which may contain other arcs and curves of $S' \cap {\cal F})$.
Let $D$ be the collection of discs of $S' \cap {\cal H}^0$ which contain a component of $C$ in their boundary.
Then $(A \cup P) \cup D$ is a collection of discs properly embedded in $M$. These are parallel to
discs in $\partial M$ via balls $B_0$. The existence of such a disc in $\partial M$ is guaranteed,
in the case where a component of $C$ lies in ${\cal F}$, by the fact that $H_0 \cap ({\cal F} \cup \gamma)$ is connected for each 0-handle $H_0$
of ${\cal H}$. When a component of $C$ is an arc in ${\cal F}$, then the existence of the
disc in $\partial M$ follows from the construction of $S'$. Each of these balls $B_0$ may be given a product structure
as $D^2 \times [-1,1]$ so that $D^2 \times \{ 1 \}$ is a component of $S' \cap B_0$.

For each 1-handle $H_1 = D^2 \times [0,1]$, the discs $D^2 \times \{ 0 \}$ and $D^2 \times \{ 1 \}$ are
each divided up by the decomposition along $S'$. For each $i \in \{ 0,1 \}$, all but one 0-handle of 
$D^2 \times \{ i \} - {\rm int}(N(S'))$ has index zero. The remaining component has index
equal to the index of $D^2 \times \{ i \}$. We are assuming that the index of $D^2 \times \{ i \}$
is positive. Hence, the product structure on $H_1$ matches $(A \cup P) \cap (D^2 \times \{ 0 \})$
with $(A \cup P) \cap (D^2 \times \{ 1 \})$. We may therefore unambiguously define $(A \cup P) \cap D^2$.
Let $B_1$ be the union, over all 1-handles, of the balls $((A \cup P) \cap D^2) \times [0,1]$.
Similarly, we may find a collection $B_2$ of balls in the 2-handles of $M$ such that $B_2 \cap {\cal H}^0 = {\cal H}^2 \cap (A \cup P)$.
Again, each component of $B_1$ and $B_2$ may be given a product structure of the form
$D^2 \times [-1,1]$ so that $S' \cap (D^2 \times [-1,1]) \supseteq D^2 \times \{ 1 \}$.
Moreover, the product structures on $B_0$, $B_1$, and $B_2$ agree on their intersection.
So, $B_0 \cup B_1 \cup B_2$ is an $[-1,1]$-bundle over a surface $F$, in which the sutures
form the zero-section over $\partial F$. So, $B_0 \cup B_1 \cup B_2$ is a product sutured manifold
that forms a product region between some components of $S'$ and a subsurface
of $\partial M$. Remove these components of $S'$ and repeat. We eventually deduce that every component of $S'$ that intersects ${\cal F}$ is boundary parallel. 

We claim also that the components of $S'$ that are disjoint from ${\cal F}$ are also boundary parallel.
Their boundary curves come from components of $\partial S \cap \partial {\cal H}^0$, possibly isotoped
a little so that they become disjoint from ${\cal F}$. Since no disc of $S \cap {\cal H}^0$ has
0-handles of $M_S$ with positive index on both sides of it, we deduce that
the discs of $S' \cap {\cal H}^0$ that are disjoint from ${\cal F}$ cannot separate
components of ${\cal F}$, which proves the claim. These discs therefore separate
off components of $(M_{S'}, \gamma_{S'})$ which are balls. These balls are product sutured manifolds,
because $(M_{S'}, \gamma_{S'})$ is taut.

We have therefore shown that each component $S''$ of $S'$ satisfies (i), (ii) and (iii) above. Cutting
$M$ along $S'$ gives a copy $Y$ of $M$, together with some product sutured manifolds. We divide $S'$
into subsurfaces $S'_1, \dots, S'_n$, as follows. Define the distance of a component of $S'$ from $Y$ to be the
minimal number of intersection points between an arc and $S'$, where the arc starts on that component and
ends in $Y$. Define $S'_i$ to be the union of components with distance $i$. So, $S'_1$ is the union of
components of $S'$ incident to $Y$, and so on.

Consider $S'_1$, the union of the components of $S'$ closest to $Y$. These separate off product regions $W$
in $(M, \gamma)$ with interior disjoint from $S$. Consider any of the product discs that was used to boundary compress
$S$, and that is incident to $S'_1$. By construction, these have interior disjoint from $W$. Hence, when we reverse 
any of these boundary compressions, the resulting surface is still boundary parallel, and still separates
off a product sutured manifold. This new surface satisfies (i), (ii) and (iii) above.
Repeat this until all boundary compressions incident to $S'_1$ have been reversed.
The resulting surface $S_1$ is a union of components of $S$. Continue with $S'_2$,
and so on.

In this way, we see that each component of $S$ satisfies (i), (ii) and (iii). So, $S$ is product-separating.
$\square$

\vskip 6pt
\noindent {\caps 5.6. The behaviour of complexity under decomposition}
\vskip 6pt

For the purposes of this paper, it is critically important that when a sutured manifold $(M, \gamma)$ with a handle
structure ${\cal H}$ is decomposed along a surface $S$, the resulting handle structure ${\cal H}'$ is not any more complex
than ${\cal H}$. One of the main reasons for introducing regulated surfaces is that this holds in this case, under some mild hypotheses.
This was essentially shown in [21]. However, there the surfaces satisfied the slightly different
Conditions 1-5. (See Propositions 10.6 and 10.7 in [21]).
In this subsection, we explain how a similar result holds for regulated surfaces.

\noindent {\bf Theorem 5.2.} {\sl Let ${\cal H}$ be a positive handle structure for a sutured manifold $(M, \gamma)$.
Let 
$$(M, \gamma) \buildrel S \over \longrightarrow (M', \gamma')$$
be a taut decomposition along a regulated surface $S$, and let ${\cal H}'$ be the handle structure obtained by decomposing ${\cal H}$ along $S$.
Then $C(H_0 \cap {\cal H}') \leq C(H_0)$ for each 0-handle $H_0$ of ${\cal H}$. 
If this inequality is an equality for some 0-handle $H_0$, then $H_0 \cap {\cal H}'$ is obtained
from $H_0$ by a weakly trivial modification. If this inequality is an equality for every 0-handle $H_0$, then
$S$ is product-separating.}

Throughout this subsection, we will use the following terminology. We denote the handle structures
of $M$ and $M'$ by ${\cal H}$ and ${\cal H}'$. As usual, ${\cal H}^i$ will denote the union of
$i$-handles of ${\cal H}$, and ${\cal F}$ denotes ${\cal H}^0 \cap ({\cal H}^1 \cup {\cal H}^2)$.
We define ${\cal F}'$ similarly, but with ${\cal H}'$ in place of ${\cal H}$.

The following is a version of Lemma 10.5 in [21]. The difference is that, in Lemma 10.5 of [21],
Conditions 1, 2 and 3 were assumed to hold. However, we note that only the weaker Condition
3$'$ was used in the proof.

\noindent {\bf Lemma 5.3.} {\sl Let $S$ be a standard surface satisfying Conditions 1, 2 and 3$'$.
Let $D$ be a component of ${\cal F}'$ with a negative index 0-handle. Then there is a 1-handle of
${\cal F}$ which lies entirely in $D$.}

The following is a version of Propositions 10.6 and 10.7 in [21].

\noindent {\bf Proposition 5.4.} {\sl Let ${\cal H}$ be a positive handle structure.
Let $H_0$ be a 0-handle of ${\cal H}$, and 
let $H'_0$ denote the 0-handles $H_0 \cap {\cal H}'$. If $S$ is a standard surface satisfying
Conditions 1, 2 and 3$'$, then $C(H_0') \leq C(H_0)$. Also, if we have equality, then
$H_0'$ is obtained from $H_0$ by a weakly trivial modification.}

%More precisely, the following hold.
%Each component of $S \cap {\cal F}$ is an arc or circle parallel to an arc or circle in $\partial {\cal F}$.
%The parallelity region inherits a handle structure from ${\cal F}$ in which each 0-handle
%has index zero. For all but one 0-handle of $H'_0$, the intersection with ${\cal F}'$ consists
%only of components with index zero.}

The difference here is that only Conditions 1, 2 and 3$'$ are assumed to hold,
whereas in Proposition 10.6 in [21], Conditions 1 - 5 were hypothesised. However,
we obtain a slightly weaker conclusion. We deduce in the case where 
$C(H_0') = C(H_0)$ that $H_0'$ is obtained from $H_0$ by a weakly trivial
modification. On the other hand, in Proposition 10.6 in [21], we could deduce that
a trivial modification had been performed.

\noindent {\sl Proof.} Let $F$ be a component of ${\cal F} \cap H_0$, and let $F'$ be
$F \cap {\cal F}'$. The first paragraph in the proof of Proposition 10.6 in [21]
gives that $C_{\cal F}(F') \leq C_{\cal F}(F)$. Since this holds for each component
$F$ of ${\cal F} \cap H_0$, we deduce that $C_{\cal F}(H'_0) \leq C_{\cal F}(H_0)$.
Note also that $n(H_0) = 1$ and $n(H'_0) \geq 1$. Therefore, $C(H'_0) \leq C(H_0)$.

Suppose that this is an equality. This implies that if $F$ is any component of
${\cal F} \cap H_0$, and $F' = F \cap {\cal F}'$, then $C_{\cal F}(F') = C_{\cal F}(F)$.
As argued in the first paragraph of the proof of Proposition 10.6 in [21], one component $X$ of $F'$ has
$C_1(X) = C_1(F')$ and $I(X) = I(F)$. As argued in the second paragraph of the proof of 
Proposition 10.6 in [21], $|\partial X| = |\partial F|$, and so $F'$ is obtained from $F$ 
by cutting along arcs and circles which are boundary parallel in $F$. Each component of
$F' - X$ has index zero. Moreover, every 0-handle of $F'$ not lying in $X$ has index zero.
We therefore have verified conditions (i) and (ii) in the definition of a weakly trivial modification.

Therefore, each component $F$ of ${\cal F} \cap H_0$ gives rise to a corresponding component
of ${\cal F}' \cap H_0$ also with positive index. All of these components of ${\cal F}'$ must lie in
the same component of $H'_0$, as otherwise $n(H'_0) > 1 = n(H_0)$. Thus, we have 
verify condition (iii) in the definition of a weakly trivial modification.
$\square$

\noindent {\sl Proof of Theorem 5.2.} According to Proposition 5.4, decomposition along
the regulated surface does not increase the complexity of any 0-handle. Also, if the
complexity of a 0-handle is unchanged, then a weakly trivial modification is performed there.
If a weakly trivial modification occurs at every 0-handle, then according to Proposition 5.1,
$S$ is product-separating. $\square$

\vskip 6pt
\noindent {\caps 5.7. The complexity and weight of standard surfaces}
\vskip 6pt

Many arguments in this field aim to modify surfaces until they are of some required type.
For example, one might first make a surface standard and then perform further modifications
to make it normal. In this second stage, one typically aims to reduce some measure of complexity 
for the surface.

In this paper, we define the {\sl complexity} of a standard surface $S$ to be the ordered pair of integers
$(|S \cap {\cal H}^2|, |\partial S \cap {\cal H}^1|)$. We compare complexities in the usual way,
using lexicographical ordering. This notion of complexity is one of the most commonly used ones.
In particular, it is the same as the one employed by the author in [21].

The first integer in the pair, $|S \cap {\cal H}^2|$, is called the {\sl weight} of the surface,
and is also denoted $w(S)$. Similarly, when $S$ is a normal surface properly embedded in a
triangulated 3-manifold, its weight $w(S)$ is its number of points of intersection with the 1-skeleton.

\vskip 6pt
\noindent {\caps 5.8. The existence of regulated surfaces}
\vskip 6pt

In Section 9 of [21], it was established that a standard decomposing surface may be often
upgraded to a regulated one. In this and the next subsection,
we will prove a version of this result, which is as follows.

\noindent {\bf Theorem 5.5.} {\sl Let
$$(M,  \gamma) \buildrel S \over \longrightarrow (M', \gamma')$$
be an allowable decomposition between taut decorated sutured manifolds
that extends to an allowable hierarchy.
Suppose that $E(M, \gamma)$ is atoroidal and that no component of 
$E(M, \gamma)$ is a Seifert fibre space, other than a solid torus or a copy
of $T^2 \times I$. Assume that $(M, \gamma)$ is not a product sutured manifold.
Let ${\cal H}$ be a positive handle structure for $(M, \gamma)$. Suppose that $S$ is standard in ${\cal H}$.
Then there is another allowable sutured manifold decomposition
$$(M,  \gamma) \buildrel S' \over \longrightarrow (M'', \gamma'')$$
such that 
\item{(i)} $S'$ is a regulated surface in ${\cal H}$;
\item{(ii)} $S'$ is connected;
\item{(iii)} $S'$ is not product-separating;
\item{(iv)} $S'$ extends to an allowable hierarchy;
\item{(v)} $(M'', \gamma'')$ is taut;
\item{(vi)} if $S$ is non-separating, then $S'$ is also.

}

This is proved using the methods developed in Section 9 of [21]. There, various modifications were made to the standard surface
$S$. These modifications do not increase the complexity of the surface at any stage. See the discussion after Lemma 9.9 in [21]
where this assertion is explicitly made. Indeed the complexity of $S$ is defined the way it is in [21] to ensure that it is decreased by these
modifications, and hence that these modifications are guaranteed to terminate. The modifications are as follows.

\noindent {\sl Modification 1.} Tubing along an arc $\alpha$ 

This was described in detail in Section 4.5. We will always require that the two arcs of $\alpha \times \{ - 1, 1 \}$
do not lie in the same trivial boundary curve of the resulting surface $S'$. Note that, according to Lemma 4.5, 
$S'$ extends to an allowable hierarchy.

\noindent {\sl Remark 5.6.} Tubing along an arc was a key tool in [21] for arranging a decomposing surface to be regulated.
In short, if a decomposing surface $S$ fails to satisfy Condition 2, then one can apply this modification to reduce $|\partial S \cap {\cal H}^1|$
without increasing $|S \cap {\cal H}^2|$. Thus, if we assume that the complexity of $S$ is minimal, 
then $S$ satisfies Condition 2. Unfortunately, though, tubing along
an arc $\alpha$ can introduce trivial boundary curves. This occurs, for instance, when $\alpha$ is parallel to a sub-arc
of $\partial S$. It is for this reason that we have had to introduce the machinery of decorated sutured manifolds in this paper.

\noindent {\sl Modification 2.} Slicing under a non-trivial annulus.

Let $A$ be an annulus in $R_\pm(M)$ which has non-trivial boundary curves and has $A \cap S = \partial A$. Suppose that the
orientation on $S$ agrees with orientation on $A$. Then we can construct a new surface $S'$ by attaching $A$ to $S$,
and then isotoping so that it becomes properly embedded. Let $(M_{S}, \gamma_{S})$ and $(M_{S'}, \gamma_{S'})$
be the sutured manifolds obtained from $(M, \gamma)$ by decomposing along $S$ and $S'$ respectively.
Note that there is a non-trivial product annulus $P$ properly embedded in $(M_{S'}, \gamma_{S'})$, with one boundary component
being a core curve of $A$, and the other boundary component being a core curve of the copy of $A$ in $S'$. Moreover,
decomposing $(M_{S'}, \gamma_{S'})$ along $P$ gives a decorated sutured manifold homeomorphic to $(M_S, \gamma_S)$.
Therefore, by Proposition 4.9, $S'$ extends to an allowable hierarchy if and only if $S$ does, under the hypotheses of Theorem 5.5.

\noindent {\sl Modification 3.} Sliding across $\gamma$.

We will not describe this move in detail here, because we will avoid using it in this paper. The reason is that it creates
new intersection points between $\partial S$ and $\gamma$. In [21], this was fine, but in this paper, $\partial S$ is not
permitted to intersect a u-suture. This is why we avoid this move.

\noindent {\sl Modification 4.} Slicing under a disc of contact.

This modification was described in detail in Section 4.3. According to Lemma 4.3, the resulting surface $S'$
also extends to an allowable hierarchy.

\noindent {\sl Modification 5.} Boundary compressing along a product disc disjoint from $\gamma$,
and then possibly removing a planar component.

This is described in Section 4.6. If the surface obtained by boundary-compressing has a component
that is planar, disjoint from $\gamma$ and with all but at most one boundary component trivial, then this component
is removed. By Proposition 4.7,
if $S'$ is obtained from $S$ by this process, and $S$ extends to an allowable hierarchy, then so does $S'$.

\noindent {\sl Modification 6.} Removal of a product region.

Suppose that a component $S_1$ of $S$ is parallel into $R_\pm(M)$, and that the interior of the product region between
$S_1$ and the subsurface of $R_\pm(M)$ is disjoint from $S$. This modification is the removal of $S_1$ from $S$. Let $S'$ be the resulting surface.
Note that the sutured manifold obtained by decomposing along $S$ is equal to the disjoint union of the manifold obtained 
decomposing along $S$ and a product sutured manifold. Hence, $S'$ extends to an allowable hierarchy if and only if $S$ does.

\noindent {\bf Lemma 5.7.} {\sl If a surface is non-separating in $M$, then this remains the case for at least
one component of the surface that results from Modifications 1, 2, 4, 5 and 6.}

\noindent {\sl Proof.} If one starts with a surface $S$, then each of these modifications does not
change $[S,\partial S] \in H_2(M, \partial M)$. If $S$ is non-separating, then $[S, \partial S]$
is non-trivial and so the resulting surface $S'$ also has $[S', \partial S']$ non-trivial.
Hence, some component is homologically non-trivial, and therefore non-separating. $\square$

\noindent {\bf Proposition 5.8.} {\sl Suppose that $(M, \gamma)$ contains no non-separating product disc or non-separating product annulus.
Let $S$ be a connected surface in $(M, \gamma)$ that is not product-separating and that extends to an allowable hierarchy. Then this remains the case
for at least one component of the surface that results from Modifications 1, 2, 4, 5 and 6.}

\noindent {\sl Proof.} We may assume that $M$ is connected, by focusing on the component containing $S$.
Note also that the statement of the lemma is empty for Modification 6, because $S$ is connected
and not product-separating, and so Modification 6 cannot be applied to it.

If $S$ is non-separating, the lemma follows from Lemma 5.7. So, suppose that $S$ is separating.
We may also assume that each component of $S'$ is separating, because a non-separating surface
is automatically not product-separating.

We now consider each of the modifications in turn. Suppose that $S'$ is obtained from $S$ by Modification 1.
As explained in Section 4.5, there is a decomposition
$$(M_{S'}, \gamma_{S'}) \buildrel D \over \longrightarrow (M_S, \gamma_S)$$
where $D$ is a product disc. We are assuming that the two arcs $\alpha \times \{ -1,1\}$
do not lie in the same trivial boundary curves of $S'$. Hence, there are a number of
other possible cases to consider.

Suppose first that at least one arc of $\alpha \times \{ -1, 1\}$ lies in a trivial boundary curve of $S'$.
We will focus on the case where exactly one arc of $\alpha \times \{ -1, 1\}$ lies in a trivial curve of $\partial S'$,
as the other case is very similar.
Then $E(S')$ is homeomorphic to $E(S)$. If $S'$ is product-separating, then $E(S')$ separates
off a product manifold homeomorphic to $E(S') \times I$ in which each attached 2-handle is vertical
and intersects $E(S')$.
The manifold $E(M_{S'}, \gamma_{S'})$ is homeomorphic to $E(M_S, \gamma_S)$. In fact, $E(M_{S'}, \gamma_{S'})$ is obtained
from $E(M_S, \gamma_S)$ by attaching a copy of $D^2 \times I$ along a sub-arc of
$\gamma_S$. This region $D^2 \times I$ is made up of attached 2-handles and $P \times I$, where
$P$ is the trivialising planar surface for the trivial boundary curve of $\partial S'$
incident to $\alpha \times \{ -1,1\}$. So, $E(M_S, \gamma_S)$ contains a component
homeomorphic to $E(S) \times I$ and the attached 2-handles
are vertical in this product structure and intersect $E(S)$. So, $S$ is product-separating.

Suppose now that neither arc of $\alpha \times \{ -1,1\}$ lies in a trivial boundary curve of $S'$.
Then there is a decomposition
$$E(M_{S'}, \gamma_{S'}) \buildrel D \over \longrightarrow E(M_S, \gamma_S).$$
One component of $E(M_{S'}, \gamma_{S'})$ is, by assumption, a copy of $E(S') \times I$.
If this component contains $D$, then we decompose along it to obtain a component
of $E(M_S, \gamma_S)$ that is a copy of $E(S) \times I$.
Hence, in this case, $S$ is product-separating. So, suppose that the product manifold $F \times I$ does not
contains $D$. Then it remains a product manifold component of $E(M_S, \gamma_S)$. However,
it may not be of the required form: it may not be homeomorphic to $E(S) \times I$.
So, we use an alternative argument. Consider an arc of $\alpha \times \{ -1,1\}$. This is
part of the boundary of $S'$. Consider the arc in $R_\pm(E(M_S, \gamma_S))$ that is
the copy of this arc of $\alpha \times \{ -1,1\}$, and that is not part of a suture. Say
it is lies in $R_+(E(M_S, \gamma_S))$. This arc lies
in the product sutured manifold $F \times I$. Since each component of ${\rm cl}(R_+(F \times I) - E(S'))$
intersects $\gamma$ by Lemma 4.10, we may pick an arc running in $R_+(M)$ from $\alpha \times \{ -1 ,1\}$
to $\gamma$. We do this for each of the two components of $\alpha \times \{ -1, 1\}$. 
Join these two arcs by an arc running over the tube. The result is an arc $\beta$ in the product
sutured manifold. This is part of a product disc in $F \times I$. Now enlarge this product disc across the tube
to form a product disc in $(M, \gamma)$. It is non-separating, since one may find a closed
curve intersecting it once, consisting of an arc in the connected surface $S$ and an arc running across the tube.
This contradicts a hypothesis of the lemma.

The argument in the case of Modification 2 is similar but simpler. In this case, there is a decomposition 
$$E(M_{S'}, \gamma_{S'}) \buildrel P \over \longrightarrow E(M_S, \gamma_S),$$
where $P$ is a product annulus. By assumption, some component of $E(M_{S'}, \gamma_{S'})$
is a product sutured manifold. If $P$ lies in it, then we deduce that $S$ is also product-separating.
On the other hand, if $P$ does not lie in the product manifold, then the union of $P$ and a vertical annulus in the product structure
forms a product annulus in $(M, \gamma)$. It is non-separating, because one may find a closed
curve that intersects it once, as in the previous case. This contradicts one of the hypotheses of
the lemma.

Modification 4 is particularly easy to handle. If $S'$ is obtained from $S$ by slicing under
a disc of contact, then $E(S')$ and $E(S)$ are homeomorphic, as are $E(M_S, \gamma_S)$
and $E(M_{S'}, \gamma_{S'})$. Moreover, in the latter case, this homeomorphism sends the attached 2-handles
to the attached 2-handles. So, $S'$ is product-separating if and only if $S$ is.

Now consider the case of Modification 5. The surface $S'$ may be connected or disconnected.
Suppose first that $S'$ is connected. Suppose also that the tubing arc $\alpha$ is disjoint from the product sutured
manifold. There are a number of possibilities. If neither curve of $\alpha \times \{ -1,1 \}$ lies in
a trivial boundary curve of $S$, then we obtain a decomposition
$$(E(S) \times [-1,1], \partial E(S) \times \{ 0 \}) \buildrel D \over \longrightarrow (E(S') \times [-1,1], \partial E(S) \times \{ 0 \}).$$
Hence, $S$ is product-separating. If just one curve of $\alpha \times \{ - 1,1 \}$ lies in a trivial boundary curve of $S$,
then $E(S)$ is homeomorphic to $E(S')$, and $E(M_S, \gamma_S)$ is homeomorphic to $E(M_{S'}, \gamma_{S'})$.
Hence, as argued in the case of Modification 1, $E(M_S, \gamma_S)$ contains a component homeomorphic
to $E(S) \times I$ and any attached 2-handles are vertical in this product structure.
So, again $S$ is product-separating. If the arcs of $\alpha \times \{ -1,1 \}$ lie
in distinct trivial boundary curves of $S$, then again $E(S)$ is homeomorphic to $E(S')$, and 
$E(M_S, \gamma_S)$ is homeomorphic to $E(M_{S'}, \gamma_{S'})$, and so $S$ is again
product-separating. Finally, suppose that the arcs $\alpha \times \{ -1,1 \}$ lie in the same trivial
boundary curve of $S$. Then we obtain $E(S)$ from $E(S')$ by attaching an annulus,
and $E(S') \times [-1,1]$ is obtained from a component of $E(M_S, \gamma_S)$ 
by decomposing along a product annulus. Therefore, this component of $E(M_S,\gamma_S)$
is homeomorphic to $E(S) \times I$, as required.

Now consider the case where $S'$ is connected, and the tubing arc $\alpha$ lies in the product sutured
manifold. Then, $\alpha$ lies in the closure of a component of  $R_\pm(M) - \partial S'$. By Lemma 4.10, this is a disc that intersects
$\partial S'$ in a single arc or an annulus that intersects $\partial S'$ in a simple closed curve. Therefore, exactly one arc of $\alpha \times \{ -1, 1\}$ 
lies in a trivial component of $\partial S$.
We have already seen in this case that $E(S)$ is homeomorphic to $E(S')$, and $E(M_S, \gamma_S)$ 
is homeomorphic to $E(M_{S'}, \gamma_{S'})$. Hence, $S$ is product-separating.

Now suppose that $S'$ is disconnected with components $S'_1$ and $S'_2$, that are both product-separating.
If neither of the resulting product manifolds contains the tubing arc, then the argument above applies.
So, suppose that at least one of the product manifolds, $E(S'_1) \times I$, say,
contains the tubing arc. Then $\alpha$ lies in the closure of a component of 
$R_\pm(M) - \partial S'_1$. By Lemma 4.10, this is a disc $D$ that intersects
$\partial S'_1$ in a single arc or an annulus that intersects $\partial S'_1$ in a simple closed curve. 
Consider the disc case first. The other end of $\alpha$ lies in $S'_2$. This
creates a suture of $\gamma_{S'}$ which lies in $D$. Since decomposition along $S'$
is taut, we deduce that $E(S'_2)$ is a disc, which separates off a product ball. 
Hence, there is a component of $E(M_S, \gamma_S)$ homeomorphic to $E(S'_1) \times I$,
and this homeomorphism takes 2-handles to 2-handles. Also, $E(S)$ is homeomorphic to $E(S'_1)$.
So, $S$ is product-separating.

Suppose now that $\alpha$ lies in an annulus $A$ that is the closure of a component of $R_\pm(M) - \partial S'_1$. 
We may assume that the intersection
between $\partial S'_2$ and $A$ does contain any arcs. For if it did, an outermost arc in $A$ would separate off
a disc, and we could then argue as in the previous case. Therefore, we may assume that $\partial S'_2$ intersects $A$
in at least one core curve. Hence, decomposing $(M, \gamma)$ along $S'_2$ gives an annular
component of $R_\pm$. Suppose that this lies in the product sutured manifold that $E(S'_2)$ separates off.
Then $E(S'_2)$ is an annulus. Now, the tubing arc $\alpha$ lies in an annular component of $R_\pm(M) - \partial S'$.
Hence, the two arcs $\alpha \times \{ - 1,1 \}$ lie in the same trivial curve of $\partial S$. This bounds
a disc of contact. So, we deduce that $E(S)$ is homeomorphic to $E(S'_1)$, and $E(M_S, \gamma_S)$
is homeomorphic to $E(M_{S_1'}, \gamma_{S_1'})$. Therefore, $S$ must have been product-separating.

There is one final case that we must consider. Suppose that $S'$ is disconnected with
components $S'_1$ and $S'_2$, both of which are product-separating. Suppose that
both of the product manifolds contain the tubing arc $\alpha$. Assume also that
$\alpha$ lies in an annular component of $R_\pm(M) - S'_1$ and an annular
component of $R_\pm(M) - S'_2$. Then, $\alpha$ lies in an annular components $A$
of $R_\pm(M)$, and it runs between a component of $\partial S'_1$ and a component of
$\partial S'_2$ that are essential curves in this annulus. Let $F_1 \times I$ be the product
manifold that $E(S'_1)$ separates off. Then $E(S'_2)$ lies in $F_1 \times I$. It is an
incompressible surface in this product, and hence vertical or parallel to a subsurface
of $F_1 \times \partial I$. In the vertical case, it is actually parallel to a regular
neighbourhood of a component of $\partial A$. In the boundary parallel case, it
separates off a product submanifold. This has an annular component of $R_\pm$,
that is a subset of $A$. Hence, $E(S_2')$ must be an annulus. It is therefore again parallel
to a regular neighbourhood of a component of $\partial A$. Now, when we tube along $\alpha$
to form $S$, we create a curve of $\partial S$ that bounds a disc of contact in $A$. Hence,
$E(S)$ is actually homeomorphic to $E(S'_1)$. Also, a component of $(M_S, \gamma_S)$
is homeomorphic to $F_1 \times I$. So, $S$ is product-separating. $\square$

\vskip 6pt
\noindent {\caps 5.9. Ensuring Conditions 1, 2, 3$'$, 4 and 5$'$}
\vskip 6pt

We follow the procedure given in Section 9 of [21], where Conditions 1 - 5 were introduced.
There, it was shown that if the surface $S$ did not satisfy these conditions, then
Modifications 1 - 6 could be performed to it which reduced its complexity.
These modifications were performed in a specific order, and we follow a similar
strategy here.

First, if a surface fails to satisfy Condition 1, 2 or 3$'$, then moves described in Section 9 of [21] are performed
which reduces its complexity. We note that when Condition $3'$ is violated, then Modification 3
is not required. Note also that Modification 4 might be required when making the new surface
standard, as in Lemma 9.9 in [21]. We therefore obtain the following lemma.

\noindent {\bf Lemma 5.9.} {\sl 
Let
$$(M, \gamma) \buildrel S \over \longrightarrow (M_S, \gamma_S)$$
be an allowable decomposition between taut decorated sutured manifolds.
Suppose that $S$ is standard in the handle structure ${\cal H}$, but fails to satisfy one of Conditions 1, 2 or 3$'$.
Then a sequence of Modifications 1, 2 and 4 can be applied to create a standard surface
with smaller complexity.}

It is shown in [21]  that if $S$ does not satisfy Conditions 4 or 5$'$, then Modifications 4, 5 or 6
may be applied to take it to a standard surface with smaller complexity. So, we obtain the
following.

\noindent {\bf Lemma 5.10.} {\sl 
Let $S$ be a standard surface in the handle structure ${\cal H}$ and suppose that
$$(M, \gamma) \buildrel S \over \longrightarrow (M_S, \gamma_S)$$
is an allowable decomposition between taut decorated sutured manifolds.
Then a sequence of Modifications 1, 2, 4, 5 and 6 
can be applied to create a regulated surface with no greater complexity.}

\noindent {\sl Proof of Theorem 5.5.}
By Proposition 4.11, $(M, \gamma)$ admits an allowable hierarchy, where the first
surface $S$ is connected and is not product-separating. Suppose first $(M,\gamma)$ admits
such a hierarchy where the first surface is non-separating. By Lemma 5.10, we may
apply a sequence of Modifications 1, 2, 4, 5 and 6 
to create a regulated surface $S''$ with smaller complexity.
This also extends to an allowable hierarchy.
By Lemma 5.7, at least one component of $S''$ is non-separating. Call this $S'$.
Then $S'$ is automatically not product-separating.

Suppose now that $(M, \gamma)$ does not admit an allowable hierarchy where the first
surface $S$ is connected and non-separating. In particular, $(M,\gamma)$ does not
contain a non-separating product disc or a non-separating product annulus by Lemma 4.4 and Proposition 4.9.
Then, by Lemma 5.10, a sequence of Modifications 1, 2, 4, 5 and 6 
can be applied to create a regulated surface with no greater complexity. 
This also extends to an allowable hierarchy. At each stage,
one of the components of the surface is not product-separating, by Proposition 5.8.
At each stage, we focus on this component, and discard the rest. At the final stage,
we end with the required surface $S'$.
$\square$

\vskip18pt
\centerline {\caps 6. New handle structures from old ones}
\vskip 6pt

\vskip 6pt
\noindent {\caps 6.1. Simplification modifications}
\vskip 6pt

According to Theorem 5.5, if a taut sutured manifold admits an allowable hierarchy, then one may find such a hierarchy
where the first decomposing surface is regulated. However, the theorem has various technical hypotheses, including that the initial handle structure
${\cal H}$ must be positive. Recall from Section 5.1 that this means that each 0-handle of ${\cal F}$ must have positive index and, 
for each 0-handle $H_0$ of ${\cal H}$, $H_0 \cap ({\cal F} \cup \gamma)$ is connected.
In this subsection, we explain some of the modifications that one can make to a handle structure of a sutured manifold to ensure that it becomes positive. 
These procedures were first introduced in Sections 7 and 8 of [21].

\noindent {\sl Procedure 1.} Slicing a 0-handle along a disc.

Suppose that $D$ is a disc properly embedded in some 0-handle $H_0$, and that $\partial D$ lies in $R_\pm(M)$, is
disjoint from ${\cal F}$, and separates ${\cal F} \cap H_0$.  Assuming that $R_\pm(M)$ is incompressible, then
$\partial D$ bounds a disc $D'$ in $R_\pm(M)$. Then Procedure 1 is decomposition along $D$, oriented so that
a suture appears along $D \cap D'$. If $(M, \gamma)$ is taut, then this decomposition is taut.

\vfill\eject
\noindent {\sl Procedure 2.} Collapsing a 2-handle and a 1-handle disjoint from $\gamma$.

Suppose that $H_1$ is a 1-handle of $M$ disjoint from $\gamma$ and that intersects ${\cal H}^2$
in a single disc. Let $H_2$ be the 2-handle containing this disc. Then Procedure 2 is the removal of $H_1$
and $H_2$. It is also the enlargement of ${\cal H}^3$ if these handles are incident to ${\cal H}^3$.

\noindent {\sl Procedure 3.} Collapsing a 2-handle and a 1-handle containing an arc of $\gamma$.

Let $H_1$ be a 1-handle of $M$ that intersects ${\cal H}^2$ in a single disc and intersects
$\gamma$ in a single arc. Let $H_2$ be the 2-handle incident to $H_1$. Then Procedure 3
is the removal of $H_1$ and $H_2$. The arc $\gamma \cap H_1$ is replaced by an arc
that runs along the part of the attaching circle of $H_2$ that is disjoint from $H_1$, and
also through the two attaching discs of $H_1$.

\noindent {\sl Procedure 4.} Decomposition along a product disc, then sliding $\gamma$.

Let $H_1$ be a 1-handle disjoint from ${\cal H}^2$ and that intersects $\gamma$ in two arcs.
Then Procedure 4 is the removal of $H_1$. For each of the two discs that form $H_1 \cap {\cal H}^0$,
the two points of intersection between this disc and $\gamma$ are joined by new sutures.
This is achieved by a decomposition along the product disc which is the co-core of $H_1$.
In this paper, we will only perform this move when neither of the arcs of $\gamma \cap H_1$
lie in a u-suture. This guarantees that decomposition along the product disc is allowable.

\noindent {\sl Procedure 5.} Collapsing a 3-ball disjoint from $\gamma$.

Suppose that a component of $M$ is a 3-ball disjoint from $\gamma$, and consists of two 0-handles
joined by a 1-handle. Then Procedure 5 is the removal of this 1-handle and one of the 0-handles.

\noindent {\sl Procedure 6.} Collapsing a 2-handle and 3-handle.

Let $H_2$ be a 2-handle that intersects ${\cal H}^3$ in a single disc. Let $H_3$ be the incident 3-handle.
Then Procedure 6 is the removal of $H_2$ and $H_3$.

\vskip 6pt
\noindent {\sl Procedure 7.} Introduction of a 2-handle.

Suppose that $M'$ is a subset of $M$ with the following properties:
\item{(i)} $M'$ is a union of handles of ${\cal H}$;
\item{(ii)} if any $j$-handle of ${\cal H}$ is incident to an $i$-handle of $M'$, where $j > i$, then the $j$-handle also lies in $M'$;
\item{(iii)} $M'$ is homeomorphic to $D^2 \times I$ (and we will make this identification);
\item{(iv)} $(D^2 \times I) \cap \partial M = D^2 \times \partial I$;
\item{(v)} $M'$ is disjoint from $\gamma$;
\item{(vi)} the only handles of $M'$ incident to $\partial D^2 \times I$ are 1-handles and 2-handles;
\item{(vii)} each 1-handle of $M'$ incident to $\partial D^2 \times I$ intersects the 2-handles in exactly two discs.

\noindent Then in Procedure 7, $M'$ is removed from $M$, and is replaced by a single 2-handle. Note that by (i) and (ii),
${\rm cl}(M - M')$ inherits a handle structure. By (iii) and (iv), attaching the 2-handle to ${\rm cl}(M - M')$ results in a
manifold homeomorphic to $M$. 
By (vi) and (vii), the attaching locus of this 2-handle is well-defined.
By (v), $M'$ inherits a sutured manifolds structure homeomorphic to $(M, \gamma)$.

\vskip 6pt
\noindent {\caps 6.2. The universal collection of 0-handle types}
\vskip 6pt

Handle structures are very general objects and, in particular, there are infinitely many
possibilities for the way that a handle may intersect the union of its neighbours.
We will need to restrict our attention to handle structures where the 0-handles come in
finitely many possible local types, as follows.

Define a 0-handle $H_0$ as {\sl tetrahedral} if $H_0 \cap {\cal F}$ has
four 0-handles and six 1-handles, so that any two 0-handles
of ${\cal F}$ are joined by a 1-handle. (The boundary of a tetrahedral
0-handle is shown in the left of Figure 16.) 
Tetrahedral 0-handles are very common. For example, if one starts with any triangulation of a
closed 3-manifold, and one forms the dual handle structure, then each 0-handle is tetrahedral.

Define a 0-handle $H_0$ as {\sl subtetrahedral} if $H_0 \cap {\cal F}$
is obtained from the tetrahedral case by removing some handles.
These arise, for example, when one starts with a triangulation of a compact 3-manifold with boundary.
Then, the dual handle structure has an $i$-handle for each $(3-i)$-simplex that does not lie entirely in $\partial M$.
It is easy to see that each 0-handle is subtetrahedral.

We will decompose handle structures along regulated surfaces,
to form new handle structures. The following result, which is explained in Section 11 of [21] and in the proof of Theorem 1.4 of [21], 
implies that, in the resulting handle structures, there are only finitely many possible types of 0-handle.
It is proved by induction on the complexity of the handle structure, as defined in Section 5.1.

\noindent {\bf Theorem 6.1.} {\sl There is a universal computable constant $N$ and a
universal computable collection of 0-handle types $H_0^1, \dots, H_0^k$, with following property.
Suppose that one starts with a subtetrahedral 0-handle disjoint from any sutures.
Suppose then that a sequence of simplifying procedures as in Section 6.1 and decompositions along 
regulated surfaces are performed. Assume also that, whenever a decomposition along a regulated
surface is made, then the handle structure is positive. Then, 
each of resulting 0-handles is one of the universal
types $H_0^i$, and at most $N$ of these have positive index.}

We say that a handle structure of a sutured manifold is {\sl of uniform type} if each 0-handle
is one of the types given in Theorem 6.1.  An algorithm to produce all the 0-handle types in Theorem 6.1 is given in Section 11 of [21].

\vskip 18pt
\centerline {\caps 7. Boundary-regulated surfaces}
\vskip 6pt

A key feature of normal surfaces is that they can be encoded by a finite list of numbers, which form a vector.
This introduction of linear algebra is very useful. Unfortunately, it does not seem to be possible to create
such a theory for regulated surfaces. The main difficulty is that regulated surfaces are required
to have a transverse orientation, which is hard to incorporate algebraically.

Therefore, in this section, we introduce a new type of surface, which we call boundary-regulated.
These are weaker than regulated surfaces, but they have the advantage of admitting a more algebraic theory,
somewhat akin to normal surface theory.

\vskip 6pt
\noindent {\caps 7.1. Definition}
\vskip 6pt

Let $(M, \gamma)$ be a sutured manifold with a handle structure ${\cal H}$.
Let $S$ be a standard surface properly embedded in $M$. We say that $S$ is {\sl boundary-oriented} if 
each elementary disc $D$ of $S$ that intersects $\partial M$ is assigned a transverse orientation,
and, at each point in $\partial S$ where two elementary discs intersect, their transverse orientations agree.

We note that the Conditions 1, 2, $3'$, $4$ and $5'$ in Section 5.2 only depend on a transverse orientation on $S$
near $\partial S$. Hence, they continue to make sense for a boundary-oriented surface $S$,
even when $S$ does not admit a global transverse orientation. We say that a boundary-oriented standard surface
satisfying Conditions 1, 2, $3'$, $4$ and $5'$ is {\sl boundary-regulated}. Thus, the only difference between regulated
and boundary-regulated surfaces is that regulated surfaces have a global transverse orientation,
whereas boundary-regulated surfaces have a transverse orientation specified only over 
a subset of the surface. Specifically, we have the following observation.

\noindent {\bf Lemma 7.1.} {\sl A regulated surface is exactly a boundary-regulated surface with
a transverse orientation that restricts to the given transverse orientation on the elementary discs
that intersect $\partial M$.}

As the decomposing surfaces in our sutured manfiold hierarchy are necessarily transversely oriented, 
we will want to realise them as regulated surfaces.

In this section, we will develop a theory of boundary-regulated surfaces.
But first, we recall some of the main points of classical normal surface theory.

\vskip 6pt
\noindent {\caps 7.2. Normal surface vectors}
\vskip 6pt

Let ${\cal H}$ be a handle structure for a 3-manifold $M$. It is well known
that normal surfaces properly embedded in $M$ may be described as solutions to a system of linear
equations, as follows.

There is one variable for each type of elementary normal disc in the 0-handles of ${\cal H}$.
So, for a normal surface $S$, one obtains a list of non-negative integers,
which count the number of elementary normal discs of $S$ of each type in each 0-handle. This is the {\sl normal
surface vector} associated with $S$, and is denoted by $(S)$.

There is one equation for each type of elementary normal disc in the 1-handles of ${\cal H}$.
When an elementary normal disc in a 0-handle runs over a 0-handle of ${\cal F}$,
it intersects this 0-handle in a collection of arcs. These arcs extend in a unique
way to elementary normal discs in the associated 1-handle of ${\cal H}$. The {\sl matching equations}
express the fact that, for each 1-handle $H_1$ of ${\cal H}$ that is attached to 0-handles
$H_0$ and $H'_0$, and for each type of elementary normal disc $D$ in $H_1$, the numbers of elementary normal
discs in $H_0$ and $H_0'$ which touch elementary discs of type $D$ are equal.

There are some other further constraints on normal surface vectors. We say that two types of normal disc in a
0-handle of ${\cal H}$ are {\sl incompatible} if there is no normal isotopy
that makes them disjoint. Clearly, a properly embedded normal surface cannot
contain incompatible normal discs. The {\sl compatibility constraints} assert that
for each pair of incompatible normal disc types in a 0-handle of ${\cal H}$, one of the associated variables
is zero.

Haken showed [8] that there is a one-one correspondence between normal surfaces,
up normal isotopy, and non-negative integer solutions to the matching equations
that satisfy the compatibility constraints.

\vskip 6pt
\noindent {\caps 7.3. Summation of normal surfaces}
\vskip 6pt

Given that normal surfaces have such a close connection with linear algebra, it makes
sense to exploit this algebraic structure. Haken therefore introduced the following definition [8].

A properly embedded normal surface $S$ is said to be the {\sl sum} of two other normal
surfaces $S_1$ and $S_2$ if $(S) = (S_1) + (S_2)$. 

In this situation, there is a well-known topological interpretation of this summation. Indeed,
suppose merely that $(S) = v_1 + v_2$ for non-negative integral vectors $v_1$ and $v_2$ satisfying the matching
equations. Then $v_1$ and $v_2$ also satisfy the compatibility constraints, and hence 
correspond to normal surfaces $S_1$ and $S_2$. Consider the elementary normal discs of $S$ in each 0-handle of ${\cal H}$. 
Since $(S) = (S_1) + (S_2)$,
we may partition these discs into two collections, so that the number of elementary
normal discs of each type in the first (respectively, second) collection is equal to the number of elementary normal
discs of that type in $S_1$ (respectively, $S_2$). Since $S_1$ satisfies the
matching equations, it is possible to build $S_1$ starting from these elementary normal
discs by inserting discs into the 1-handles and 2-handles.
Specifically, we consider, for each 1-handle $D^1 \times D^2$, the arcs of $S_1 \cap (\partial D^1 \times D^2)$
of each type. Because $(S_1)$ satisfies the matching equations,
there are the same number of arcs of this type in each component of
$\partial D^1 \times D^2$. Hence, we may insert discs into $D^1 \times D^2$
interpolating between these arcs. We may do this for each disc type in the 1-handles.
We may similarly insert discs into the 2-handles.
The result is not quite a standard surface, because it does not respect the product 
structure on the 1-handles and 2-handles. Specifically, for any 1-handle $D^1 \times D^2$
and each point $p$ in $D^1$,
the intersection between $\{ p \} \times D^2$ is a collection of properly embedded arcs,
but these arcs vary as $p$ varies. See Figure 18. Clearly, one may ambient isotope
this surface to form a normal surface, normally isotopic to $S_1$, but we will not do so.

We may perform the same construction for $S_2$. Then $S_1$ and $S_2$ may intersect in the
1-handles and 2-handles, and after a small isotopy, this intersection will be a collection of
simple closed curves and properly embedded arcs. These are called {\sl trace curves}. 

\vskip 18pt
\centerline{
\epsfxsize=4in
\epsfbox{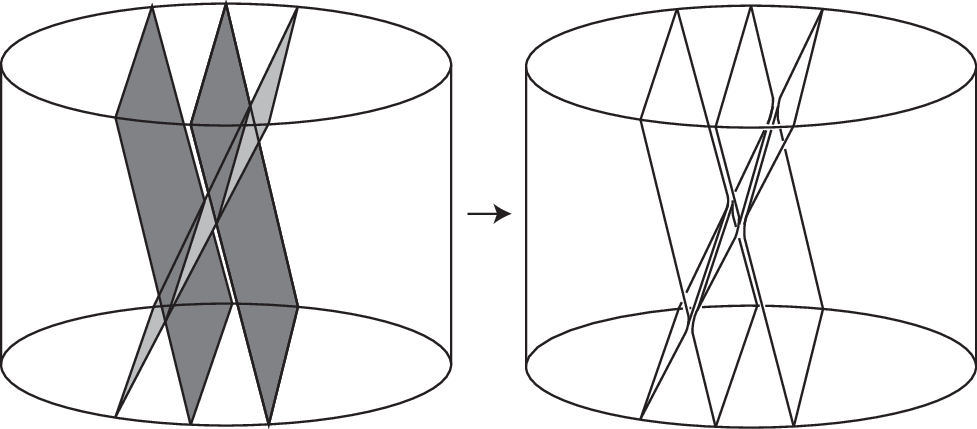}
}
\vskip 6pt
\centerline{Figure 18: Summation of normal surfaces within a 1-handle}

Each of these trace curves may be resolved in two possible ways. One way is called
a {\sl regular switch} if, when we perform the resolution along such a trace curve,
the resulting discs in the 1-handles are still normally isotopic to standard discs. When the trace curves are
all resolved in this way, the result is a surface normally isotopic to $S$. The other way
of resolving a trace curve is called an {\sl irregular switch}. If a regular switch is
performed along all but one trace curve, and an irregular switch is performed there, then the
resulting surface is not normal. In fact, there is an ambient isotopy that reduces its complexity.

If one removes an open regular neighbourhood of the trace curves from $S_1$ and $S_2$,
the components of the resulting compact surface are called {\sl patches}.

The trace curves $S_1 \cap S_2$ give rise to discs, annuli and M\"obius bands properly embedded in the exterior of $S$,
as follows. Since $S$ is obtained from $S_1 \cup S_2$ by resolving the intersections $S_1 \cap S_2$, an
$I$-bundle over $S_1 \cap S_2$ lies in the manifold, with the $\partial I$-bundle forming precisely the intersection with $S$.
This $I$-bundle forms the {\sl trace discs, annuli and M\"obius bands}.

The phrases {\sl regular switch}, {\sl irregular switch}, {\sl patch} and {\sl trace curves} are standard terminology
in normal surface theory (as in Section 1 of [13] and Sections 3.3.5 and 4.1.2 in [23]).

\vskip 6pt
\noindent {\caps 7.4. A vector for boundary-regulated surfaces}
\vskip 6pt

We now develop an analogue of the results in the previous two subsections for
boundary-regulated surfaces. 

For each type of elementary disc $D$ in a 0-handle of ${\cal H}$, we introduce either one or two variables,
depending on whether $D$ is disjoint from the boundary of $M$. When $D$ intersects $\partial M$,
the two different variables corresponding to $D$ relate to the two possible transverse orientations
on $D$.

When $S$ is a boundary-regulated surface, the associated {\sl boundary-regulated vector}
$(S)_{\partial r}$ counts the number of elementary  discs of each type, except
that when elementary discs intersect $\partial M$, their orientations are taken into account,
and therefore one obtains two non-negative integers for this disc type.

The boundary-regulated vector satisfies a collection of equations, which we again term
the {\sl matching equations}. These fall into two types. Equations of the first type are
just the classical matching equations, viewing $S$ as an unoriented standard surface. 
Equations of the second type are concerned with the transverse orientations on the intersection arcs between $\partial S$ and the handles.
For each component of ${\cal H}^1 \cap \partial M$, the curves $\partial S$ intersect
this component in a collection of parallel arcs. These arcs have two possible transverse
orientations. So, for this component of ${\cal H}^1 \cap \partial M$, we obtain two
equations, one for each of these two transverse orientations, as follows. In the two
adjacent 0-handles of ${\cal H}$, the elementary discs of $S$ that are incident to this
component of ${\cal H}^1 \cap \partial M$ are transversely oriented.
So, the two equations assert that the number of arcs with some transverse orientation,
arising from elementary discs in each of the two adjacent 0-handles, agree.

There are also {\sl compatibility constraints}, which are defined as follows. Again, these come in two types. The constraints of the first type
are just like those from classical normal surface theory. They assert that, for
each pair of  elementary disc types in a 0-handle, only one of these two types can occur
if two discs of these types inevitably intersect. The second type of compatibility constraint
is concerned with transverse orientations. For each 0-handle of ${\cal H}$ and for each pair of transversely oriented elementary disc types
that intersect $\partial M$ and that lie within that 0-handle, the co-ordinate of one of these is forced to be zero if
any two disjoint discs of these types inevitably give rise to a tubing arc $\alpha$, as in Condition 2 of Section 5.2.

%\noindent {\bf Proposition.} {\sl There is a one-one correspondence between 
%boundary-regulated normal surfaces, up to normal isotopy, and non-negative
%integer solutions to the boundary-regulated matching equations satisying
%the compatibility constraints.}

Let $v$ be a vector with non-negative integer entries satisfying the boundary-regulated
matching equations and compatibility constraints. Then we may construct a boundary-regulated
surface $S$ with $(S)_{\partial r} = v$, as follows. 

For each elementary disc type in the 0-handles, insert into the 0-handle as many elementary discs of that
type as specified by the vector $v$. At this stage, we do not attempt to orient them. Since the vector
$v$ satisfies the first type of compatibility constraint, there is a normal isotopy that makes any two
of these elementary discs disjoint. Moreover, all of the discs may be made disjoint simultaneously, and there is a unique
position, up to normal isotopy, for their union. We now impose the transverse orientations on the elementary
discs in the 0-handles that intersect $\partial M$, as specified by the vector $v$. For each unoriented disc type,
there is at most one way of doing this without creating a tubing arc. More precisely, there can be at most one switch
of transverse orientations between parallel discs of the same type. Moreover, the switch can occur in just one way.
For example, suppose that the discs of type $D$ run over $R_-(M) \cap {\cal H}^1$. Then, at a switch of orientations
between two discs of this type, their transverse orientations point away from each other. Since $v$ satisfies the
compatibility constraints, there is therefore a unique way of imposing these transverse orientations without
creating a tubing arc.

We now insert the elementary discs into the 1-handles. The method of doing this is described in Section 7.3.
Since $v$ satisfies the unoriented matching equations, one may insert unoriented discs into each 1-handle $D^1 \times D^2$
to interpolate between the elementary discs in the incident 0-handles. As in Section 7.3, each of these discs
intersects $\{ p \} \times D^2$ in a properly embedded arc, for each point $p$ in $D^1$. If any of these elementary
discs in a 1-handle intersects $\partial M$, then a transverse orientation may be imposed upon it, so that, near $\partial M$,
this transverse orientation agrees with the
transverse orientations on the two elementary discs in the 0-handles to which it is incident. This is because 
$v$ satisfies the second type of equation, which is concerned with the transverse orientation on arcs
in ${\cal H}^1 \cap \partial M$.

Finally, elementary discs may be inserted into the 2-handles. For each 2-handle $D^2 \times D^1$,
the surface that has been constructed so far intersects the 2-handle in a collection of circles, each isotopic
to $D^2 \times \{ p \}$ for some point $p$ in $D^1$. Hence, these curves may be extended to properly embedded
disjoint discs in the 2-handle. No transverse orientations are imposed upon these discs, because they are
disjoint from $\partial M$.

This procedure creates a surface $S$. As explained in Section 7.3, this is not yet boundary-regulated, simply
because it is not standard. However, there is a normal isotopy, after which $S$ respects the product structures
of the 1-handles and 2-handles, and is therefore a standard surface. This is the boundary-regulated surface
that is required. We have therefore proved the following lemma.

\noindent {\bf Lemma 7.2.} {\sl Let $v$ be a non-negative integral solution to the boundary-regulated equations,
satisfying the compatibility constraints. Then there is some properly embedded, boundary-regulated surface
$S$ such that $(S)_{\partial r} = v$.}

\vfill\eject
\noindent {\caps 7.5. Summation of boundary-regulated surfaces}
\vskip 6pt

In this subsection, we consider the case where $S$ is a boundary-regulated surface, and its vector
$(S)_{\partial r}$ can be expressed as a sum $v_1 + v_2$, where $v_1$ and $v_2$ are non-negative
and integral and satisfy the boundary-regulated matching equations. Our aim is first to show that
there are boundary regulated surfaces $S_1$ and $S_2$ such that $v_1 = (S_1)_{\partial r}$ and
$v_2 = (S_2)_{\partial r}$. Then we will show how this summation may be interpreted topologically.

Note that $v_1$ and $v_2$ also satisfy the boundary-regulated compatibility constraints, for the following reasons.
If $v_1$, say, contains two disc types in a 0-handle that inevitably intersect, then the same is true for $(S)_{\partial r}$,
which is impossible. Suppose that there is a tubing arc $\alpha$ joining two elementary discs that are present in $v_1$.
Then these discs are also present in $S$. So, some subset of $\alpha$ forms a tubing arc for $S$,
which contradicts the fact that it is boundary-regulated. Hence, by Lemma 7.2, $v_1$ and $v_2$ do correspond to 
boundary-regulated surfaces $S_1$ and $S_2$. Of course,
$(S)_{\partial r} = (S_1)_{\partial r} + (S_2)_{\partial r}$.

We wish to develop a geometric interpretation for this summation. Again
consider the elementary discs of $S$ in each 0-handle. Since $(S)_{\partial r} = (S_1)_{\partial r} + (S_2)_{\partial r}$,
we may partition these discs into two collections, so that the number of elementary
discs of each type in the first (respectively, second) collection is equal to the number of elementary normal
discs of that type in $S_1$ (respectively, $S_2$). Since $S_1$ satisfies the
matching equations, it is possible to build $S_1$ starting from these elementary 
discs by inserting discs into the 1-handles and 2-handles. We build $S_2$ similarly.

Note that, exactly as in Section 7.3, $S$ is obtained, up to normal isotopy, by cutting $S_1 \cup S_2$ along $S_1 \cap S_2$
and then resolving these intersections. We borrow the terminology from the normal case, by using the
phrases {\sl regular switch}, {\sl irregular switch}, {\sl trace curves} and {\sl patches} in the obvious way.

Note that the boundary-regulated matching equations take account of the transverse orientations of the elementary discs
that intersect $\partial M$. Hence, we deduce that, in each component $D$ of $\partial M \cap {\cal H}^1$,
the arcs $D \cap S_1$  and $D \cap S_2$ are compatibly oriented whenever they intersect. (See Figure 19). 
We therefore deduce the following result.

\vskip 6pt
\centerline{
\epsfxsize=3.6in
\epsfbox{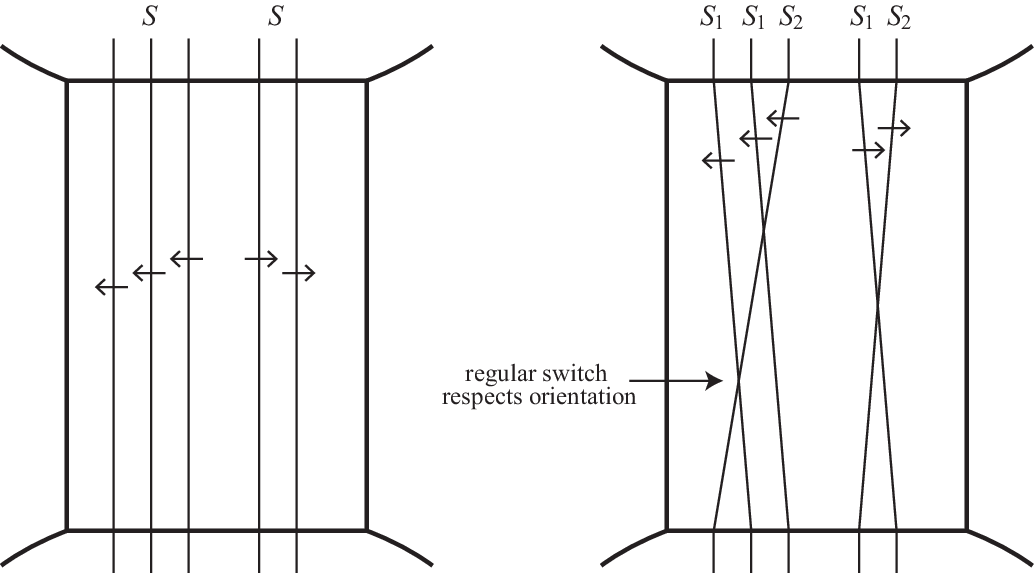}
}
\vskip 6pt
\centerline{Figure 19: A component of ${\cal H}^1 \cap \partial M$}

\noindent {\bf Lemma 7.3.} {\sl Let $S$ be a properly embedded, boundary-regulated surface. 
Let $(S)_{\partial r}$ be its boundary-regulated vector.
Suppose that $(S)_{\partial r} = (S_1)_{\partial r} + (S_2)_{\partial r}$ for boundary-regulated
surfaces $S_1$ and $S_2$. Then for each properly embedded arc of $S_1 \cap S_2$, 
the regular switch along this arc respects the transverse orientations of $\partial S_1$ and $\partial S_2$ near that arc.}

Note that the lemma only refers to {\sl arcs} of $S_1 \cap S_2$, not simple closed curves.

A properly embedded boundary-regulated surface $S$ is {\sl fundamental} if it cannot
be written as a sum of two non-empty boundary-regulated surfaces.

\vfill\eject
\noindent{\caps 7.6. Avoiding trace annuli with a trivial boundary curve}
\vskip 6pt

A large part of Section 4 was concerned with decompositions along annuli disjoint from the sutures.
One of the reasons for this is that trace annuli are of this form. However, it will be important
that these annuli have no trivial boundary curves. In the next few sections, we will prove the following proposition,
which implies that trivial boundary curves may be avoided.

\noindent {\bf Proposition 7.4.} {\sl Let ${\cal H}$ be a handle decomposition of a compact orientable
sutured manifold $(M, \gamma)$. Let $S$ be a connected regulated surface properly embedded in $M$. Suppose that
$S$ extends to an allowable hierarchy and is not product-separating. Assume that at least one of the following holds:
\item{(i)} $S$ is non-separating in the component of $M$ that contains it, or
\item{(ii)} $(M, \gamma)$ contains a non-separating product disc or a non-separating product annulus.

\noindent Suppose also that $S$ has smallest
possible complexity among connected standard surfaces with these properties. Let $(M_S, \gamma_S)$ be the result of decomposing
along $S$. Suppose that $S$ 
is a sum of non-empty boundary-regulated surfaces $S_1$ and $S_2$. Suppose also that $S_1 \cap S_2$ has fewest components,
among all ways of expressing $S$ as such a sum. Then no trace annulus for this summation has a boundary curve that is trivial in $(M_S, \gamma_S)$.}

\vskip 6pt
\noindent {\caps 7.7. Generalised summation}
\vskip 6pt

We say that a properly embedded surface $S$ in an orientable 3-manifold is a {\sl generalised sum} of two properly embedded surfaces 
$S_1$ and $S_2$ if $S_1$ and $S_2$ are in general position, and hence $S_1 \cap S_2$ is a collection of properly embedded
arcs and simple closed curves, and $S$ is obtained from $S_1 \cup S_2$ by resolving these intersections in some way.
We write $S = S_1 + S_2$, although there is more than one surface that can be obtained from $S_1$ and $S_2$ by
generalised summation.

Generalised summation arises in several possible ways. Summation of normal surfaces, as described in Section 7.3,
is a type of generalised summation. Similarly, summation of boundary-regulated surfaces is also. Another natural
method of performing generalised summation is when $S_1$ and $S_2$ are oriented and the resolution of $S_1 \cap S_2$
occurs in the way that respects this orientation.

We again borrow terminology from normal surface theory, by speaking of patches, trace curves, and trace discs,
annuli and M\"obius bands. A trace curve is {\sl two-sided} if the $I$-bundle over it with $(\partial I)$-bundle in $S$
is orientable; otherwise it is {\sl one-sided}. Thus, trace annuli and trace discs arise from two-sided
trace curves, and trace M\"obius bands arise from one-sided trace curves.

Consider a trace arc or curve $a$ of $S_1 \cap S_2$. Suppose the four patches incident to $a$ are all distinct.
We say that two of these patches are {\sl opposing along $a$} if one lies in $S_1$, the other lies in $S_2$, and the two
copies of these patches in $S$ do not intersect along $a$. Thus, the four patches incident to $a$ decompose into
two opposing pairs. (See Figure 20.)

\vskip 12pt
\centerline{
\epsfxsize=4in
\epsfbox{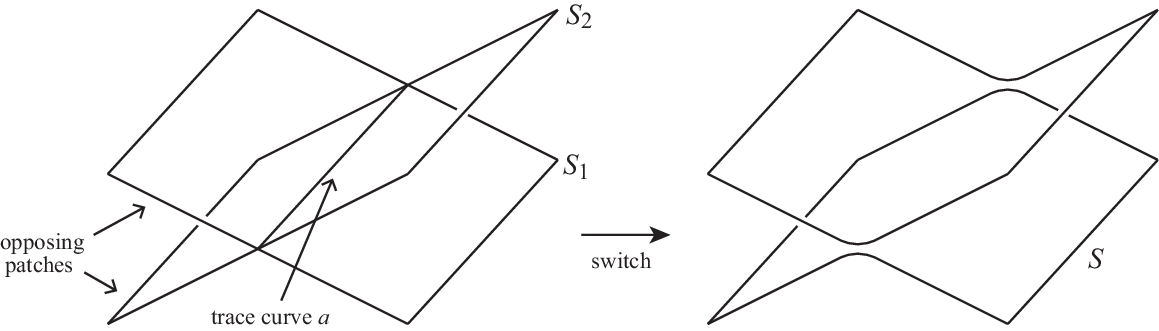}
}
\vskip 6pt
\centerline{Figure 20: Patches opposing along a trace curve $a$}

\vfill\eject
\noindent {\caps 7.8. Reducible summation}
\vskip 6pt

The summation $S = S_1 + S_2$ is {\sl reducible} if there are subsurfaces $P_1 \subset S_1$ and $P_2 \subset S_2$
which are each a union of patches, such that no two patches of $P_1$ are incident in $S_1$, no two patches of $P_2$
are incident in $S_2$, and along every arc and curve of $\partial P_1 - \partial S_1$ and along every arc and curve
of $\partial P_2 - \partial S_2$, $P_1$ and $P_2$ are opposed.

Suppose that the summation is reducible. Then one may attach $P_2$ to $S_1 - P_1$, and attach $P_1$ to $S_2 - P_2$,
giving two surfaces $S'_1$ and $S'_2$. Then $S = S'_1 + S'_2$, but $|S_1' \cap S'_2| < |S_1 \cap S_2|$.

Note that if the summation $S_1 + S_2$ arose from the summation of normal or boundary-regulated surfaces, then
so does the summation $S'_1 + S'_2$. So, if one considers an expression of $S$ as a normal or boundary-regulated
sum $S_1 + S_2$ of non-empty surfaces, and $|S_1 \cap S_2|$ is minimised, then the summation is not reducible.

\vskip 6pt
\noindent {\caps 7.9. Alternative summands}
\vskip 6pt

Let $S$ be a generalised sum $S_1 + S_2$. Then an {\sl alternative summand} for $S$ is a surface
$F$ that is a union of patches of $S_1 \cup S_2$ such that, for each trace arc and curve of $S_1 \cap S_2$,
exactly two of the four patches emanating from this curve lie in $F$, and these are not opposing.

The union of the patches that do not lie in $F$ forms another surface $F'$. Note that $S$ is also the
generalised sum $F + F'$. Moreover, each of the trace arcs and curves for $F+F'$ is a trace
arc or curve for $S_1 + S_2$, and the way that these are resolved is the same. Hence if $S = S_1 + S_2$
is a normal (or boundary-regulated) sum of normal (or boundary-regulated) surfaces, then $F$ and $F'$
are normal (or boundary-regulated) and the summation is normal (or boundary-regulated).

\vskip 6pt
\noindent {\caps 7.10. Irregular switches}
\vskip 6pt

Let $S$ be a generalised sum $S_1 + S_2$. Then a surface is obtained from $S_1 \cup S_2$ by
making some {\sl irregular switches} if it is obtained from $S_1 \cup S_2$ by resolving the
intersections in some way, and at least one of these arcs or curves, the way that it resolved
is different from that of $S$.

\vskip 6pt
\noindent {\caps 7.11. Tubing and summation}
\vskip 6pt

Let $S$ be a generalised sum $S_1 + S_2$. Suppose that $S_1$ and $S_2$ have been given a transverse
orientation in a regular neighbourhood of $\partial M$ and that the summation respects this orientation.
Then $S$ inherits a well-defined transverse orientation in a regular neighbourhood of $\partial M$.
The main case where this hypothesis holds is in the case of summation of boundary-regulated surfaces.

Let $\alpha$ be the closure of an arc component of
$\partial S_1 - \partial S_2$. Suppose that $\alpha$ is a tubing arc for $S_2$. Thus, $\alpha$ is disjoint from the sutures $\gamma$,
and the transverse orientations of $S_2$ at the endpoints of $\alpha$ and the transverse orientation of $R_\pm(M)$
near $\alpha$ all agree. 

Let $S_2^+$ be the result of attaching a tube to $S_2$ along $\alpha$. Then $S_1 \cap S_2^+$ is obtained
from $S_1 \cap S_2$ by modifying it near $\alpha$, in the following way. The arcs of $S_1 \cap S_2$ at the endpoints
of $\alpha$ are truncated just before $\partial M$, and then a parallel copy of $\alpha$ is attached to them.

One may perform a corresponding
sum $S_1 + S_2^+$. For the arcs and curves of $S_1 \cap S_2^+$ that are equal to an arc or
curve of $S_1 \cap S_2$, the way that these are resolved is unchanged. For the new
arc or curve of $S_1 \cap S_2^+$, the resolution is the one that agrees with 
the resolution of the arcs of $S_1 \cap S_2$ at the endpoints of $\alpha$.
The resulting surface $S^+$, which is $S_1 + S_2^+$, is obtained from $S$ by attaching a tube that runs parallel
to $\alpha$.

\vskip 18pt
\centerline{
\epsfxsize=6in
\epsfbox{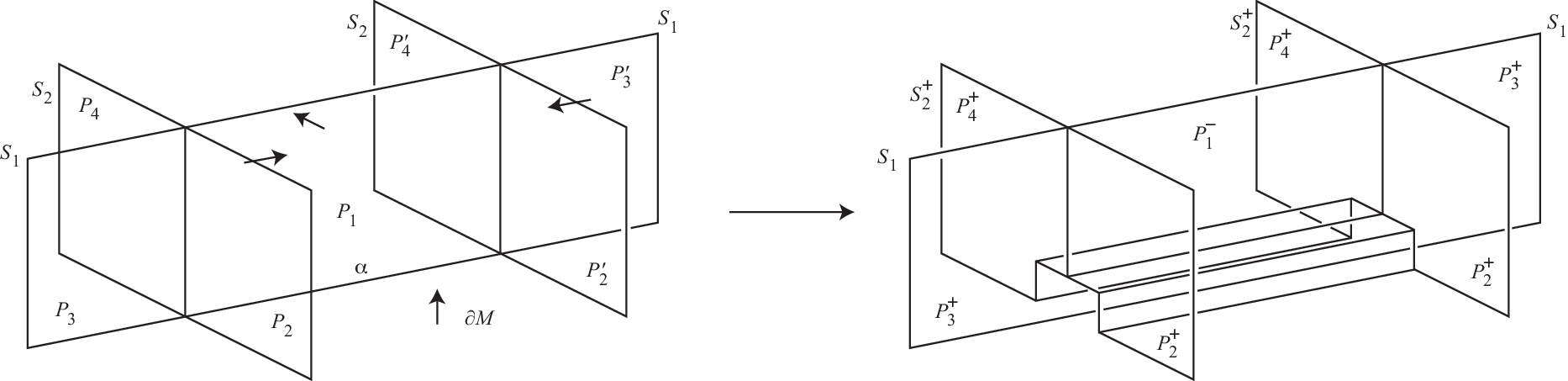}
}
\vskip 6pt
\centerline{Figure 21: Attaching a tube and summation}

Note that $S_1$ and $S_2^+$ also have a well-defined transverse orientation in a regular neighbourhood
of $\partial M$ and the summation $S_1 + S_2^+$ respects this orientation.

\noindent {\bf Lemma 7.5.} {\sl If the summation $S_1 + S_2^+$ is reducible, then so was
the summation $S_1 + S_2$.}

\noindent {\sl Proof.} Figure 21 shows the patches of $S_1 + S_2$ near $\alpha$.
These are labelled $P_1$, $P_2$, $P_3$, $P_4$, $P_2'$, $P_3'$ and $P_4'$.
It is possible that some of these patches may actually be the same patch.
Up to changing all transverse orientations near $\partial M$ and up to symmetry, we may assume that the transverse
orientations are also as shown in Figure 21. Thus, the following form opposing pairs:
$P_1$ and $P_2$; $P_3$ and $P_4$; $P_1$ and $P_2'$; $P'_3$ and $P'_4$.

The patches of $S_1 + S_2^+$ can also be read off from Figure 21. One patch $P_3^+$ is obtained from $P_3 \cup P_3'$
by attaching a band along $\alpha$. This opposes a similar patch $P_4^+$ obtained from $P_4 \cup P_4'$. There is a
patch $P_1^-$ obtained from $P_1$ by removing a regular neighbourhood of $\alpha$. This opposes
a patch $P_2^+$ obtained from $P_2 \cup P_2'$ obtained by attach a band along $\alpha$. 

Suppose that $S_1 + S_2^+$ is reducible. If none of the patches near $\alpha$ form the specified
subsurfaces of $S_1$ and $S_2^+$, then the corresponding subsurfaces of $S_1$ and $S_2$
give that $S_1 + S_2$ is reducible. If the patch $P_2^+$ forms part of one of the subsurfaces,
then $P_1^-$ must form part of the other subsurface, and then we obtain similar subsurfaces
of $S_1$ and $S_2$ giving reducibility there. A similar argument holds when $P_3^+$
and $P_4^+$ form part of the subsurfaces. $\square$

\noindent {\bf Lemma 7.6.} {\sl If the surface $S^+$ has an alternative summand $F^+$, then
$S$ has a corresponding alternative summand $F$. Moreover, $F$ is either isotopic to $F^+$ or is obtained from $F^+$
by performing a boundary-compression along a product disc disjoint from $\gamma$.}

\noindent {\sl Proof.} Consider the right-hand diagram in Figure 21. Along the arc or curve of $S_1 \cap S_2^+$ shown,
exactly two patches lie in $F^+$. The possibilities are $P_1^-$ and $P_3^+$; $P_2^+$ and $P_4^+$; $P_2^+$ and $P_3^+$;
$P_1^-$ and $P_4^+$. In each case, we may define $F$ to be the union of the corresponding patches for $S_1 \cup S_2$.
For example, in the former case, we take $P_1$, $P_3$ and $P_3'$ together with all the other patches of $F^+$.
In this case, $F$ is equal to $F^+$. But in two of the other three cases, $F$ is obtained from $F^+$ by boundary-compressing
along the obvious product disc disjoint from $\gamma$. $\square$

A very similar argument gives the following result. Its proof is omitted.

\noindent {\bf Lemma 7.7.} {\sl If the surface $T^+$ is obtained from $S^+$ by performing irregular switches,
then there is a corresponding surface $T$ obtained from $S$ by performing irregular switches. If $T^+$ and $S^+$
can be transversely oriented so that the transverse orientations agree on each patch of $S^+$ and $T^+$ and agree with the given transverse
orientations near $\partial M$, then $T$ and $S$ can also
be transversely oriented so that the transverse orientations agree on each patch of $S$ and $T$ and agree with the given transverse
orientations near $\partial M$.
Moreover, in this case, $T$ is obtained from $T^+$ by boundary-compressing along a product disc
disjoint from $\gamma$.}

\vfill\eject
\noindent {\caps 7.12. Trivial patches}
\vskip 6pt

Let $S = S_1 + S_2$ be a generalised summation. Then a patch $P$ is {\sl trivial}
if it is a planar surface, with one boundary curve forming a trace curve of $S_1 \cap S_2$
and the remaining boundary curves being trivial boundary curves of $S$.

\vskip 12pt
\centerline{
\epsfxsize=3.1in
\epsfbox{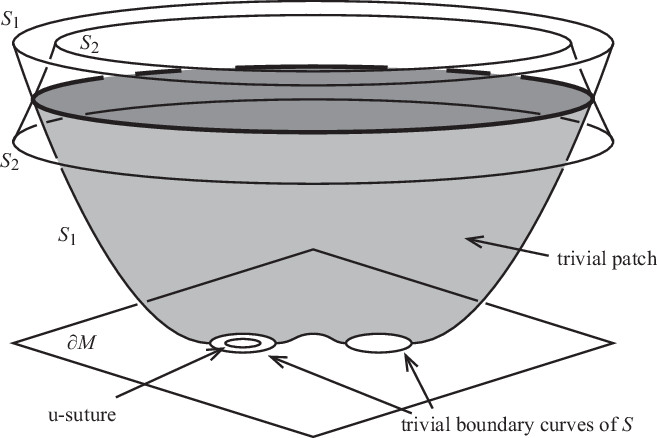}
}
\vskip 6pt
\centerline{Figure 22: Trivial patch}

Suppose that
$$(M, \gamma) \buildrel S \over \longrightarrow (M_S, \gamma_S)$$
is a taut allowable decomposition between decorated sutured manifolds. Suppose also that $S$ is connected.
One would like to say that if there is a trace annulus $A$ with a boundary curve that is trivial in 
$(M_S, \gamma_S)$, then there has to be a trivial patch. However, this
need not obviously be the case. Since the boundary curve is trivial, it bounds a 
disc in $E(M_S, \gamma_S)$. The restriction of this disc to $M_S$ is a planar surface,
and the intersection between this and $S$ is a planar surface $\tilde P$.
The complication is that $\tilde P$ may contain trace curves and trace arcs in its interior.
The possible presence of trace arcs in $\tilde P$ is problematic, because
an `innermost' patch $P$ in $\tilde P$ does not necessarily form a trivial patch.
An example is given in Figure 23.

\vskip 12pt
\centerline{
\epsfxsize=1.8in
\epsfbox{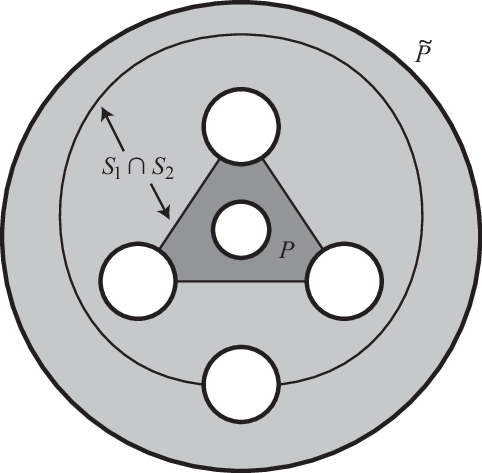}
}
\vskip 6pt
\centerline{Figure 23: The `innermost' planar surface $P$ is not a trivial patch}

\vskip 6pt
\noindent {\caps 7.13. Proof of Proposition 7.4.}
\vskip 6pt

We now embark on this proof. We use the terminology of the proposition. We suppose for contradiction
that $(M_S, \gamma_S)$ contains a trace annulus with a trivial boundary curve.

We are assuming that $S$ is a sum of non-empty boundary-regulated surfaces $S_1$ and $S_2$, and where
$|S_1 \cap S_2|$ is minimal. Therefore, the summation is not reducible.

The plan of the proof is to follow the argument of Lemma 2.1 in [13] where, in the case of normal summation,
disc patches in the interior of $M$ were excluded. However, the difficulties described in Section 7.12 relating
to innermost patches cause complications. To circumvent these, we introduce a new surface $S'$ which will
be a generalised sum $S_1 + S'_2$. The problematic innermost patches, as described in Figure 23, will not arise
in this summation. We can then apply the argument of Lemma 2.1 in [13]. We then use the lemmas in
Section 7.11 to provide information about the original summation $S = S_1 + S_2$. This will lead to
a contradiction.

The new properly embedded surface $S'$ will be obtained from $S$ by attaching tubes
along some arcs of $\partial S_1 - \partial S_2$. These will be tubing
arcs. Because of Lemma 7.3, the summation at the endpoints of these arcs satisfies the conditions
of Section 7.11. Hence, this new surface $S'$ will be a generalised summation $S_1 + S'_2$.
By Lemma 7.5, this new summation will not be reducible. By Lemma 4.5, $S'$ will extend to
an allowable hierarchy.
Also, by Lemma 5.7 or Proposition 5.8, $S'$ will not be product-separating.
We denote the manifold obtained by decomposing along $S'$ by $(M_{S'}, \gamma_{S'})$.
Since this manifold is at the start of an allowable hierarchy, it is taut by Corollary 3.7,
its canonical extension is taut by Lemma 3.3, and $E(S')$ is taut by Lemma 3.3.

This surface $S'$ is constructed as follows. Consider any trivial curve $c$ of $\partial S$ that is neither
a component of $\partial S_1$ nor of $\partial S_2$. Thus, $\partial S$ is comprised of arcs of $\partial S_1 - \partial S_2$
and arcs of $\partial S_2 - \partial S_1$ attached to each other. Since $c$ is trivial, it bounds
a disc in $\partial E(M, \gamma)$ disjoint from the sutures. By choosing $c$ appropriately,
we may assume that if any curve of $\partial S$ lies in the interior of this disc, then it is
a component of $\partial S_1$ or $\partial S_2$. This implies that every arc of $\partial S_1 - \partial S_2$
and $\partial S_2 - \partial S_1$ in $c$ is a tubing arc. Pick all such arcs of $\partial S_1 - \partial S_2$
in $c$, and attach a tube to $S_2$ along them. The new surface is a generalised sum.
Continuing in this way as far as possible, we end with a surface $S'$ which is a generalised sum
of surfaces $S_1$ and $S'_2$. These have the property that any trivial curve of $\partial S'$
is a component of $\partial S_1$ or $\partial S'_2$.

This produces new trace curves, all of which are two-sided, and hence new trace annuli. Some
of their boundary curves may be trivial in $(M_{S'}, \gamma_{S'})$. But the trivial boundary curves in $(M_S, \gamma_S)$
remain trivial boundary curves in $(M_{S'}, \gamma_{S'})$, by Lemma 4.6.
So, there is certainly at least one trace annulus with a trivial boundary curve in $(M_{S'}, \gamma_{S'})$.
This boundary curve bounds a disc in $E(M_{S'}, \gamma_{S'})$. The restriction of this disc to $M_{S'}$ is a planar surface,
and the intersection between this and $S'$ is a planar surface $\tilde P$. This is a union of patches.
Any boundary component of $\tilde P$ other than the one lying in the trace annulus
is a trivial curve of $\partial S'$. Hence, it lies entirely in $\partial S_1$
or $\partial S'_2$. Therefore, there are only trace curves in $\tilde P$, and no
trace arcs. Consider one that is innermost in $\tilde P$. This encloses
a patch $P_1$. This extends to a disc $D_1$ in $\partial E(M_{S'}, \gamma_{S'})$
disjoint from the sutures. Now consider the trace annulus $A$ incident to $P_1$. The other boundary component of $A$ is 
also a trivial curve in $\partial M_{S'}$. This is because $E(M_{S'}, \gamma_{S'})$ is
taut by Lemma 3.3. Hence, this component of $\partial A$ also bounds a disc $\tilde D_1$ in $\partial E(M_{S'}, \gamma_{S'})$ disjoint
from the sutures. The intersection between this disc and $S'$ is a planar surface $\tilde P_1$
that is a union of patches. We say that $\tilde P_1$ is {\sl associated} with the patch $P_1$.

There are several cases to consider. Suppose first that the patch in $\tilde P_1$ incident to $\partial A$
is not opposing $P_1$ along the trace curve. If $\tilde D_1$ is disjoint from $D_1$, then $\tilde D_1 \cup D_1$
forms a 2-sphere, which, after a small isotopy, can be made disjoint from 
$E(S')$ and which separates $E(S')$. This contradicts the fact $E(M_{S'})$ is irreducible
or that $E(S')$ is incompressible or that $E(S')$ has no 2-sphere components. So, suppose now that $\tilde D_1$ and $D_1$
are not disjoint. Since $P_1$ is a patch, and $\tilde P_1$ is made from a union of patches, we deduce that $P_1 \subset \tilde P_1$.
In this case, we remove $\tilde P_1$ from $S'$ and replace it by $P_1$. Let $S''$ denote the new surface.
This is obtained from $S_1 \cup S'_2$ by performing an irregular switch along $\partial D_1$, giving
a surface $F'$ and then removing the punctured torus formed from $\tilde P_1 - P_1$. Note that $S''$
is obtained from $S'$ by slicing under discs of contact and then performing an isotopy.
Hence, $S''$ also extends to an allowable hierarchy by Lemma 4.3. If $S$ was non-separating,
then so is $S''$. The surface $F'$
can be transversely oriented so that the patches of $F'$ and the patches of $S'$ have the
same transverse orientations. So, by Lemma 7.7, there is a corresponding surface $F$ obtained
from $S_1 \cup S_2$ by performing an irregular switch. It is isotopic to a standard
surface with complexity strictly less than that of $S$. By Lemma 7.7, this is obtained
from $F'$ by possibly performing some boundary compressions along product discs disjoint
from $\gamma$. The surface $S''$ is obtained from $F'$ by removing a component.
This component corresponds to a component or components of $F$. Remove these
to give a surface $S'''$. Then $S'''$ is obtained from $S''$ by possibly performing
boundary compressions along product discs disjoint from $\gamma$. Hence, by Proposition 4.7,
$S'''$ extends to an allowable hierarchy. Since it is a union of components
of $F$, its complexity is at most that of $F$, and this is (after an isotopy) strictly
less than that of $S$. If $S''$ was non-separating, then so is at least one component of
$S'''$. Thus, we contradict our assumption that $S$ had minimal complexity.

Now consider the case where $\tilde P_1$ and $P_1$ are opposing. Then $\tilde P_1$ cannot be a single patch
because the summation $S_1 + S'_2$ would then be reducible, contradicting Lemma 7.5. Hence, $\tilde P_1$ is a union of patches.
Again, we may find one that is innermost in $\tilde P_1$. Call this patch $P_2$. Since there are not trace
arcs in $\tilde P_1$, $P_2$ is a trivial patch. Continuing in this way, we obtain
trivial patches $P_1, \dots, P_n$ and associated planar surfaces $\tilde P_1, \dots, \tilde P_n$. We stop the process
when $\tilde P_n$ contains one of the previous trivial patches in its interior. By relabelling, we may
assume that this trivial patch is $P_1$.

We can view the surface $\tilde P_i - P_{i+1}$ (where the index $i$ is modulo $n$) as a punctured annulus.
As explained in the proof of Lemma 2.1 of [13], these punctured annuli combine to form a surface $T'$,
which is a punctured torus. This is an alternative summand for $S'$.
Thus $S' = T' + F'$, for some properly embedded surface $F'$.
We can now fill in each of these punctures of $T'$, by attaching discs of contact,
and the result is a torus $\overline{T}'$. We can fill in the corresponding discs
of contact for $S'$, giving a surface $\overline{S}'$, which also extends to an allowable hierarchy by Lemma 4.3.
It is not product-separating by Proposition 5.8.
Then $\overline{S}' = \overline{T}' + F'$. Hence, $F'$ is in fact ambient isotopic to $\overline{S}'$.
Therefore, $F'$ also extends to an allowable hierarchy with the same reduced length and is
not product-separating.
Since $F'$ is an alternative summand for $S'$, Lemma 7.6 implies that $S$ also has a corresponding
alternative summand $F$. By Lemma 7.6,  this is obtained from $F'$ by performing boundary-compressions
along product discs disjoint form $\gamma$. So by Proposition 4.7, $F$ also extends to an allowable hierarchy.
By Proposition 5.8, $F$ is not product-separating. If $S$ was non-separating, so too is at least one
component of $F$. Now the summation $S_1 + S_2$ was boundary-regulated,
and so the alternative summand $F$ is also boundary-regulated. Since it is a boundary-regulated summand,
its complexity is strictly less than that of $S$. But this contradicts the assumption that $S$
had minimal complexity. $\square$

The surface $S'$ constructed in the above proof will be useful to us later, and so we record some of its properties now.

\noindent {\bf Addendum 7.8.} {\sl Let $S$, $S_1$ and $S_2$ be as in Proposition 7.4. Then there is a surface $S'$ properly embedded
in $(M, \gamma)$ with the following properties:
\item{(i)} $S'$ is a generalised sum of $S_1$ and a surface $S'_2$;
\item{(ii)} when $S_1$ and $S_2$ have transverse orientations and $S$ is the oriented double-curve sum of $S_1$ and
$S_2$, then $S'_2$ also has a transverse orientation and $S'$ is the oriented double-curve sum of $S_1$ and $S'_2$;
\item{(iii)} no trace annulus for this summation has a trivial boundary curve in the sutured manifold $(M_{S'}, \gamma_{S'})$ that is
obtained by decomposing along $S'$;
\item{(iv)} $S'$ is obtained from $S$ by tubing along arcs;
\item{(v)} every trivial curve of $\partial S'$ is a component of $\partial S_1$ or $\partial S'_2$.

}

\vskip 6pt
\noindent {\caps 7.14. Decomposition along fundamental surfaces}
\vskip 6pt

We now come to a central result of this paper.

\noindent {\bf Theorem 7.9.} {\sl Let $(M, \gamma)$ be a taut decorated sutured manifold with a positive handle structure ${\cal H}$.
Suppose that $(M, \gamma)$ admits an allowable hierarchy.  Suppose also that $E(M, \gamma)$ is atoroidal, 
and that no component of $E(M, \gamma)$ is a Seifert fibre space other than a solid torus or a copy of $T^2 \times I$.
Suppose that no component of $(M, \gamma)$ has boundary a single torus with no sutures. 
In addition, suppose that $(M, \gamma)$ is not a product sutured manifold.
Then it admits a taut allowable sutured manifold decomposition along a regulated surface $S$, that is
fundamental as a boundary-regulated surface, that extends to an allowable
hierarchy, and such that decomposition along $S$ reduces the complexity of the handle structure.}

\noindent {\sl Proof.} The proof divides into two cases.

\noindent {\sl Case 1.} $(M, \gamma)$ admits an allowable hierarchy where the first surface $S$ is connected and non-separating.

We may assume that $S$ is incompressible
by Lemma 4.2, and so by Lemma 4.9 in [21], it can be placed in standard form.
By Theorem 5.5, we may further assume that $S$ is regulated, as long as we drop
the requirement that $S$ is incompressible. We pick $S$ to have minimal
complexity among regulated non-separating surfaces that extend to an 
allowable hierarchy. Let $(M_S, \gamma_S)$ be obtained from $(M, \gamma)$ by 
decomposing along $S$.

Suppose that $S$ is not fundamental as a boundary-regulated surface. Then it can be written as a sum
of non-empty boundary-regulated surfaces $S_1$ and $S_2$. 
We may perform a normal isotopy to $S_1$ and $S_2$, leaving $\gamma$ invariant, so that they intersect
in a collection of simple closed curves and properly embedded arcs. We may assume that $S_1$
and $S_2$ have been chosen so that the number of these curves and arcs is minimal. 

Each patch of this summation inherits a transverse orientation from $S$.

\noindent {\sl Case 1A.} These orientations induce transverse orientations of $S_1$ and $S_2$.

Then $S$ is the oriented double-curve sum of $S_1$ and $S_2$. Since the summation respects the
transverse orientations on $S_1$ and $S_2$, we deduce that $[S, \partial S] = [S_1, \partial S_1] + [S_2, \partial S_2]
\in H_2(M, \partial M)$. Therefore, at least one of $[S_1, \partial S_1]$ and $[S_2, \partial S_2]$ is non-trivial,
say the former. We will show that $S_1$ extends to an allowable hierarchy. 
Since $S_1$ is oriented and boundary-regulated, it is regulated. But we have
then found a non-separating, regulated surface $S_1$ that extends to an allowable
hierarchy, but with smaller complexity than that of $S$. This contradicts our
minimality assumption.

Note first that $S_1$ is connected. For if it is not, then it has a homologically non-trivial component
$F$. We can then instead consider the summation $S = F + (S_2 + (S_1 - F))$,
where $S_2 + (S_1 - F)$ is the oriented double-curve sum of $S_2$ and $S_1 - F$. There are fewer curves of intersection between 
$F$ and $S_2 + (S_1 - F)$ than between $S_1$ and $S_2$, contradicting our minimality assumption.

By Addendum 7.8, we have a surface $S'$ as in the statement there. This is the oriented double-curve
sum of $S_1$ and a surface $S'_2$. We claim that the trivial boundary curves of $S'$ are precisely
the trivial curves of $\partial S_1$ and $\partial S'_2$. Certainly, every trivial curve of $\partial S'$ arises
as a trivial curve of $\partial S_1$ or $\partial S'_2$, by (v) of Addendum 7.8. To prove the claim,
we must show that every trivial curve of $\partial S_1$ and $\partial S'_2$ is preserved as a curve of
$\partial S'$. We will show this for curves of $\partial S_1$, as the other case is similar.
Consider a trivial curve $C$ of $\partial S_1$ and its trivialising planar surface $F$.
The intersection between $F$ and $\partial S'_2$ is a collection of simple closed curves and arcs.
Consider just the arcs, which all start and end on $C$, because the remaining components
of $\partial F$ are u-sutures of $\gamma$. Transversely orient the arcs according to the transverse
orientation of $S'_2$. We may construct a tree, where each vertex corresponds to a
component of the complement of these arcs. This has a sink vertex and a source vertex.
One of these corresponds to a trivialising planar surface for a component of $\partial S'$.
This contradicts (v), and so the claim is proved.

Let $P$ be the trace discs, annuli and M\"obius bands arising from the summation $S' = S_1 + S'_2$. In fact,
because the summation respects the transverse orientations of $S_1$ and $S'_2$,
$P$ has no M\"obius band components, because incident to each component of $P$,
there are parts of $S'$ pointing away from $P$ and parts of $S'$ pointing towards $P$.
Give $P$ some transverse orientation. These then form a collection of product discs and product annuli.
By Addendum 7.8, the annuli are all non-trivial in $(M_{S'}, \gamma_{S'})$. 
We have a commutative diagram of sutured manifold decompositions
$$

\matrix{
(M, \gamma) & \buildrel S' \over \longrightarrow & (M_{S'}, \gamma_{S'}) \cr
\Big\downarrow {\scriptstyle S_1} && \Big\downarrow {\scriptstyle P} \cr
(M_1, \gamma_1) & \buildrel S'_2 - {\rm int}(N(S_1)) \over \longrightarrow & (M_{12}, \gamma_{12}). \cr}
$$

We claim that any trivial boundary curve $C$ of $S'_2 - {\rm int}(N(S_1))$ in $M_1$ is disjoint from $S_1$.
Suppose not and consider the trivialising planar surface $F$ bounded by $C$. This is divided into subsurfaces
$F \cap R_\pm(M)$ and $F \cap S_1$. Any component of $F \cap R_\pm(M)$ incident to $C$ gives rise to
a trivial boundary curve of $S'$ that is not a component of $\partial S_1$ or $\partial S'_2$.
This contradicts (v) in Addendum 7.8. Thus if $C$ has non-empty intersection with $S_1$,
then it must lie entirely in $S_1$. But it then gives rise to
a trace annulus in $(M_{S'}, \gamma_{S'})$ with a trivial boundary curve, which contradicts Addendum 7.8 (iii).
This proves the claim.

We claim that all the decompositions in the above commutative diagram are allowable. Since $S'$ is obtained from $S$ by tubing along arcs,
decomposition along $S'$ is allowable, by Lemma 4.5. The annuli in $P$ are all allowable.
Moreover, the product discs in $P$ are allowable, because they are disjoint from the u-sutures of $\gamma_{S'}$.
This is because each such product disc lies within a regular neighbourhood of a trace arc and so  
is clearly disjoint from the u-sutures of $\gamma$. It is also disjoint from the trivial
boundary curves of $S'$, as each trivial boundary curve of $S'$ is a component of $\partial S_1$ or $\partial S'_2$.
Note also that $\partial S_1$ is disjoint from the u-sutures of $\gamma$. Also, any trivial boundary curve
of $S_1$ is a trivial boundary curve of $S$, by the above claim. Hence, $S'_2 - {\rm int}(N(S_1))$
is disjoint from these u-sutures of $\gamma_1$. Thus, we have shown that all of these surfaces
are disjoint from the u-sutures. Finally, any trivial boundary curve of $S_1$ or $S'_2 - {\rm int}(N(S_1))$
is a trivial curve of $\partial S'$, and so its trivialising planar surface is correctly oriented.
This proves the claim.

Thus, the above diagram consists of decorated sutured manifolds. We show that the decoration on 
$(M_{12}, \gamma_{12})$ is the same, no matter which way one goes round the digram.
If we decompose along $S'$ then $P$, the u-sutures come from u-sutures of $\gamma$ and
trivial curves of $\partial S'$. Note that decomposition along $P$ does not create any u-sutures
because every annular component of $P$ has non-trivial boundary curves. If we decompose along $S_1$ then 
$S_2' - {\rm int}(N(S_1))$, then the u-sutures come from the u-sutures of $\gamma$
and the trivial boundary curves of $\partial S_1$ and $\partial S'_2$, which are disjoint.
Thus, the u-sutures of $\gamma_{12}$ are the same in both cases.

By (iv) in Addendum 7.8, $S'$ is obtained from $S$ by tubing along arcs, and so by Lemma 4.5, $S'$
extends to an allowable hierarchy. We want to show that $(M_{12}, \gamma_{12})$ admits an
allowable hierarchy. To do this, we will use Lemma 4.4 and Proposition 4.9. We need to verify the
hypotheses of Proposition 4.9. We are assuming that $E(M,\gamma)$ is atoroidal and no component of 
$E(M,\gamma)$ is a Seifert fibre space other than a solid torus  or a copy of $T^2 \times I$.
By Lemma 3.11, $E(M_{S'}, \gamma_{S'})$ also has these properties. Each decomposition along a 
product annulus in $P$ cannot create a solid torus disjoint from the sutures, since each product 
annulus has one boundary component in $R_-(M_{S'}, \gamma_{S'})$ and one boundary component in
$R_+(M_{S'}, \gamma_{S'})$. So, by Proposition 4.9 and Lemma 4.4, $(M_{12}, \gamma_{12})$ admits an
allowable hierarchy.

However, this does not yet prove that $S_1$ extends to an allowable hierarchy,  because
it is not clear that $S_1$ and $S'_2 - {\rm int}(N(S_1))$ can be part of an allowable hierarchy,
since they may have components that are planar, disjoint from the sutures and where all
but at most one boundary curve is trivial. In fact, $S_1$ cannot have such a component.
Since $S_1$ is connected, this component would be all of $S_1$.
Its canonical extension would have to be a boundary parallel disc in $E(M, \gamma)$
by the tautness of $E(M, \gamma)$. Hence, it would be homologically trivial in $M$,
which is a contradiction.

Now, $S'_2  - {\rm int}(N(S_1))$ may possibly have a planar component $F$ disjoint from
the sutures and where all but at most one boundary curve is trivial. We now show how
to remove such a component $F$.

Since $(M_{12}, \gamma_{12})$ admits an allowable hierarchy, $E(M_{12}, \gamma_{12})$ is taut by Lemma 3.3.
We have a decomposition
$$E(M_1, \gamma_1) \buildrel E(S'_2 - {\rm int}(N(S_1))) \over \longrightarrow E(M_{12}, \gamma_{12}),$$
which we will show is taut using Theorem 2.3. No boundary component of $E(S'_2 - {\rm int}(N(S_1)))$
bounds a disc in $E(M_1, \gamma_1)$ disjoint from the sutures. Also, no component of $E(M_1, \gamma_1)$ can be
a solid torus with no sutures. This is because this would be a component of $(M_1, \gamma_1)$,
and this would imply that $(M, \gamma)$ is not taut or that $S_1$ is homologically trivial.
Thus, Theorem 2.3 implies that the above decomposition is taut and that $E(S'_2 - {\rm int}(N(S_1)))$ is taut. 
This implies that $E(F)$ is a boundary-parallel
disc in $E(M_1, \gamma_1)$. By choosing $F$ appropriately, we can assume that no component
of $S'_2 - {\rm int}(N(S_1))$ lies in the product region that $E(F)$ separates off.
Hence, $E(F)$ separates off a ball component of $E(M_{12}, \gamma_{12})$.
This gives a pre-ball component of $(M_{12}, \gamma_{12})$, by Lemma 3.9. Hence, $F$ is product-separating in $(M_1, \gamma_1)$.
We remove $F$, as follows. We first slice under any discs of contact, creating a surface $\overline{F}$. The effect of this on 
$(M_{12}, \gamma_{12})$ is to remove some u-sutures and attach 2-handles along them.
Hence, the resulting sutured manifold still admits an allowable hierarchy.
But then $\overline{F}$ is then boundary-parallel, by Lemma 4.10, separating off a pre-ball.
Hence, we can remove it, and again the resulting sutured manifold
admits an allowable hierarchy.

Thus, $S_1$ does extend to an allowable hierarchy, as required.

\vskip 6pt
\noindent {\sl Case 1B.} The orientations on the patches do not induce transverse orientations of $S_1$ and $S_2$.
\vskip 6pt

Then, for some trace curve $C$ of $S_1 \cap S_2$, the orientations on $S_1$ or $S_2$ near that curve
do not agree. By Lemma 7.3, $C$ must be a simple closed curve, rather than an arc.
So, there is an annulus or M\"obius band $A$ embedded in a regular neighbourhood of this curve such that $A \cap S = \partial A$.
The transverse orientations on $S$ near $\partial A$ either all point towards $A$ or all point away from $A$.
Let $S'$ be the surface that results from performing an annular swap to $S$, using the annulus or M\"obius band
$A$. In other words, we perform an irregular switch along $C$.
Orient $S'$ so that the orientations of $S'$ and $S$ agree on their intersection.
If $A$ is an annulus, no component of $\partial A$ bounds a disc in $E(S)$, by Proposition 7.4.
If $A$ is a M\"obius band, $\partial A$ also cannot bound a disc in $E(S)$, since this would give
rise to an ${\Bbb RP}^3$ summand of $M$. So, by Proposition 4.8, either decomposition along $S'$ extends to an allowable
hierarchy or decomposition along $S' - S'_1$ does, where $S'_1$ is a component
of $S'$ that separates off a solid torus.
Therefore, since $S$ was non-separating, so too is at least one component of the new decomposing surface.
Note that $S$ and $S'$ have the same complexity,
but there is an ambient isotopy, leaving $\gamma$ invariant, that makes $S'$ standard, and reduces its complexity. This contradicts the fact
that $S$ had minimal complexity.

\noindent {\sl Case 2.} For every allowable hierarchy for $(M, \gamma)$ using connected decomposing surfaces, the first decomposing
surface is separating.

Note that this case may arise in several ways. For example, if $M$ is the union of two non-trivial knot exteriors, then the first
decomposing surface might have to be their common torus. Of course, though, in this situation, $E(M, \gamma)$ would
be toroidal, which we are assuming is not the case. The real situation where Case 2 might hold arises
in the proof of Propositions 4.8 and 4.9. Many of the annuli considered there are separating, and so
it seems unavoidable to consider this possibility.

By Proposition 4.11, we may assume that the first decomposing surface is not product-separating.
As in Case 1, we may arrange, using Lemma 4.2 and Theorem 5.5, for $S$ also
to be regulated. We pick $S$ to have minimal
complexity among connected regulated surfaces that extend to an 
allowable hierarchy, and that are not product-separating. 
Suppose that $S$ is not fundamental as a boundary-regulated surface. 
Then it can be written as a sum of non-empty boundary-regulated surfaces $S_1$ and $S_2$. 
Again, as in Case 1, we may assume that this summation has no trace annulus with a trivial boundary curve.
Each patch inherits a transverse orientation from $S$.

\noindent {\sl Case 2A.} These orientations induce transverse orientations of $S_1$ and $S_2$.

As explained in Case 1A, both $S_1$ and $S_2$ extend to allowable hierarchies.
We claim that at least one of $[S_1, \partial S_1]$ and $[S_2, \partial S_2]$ is non-trivial in $H_2(M, \partial M)$.
This will contradict the hypothesis that there is no connected non-separating surface in $M$
that extends to an allowable hierarchy.

Suppose that, contrary to the claim, both $S_1$ and $S_2$ are trivial in $H_2(M, \partial M)$.
Pick a point $p$ in $M$ disjoint from $S_1 \cup S_2$. We define an integer associated
to each component of $M - (S_1 \cup S_2)$, as follows. Pick an arc that starts at $p$
and ends within some component $R$ of $M - (S_1 \cup S_2)$. Make this arc transverse
to $S_1$ and $S_2$, and define the value of $R$ to be the signed intersection number
between this arc and $S_1 \cup S_2$. This is well-defined because both $S_1$
and $S_2$ are homologically trivial. Since $S$ is homologically trivial, one can also
define integers associated with the components of $M - S$. These are just the
same as those for $M - (S_1 \cup S_2)$, except near $S_1 \cap S_2$,
where components of $M - (S_1 \cup S_2)$ with the same value are glued together to form
$M - S$.
Now, as one crosses $S_1$ or $S_2$, these integers change by $1$.
So, near an arc or curve of $S_1 \cap S_2$, these integers take three different values.
But this is a contradiction, because, as $S$ is connected, $M - S$ has only two components.
Therefore, we deduce that $S_1 \cap S_2$ is empty. Hence, $S$ is the disjoint union
of $S_1$ and $S_2$, and so is disconnected. This is contrary to hypothesis.

\noindent {\sl Case 2B.} The orientations on the patches do not induce transverse
orientations for $S_1$ and $S_2$.

We would like to apply the same argument as in Case 1B, but there are some added
complications. Again, there is some simple closed
curve $C$ of $S_1 \cap S_2$, such that the orientations on $S_1$ or $S_2$ near that curve
do not agree. 
So, there is an annulus or M\"obius band $A$ embedded in a regular neighbourhood of this curve such that $A \cap S = \partial A$.
The transverse orientations on $S$ near $\partial A$ either all point towards $A$ or all point away from $A$.
Let $S'$ be the surface that results from performing an annular swap to $S$, using $A$. 
In other words, we perform an irregular switch along $C$. Let $A'$ be the annulus
or M\"obius band such that applying an annular swap to $S'$ along $A'$ gives back $S$.
Orient $S'$ so that the orientations of $S'$ and $S$ agree on their intersection.
If some component of $S'$ separates off a solid torus with no sutures, call this
$S'_1$. Otherwise set $S'_1 = \emptyset$. Let $S'' = S' - S'_1$.
Then, by Proposition 4.8, decomposition along $S''$ extends an allowable hierarchy.
The difficulty is that each component of the new decomposing surface $S''$ may be product-separating.
If this is not the case, then we have found a decomposing surface with smaller complexity than $S$,
that is not product-separating and that extends to an allowable hierarchy. This contradicts
our minimality assumption about $S$.

So, suppose that $S''$ is product separating. Our aim is to reach a contradiction.

We slice under any discs of contact for $S'$. This does not change the fact that $S''$
(which is a union of components of $S'$) is product-separating. But by Lemma 4.10, each component of 
$S''$ is now boundary-parallel. Note that $A'$ still satisfies $A' \cap S' = \partial A'$.
When we perform an annular swap along $A'$, we get a surface that is obtained from $S$
by slicing under discs of contact.

\noindent {\sl Case 2B(i).} $S'' = S'$ and $S'$ is connected.

Let $F \times I$ be the product sutured manifold that $S'$ separates off.  Suppose first that the interior of $A'$ is disjoint from $F \times I$.
If $A'$ separates the component of $(M_{S'}, \gamma_{S'})$ that contains it, then when we
apply the annular swap along $A'$, the resulting surface $S$ is disconnected, which is contrary to assumption.
So, $A'$ is non-separating. Using the product structure on $F \times I$,
extend $A'$ into $F \times I$, to form a properly embedded annulus in $M$, disjoint from $\gamma$.
This annulus is non-separating in $M$.
But then this annulus extends to an allowable hierarchy, by Proposition 4.9. This contradicts our assumption
that no non-separating surface extends to an allowable hierarchy.

Suppose now that $A'$ lies in $F \times I$. Then it is a non-trivial annulus with both boundary components in
the same component of $F \times \partial I$. It is therefore boundary parallel in $F \times I$.
So, when we perform the annular swap to $S'$ along $A'$, the resulting surface $S$ is disconnected,
which is contrary to assumption.

\noindent {\sl Case 2B(ii).} $S'' = S'$ and $S'$ is disconnected.

We are assuming that $S''$ is product-separating, and hence each component of $S'$ separates off a product
manifold. If these are disjoint from each other, then the interior of $A'$ is also disjoint from them,
and the argument of Case 2B(i) applies. So, suppose that the product sutured manifolds overlap.
Denote these by $F_1 \times I$ and $F_2 \times I$, which are separated off by $S'_1$ and $S'_2$ respectively.

Suppose first that these are nested: say that $F_2 \times I$ is a subset of $F_1 \times I$. Say that
$S'_1$ has transverse orientation pointing out of $F_1 \times I$. Then $R_+(F_1 \times I) - S'_1$
just consists of discs and annuli by Lemma 4.10. Since the transverse orientation on $S'_1$ points
away from $A'$, so does the transverse orientation on $S'_2$. So, it points into $F_2 \times I$.
Hence, $R_+(F_2 \times I)$ is a subset of $R_+(F_1 \times I) - S'_1$. It is an essential
subsurface of $R_+(F_1 \times I) - S'_1$, since otherwise $S'$ has a disc of contact, but we
have removed any of those. Also, $R_+(F_2 \times I)$ cannot be disc, because this implies
that a curve of $\partial A$ is trivial. Hence, $S'_2$ is an annulus, which is parallel to
an annulus in $R_+(F_1 \times I)$. So, when we perform an annular swap to $S'$ along $A'$,
we deduce that $S$ is disconnected, which is a contradiction.

Now consider the case where $F_1 \times I$ and $F_2 \times I$ are not nested. Then $S'_2$ is
an incompressible surface in $F_1  \times I$, and hence it is parallel to a subsurface $G$ of $\partial M \cap (F_1 \times I)$.
Note that $\partial S'_2$ cannot intersect any of the disc components of $R_+(F_1 \times I) - S'_1$, because
an innermost curve of intersection between one of these discs and $\partial S'_2$ would form a disc of contact,
and an outermost arc of intersection between one of these discs and $\partial S'_2$ would separate off
a disc component of $R_+(M_{S'}, \gamma_{S'})$, which would imply that $S'_2$ must be a disc.
Hence, the entirety of $G$ lies in $R_-(M)$, except possibly some annular collars on $\partial G$. Now,
$S'_2$ also separates off a product sutured manifold from $F_2 \times I$, and by assumption,
this does not lie wholly in $F_1 \times I$. Hence, its interior is disjoint from the product region between
$S'_2$ and $G$. Therefore, the union of $F_2 \times I$ and this product region is a product sutured manifold,
which is all of $(M, \gamma)$. This is contrary to assumption.

\noindent {\sl Case 2B(iii).} $S'$ has two components $S'_1$ and $S'_2$, and one of these, $S'_1$ say, separates off a solid torus 
with no sutures.

So, $S'' = S'_2$ and this is product separating. Let $V$ be the solid torus bounded by $S'_1$. 

We first consider the case where $V$ lies in the product manifold $F \times I$ separated off by $S'_2$.
Remove a regular neighbourhood of $A' \cap \partial V$ from $\partial V$ and attach the two copies of $A'$. The result is an annulus
in the product manifold. Remove the parts of this annulus that intersect $F \times \partial I$.
The result is a collection of annuli properly embedded in $F \times I$. Decomposing 
$F \times I$ along these annuli gives a product manifold. But the resulting manifold is a component
of $(M_S, \gamma_S)$. Hence, we deduce that $S$ is product-separating, which is a contradiction.

So we now suppose that $V$ is disjoint from the product manifold $F \times I$. Consider all the trace annuli
with interiors lying in $V$. Suppose that there is at least one of these. Each one separates $V$, and so there is one that is outermost in $V$. This separates off a
solid torus in $V$ with interior that is disjoint from the trace annuli and from $A$. So, we may assume that the interior of $V$
contains no trace annuli.

The solid torus $V$ cannot intersect $\partial M$ for the following reason. If it did, we could find an annulus
in $\partial V$ that runs between this intersection with $\partial M$ up to $A'$. Then extend the annulus across 
the product region that $S'_2$ separates off. We obtain an annulus
properly embedded in $M$ disjoint from $\gamma$ and also disjoint from $S$. If this annulus is separating,
then because parts of $S$ lie on both sides of it, we deduce that $S$ is disconnected, which is a contradiction.
On the other hand, if the annulus is non-separating, then we again have a contradiction to our assumption
that no non-separating surface extends to an allowable hierarchy.

We may give $V$ a Seifert fibration so that its intersection with $A'$ is a fibre. Then $V$ must have an
exceptional fibre, as otherwise $S$ is isotopic to $S'_2$ and we deduce that $S$ is product-separating.

Now we have arranged that no trace annulus has interior that intersects $V$.
However, other trace annuli may intersect $V$, because their interiors may lie outside of $V$ but their boundaries lie on $V$.
Suppose first that there is no such trace annulus. In other words, suppose that the only trace annulus
to intersect $V$ is $A$. The two parts of $S$ in $\partial V$ near $A$ lie in $S_1$ and $S_2$. So, as one goes around
$\partial V$, the surface must switch from $S_1$ to $S_2$. It cannot do so at a trace annulus, because by assumption
no others intersect $\partial V$. So, we deduce that $\partial V$ has non-empty intersection with $\partial M$,
and we have already ruled this case out.

So, we may assume that at least one trace annulus intersects $\partial V$ but has interior outside of $V$. 
This annulus therefore has interior disjoint from $V$ and from the product region $F \times I$.
Because of the orientations on $S'$, we deduce that each such trace annulus gives an
irregular switch. So, we may concentrate instead on one of these annuli. The only problem
case is when the switch gives a decomposing surface that separates off a solid torus $V'$ with no sutures.
We give $V$ and $V'$ Seifert fibrations in which the trace annuli are unions of fibres.
These fibrations must each have a singular fibre.
Since $V$ and $V'$ intersect along at least an essential annulus, we deduce that
there is a Seifert fibre space embedded in $M$ with base space a disc and with
two singular fibres. We deduce that $E(M, \gamma)$ is either toroidal or has a Seifert fibred
component which is neither a solid torus nor a copy of $T^2 \times I$. This is contrary to
assumption. $\square$

The above argument justifies the introduction of boundary-regulated surfaces. It is critical that we control the complexity of
our decomposing surfaces, by ensuring that they are fundamental. So, we are forced to consider the situation where a decomposing
surface $S$ is a sum $S_1 + S_2$. The problem is that this summation need not respect orientations; indeed $S_1$ and $S_2$ need
not even be transversely orientable. This sort of issue was faced by Tollefson and Wang [28], who showed that if the summation does not
respect orientations, then one should perform an irregular switch along a trace curve and then simplify the resulting surface $S'$ by
reducing its weight. The difficulty here is that it is not at all clear that the decomposition along $S'$ extends to an allowable hierarchy
or is just taut. If the trace curve is a simple closed curve,
then Proposition 7.4 and Proposition 4.8 give this. But this argument breaks down when a trace curve is an arc.
Fortunately, it is at this point that Lemma 7.3 intervenes. This asserts that when a surface is a sum, as a boundary-regulated
surface, then the trace curves that intersect $\partial M$ must respect the orientations on the surfaces.

\vskip 18pt
\centerline {\caps 8. Bounding the complexity of normal and boundary-regulated surfaces}
\vskip 6pt

We have been encoding normal surfaces and
boundary-regulated surfaces using vectors, which count the number
of elementary discs of each type in the surface. In this section, we find upper bounds on the weights of these surfaces,
by exploiting algebraic methods. Much of the material here is fairly well known, and goes
back to Haken [8] and Hass-Lagarias [9].

We consider the following general set-up. Suppose that $A$ is an $m \times n$
matrix with integer entries. 
We will examine solutions to the system
$$Ax = 0 \eqno{(1)}$$
where $x = (x_1, \dots, x_n)^T \in {\Bbb R}^n$, subject also to the inequalities
$$x_i \geq 0 \quad \hbox{ for } i = 1, \dots, n. \eqno{(2)}$$
Consider also a system of constraints, each of which is of the form
$$x_i =0 \quad \hbox{ or } \quad x_j = 0 \eqno{(3)}$$
where $1 \leq i , j \leq n$. We call (3) the {\sl compatibility constraints}.

We say that (1), (2) and (3) form a {\sl system} $\Sigma$.
We say that a solution $x \in {\Bbb R}^n$ to this system is {\sl integral} if all its co-ordinates are integers.
An integral solution $x$ is {\sl fundamental} if it cannot be written as $x = y + z$ where $y$ and $z$ are also integral
solutions to the system, neither of which is zero.

An integral solution $x$ to $\Sigma$ is a {\sl vertex} solution if
\item{(i)} whenever $k x = y + z$, for some positive real number $k$ and some solutions
$y$ and $z$, then both $y$ and $z$ are multiples of $x$; and
\item{(ii)} $x \not= k y$ for some integer $k \geq 2$ and some integral solution $y$.

One can interpret vertex solutions another way. It is well known that the set of real solutions to the system $\Sigma$ forms a union $P$ of finitely many convex polytopes.
Each polytope is obtained as follows. For each of the compatibility constraints, which asserts that $x_i = 0$ or
$x_j =0$ for some pair of integers $i$ and $j$ between $1$ and $n$, one makes a choice to fix $x_i$ to be
zero or to fix $x_j$ to be zero. Once this choice is made, the resulting system is a collection of equations $Ax = 0$
together with inequalities $x_i \geq 0$ for all $i$, and equalities of the form $x_i = 0$ for certain values of  $i$.
The solution set of this restricted system is invariant under dilations based at the origin.  Hence it forms a cone,
with cone point being the origin. Each line through the origin in $P$ intersects the plane $x_1 + \dots + x_n = 1$
in exactly one point. If we consider the subset of $P$ that satisfies this extra normalising constraint, then
we obtain a union of finitely many compact convex polytopes. A vertex solution is exactly a vertex of one of these
polytopes, rescaled minimally so that its coefficients are integers. 

\noindent {\bf Theorem 8.1.} {\sl Suppose that each row of $A$ has $\ell^2$ norm at most $k$. Then we have the following bounds.
\item{(1)} For any vertex solution $x$ to the system $\Sigma$, each co-ordinate of $x$ is at most $n^{1/2} k^{n-1}$.
\item{(2)} For any fundamental solution $x$ to the system $\Sigma$, each co-ordinate of $x$ is at most $n^{3/2} k^{n-1}$.

}

\noindent {\sl Proof.} (1) There is some positive real number $k$ and some vector $y = (y_1, \dots, y_n)^T$ satisfying
$y = kx$, such that $y$ is the unique solution to the following system of linear equations. The equations include
the equations $Ay = 0$. They also include extra equalities of the form $y_i = 0$ for some $i$
between $1$ and $n$. Finally, there is the normalising equation 
$$y_1 + \dots + y_n = 1.$$
We can write this system as 
$$\widetilde{A} y = \widetilde{b}$$
where the first $m$ rows of $\widetilde{A}$ are $A$. The vector $\widetilde{b}$ has all co-ordinates zero except the last,
which is $1$. The solution $y$ is unique and hence the kernel of $\widetilde{A}$ is zero. This implies that the rank of $\widetilde{A}$ is $n$.
Hence, there is some collection of $n$ rows of $\widetilde{A}$ with rank $n$. Let $\widehat{A}$ be this $n \times n$
matrix, and let $\widehat{b}$ be the corresponding $n$ entries of $\widetilde{b}$. Then $y$ is the unique solution to
$$\widehat{A} y = \widehat{b}.$$
The $n \times n$ matrix $\widehat{A}$ is invertible, and so
$$y = (\widehat{A})^{-1} \widehat{b}.$$
Now $(\widehat{A})^{-1} = ({\rm det}(\widehat{A}))^{-1} {\rm adj}(\widehat{A})$, where
${\rm adj}(\widehat{A})$ is the adjugate matrix for $\widehat{A}$. Hence, ${\rm det}(\widehat{A}) y = {\rm adj}(\widehat{A}) \widehat{b}$ 
is an integral vector. Therefore, $x$ is obtained from ${\rm det}(\widehat{A}) y$ by scaling by a factor at most 1.
We may therefore bound the entries of $x$ by bounding the entries of ${\rm det}(\widehat{A}) y$.
Each entry of ${\rm adj}(\widehat{A})$ is equal to the determinant of an $(n-1)\times (n-1)$ minor of $\widehat{A}$.
The modulus of this determinant is at most the product of the $\ell^2$ norms of its rows,
and this is at most $n^{1/2} k^{n-1}$, for the following reason. Each row of $\widehat{A}$ is either a row of $A$, a row consisting of a single $1$
and the remaining entries zero, or consists entirely of $1$s, in the case where it is the normalising
equation. So, each co-ordinate of ${\rm adj}(\widehat{A}) \widehat{b}$ has modulus at most $n^{1/2} k^{n-1}$.

(2) As described above, the real solutions to the system $\Sigma$ and the normalising condition $x_1 + \dots + x_n = 1$
form a union of convex compact polytopes. The vertices of the polytopes are multiples of the vertex solutions.
Hence, any solution to $\Sigma$ can be written as a linear combination $\sum \lambda_i v_i$,
where $\lambda_i \geq 0$ for each $i$, and each $v_i$ is a vertex solution. One may choose this expression so that
at most $n$ of the $\lambda_i$ are non-zero. It is clear that if $x = \sum \lambda_i v_i$ is fundamental, then 
each $\lambda_i < 1$, since otherwise $x$ has some $v_i$ as a summand. Therefore, we deduce
that each co-ordinate of $x$ has modulus at most $n$ times the maximal size of a co-ordinate for a vertex solution.
So, by (1), each co-ordinate of $x$ has modulus at most $n^{3/2} k^{n-1}$. $\square$

\noindent {\bf Theorem 8.2.} {\sl There is a universal computable constant $c$ with the following property.
Let $(M,\gamma)$ be a sutured manifold, and let ${\cal H}$ be a handle structure of uniform type.
Then, for each fundamental
boundary-regulated surface $S$ properly embedded in $M$, the weight of $S$
is at most $c^h$, where $h$ is the number of 0-handles of ${\cal H}$.}

\noindent {\sl Proof.} Each boundary-regulated surface can be described as an integral solution
to a system $\Sigma$ as above. The number of variables $n$ is equal to the number of disc types in the 0-handles.
Since ${\cal H}$ is of uniform type, the number of disc types $n$ is at most $k_1 h$, for some universal
computable constant $k_1$. Also, there is some universal computable
constant $k_2$ which is an upper bound for the $\ell_2$ norm of each row of the matrix $A$. So, by (2) of
Theorem 8.1, each co-ordinate of a fundamental boundary-regulated surface has modulus at most $n^{3/2} {k_2}^{n-1}$.
Furthermore, each elementary disc intersects at most $k_3$ 2-handles, for some universal computable constant $k_3$. 
So, the weight of the surface is at most $k_3 n^{3/2} {k_2}^{n-1}$. Since $n^{3/2} < 2^n$ for all $n \geq 2$,
we get an upper bound of $c^h$, for some universal computable constant $c$. $\square$

The above result is sufficient for most surfaces in our hierarchy. However, the first surface is different, because
its homology class is fixed. It is the Thurston norm of this homology class that is being certified. In order to
deal with this, we use some work of Tollefson and Wang [28]. They used triangulations rather than
handle structures, and so we will do the same here. 

Central to the approach of Tollefson and Wang is the notion of a lw-taut surface.
They say that a normal surface $S$ properly embedded in a compact triangulated 3-manifold $M$
is {\sl lw-taut} if 
\item{(i)} it has smallest Thurston complexity in its class in $H_2(M, \partial M)$;
\item{(ii)} no collection of components of $S$ is homologically trivial; and
\item{(iii)} it has smallest possible weight among normal surfaces satisfying (i) and (ii) in the same homology class of $H_2(M, \partial M)$.

\noindent Note that there is an important distinction in the way that Tollefson and Wang use the word `taut' compared with
its use by Scharlemann [25] and Gabai [3]. Tollefson and Wang require a `taut' surface to have minimal Thurston complexity in its class
in $H_2(M, \partial M)$, whereas Scharlemann and Gabai require Thurston complexity only to be minimised in the class in $H_2(M, N(\partial S))$.
Fortunately when $\partial M$ is a collection of tori, these concepts are closely related, using Lemma 2.1.

Our result bounds the weight of the first surface in the hierarchy, as follows.

\noindent {\bf Theorem 8.3.} {\sl Let $M$ be a compact orientable irreducible 3-manifold with boundary a
(possibly empty) collection of tori. Let ${\cal T}$ be a triangulation of $M$. Let $\phi$ be a simplicial 1-cocycle. 
Then there is a compact oriented lw-taut normal surface $S$ such that
$[S,\partial S] $ is Poincar\'e dual to $[\phi]$, satisfying $w(S) \leq k^{t^2} ||\phi||_1$, where $k$ is
a universal computable constant and $t$ is the number of tetrahedra of ${\cal T}$.
Here, $||\phi||_1$ is the $\ell^1$ norm of $\phi$.}

Tollefson and Wang consider the non-negative real solutions to the normal surface matching equations and compatibility
constraints. These form a union of convex polytopes. The boundary of any such polytope has a natural decomposition into faces.
The polytope itself is also called a face. A solution lying in such a face is said to be {\sl carried} by the face,
and when the solution corresponds to a normal surface $S$, we also say that $S$ is {\sl carried} by the face.
There is a unique face of minimal dimension that carries a given solution.

The following result of Tollefson and Wang is a central part of their work (see Theorem 3.3 in [28]).

\noindent {\bf Theorem 8.4.} {\sl Let $M$ be a compact orientable triangulated 3-manifold.
Let $S$ be an oriented, lw-taut normal surface, and let $C_S$ be the minimal face
carrying $S$. Then every surface carried by $C_S$ is lw-taut. Moreover, every surface carried by $C_S$ inherits a
well-defined orientation, and for any two surfaces $G$ and $H$ carried by $C_S$, the homology class in $H_2(M, \partial M)$ of the normal sum
$G + H$ is equal to the sum of the homology classes of $G$ and $H$.}

\noindent {\sl Proof of Theorem 8.3.} Let $S$ be an oriented, lw-taut surface such that $[S, \partial S]$ is Poincar\'e dual to $[\phi]$, 
and let $C_S$ be the minimal face
carrying $S$. This face is the convex hull of the rays through finitely many vertex surfaces $S_1, \dots, S_m$. 
By Theorem 8.1 (1), each co-ordinate of each $S_i$
has modulus at most $(c_1)^t$, for some universal computable constant $c_1$. 
Now, any element of $C_S$ can be written as a linear combination $\sum_i \lambda_i S_i$, where
each $\lambda_i \geq 0$. We may in fact ensure that there is no linear dependence between
the homology classes of the $S_i$ for which $\lambda_i$ is non-zero, as follows. Say that $\sum \mu_i [S_i, \partial S_i] = 0$,
where not all the $\mu_i$ are zero. Then we
add a small multiple of $\sum \mu_i S_i$ to $\sum_i \lambda_i S_i$. By Theorem 8.4, the effect on
the homology class is to add a multiple of $\sum \mu_i [S_i, \partial S_i]$, which is zero.
By choosing this multiple of $\sum \mu_i S_i$ correctly, we may arrange that in $\sum_i (\lambda_i + \mu_i) S_i$, all the
coefficients $\lambda_i + \mu_i$ are non-negative and at least one more coefficient is zero than is the case for $\sum \lambda_i S_i$. 
In this way, we may decrease the number of
non-zero coefficients. When this is minimal, the $S_i$ with non-zero coefficients are homologically linearly independent.
When we apply this procedure to $S$, we get $S = \sum_i \lambda_i S_i$ with exactly $r$ non-zero terms, say,
where $r \leq {\rm dim}(H_2(M,\partial M; {\Bbb Q})) \leq 2t$. By relabelling, we may ensure that $\lambda_1, \dots, \lambda_r$
are non-zero.

Pick a maximal tree in the 1-skeleton of ${\cal T}$. Then each edge not in the tree determines a loop
in $M$ that starts at some basepoint, travels in the tree to the start of the edge, runs along the edge, and then back through
the tree to the basepoint. These loops form a generating set for $\pi_1(M)$ and hence a generating set for $H_1(M)$.
So, some subset $\ell_1, \dots, \ell_d$ forms a basis for $H_1(M; {\Bbb Q})$.
Therefore, an element of $H_2(M, \partial M)$ is determined by its signed intersection numbers with these loops $\ell_1, \dots, \ell_d$.
In the subspace of $H_2(M, \partial M)$ spanned by the surfaces in $S_1, \dots, S_r$, we may determine the class of 
the surface $\sum_i \lambda_i S_i$ by examining its signed intersection number 
with just $r$ loops in $\ell_1, \dots, \ell_d$. By re-ordering, we may assume that these loops are $\ell_1, \dots, \ell_r$.

Form an $r \times r$ matrix $P$, with $(i,j)$ entry equal to the signed intersection number between $\ell_i$ and $S_j$. Similarly
form an $e \times r$ matrix $Q$, with $(k,j)$ entry equal to the number of intersection points between an edge $e_k$ and $S_j$,
where $e$ is the number of edges of the triangulation.
Hence, for a surface $\sum_{i=1}^r \lambda_i S_i$ carried by $C_S$, the unsigned intersection numbers with the edges are the entries of
$Q \lambda$, where $\lambda = (\lambda_1, \dots, \lambda_r)^T$. Similarly, the
signed intersection numbers with the loops $\ell_1, \dots, \ell_r$ are $P \lambda$, since by Tollefson and Wang's Theorem 8.4, the homology class
of $\sum \lambda_i S_i$ is $\sum \lambda_i [S_i, \partial S_i]$. Now, the signed intersection number of $S$ with the loops $\ell_i$ is the
evaluation $\phi(\ell_i)$, where $\phi$ is the given simplicial 1-cocycle. Let $\mu$ be the vector with entries $\phi(\ell_i)$. 
Hence, the unsigned intersection numbers of $S$ with the edges of ${\cal T}$ are
given by $QP^{-1} \mu$. 

Each co-ordinate of the normal surface vector $(S_j)$ is, by Theorem 8.1(1), at most $(c_1)^t$. 
So, each entry of $P$ has modulus at most $20t(c_1)^t$, since each elementary normal disc intersects the 1-skeleton
in at most 4 points and there are at most $5t$ types of elementary disc in $S_i$.
Similarly, each entry of $Q$ is at most $20t(c_1)^t$.
Therefore, the $\ell^2$ norm of each row of $P$ is at most $r^{1/2} 20t(c_1)^t$.
Hence, each entry of $P^{-1}$ has modulus at most $(r^{1/2} 20t(c_1)^t)^r$. Therefore, each entry of 
$QP^{-1} \mu$ has modulus at most $20t(c_1)^t r^2(r^{1/2} 20 t(c_1)^t)^r ||\phi||_1$. This is at most $k^{t^2}||\phi||_1$,
for some universal computable constant $k$. $\square$

Theorem 8.3 is stated in terms of the $\ell^1$ norm of a simplicial 1-cocycle. We can also provide a bound on this norm, as follows.

\noindent {\bf Lemma 8.5.} {\sl Let $M$ be a compact orientable 3-manifold with a triangulation ${\cal T}$ that has $t$ tetrahedra. If $b_1(M)> 0$,
there is an integral simplicial 1-cocycle $\phi$ representing a non-trivial element of $H^1(M)$ and satisfying $||\phi||_1 \leq (12t)^{3/2}(3 \sqrt{2})^{12t - 1}$.}

\noindent {\sl Proof.} Let $T$ be a maximal tree in the 1-skeleton of ${\cal T}$. Any class in $H^1(M)$ may be represented by a 1-cocycle that is zero on each edge of $T$. Conversely, any 1-cocycle that is zero on each edge in $T$ but not identically zero represents a non-trivial element of $H^1(M)$. Thus, the condition that $b_1(M)>0$ is exactly that there is such a 1-cocycle. This can be rephrased using a system of equations. For convenience, we introduce two variables for each edge not in $T$, each corresponding to an orientation on the edge. Thus, there are at most $12t$ variables. We also introduce a compatibility constraint requiring that, for each edge, at least one of these variables is zero. Thus, a simplicial 1-cochain that is zero on each edge of $T$ corresponds to an assignment of non-negative integers to these variables satisfying the compatibility constraints. This is a cocycle if and only if these variables satisfy a system of equations, with one equation for each face of the triangulation.  Each row of the matrix encoding this system has $\ell^2$ norm at most $3 \sqrt{2}$. We are assuming that $b_1(M)>0$ and so this system as a non-zero integral solution. Hence by Theorem 8.1 (1), it has a non-zero integral solution for which each co-ordinate is at most $(12t)^{1/2} (3 \sqrt{2})^{12t-1}$. This gives a cocycle $\phi$ with $||\phi||_1 \leq (12t)^{3/2} (3 \sqrt{2})^{12t-1}$. $\square$

\vfill\eject
\centerline{\caps 9. Determining the components of the parallelity bundle}
\vskip 6pt

We saw in Theorem 7.9 that when $(M, \gamma)$ is a taut sutured manifold
with a positive handle structure ${\cal H}$ satisfying some natural conditions, then a taut sutured manifold decomposition may be
performed along a regulated surface $S$ that is fundamental as a boundary-regulated surface.
Hence, by Theorem 8.2, there is an upper bound
to the number of elementary discs of $S$, that is an exponential function of the number of
0-handles of ${\cal H}$. Examples due to Hass, Snoeyink and Thurston [11] demonstrate that
one cannot, in general, find a better upper bound than this. This is potentially
problematic, because the sutured manifold $(M', \gamma')$ obtained by decomposing along $S$ inherits a handle
structure ${\cal H}'$ that may have many more handles than ${\cal H}$. Then, continuing this down
the hierarchy, the resulting handle structures may become so complex that they cannot be
analysed efficiently and so are not suitable for an NP algorithm. 

Fortunately, there is a method to get around this problem. Although
$S$ may have a large number of elementary discs, there is a linear upper bound (as a function
of the number of handles of ${\cal H}$) on the number of elementary disc types of $S$. Between
two discs of the same type, there is a subset of $M'$ homeomorphic to $D^2 \times I$.
These pieces patch together to form an $I$-bundle embedded in $M'$ known as the
{\sl parallelity bundle}. There is a linear upper bound on the number of handles of ${\cal H}'$
that are not part of this bundle. Thus, when $S$ has very large weight, almost all the handles
of $M'$ lie within this bundle. The solution, therefore, is to remove or modify
this bundle. The parts that are $I$-bundles over discs are replaced by 2-handles. The
parts that are $I$-bundles over more complicated surfaces are typically removed by decomposing
along the annuli and discs that separate them from the remainder of $M'$. In order to achieve this algorithmically,
it is important to be able determine the components of the parallelity bundle and to decide efficiently
which are $I$-bundles over discs. 

In this section, we explain how to do this. In Sections 9.1 and 9.2, we give the precise definition of the
parallelity bundle. In Section 9.3, we describe an algorithm, due to Agol, Hass and Thurston [2],
that was originally designed to determine the Euler characteristic of the components of a normal surface. 
In Section 9.4, we give our main application, which is a method for determining
the components of the parallelity bundle of the exterior of a normal or regulated surface, and for deciding which
of these components are $I$-bundles over discs. In Sections 9.5 and 9.6, we give further applications of
the Agol-Hass-Thurston algorithm. The first of these is concerned with `normal' 1-manifolds in a surface,
and it gives an efficient method for determining the components of the complement of such a 1-manifold.
In Section 9.6, we give a method for determining whether a boundary-regulated surface can be transversely oriented.

\vskip 6pt
\noindent {\caps 9.1. The parallelity bundle for sutured manifolds}
\vskip 6pt

Let ${\cal H}$ be a handle structure for the sutured manifold $(M, \gamma)$. A handle $H$ of ${\cal H}$
is a {\sl parallelity handle} if it satisfies one of the following:
\item{(i)} it is a 2-handle that is disjoint from 3-handles;
\item{(ii)} it is a 1-handle that intersects ${\cal H}^2 \cup \gamma$ in precisely two components and that is disjoint from the 3-handles; or
\item{(iii)} it is a 0-handle $H$ such that $H \cap ({\cal F} \cup \gamma)$ is connected, every
0-handle of $H \cap {\cal F}^0$ has index 0 and $H$ is disjoint from the 3-handles.

A product structure $D^2 \times I$ may be imposed on any parallelity handle $H$ so that
\item{(i)} the intersection between $H$ and
any other handle is of the form $\beta \times I$, for a subset $\beta$ of $\partial D^2$, and
\item{(ii)} each component of $H \cap \gamma$ is of the form $\beta \times \{ 1/2 \}$ for
some subset $\beta$ of $\partial D^2$.

These product structures may all be chosen so that the $I$-bundles agree on the intersection between any collection
of parallelity handles. Hence, they form an $I$-bundle ${\cal B}$ over a surface $F$, known as the {\sl base surface}. 
(See the proof of Lemma 5.3 in [22] for example.) This $I$-bundle is called the {\sl parallelity bundle}.
The {\sl horizontal boundary} $\partial_h{\cal B}$ is the $\partial I$-bundle over $F$. The {\sl vertical boundary} $\partial_v {\cal B}$
is the $I$-bundle over $\partial F$. The intersection between ${\cal B}$ and ${\rm cl}(M - {\cal B})$ is 
a subset of $\partial_v {\cal B}$ consisting of annuli disjoint from $\gamma$ and product discs.

Note that this notion of the parallelity bundle is related to, but not the same as, the notion of an `amalgam' presented
in [21]. An amalgam was also an $I$-bundle over a surface. However, away from its vertical boundary, it was allowed
to contain handles that were not parallelity handles. Another key difference was that its vertical boundary was required
to have interior disjoint from $\partial M$, and hence be annuli properly embedded in $M$.

\vskip 6pt
\noindent {\caps 9.2. The parallelity bundle for pairs}
\vskip 6pt

There is also a version of the parallelity bundle in the case where the 3-manifold has a specified subsurface $S$ in its boundary,
which is as follows.

Let ${\cal H}$ be a handle structure for the compact orientable 3-manifold $M$. Let $S$ be a compact
subsurface of $\partial M$, such that $\partial S$ is disjoint from the 2-handles and respects the product structure on the
1-handles of ${\cal H}$. In this case we say that ${\cal H}$ is
a {\sl handle structure} for the pair $(M, S)$.

A handle $H$ of ${\cal H}$ is a {\sl parallelity handle} if it admits a product structure $D^2 \times I$ such that
\item{(i)} $D^2 \times \partial I = H \cap S$;
\item{(ii)} each component of ${\cal F}^0 \cap H$ and ${\cal F}^1 \cap H$ is
$\beta \times I$, for a subset $\beta$ of $\partial D^2$.

Again, we will view this as an $I$-bundle over $D^2$, and again, these bundle structures may be
chosen so that they agree on the intersection between handles. So, their union is an $I$-bundle over
a surface. This is the {\sl parallelity bundle} for the pair $(M,S)$.

The main case when we consider pairs in this way is the following situation. Suppose that $M$ is
a compact orientable 3-manifold with a triangulation ${\cal T}$. Suppose that $S$ is a properly
embedded orientable normal surface in $M$. Then cutting $M$ along $S$ gives a 3-manifold $M'$.
It has a handle structure ${\cal H}'$, which arises by first dualising ${\cal T}$ to form a handle structure for $M$,
and then decomposing this along $S$. Inside $\partial M'$, the two copies of $S$ form a subsurface $S'$.
Then we will frequently consider the parallelity bundle for the pair $(M', S')$. In this case,
parallelity handles arise precisely in the space between two elementary normal discs of $S$.

\vskip  6pt
\noindent {\caps 9.3. The extended orbit counting algorithm of Agol, Hass and Thurston}
\vskip 6pt

In this subsection, we recall the algorithm of Agol, Hass and Thurston [2].
The setting for the algorithm is the natural numbers $\{ 1, 2, \dots, N \}$, which we denote by $[1,N]$. For integers
$a \leq b \in [1,N]$, the set of integers lying between $a$ and $b$, including $a$ and $b$, is denoted $[a,b]$.

The input to the algorithm is the following data:
\item{(i)} positive integers $N$, $k$, $d$ and $m$;
\item{(ii)} a collection of linear bijections $g_i \colon [a_i, b_i] \rightarrow [c_i, d_i]$, where $a_i, b_i, c_i, d_i \in [1,N]$,
and where $1 \leq i \leq k$;
\item{(iii)} a function $z \colon [1,N] \rightarrow {\Bbb Z}^d$, with the property that fewer than $m$ integers $j$ satisfy
$z(j) \not= z(j+1)$.

The linear bijections are called {\sl pairings}. Note that because each $g_i$ is a linear bijection, it sends the endpoints of
the interval $[a_i, b_i]$ to the endpoints of $[c_i, d_i]$. It is {\sl orientation preserving} if $g_i(a_i) = c_i$ and $g_i(b_i) = d_i$.
Otherwise it is {\sl orientation reversing}. Obviously, an orientation preserving pairing is just a translation.

Two integers $x$ and $y$ in $[1,N]$ are said to be in the same {\sl orbit} if there is a sequence of these bijections
$g_{i_1}, \dots, g_{i_s}$ such that $g_{i_s}^{\pm 1} \dots g_{i_1}^{\pm 1}(x) = y$. Clearly, this forms an equivalence
relation and the equivalence classes are the orbits.

The {\sl weight} of an orbit is the sum, over all elements $x$ in the orbit, of $z(x)$.

Let $D$ be the maximal value, over all integers $x \in [1,N]$ of $||z(x)||_1$, where the $\ell^1$ norm on
${\Bbb Z}^d$ is used.

The following theorem is due to Agol, Hass and Thurston [2].

\noindent {\bf Theorem 9.1.} {\sl There is an algorithm, with running time that is bounded above by 
a polynomial in $kmd (\log D)(\log N)$,
that produces a list of all orbits with their weights.}

It is the $\log N$ that is critical here.
It is perhaps helpful here to consider a simplified situation, where $d = 1$ and the weight function $z$ takes the constant value $1$.
Then the weight of an orbit is just the number of elements in it. Thus, the algorithm in this situation produces a list of orbits,
together with the number of elements in each orbit.

However, it is useful to allow more complicated weight functions. For example, suppose that we wanted to determine
the number of elements in the orbit of $x$, for some given $x \in [1,N]$. Then we could set $d = 2$, and assign every element
to have weight $(1,0)$, except for $x$ which has weight $(1,1)$. Hence, the weight of an orbit is either $(r,0)$ or $(r,1)$,
where $r$ is the number of elements in the orbit, and where the second co-ordinate is $1$ if and only if the orbit contains $x$.
Hence, by examining the unique orbit of weight $(r,1)$, we can determine the number elements in the orbit of $x$.

Agol, Hass and Thurston used their algorithm to determine the Euler characteristic of the components of a normal surface $S$
in a compact triangulated 3-manifold $M$, as follows. The input to the algorithm is the triangulation ${\cal T}$ of $M$
and the normal surface vector $(S)$. The integer $N$ is $w(S)$, the weight of $S$. We think of $[1,N]$ as the points of
intersection between $S$ and the 1-skeleton of ${\cal T}$. More precisely, each edge $e$ of ${\cal T}$ is oriented in some way, and the points $S \cap e$
are identified with some interval $[a,b]$ in $[1,N]$, so that the linear ordering of the points
along $e$ agrees with the linear order on $[a,b]$. The bijections $g_i$ arise from the faces 
of the triangulation. For each face $f$, the arcs $S \cap f$ come in at most three types. The arcs
of each type run from the 1-skeleton to the 1-skeleton, and so specify a linear bijection from a
sub-interval of $[1,N]$ to another such sub-interval.

It fairly clear that there is a one-one correspondence between the components of $S$ and the
orbits, for the following reason. Each component
of $S$ contains a point of intersection $p$ with the 1-skeleton, and this corresponds to an integer $x$ in
$[1,N]$. For any other point of intersection with the 1-skeleton, corresponding to an integer $y$,
this lies in the same component of $S$ as $p$ if and only if there is a path joining it to $p$ that
lies in the intersection between $S$ and the 2-skeleton of ${\cal T}$. 
This is equivalent to the existence of an equation $g_{i_s}^{\pm 1} \dots g_{i_1}^{\pm 1}(x) = y$,
for some bijections $g_{i_1}, \dots, g_{i_s}$.

By choosing the weight function $z$ appropriately, we can count the number of times components of $S$
intersect each edge $e$ of the triangulation. The Euler characteristic of a component $S'$ is a linear function
of the integers $|S' \cap e|$, as $e$ runs over all edges of the triangulation. Hence, this Euler characteristic
can be calculated in time that is a polynomial function of $t \log w(S)$, where $t$ is the number of 
tetrahedra of ${\cal T}$.

\vskip 6pt
\noindent {\caps 9.4. Determining the component types of the parallelity bundle}
\vskip 6pt

We now come to our main applications of Theorem 9.1.

\noindent {\bf Theorem 9.2.} {\sl There is an algorithm that takes, as its input, 
\item{(i)} a handle structure ${\cal H}$, of uniform type, with $h$ handles, for a sutured manifold $(M, \gamma)$;
\item{(ii)} a boundary-regulated vector $(S)_{\partial r}$ for a regulated surface $S$;
\item{(iii)} the transverse orientation on the elementary discs of $S$ that do not lie between
two normally parallel elementary discs;

\noindent and provides as its output, the following data. If $(M', \gamma')$ is the sutured manifold
that results from decomposing along $S$, and ${\cal B}$ is the parallelity bundle for the handle structure ${\cal H}'$
that it inherits, then the algorithm produces a handle structure for ${\rm cl}(M' - {\cal B})$, specifying which parts of its
boundary lie in $R_-(M')$ and which parts lie in $R_+(M')$, and, for each
component $B$ of ${\cal B}$ that does not lie between normally parallel closed components of $S$, it determines:
\item{(i)} the genus and number of boundary components of its base surface;
\item{(ii)} whether $B$ is a product or twisted $I$-bundle;
\item{(iii)} whether $B$ is a component of $M'$ that is a product sutured manifold; and
\item{(iv)} the way that $\partial_v B$ and ${\rm cl}(M' - {\cal B})$ intersect.

\noindent This algorithm runs in time that is bounded by a polynomial in $h \log(w(S))$.

}

We now explain in a bit more detail what is meant by `the way that $\partial_v B$ and ${\rm cl}(M' - {\cal B})$ intersect'.
The intersection between $\partial_v B$ and ${\rm cl}(M' - {\cal B})$ is a union of squares, where
each square is a component of intersection between a handle of $B$ and a handle of ${\rm cl}(M' - {\cal B})$.
Each square is a union of fibres in the $I$-bundle structure on $B$. Thus, when the algorithm provides
`the way that $\partial_v B$ and ${\rm cl}(M' - {\cal B})$ intersect', it gives, for each component $A$ of $\partial_v B$,
the squares forming the intersection between $A$ and the handles of ${\rm cl}(M' - {\cal B})$. It provides
these in the order in which they appear as one travels around $A$ in some direction.
% There are two ways of encircling the annulus
% $A$ and in the case where $B$ is a product $I$-bundle $F \times I$, one must distinguish between them. Thus,
% the algorithm picks an orientation on $F$ which determines an orientation on each component of $\partial F$ and
%thus a way of encircling each annulus $A$ in $\partial F \times I$. The squares  forming the intersection between $A$ and the handles of ${\rm cl}(M' - {\cal B})$
%are provided in this order. 
Note that $A$ need
not be properly embedded in $M'$, since its interior may intersect $\partial M'$. Each component of 
intersection between ${\rm int}(A)$ and $\partial M'$ is a union of $I$-fibres with their endpoints
removed. The algorithm also records when one meets such a component of ${\rm int}(A) \cap \partial M'$
as one travels around $A$, together with the way that it is attached to the incident squares in 
$A \cap {\rm cl}(M' - {\cal B})$. Also, in the case where $B$ is a product $I$-bundle, the algorithm
labels the components of $\partial_h B$ by $F_1$ and $F_2$ and specifies, for each component of $\partial_v B \cap {\rm cl}(M' - {\cal B})$, 
its arcs or curves of intersection with $F_1$ and $F_2$.

Note that the output provided by Theorem 9.2 is enough to specify the topological type of $M'$.
Specifically, the handle structure on ${\rm cl}(M' - {\cal B})$ of course determines its topology,
and the information provided about each component $B$ of ${\cal B}$ is enough to determine
its topology. Furthermore, there is an essentially unique way to attach $B$ onto ${\rm cl}(M' - {\cal B})$
in a way that makes the resulting manifold orientable and that is compatible with output provided
in (iv).

The reason that we have to exclude components of ${\cal B}$ that lie between normally parallel closed 
components of $S$ is as follows. Such a component is homeomorphic to a product. However,
to determine whether this is a product sutured manifold, we would need to know the transverse orientation
on the adjacent components of $S$, but this information might not be provided as part of the input.

We will also obtain the following variant, which deals with pairs.

\noindent {\bf Theorem 9.3.} {\sl There is an algorithm that takes, as its input, 
\item{(i)} a triangulation ${\cal T}$, with $t$ tetrahedra, for a compact orientable manifold $M$;
\item{(ii)} a vector $(S)$ for an orientable normal surface $S$;

\noindent and provides as its output, the following data. If $M'$ is the manifold
that results from decomposing along $S$, and $S'$ is the two copies of $S$ in $\partial M'$,
and ${\cal B}$ is the parallelity bundle for the pair $(M',S')$ with its induced handle structure, then
the algorithm produces a handle structure for ${\rm cl}(M' - {\cal B})$ and, for each component $B$ of ${\cal B}$, it determines:
\item{(i)} the genus and number of boundary components of its base surface;
\item{(ii)} whether $B$ is a product or twisted $I$-bundle; and
\item{(iii)} the way that $\partial_v B$ and ${\rm cl}(M' - {\cal B})$ intersect.

\noindent It runs in time that is bounded by a polynomial in $t \log(w(S))$.}

\noindent {\sl Proof.} We will focus on the proof of Theorem 9.2. The proof of
Theorem 9.3 is very similar, and requires only very minor modifications.

Let ${\cal B}_S$ be the parallelity bundle for the pair $(M', \partial N(S) \cap M')$. Let ${\cal B}_-$
be ${\cal B} \cap {\cal B}_S$. This is composed of handles of $M'$ that are parallelity handles in the
sense of Section 9.1 and that lie between parallel elementary discs of $S$. We will show how to
determine the handle structure of ${\rm cl}(M' - {\cal B}_-)$, together with (i), (ii), (iii) and (iv) of the theorem
for each component $B$ of ${\cal B}_-$. From this, it is straightforward to obtain
the required information about ${\rm cl}(M - {\cal B})$ and each component of ${\cal B}$. This
is because the number of handles of ${\rm cl}(M - {\cal B}_-)$ is bounded above by a linear function
of $h$. Hence, one can determine, in polynomial time, those handles of ${\rm cl}(M - {\cal B}_-)$ that
lie in ${\cal B}$. One can then attach these handles onto ${\cal B}_-$, and using the information
that we already have about ${\cal B}_-$, we can deduce what we need to know about 
${\rm cl}(M - {\cal B})$ and each component of ${\cal B}$.

The first stage in the procedure is to construct the handle structure of 
${\rm cl}(M' - {\cal B}_-)$. Each handle $H$ of ${\cal H}$ is divided up into handles of ${\cal H}'$.
It is a straightforward matter to determine those handles of $H \cap {\cal H}'$ that are not parallelity handles.
This is done by examining the entries of $(S)_{\partial r}$ corresponding to the elementary disc
types in $H$ and noting which of these are zero and which are non-zero.
The number of handles of $H \cap {\cal H}'$ that are
not parallelity handles is bounded above by a universal computable constant, for the following reason.
Since ${\cal H}$ is of uniform type and $S$ is boundary-regulated, there is a universal computable
upper bound on the number of elementary disc types of $S$ in $H$. If two adjacent elementary discs
of the same type are disjoint from $\partial M$, then the region between them becomes a parallelity handle
of $H \cap {\cal H}'$. If two adjacent elementary discs
of the same type intersect $\partial M$ and are coherently oriented, then also 
the region between them becomes a parallelity handle
of $H \cap {\cal H}'$. Because $S$ is boundary-regulated, at most one pair of adjacent
discs of the same type intersecting $\partial M$ can be incoherently oriented.

Consider a handle of $H'$ of ${\cal H}'$ that is not a parallelity handle. It is adjacent to a bounded number of other handles
of ${\cal H}'$. Some of these lie in ${\cal B}_-$, and so will form pieces of $\partial_v {\cal B}_- \cap H'$.
The remainder form handles of ${\rm cl}(M' - {\cal B}_-)$.

Thus, it is clear that a handle structure for ${\rm cl}(M' - {\cal B}_-)$ can be constructed in time that is bounded above by
a polynomial function of $h$. Moreover, the components of $\partial_v {\cal B}_- \cap {\rm cl}(M' - {\cal B}_-)$ can be
located. In addition, because of the transverse orientations provided in (iii), we can determine which parts of
the boundary of ${\rm cl}(M' - {\cal B}_-)$ lie in $R_-(M')$ and which parts lie in $R_+(M')$.

Since ${\cal H}$ is of uniform type, there is a universal computable upper bound $c$ to the number
of types of elementary discs of $S$ in each handle of ${\cal H}$. Say that there are $n \leq c h$ elementary
disc types of $S$. For each integer $i$ between 1 and $n$, let $x_i$ be the number of elementary discs of
that type. When a disc type intersects $\partial M$, we distinguish between the two possible transverse orientations on it.
So, when the elementary disc type lies in a 0-handle of ${\cal H}$, then $x_i$ is a co-ordinate
of  $(S)_{\partial r}$.

For each $i$ between 1 and $n$, set $N_i$ to be $\max \{ 0, 2(x_i-1) \}$. So, $N_i$ is the number of components of intersection between $\partial_h {\cal B}_-$
and the parallelity handles of $M'$ that lie between parallel elementary discs of $S$ of the $i$th type. 
Let $N =  \sum_{i=1}^n N_i$ be the number of elementary discs that make up $\partial_h {\cal B}_-$.

We now define the bijections $g_i$. Consider any $q$-handle $H_q$ of ${\cal H}$ that is attached to some $p$-handle $H_p$ where
$p < q$. Consider any elementary disc types in $H_p$ and $H_q$ that intersect. These correspond to intervals $[a_j, b_j]$ and $[c_j, d_j]$ in $[1, N]$.  
Thus, we get a bijection $[a_j, b_j] \rightarrow [c_j, d_j]$.
We consider all such bijections, as we run over all pairs of incident handles $H_p$ and $H_q$ and
all elementary disc types in $H_p$ and $H_q$. The number of these bijections is at most a constant times $h^2$.

We now define the weight function $z \colon [1,N] \rightarrow {\Bbb Z}^d$. The integer $d$ is equal to
$2|\partial_v{\cal B}_- \cap {\rm cl}(M' - {\cal B}_-)| + 4$. The first three co-ordinates are integers used for counting. 
The fourth co-ordinate is used to detect whether a component of ${\cal B}_-$ is incident to $\partial M$. 
The remaining
co-ordinates correspond to the components of $\partial_v {\cal B}_- \cap \partial_h {\cal B}_- \cap {\rm cl}(M' - {\cal B}_-)$. Each integer $j$ in $[1,N]$ 
corresponds to an elementary disc $E$ of $\partial_h {\cal B}_-$. Then $z(j)$ is a $d$-tuple. The first co-ordinate is set to be $1$ if
$E$ lies in a 0-handle of ${\cal H}$. The second is set to be $1$ if $E$ is in a 1-handle. The third co-ordinate is $1$
if $E$ is in a 2-handle. The fourth co-ordinate is set to be $1$ if $E$ intersects $\partial M$. 
The remaining
co-ordinates are set to be zero, unless the disc $E$ is incident to $\partial_v{\cal B}_-$. In this case,
the relevant co-ordinate of $z(j)$ is set to $1$ if $E$ runs over the relevant
component of $\partial_{v} {\cal B}_- \cap \partial_h {\cal B}_- \cap {\rm cl}(M' - {\cal B}_-)$; otherwise, it is zero.

It is clear that the orbits correspond to the components of $\partial_h {\cal B}_-$. Moreover, the first three co-ordinates
of the weight function count the number of handles of each index of $\partial_h {\cal B}_-$. Therefore, they can be used
to compute the Euler characteristic of this component. 
The fourth co-ordinate can be used
to determine whether a component of $\partial_h {\cal B}_-$ intersects $\partial M$. The remaining co-ordinates can be used
to determine whether a component of $\partial_h{\cal B}_-$ contains one of
the components of $\partial_v {\cal B}_-  \cap \partial_h {\cal B}_- \cap {\rm cl}(M' - {\cal B}_-)$. Theorem 9.1 therefore produces a list of components of
$\partial_h {\cal B}_-$, and for each such component $F$, it gives the Euler characteristic $\chi(F)$ and the components of
intersection between $\partial F$ and $\partial_{v} {\cal B}_- \cap \partial_h {\cal B}_- \cap {\rm cl}(M' - {\cal B}_-)$, and
whether $F$ intersects $\partial M$. Note that if a component of ${\cal B}_-$ is not a component of $M'$, then
we can determine whether it is a twisted $I$-bundle or a product, by examining a component $A$ of $\partial_v {\cal B}_- \cap {\rm cl}(M' - {\cal B}_-)$ to which it
is incident and seeing whether the two components of $A \cap \partial_h {\cal B}_-$ lie in the same
component of $\partial_h {\cal B}_-$. On the other hand, if a component $B$ of ${\cal B}_-$ is an entire component of
$M'$, then we can determine whether or not its horizontal boundary intersects $\partial M$.
If $B$ intersects $\partial M$, then it is a product sutured manifold, because a parallelity handle
of ${\cal B}_-$ incident to $\partial M$ must intersect $\gamma'$. But if $B$ is disjoint from
$\partial M$, then it is an $I$-bundle over a closed surface. As its horizontal boundary lies in the copies of $S$ in $\partial M'$, 
we deduce that $B$ either lies between normally parallel closed components of $S$ or is a thin regular neighbourhood
of a closed non-orientable surface.

Thus, if $B$ is any component of ${\cal B}_-$ that is also a component of $M'$ and an $I$-bundle over a closed surface, we need to be able to decide
whether $B$ is a product $I$-bundle or a twisted $I$-bundle. This is acheived as follows. We use the Agol-Hass-Thurston algorithm to determine the
normal vectors of the components of $S$. For each closed component $S'$, we cut along it and determine whether any component of the resulting manifold
is a component of the parallelity bundle. Such components are exactly the components of ${\cal B}_-$ that are a component of $M'$ and a twisted $I$-bundle 
over a closed surface.

The running time for the algorithm of Theorem 9.1 is bounded above by a polynomial function of
$kmd (\log D)(\log N)$. Here, $k$ is the number of pairings, which is at most a constant times $h^2$.
The integer $D$ is at most $2|\partial_v{\cal B}_- \cap {\rm cl}(M' - {\cal B}_-)| + 4$, which is bounded above by a linear function of $h$.
The integer $m$ is the number of values of $j$ such that $z(j) \not= z(j+1)$. These values of $j$
can occur at the outermost elementary discs of the same type, and when an elementary normal disc or one of
its neighbours intersects $\partial_v {\cal B}_- \cap {\rm cl}(M' - {\cal B}_-)$. So, again, there is an upper bound on the number of
such values of $j$, which is a polynomial function of $h$. Finally, $N$ is at most a linear function of $w(S)$ and $h$.
So, the running time is at most a polynomial function of $h \log (w(S))$.

We also need to compute, for each component $A$ of $\partial_v {\cal B}_-$, the components of intersection between
$A$ and the handles of ${\rm cl}(M' - {\cal B}_-)$ in the order in which they appear as one travels along $A$. Within a
single component of $A \cap {\rm cl}(M' - {\cal B}_-)$, this can be read straightforwardly from the handle structure
on ${\rm cl}(M' - {\cal B}_-)$ (up to reversing the direction of travel along $A$). 
However, some work needs to be done to go from one component of $A \cap {\rm cl}(M' - {\cal B}_-)$
to the next one. Again we use the Agol-Hass-Thurston algorithm. 

For each elementary disc type that appears in $S$, consider all the parallelity handles of ${\cal B}_-$ that are incident to discs of that type
and let $N'_i$ be the number of components of intersection between the interior of their vertical boundary and $\partial M'$.
Let $N'$ be $\sum_{i=1}^n 2 N'_i$.
So, $[1,N']$ records the components of intersection between $(\partial_v {\cal B}_- \cap \partial_h {\cal B}_-) - {\rm cl}(M' - {\cal B}_-)$
and the parallelity handles of ${\cal B}_-$.

We can set up bijections that record when one such component is incident to another such component.
Thus, the orbits correspond to the components of $(\partial_v {\cal B}_- \cap \partial_h {\cal B}_-) - {\rm cl}(M' - {\cal B}_-)$.
We can compute where these are attached to ${\rm cl}(M' - {\cal B}_-)$ as follows.
Let $d'$ be the number of components of $\partial_v {\cal B}_- \cap \partial_h {\cal B}_- \cap {\rm cl}(M' - {\cal B}_-)$. \break
We define a weight function $[1,N'] \rightarrow {\Bbb Z}^{d'}$ that records when a component of intersection 
between \break $(\partial_v {\cal B}_- \cap \partial_h {\cal B}_-) - {\rm cl}(M' - {\cal B}_-)$
and a parallelity handle of $M'$ is incident to a component of $\partial_v {\cal B}_- \cap \partial_h {\cal B}_- \cap {\rm cl}(M' - {\cal B}_-)$.
Thus, when an orbit corresponds to a component $\alpha$ of $(\partial_v {\cal B}_- \cap \partial_h {\cal B}_-) - {\rm cl}(M' - {\cal B}_-)$,
the total weight of $\alpha$ gives the components (or component) of $\partial_v {\cal B}_- \cap \partial_h {\cal B}_- \cap {\rm cl}(M' - {\cal B}_-)$
to which $\alpha$ is attached.

Thus, using the Agol-Hass-Thurston algorithm, we can compute  for each component $A$ of $\partial_v {\cal B}_-$, the components of intersection between
$A$ and the handles of ${\rm cl}(M' - {\cal B}_-)$ in the order in which they appear as one travels along $A$.
Again, the running time is a polynomial function of $h \log (w(S))$, as required. 
$\square$

We remark that one can easily establish an alternative version of Theorem 9.3. In this variation, one is not given a triangulation of $M$ with $t$ tetrahedra, but instead a handle structure of $M$ of uniform type, with $t$ handles. The conclusion of the theorem remains unchanged.

In the remainder of this section, we give two further applications of the Agol-Hass-Thurston algorithm.

\vskip 6pt
\noindent {\caps 9.5. Cutting a surface along a normal 1-manifold}
\vskip 6pt

In this section, we consider a handle structure ${\cal H}$ for a compact surface $F$. As usual, we let
${\cal H}^i$ denote the union of the $i$-handles, for $i = 0$, $1$ and $2$.

We say that a 1-manifold $C$ properly embedded in $F$ is {\sl normal} if
\item{(i)} it misses the 2-handles of ${\cal H}$;
\item{(ii)} it respects the product structure on each 1-handle;
\item{(iii)} its intersection with each 0-handle is a collection of properly embedded arcs;
\item{(iv)} for each 0-handle $H_0$, no arc of $C \cap H_0$ has both endpoints on the
same component of $H_0 \cap {\cal H}^1$ or on the same component of $\partial H_0 - {\cal H}^1$.

It is clear that normal 1-manifolds may be encoded algebraically, as follows.
The intersection between a normal 1-manifold $C$ and a 0-handle $H_0$ is a collection
of arcs, which come in finitely many types. One can form a vector $(C)$ which counts the
number of arcs of each type. These satisfy a collection of {\sl matching equations},
one for each 1-handle. There are also {\sl compatibility constraints} which prevent
distinct arc types within a 0-handle from coexisting if they inevitably intersect each other.

The {\sl weight} of a normal curve $C$ is $|C \cap {\cal H}^1|$.

We say that the handle structure ${\cal H}$ has {\sl uniform type} if, for each 0-handle $H_0$,
the number of components of $H_0 \cap {\cal H}^1$ is bounded above by some universal
computable constant.

\noindent {\bf Theorem 9.4.} {\sl There is an algorithm that takes, as its input, 
\item{(i)} a handle structure ${\cal H}$, of uniform type, with $h$ handles, for a compact surface $F$;
\item{(ii)} a vector $(C)$ for a normal 1-manifold $C$,

\noindent and provides as its output, the following data for each component $F'$ of $F - {\rm int}(N(C))$:
\item{(1)} the genus of $F'$;
\item{(2)} the number of boundary components of $F'$;
\item{(3)} the location of the components of $F' \cap \partial F$;
\item{(4)} the vector $(C')$ for each component $C'$ of ${\rm cl}(\partial F' - \partial F)$.

The algorithm runs in time that is bounded by a polynomial in $h \log (w(C))$. This polynomial depends
on the implied constant in the definition of uniform.

}

The proof is an easier version of Theorem 9.2, and so we omit it.

\vfill\eject
\noindent {\bf Corollary 9.5.} {\sl There is an algorithm that takes, as its input, 
\item{(i)} a handle structure ${\cal H}$, of uniform type, with $h$ handles, for a decorated sutured manifold $(M, \gamma)$, such that
for each 0-handle $H_0$ of ${\cal H}$, $H_0 \cap ({\cal F} \cup \gamma)$ is connected;
\item{(ii)} a vector $(S)_{\partial r}$ for a boundary-regulated surface $S$,

\noindent and determines whether any component of $\partial S$ bounds a disc disjoint from $\gamma$.
It also determines those components of $\partial S$ disjoint from $\gamma$ that are trivial, and for each such component, the
algorithm determines the side of $\partial S$ on which the trivialising planar surface lies.
The algorithm runs in time that is bounded by a polynomial in $h \log (w(S))$.

}

\noindent {\sl Proof.} The handle structure ${\cal H}$ induces a handle structure on $R_\pm(M)$.
Since ${\cal H}$ is of uniform type, there is a universal computable upper bound to the number of
times that a 0-handle of $R_\pm(M)$ can intersect the 1-handles. Therefore, we may declare
that the handle structure on $R_\pm(M)$ is of uniform type.

The intersection between $\partial S$ and $R_\pm(M)$ is a normal 1-manifold $C$ in this handle structure,
possibly after removing components of $\partial S \cap R_\pm(M)$ that lie within a single 0-handle $H_0$
of ${\cal H}$ and have endpoints on the same component of $H_0 \cap \gamma$.
The vector $(C)$ is a simple linear function of $(S)_{\partial r}$. Using Theorem 9.4, we can determine whether some
component $F$ of $R_\pm(M) - {\rm int}(N(C))$ is a disc with boundary disjoint from $N(\gamma)$.
Hence, we can decide whether any component of $\partial S$ bounds a disc disjoint from $\gamma$.

We can determine the trivial curves of $C$ as follows. Applying Theorem 9.4, we build a graph.
It has a vertex for each component of $R_\pm(M) - {\rm int}(N(C))$ that does not lie
between two normally parallel curves of $C$. For a vertex corresponding to a component $F'$ of $R_\pm(M) - {\rm int}(N(C))$,
the edges emanating from this vertex correspond to the components of ${\rm cl}(\partial F' - \gamma)$.
When an edge starts at one vertex and ends at another, the corresponding normal arcs
of $\partial S \cap R_\pm(M)$ are normally parallel. Each edge of the graph therefore corresponds to a maximal normally
parallel collection of arcs or simple closed curves. We call these {\sl arc-type} and {\sl curve-type} edges
respectively.

We now can characterise the trivial curves of $C$. Each corresponds to a curve-type edge of this graph
that separates off a subgraph that is a tree and that corresponds 
to a subsurface $F'$ in which all the vertices have genus zero, all the edges in this
subgraph are of curve-type and all the components of $F' \cap \gamma$ are u-sutures. This can readily
be determined using Theorem 9.4. $\square$

\vskip 6pt
\noindent {\caps 9.6. Determining whether a surface is transversely oriented}
\vskip 6pt

Although boundary-regulated surfaces are convenient to describe by means of their vectors,
it is regulated surfaces that we actually need to decompose along. The following result can be used
to determine whether a boundary-regulated surface is regulated.

\noindent {\bf Theorem 9.6.} {\sl There is an algorithm that takes, as its input
\item{(i)} a handle structure ${\cal H}$, of uniform type, with $h$ handles, for a sutured manifold $(M, \gamma)$;
\item{(ii)} a vector $(S)_{\partial r}$ for a boundary-regulated surface $S$;
\item{(iii)} a transverse orientation on a collection of at most $r$ elementary discs of $S$;
\item{(iv)} a specified 2-handle of ${\cal H}$ with a transverse orientation on its cocore $e$;

\noindent and determines whether there is a transverse orientation on $S$ that is compatible
with the transverse orientation on the discs in (iii) and
the given transverse orientation on the discs that intersect
$\partial M$. It also computes the signed intersection number with between $e$ and $S$
with this transverse orientation. The algorithm runs in time bounded above by a polynomial in $hr \log (w(S))$.

}

Using essentially the same proof, we can also determine whether a normal surface is orientable.

\vfill\eject
\noindent {\bf Theorem 9.7.} {\sl 
There is an algorithm that takes, as its input
\item{(i)} a triangulation ${\cal T}$ with $t$ tetrahedra for a compact orientable $3$-manifold $M$;
\item{(ii)} a vector $(S)$ for a normal surface $S$;
\item{(iii)} a transverse orientation on a collection of at most $r$ elementary discs of $S$;
\item{(iv)} a specified edge $e$ of ${\cal T}$;

\noindent and determines whether there is a transverse orientation on $S$ that is compatible
with the transverse orientation on the discs in (iii). It also computes the signed intersection number with between $e$ and $S$
with this transverse orientation. The algorithm runs in time bounded above by a polynomial in $tr \log (w(S))$.

}

We will focus on the proof of Theorem 9.6, as the proof of Theorem 9.7 is essentially the same.

\noindent {\sl Proof.} From the vector $(S)_{\partial r}$, one may produce a list of components of $S$.
Indeed, this was one of the original applications of the Agol-Hass-Thurston algorithm.
These are produced as a collection of vectors. It suffices to check, for each such component,
that it admits a transverse orientation that is compatible with the given orientations and to compute its
signed intersection number with $e$. Thus, we may assume that $S$ is connected.

We first determine whether the connected surface $S$ is transversely orientable, but with
no regard as to whether this is compatible with the pre-assigned transverse orientations on various elementary discs.
To do this, one forms the orientable double cover, and one checks whether it has two or one
component, as follows.

From the vector $(S)_{\partial r}$, one may readily compute the number of elementary discs
of $S$ of each type. For elementary discs in a 0-handle of ${\cal H}$, this is just a co-ordinate of $(S)_{\partial r}$,
or a sum of two co-ordinates.
For elementary discs in 1-handles and 2-handles, this is a simple linear function of the vector $(S)_{\partial r}$.
Let $N$ be twice the number of elementary discs of ${\cal H}$. We view $[1,N]$ as divided into two halves.
Integers in the first half $[1,N/2]$ correspond to the elementary discs of $S$ with some transverse orientation,
and integers in the second half correspond to the same discs but with the opposite orientation.
Elementary discs of the same type correspond to sub-intervals of $[1,N/2]$ or $[N/2+1, N]$.
The number of such subintervals is at most $k_1 h$, for some universal computable $k_1$.
Whenever a $p$-handle and $q$-handle of ${\cal H}$ are incident, the elementary normal
discs in these handles patch together. Hence, we obtain bijections between sub-intervals of $[1,N]$.
These are chosen so that they respect the transverse orientations on the discs.
The number $k$ of such pairings is at most $(k_2 h)^2$, for some universal computable $k_2$.

Clearly, $S$ induces either one or two orbits, depending on whether $S$
is unorientable or orientable. However, we also need to check whether it has
a transverse orientation that is compatible with the given ones on various elementary discs.
To do this, we introduce a weight function $z \colon [1,N] \rightarrow {\Bbb Z}^d$, where $d$ is $r$ plus twice
the number of elementary disc types that intersect $\partial M$. For an integer $j \in [1,N]$
corresponding to an elementary disc, and integer $i$ satisfying $1 < i \leq d$, the $i$th co-ordinate of $z(j)$ is an `error' co-ordinate.
It is set to $1$ if the orientation on the disc corresponding to $j$ is contrary to the orientation specified by
the fact that it is boundary-regulated, or by the transverse orientation given in (iii). Otherwise, this co-ordinate is set to $0$.
The first co-ordinate counts the signed intersection number between the elementary disc and the arc $e$.
Thus, it is zero unless the disc lies in the same 2-handle as $e$, in which case it is $\pm 1$.
So, the number of values of $j$ such that $z(j) \not= z(j+1)$ is at most $2r$ plus twice the number of disc types.

Theorem 9.1 produces a list of orbits and their weights. There is a transverse orientation on the connected surface $S$
that is compatible with the given ones if and only if there are two orbits and at least one of these orbits has weight $(q, 0, 0, \dots, 0)$.
In this case, $q$ is the signed intersection number between $S$ and $e$.
$\square$

\vfill\eject
\centerline {\caps 10. An algorithm to simplify handle structures}
\vskip 6pt

In Theorem 9.2, we showed that, when a taut sutured manifold decomposition is performed along a regulated surface,
the components and structure of the parallelity bundle of the resulting sutured manifold $(M', \gamma')$ may be
efficiently determined. In this section, we utilise this information to produce a simplified handle structure for
a sutured manifold obtained from $(M', \gamma')$.

\noindent {\bf Theorem 10.1.} {\sl Let ${\cal H}$ be a positive handle structure of a connected decorated sutured manifold $(M, \gamma)$,
of uniform type, with $h$ handles. Let $S$ be a regulated surface properly embedded in $M$. Suppose that any normally parallel closed
components of $S$ are coherently oriented.
Let $(M', \gamma')$ be obtained by performing an allowable decomposition along $S$. 
Suppose that no component of $M'$ has boundary a single torus with no sutures.
Then, there is an algorithm that takes as its
input ${\cal H}$, $(S)_{\partial r}$ and the transverse orientations on the elementary discs of $S$ that do not lie between elementary discs of the same type,
and either correctly
asserts that $E(M', \gamma')$ is not taut, or supplies
a handle structure ${\cal H}'''$ for a decorated sutured manifold $(M''', \gamma''')$ with the following properties:
\item{(i)} there are decompositions
$$E(M', \gamma') \buildrel A \over \longrightarrow E(M'', \gamma'') \sharp Y  \buildrel D \over \longrightarrow E(M''', \gamma''') \sharp Y$$
for some closed 3-manifold $Y$, and
where $A$ is a collection of product discs and oriented non-trivial annuli disjoint from $\gamma'$ and $D$ is a collection of product discs; 
this decomposition is taut if and only if $E(M', \gamma')$ is taut;
\item{(ii)} if $E(M', \gamma')$ is taut, then it can be certified in polynomial time as a function of $h$ that $Y$ is homeomorphic to
the 3-sphere; 
\item{(iii)} if $E(M', \gamma')$ is atoroidal and no component of $E(M',\gamma')$ is a Seifert fibre space other than a solid torus 
or a copy of $T^2 \times I$, and $(M', \gamma')$ admits an allowable hierarchy, then so does $(M''', \gamma''')$;
\item{(iv)} ${\cal H}'''$ is of uniform type;
\item{(v)} each 0-handle of ${\cal H}'''$ lies inside a 0-handle of ${\cal H}$;
\item{(vi)} for each 0-handle $H$ of ${\cal H}$, the complexity of $H \cap {\cal H}'''$ is at most that of $H$,
and in the case of equality, $H \cap {\cal H}'''$ contains a single 0-handle isotopic to $H$;
\item{(vii)} ${\cal H}'''$ has no parallelity handles other than 2-handles;
\item{(viii)} ${\cal H}'''$ is positive.

\noindent The algorithm runs in time that is bounded by a polynomial in $h \log(w(S))$.

}

The idea behind this theorem is simply that if one were to decompose ${\cal H}$ by cutting along $S$, then the resulting
handle structure ${\cal H}'$ might have many parallelity handles. Our goal is to remove these parallelity handles,
primarily by decomposing along the vertical boundary components of the parallelity bundle ${\cal B}$ for ${\cal H}'$, or
by replacing them by 2-handles. The exact procedure is slightly delicate, and requires modifications to
be made in the right order. We follow the method given in Section 8 of [21]. 

\vskip 6pt
\noindent {\caps 10.1. Decomposition along product discs} 
\vskip 6pt

Some component $B$ of the parallelity bundle ${\cal B}$ may have vertical boundary containing a
properly embedded product disc. This happens when a component of $\partial_v B$ neither is a properly embedded annulus
nor lies wholly in $\partial M$.

If any of these product discs intersects a u-suture of $\gamma'$, we attach a 2-handle along the u-suture,
and enlarge ${\cal B}$ by including this 2-handle in the parallelity bundle. This may
increase the complexity of the handle structure, but this is only a temporary problem,
because all enlarged components of ${\cal B}$ will be later removed from the manifold.

So suppose that a component $B$ of ${\cal B}$ has vertical boundary that contains a properly embedded product disc, 
even after $B$ has possibly been enlarged. In this case, we may remove $B$, as follows.
The product discs and annuli in ${\rm cl}(\partial_v B - \partial M)$ become discs and annuli in the boundary
of the new manifold, each of which contains a new arc of a suture. The reason that we
may remove $B$ in this way is as follows. 

We can decompose along any product discs in ${\rm cl}(\partial_v B - \partial M)$, since
these are allowable. If $B$ becomes detached from the remainder of $M'$, it is a product
sutured manifold, and hence may be decomposed along product annuli and discs to a taut 3-ball.
On the other hand, if $B$ is still
attached to the remainder of $M'$, then we may also decompose along vertical product discs in $B$, until $B$
becomes a collection of collar neighbourhoods of the annuli in ${\rm cl}(\partial_v B - \partial M)$.
We can then remove these collar neighbourhoods, without changing the homeomorphism
type of the manifold.

Thus, in this way, we may assume that every component of $\partial_v {\cal B}$ either is a properly embedded
annulus in $M'$ or lies wholly in $\partial M'$. All future modifications that we make will preserve this property.

\vskip 6pt
\noindent {\caps 10.2. The boundary graph}
\vskip 6pt

In order to be able 
to decide which vertical boundary components of the parallelity bundle to cut along, we introduce
two graphs. The first of these, the boundary graph, encodes the way that $\partial M'$ is cut into pieces
by the parallelity bundle. The second graph, the connectivity graph, stores information about the 
way that the pieces of $M' - {\cal B}$ and ${\cal B}$ are connected together.

Let ${\cal H}'$ be the handle structure obtained by decomposing ${\cal H}$ along $S$.
Let ${\cal B}'$ be the union of the components of its parallelity bundle that are not 2-handles.

We define the {\sl boundary graph} $X_\partial$ which encodes information about the boundary of $M'$. It has two types of vertices: a vertex for each component
of $\partial E(M', \gamma') - {\cal B}'$ (known as $G$-vertices, where `$G$' stands for `guts'), 
and a vertex for each component of $\partial_h{\cal B}'$ (known as $B$-vertices, where `$B$' stands for `bundle'). 
It has an edge for each boundary component of $\partial_v {\cal B}'$. Note that every component of $\partial_v {\cal B}'$ is an annulus,
and so the edges of $X_\partial$ come in pairs, where edges in a pair correspond to boundary components of the same annulus.
Since $\partial_v {\cal B}'$ separates ${\cal B}'$ from the remainder of $E(M', \gamma')$, each edge of $X_\partial$ joins a $B$-vertex
to a $G$-vertex. We orient the edge, so that it points from the $B$-vertex to the $G$-vertex

We give each vertex of $X_\partial$ an integer value, which we call its {\sl $\chi$-value}. This is the Euler characteristic of the
corresponding component of ${\rm cl}(\partial E(M', \gamma') - {\cal B}')$ or $\partial_h {\cal B}'$.

Each $B$-vertex is assigned the label $R_-$ or $R_+$, according to whether the corresponding component of $\partial_h {\cal B}$
lies in $R_-(M')$ or $R_+(M')$. Each edge is also assigned the label $R_-$ or $R_+$, according to whether the
corresponding curve lies in $R_-(M')$ or $R_+(M')$. This label is the same as that of the $B$-vertex that it emanates from.

Clearly, $X_\partial$ and its $\chi$-values are constructible in time that is bounded above by a polynomial function of $h \log (w(S))$,
for the following reasons. The $G$-vertices and edges can be readily
constructed from the handle structure on ${\rm cl}(M' - {\cal B}')$ given by Theorem 9.2.
The $B$-vertices, the way that
they are incident to the edges of $X_\partial$ and their $\chi$-values are also provided by the algorithm in Theorem 9.2.

It is perhaps convenient to view $X_\partial$ as a graph embedded in $\partial M'$. 

%A simple example is shown in Figure 24, where the $\chi$-values of vertices are also labelled.

%\vskip 18pt
%\centerline{
%\includegraphics[width=3in]{boundary-graph-example.eps}
%}
%\vskip 6pt
%\centerline{Figure 24: An example of the boundary graph}

\vskip 6pt
\noindent {\caps 10.3. The connectivity graph}
\vskip 6pt

The {\sl connectivity graph} $X_c$ has an edge for each component of $\partial_v {\cal B}'$.
It has a vertex for each component of $M' - \partial_v {\cal B}'$. These are clearly of two types:
components of ${\cal B}'$ and components of ${\rm cl}(M' - {\cal B}')$. We orient
the edges so that they point from the former vertices to the latter ones.

Again, the graph $X_c$ may be readily constructed using the data provided by Theorem 9.2.

Note that there is a morphism of graphs $X_\partial \rightarrow X_c$. This sends edges to edges, by
sending each boundary component of $\partial_v {\cal B}'$ to the component of
$\partial_v {\cal B}'$ that contains it. It sends vertices to vertices in a similar way.

\vskip 6pt
\noindent {\caps 10.4. The algorithm}
\vskip 6pt

We now explain how to modify $M'$ and its handle structure. As explained in Section 10.1,
the first stage is to attach a 2-handle along any u-suture that intersects ${\cal B}$,
and enlarge ${\cal B}$ to include these 2-handles.
We then remove from $M'$ any component $B$ of ${\cal B}$ that contains
a product disc in $\partial_v B$. (This is the procedure described in Section 10.2.)
We then construct the boundary graph $X_\partial$
and the connectivity graph $X_c$.

We now modify the manifold $M'$, and at the same time modify the boundary graph
and the connectivity graph.
We follow the same procedure used in Section 8 of [21], and so, for ease of reference, we use the
same numbering of cases.

The algorithm starts by checking whether the hypothesis of the first case holds. If it does, then the
procedure given in that case is followed. This changes the handle structure and the graphs $X_c$ and $X_\partial$.
The algorithm then goes back to the beginning, and checks whether Case 1 holds in this new structure.
On the other hand, if Case 1 does not apply, then the algorithm proceeds to Case 2, and so on.
In this way, a case is considered only if all the earlier cases do not apply.

\noindent {\sl Case 1.} When there a $B$-vertex of $X_\partial$ with $\chi$-value 1.

This corresponds to component of $\partial_h {\cal B}'$ that is a disc. This therefore lies in a component
of ${\cal B}'$ that is an $I$-bundle over a disc. If its vertical boundary lies in $\partial M'$, then it
is a component of $E(M', \gamma')$ that is a taut ball. No further decompositions are needed for 
this ball, and so it can be ignored. So, we now assume that the vertical boundary of this component
of ${\cal B}'$ is a properly embedded annulus.

Suppose first that this component of ${\cal B}'$ had been enlarged when we added 2-handles along u-sutures.
Then we modify this component of ${\cal B}'$ as follows. We remove one of the added 2-handles, forming an
$I$-bundle over an annulus. We then remove this $I$-bundle, which does not change the homeomorphism type
of the manifold. In doing so, we replace the vertical boundary component that is incident to the remainder of $M'$
by a u-suture. We therefore remove the vertex of $X_c$ corresponding to this component of ${\cal B}'$ and
its incident edge. We also remove the two corresponding vertices of $X_\partial$ and their incident edges.
We combine the two $G$-vertices at the other endpoints of these edges, forming a single $G$-vertex. Its
$\chi$-value is the sum of the $\chi$-values of the two original $G$-vertices.

Suppose now the component of ${\cal B}'$ had not been enlarged when we added 2-handles along u-sutures.
Then we remove this component of ${\cal B'}$ and replace it by a 2-handle. We therefore remove
the two corresponding vertices of $X_\partial$ and the two edges that are incident to them. For each of these
edges, we add 1 to the $\chi$-value of the $G$-vertex that is incident to it. (If both edges are incident
to the same $G$-vertex, we add to $2$ to its $\chi$-value.) We also remove the corresponding vertex of $X_c$
and the edge that is incident to it.

We therefore now assume that every $B$-vertex has $\chi$-value other than 1.

\noindent {\sl Case 2.} When there is a $B$-vertex with no incident edges.

This corresponds to a component of ${\cal B}'$ that is an $I$-bundle over a surface.
If this component is a product sutured manifold (which is information provided by the algorithm
in Theorem 9.2), then there is a sequence of decompositions along product annuli and then product discs
which takes it to a taut ball. So, we replace this component this component of ${\cal B}'$ with a single
0-handle containing a single suture.

On the other hand, if this component of ${\cal B}'$ is not a product sutured manifold,
then it is an $I$-bundle over a closed surface and its horizontal boundary lies entirely in 
$R_-(M')$ or entirely in $R_+(M')$. Our procedure checks the $\chi$-value of this $B$-vertex. 
If it is zero, then this component of $M'$ is an $I$-bundle
over a torus or Klein bottle, and again there is a sequence of taut decompositions along annuli
and discs taking it to a taut ball. If the $\chi$-value of the vertex is non-zero, then this component of
$E(M', \gamma')$ is not taut, and the algorithm terminates with a declaration to this effect.

%This corresponds to a component of ${\cal B}'$ that is an $I$-bundle over a surface. If its horizontal boundary
%intersects $\partial M'$, then (by definition of the parallelity bundle), its horizontal boundary intersects
%both $R_-(M')$ and $R_+(M')$. On the other hand, if its horizontal boundary is disjoint from $\partial M'$, then
%it corresponds to parallel closed components of $S$, or possibly a closed component of $\partial S$ parallel to a component
%of $\partial M$ without sutures, or possibly just one component of $S$ that is a fibre in a fibration of $M'$ over the circle.
%Each case is either excluded by hypothesis or leads to the conclusion that this component of ${\cal B}'$
%intersects both $R_-(M')$ and $R_+(M')$. So, there is a sequence of decompositions along product annuli and then product discs
%which takes it to a taut ball. So, we replace this component this component of ${\cal B}'$ with a single
%0-handle containing a single suture.

%On the other hand, if this component of ${\cal B}'$ is incident only to $R_-(M')$ or to $R_+(M')$, then our
%procedure checks the $\chi$-value of this vertex. If it is zero, then this component of $M'$ is an $I$-bundle
%over a torus or Klein bottle, and again there is a sequence of taut decompositions along annuli
%and discs taking it to a taut ball. If the $\chi$-value of the vertex is non-zero, then this component of
%$E(M', \gamma')$ is not taut, and the algorithm terminates with a declaration to this effect.

We therefore now assume that every $B$-vertex is incident to some edge.

\noindent {\sl Case 3.} When every vertex of $X_\partial$ has non-positive $\chi$-value, and for every pair of edges in
$X_\partial$, one lies in $R_-$ and one lies in $R_+$.

This is exactly the situation where every component of $\partial_v {\cal B}'$ that does not lie in $\partial M'$ is a non-trivial product annulus in $E(M',\gamma')$, 
for the following reason.
If one of these annuli was trivial, then one of its boundary components would, by definition,
bound a disc in $R_-(E(M', \gamma'))$ or $R_+(E(M', \gamma'))$. The intersection between this disc and $\partial_v {\cal B}'$
is a collection of simple closed curves, and an innermost one bounds a disc which corresponds
to a vertex of $X_\partial$ with $\chi$-value 1. This is contrary to our hypothesis.

Thus, we decompose $(M',\gamma')$ along the annuli in $\partial_v {\cal B}'$ that do not lie in $\partial M'$. This has the effect
of removing ${\cal B}'$ from the remainder of $M'$. Then ${\cal B}'$ becomes a product sutured manifold,
which we may further decompose along product discs, resulting in a collection of taut 3-balls.
The manifold ${\rm cl}(M' - {\cal B}')$ inherits a handle structure. Each component of $\partial_v {\cal B}'$ that does not lie in $\partial M'$ 
becomes a suture in this manifold that is not a u-suture. It is shown in [21] that this new handle structure has smaller complexity that
that of $M'$. (Recall that a notion of complexity for handle structures was defined in [21]
and summarised in Section 5.1.)

We may therefore assume that some annulus of $\partial_v {\cal B}'$ does not lie in $\partial M'$ and is not a non-trivial product annulus.

\noindent {\sl Case 4.} Some $G$-vertex $v$ has positive $\chi$-value.

This corresponds to component of ${\rm cl}(\partial E(M', \gamma') - {\cal B}')$ that is a disc $D'_1$. This vertex $v$ is therefore
incident to a single edge $e$ of $X_\partial$. This edge is sent to an edge of the connectivity graph by
the map $X_\partial \rightarrow X_c$. If this edge of $X_c$ is non-separating, then the algorithm
terminates with the correct declaration that $E(M', \gamma')$ is not taut, for the following
reason. Let $A$ be the annulus of $\partial_v {\cal B}'$ corresponding 
to this edge of $X_c$. If the edge is non-separating, then the annulus $A$ is also
non-separating, and hence so is the disc $A \cup D'_1$. Hence, $E(M', \gamma')$ contains
an essential disc with boundary disjoint from $\gamma'$, and therefore $E(M' ,\gamma')$ is not taut.

Therefore, suppose that this edge of $X_c$ is separating.
Let $e'$ be the other edge in this pair in $X_\partial$. Then $e'$ is separating.
The curve in $R_\pm(M)$ corresponding to $e'$ bounds the disc $A \cup D'_1$ in $E(M',\gamma')$.
So if $E(M', \gamma')$ is taut, then this curve bounds a disc in $\partial E(M', \gamma')$
disjoint from the sutures. So, the algorithm computes the total $\chi$-values of the vertices in
each of the components of $X_\partial - e'$. It determines whether one of these components
$C$ has $\chi$-value $1$ and is disjoint from the sutures of $E(M', \gamma')$. If this applies
to neither component of $X_\partial - e'$, then the algorithm terminates with the correct declaration
that $E(M' ,\gamma')$ is not taut. If it applies to both components of $X_\partial - e'$, then this component of
$\partial E(M' , \gamma')$ is a sphere with no sutures. This is impossible unless $E(S)$ has a 2-sphere
component, in which case $E(M' , \gamma')$ is not taut. So, again the algorithm terminates.
So, suppose that exactly one component $C$ of $X_\partial - e'$ has $\chi$-value $1$ and is disjoint from the sutures of $E(M', \gamma')$.
This corresponds to a disc $D_2'$ in $R_\pm(E(M', \gamma'))$. The argument now divides into two cases.

\noindent {\sl Case 4A.} $v$ does not lie in $C$.

Then $D'_1$ and $D'_2$ are disjoint. Then the union of $D'_1$, $D'_2$ and the annulus $A$
of $\partial_v {\cal B}'$ incident to them together forms a 2-sphere. Since $A$ is separating,
so too is this 2-sphere. Using the Agol-Hass-Thurston algorithm, we can determine whether
this 2-sphere separates components of $\partial M'$. If it does, then the algorithm terminates with the
correct declaration that $E(M', \gamma')$ is not taut. On other hand, if the entirety  of $\partial M'$
lies on one side of the 2-sphere, then we obtain an expression of $E(M', \gamma')$ as a connected sum,
where one of the summands $Y$ is closed. If $E(M', \gamma')$ is taut, then $Y$ is the 3-sphere.
Using the result of Schleimer [26] and Ivanov [12] that 3-sphere recognition is in NP, then we can certify
that $Y$ is the 3-sphere, once we have a handle structure for $Y$ of uniform type and where the number
of handles is bounded above by a linear function of $h$. We can construct this handle structure as follows.
Let $Z$ be the component of ${\rm cl}(M' - A)$ containing
$D_1'$ and $D'_2$. Then $Y$ is obtained from $Z$ by attaching a 3-ball.
The parts of $Z$ lying in $M' - {\cal B'}$ already have a handle structure. The parts lying in 
${\cal B}'$ are all $I$-bundles over planar surfaces, since their horizontal boundaries lie in $D'_1 \cup D'_2$.
So we can give them a handle structure using just 1-handles and 2-handles, knowing just the number
of vertical boundary components of each component of ${\cal B'} \cap Z$ and whether or not it is product $I$-bundle.
This is all information given by Theorem 9.2. Thus, we can construct a handle structure for $Y$ of uniform type,
and where the number of handles is bounded by a linear function of $h$. If $Y$ is a 3-sphere, this can be
certified in polynomial time. If so, then $Z$ is a 3-ball, and we can therefore enlarge ${\cal B}'$ by adding the ball to it. Note that this may change
$(M', \gamma')$ since the ball was a subset of $E(M', \gamma')$ and not necessarily
$M'$. So, the ball may contain some attached 2-handles. But clearly $E(M' ,\gamma')$
remains unchanged. Moreover, when we attach these 2-handles to the u-sutures,
the resulting manifold still has an allowable hierarchy, which is just the original hierarchy.

The way that the algorithm processes this
is as follows. The annulus $A$ corresponds to an edge of the connectivity graph. Since $A$ is
separating in $M'$, removing this edge from the connectivity graph creates two components.
The edge $e$ points into one of these components, and this corresponds to a component $Z$
of $M' - A$. Once we have verified the certificate that $Z$ is a 3-ball,
we discard all handles of $M' - {\cal B}'$ that lie in $Z$, together with all components of ${\cal B}'$ lying 
in $Z$. We also discard the corresponding edges and vertices of $X_\partial$ and $X_c$. The two edges of $X_\partial$
corresponding to the two components of $A$ are also removed. The adjacent $B$-vertices simply have
1 added to their $\chi$-values.

\noindent {\sl Case 4B.} $v$ does lie in $C$.

Then the discs $D'_1$ and $D'_2$ are nested: $D'_1$ lies in $D'_2$. The union of $D_2'$, $A$ and $D_1'$ pushed
a little into the interior of $M'$ is a separating 2-sphere. Again, we determine whether it separates components of
$\partial M'$. If it does not, it forms an expression of $E(M', \gamma')$ as a connected sum with a closed summand $Y$.
Again, if $Y$ is the 3-sphere, this can be certified in polynomial time.
In this case, we remove the component of $M' - A$ that contains the component of ${\rm int}({\cal B}')$ incident to $A$. The handles of $M' - {\cal B}'$
that lie in this component are removed, as are any components of ${\cal B}'$. The connectivity graph
is modified accordingly, by removing vertices and edges. Similarly, vertices and edges of $X_\partial$
are removed. The vertex $v$ is removed. But the $G$-vertex incident to $e'$ is retained, and its
$\chi$-value is increased by $1$. As shown in Lemma 7.3 in [21], this decreases the complexity of the
handle structure.

\noindent {\sl Case 5.} Every vertex of $X_\partial$ has non-positive $\chi$-value, and for some pair
of edges of $X_\partial$, both lie in $R_-(M')$ or both lie in $R_+(M')$.

Let $e_1$ and $e_2$ be this pair of edges. Let $A$ be the annulus of $\partial_v {\cal B}'$ that
contains the corresponding curves. It is properly embedded in $M'$, since a component of 
$\partial_v {\cal B}'$ lying in $\partial M'$ must have one boundary component in $R_-(M')$
and the other boundary component in $R_+(M')$.
Let $A_1$ and $A_2$ be parallel copies of this annulus, incoherently
oriented in such a way that the region between them inherits four sutures. Isotope $A_1$
and $A_2$ a little so that they become standard surfaces.

\noindent {\sl Case 5A.} $e_1$ and $e_2$ point into the same component of $X_\partial - \{ e_1, e_2 \}$
and this has $\chi$-value $0$.

Then (as shown in Case 5B(i) of Section 8 in [21]), we decompose along the component of $A_1 \cup A_2$ that
is closest to ${\cal B}'$. 

\noindent {\sl Case 5B.} $e_1$ and $e_2$ point out of the same component of $X_\partial - \{ e_1, e_2 \}$
and this has $\chi$-value $0$.

Then we decompose along the component of $A_1 \cup A_2$ that
is furthest from ${\cal B}'$. 

Note that it is not possible for both Case 5A and Case 5B to hold because of our hypothesis that no component of $M'$
has boundary a single torus with no sutures.

\noindent {\sl Case 5C.} Neither Case 5A nor Case 5B applies.

Then we decompose along $A_1 \cup A_2$. 

In each of these cases, the decomposition of $E(M' , \gamma')$ along these annuli is taut if and only if $E(M', \gamma')$ is taut, for the
following reason. If the decomposition is taut, then, by definition $E(M', \gamma')$ is taut.
So, suppose that $E(M', \gamma')$ is taut. One cannot apply Proposition 2.5, because the decomposition is
not along product annuli. However, the argument to establish tautness is given in Case 5 of Section 8 in [21].
Briefly, there it is shown that the only way that the decomposition can fail to be taut is if the resulting
sutured manifold has a solid toral component with no sutures. Similarly, Proposition 4.9 implies that
if $(M' ,\gamma')$ admits an allowable hierarchy, then so does the manifold obtained by decomposing
along non-trivial annuli disjoint from the sutures, unless they separate off a solid torus $V$ with no sutures.
The boundary of $V$ intersects
$\partial M'$ in a collection of annuli. Each annulus gives a component of 
$X_\partial - \{ e_1, e_2 \}$ and this has $\chi$-value $0$. The various arrangements for such
a solid torus are considered in Case 5 of Section 8 in [21]. In each case,
they are avoided by the suitable choice of decomposing surface, either $A_1$, or $A_2$, or $A_1 \cup A_2$.

It is also shown in Case 5 of Section 8 in [21] that these decompositions do not increase the complexity
of any 0-handle. Moreover, if the complexity of a 0-handle is unchanged, then it is only
modified by performing a trivial modification.

At the end of this process, we have created a sutured manifold $(M''', \gamma'')$ with a handle structure ${\cal H}'''$
that satisfies all the conditions of Theorem 10.1, except possibly (viii), which requires that ${\cal H}'''$
be positive. We now explain how to guarantee this extra condition.

\vskip 6pt
\noindent {\caps 10.5. Arranging positivity}
\vskip 6pt

The procedure that we will use will decrease the number of handles in the handle structure. 
Hence, we must initially find an upper bound on the number of handles in our handle structure
for $(M''', \gamma''')$. Note the number of handles of ${\cal H}'$ that are not parallelity handles
is at most $ch$, for some universal computable constant $c$. Each of the parallelity handles, other than isolated 2-handles
incident to no other parallelity handle, is removed
in the above process. Some are replaced with 2-handles.
Hence, the number of handles in the resulting handle structure ${\cal H}'''$ for $(M''', \gamma''')$
is at most $ch$ plus the number of components of $\partial_v {\cal B}$. This is at most
$c'h$ for some universal computable constant $c'$.

The first stage in the procedure is to arrange that ${\cal H}'''$ has no 3-handles. This is achieved using Procedure 6 of Section 6.1.
Since no component of $M'''$ is closed, then if ${\cal H}'''$ contains a 3-handle, then it contains a 
2-handle that intersects ${\cal H}^3$ in a single disc. Hence, Procedure 6 may be applied until there
are no more 3-handles. None of the later modifications will introduce 3-handles, and so we
may henceforth assume that none of the handle structures that we consider contains them.

Suppose that the handle structure is not positive. Then, as explained in Lemma 7.6 of [21],
one of Procedures 1-6 of Section 6.1 can be applied, or the handle structure contains a parallelity handle
that is not an isolated 2-handle. However, there are two slight discrepancies between the procedures referred
to in Lemma 7.6 of [21] and Procedures 1-6 of Section 6.1. Firstly,
a slight variant of Procedure 1 is used in [21], but the version that we present here works just
as well; it has the advantage that one does not need to recognise 3-balls algorithmically.
Secondly, Procedure 4 may be applied here only when the relevant 1-handle $H_1$ of ${\cal H}'''$
is disjoint from the u-sutures, whereas in [21] there was no such constraint. If just one arc
of $H_1 \cap \gamma$ lies in a u-suture or the arcs of $H_1 \cap \gamma$ lie in distinct
u-sutures, then we attach a 2-handle along one of these u-sutures. The resulting manifold
still has the same canonical extension and it still admits an allowable hierarchy. However,
we have introduced a handle and so we may have increased the complexity of the handle
structure. But we can immediately apply Procedure 3, which cancels $H_1$ and the new 2-handle.
The resulting handle structure has no greater complexity than the original one, and it has
fewer handles. If both arcs of $H_1 \cap \gamma$ lie in the same u-suture, then we still attach
a 2-handle along it. Then $H_1$ becomes a parallelity handle disjoint from $\gamma$. So we apply the procedures given in Sections 10.1, 10.2 and 10.3
again, except that we do not use Theorem 9.2. We again define ${\cal B}'$ to be
the union of the parallelity handles that are not 2-handles. We calculate the boundary graph
$X_\partial$ and connectivity graph $X_c$ directly, without using Theorem 9.2. This is possible
because the handle structure has at most $c'h$ handles. Then we apply the algorithm given
in Section 10.4. 

So, we proceed initially, by applying Procedures 1-6 of Section 6.1 as many times as possible. 
If we come across a 1-handle intersecting $\gamma$ in two arcs and disjoint from the 2-handles,
but which intersects at least one u-suture, then we proceed as above.

This process may create new parallelity handles. If so, we apply the procedures given in Sections 10.1, 10.2 and 10.3
again, except that we do not use Theorem 9.2. Iterating in this way,
we eventually reach the desired handle structure ${\cal H}'''$.

\vskip 6pt
\noindent {\caps 10.6. Decomposition along normal surfaces}
\vskip 6pt

In this section, we give a variation of Theorem 10.1. Instead of decomposing a handle structure along a regulated surface,
we decompose a triangulation along a normal surface.

\noindent {\bf Theorem 10.2.} {\sl Let ${\cal T}$ be a triangulation of a connected compact orientable 3-manifold $M$ with $t$ tetrahedra.
Suppose that every boundary component of $M$ is a torus that is triangulated using just two triangles. Give $M$ a sutured manifold
structure $(M, \gamma)$ where $R_-(M) = \partial M$ and $R_+(M) = \emptyset$. Let $S$ be a taut normal surface properly embedded in $M$,
such that the intersection between $S$ and each component of $\partial M$ is a (possibly empty) collection of coherently oriented essential curves.
Let $(M', \gamma')$ be obtained by decomposing along $S$. 
Then, there is an algorithm that takes as its input ${\cal T}$ and $(S)$ and either correctly
asserts that $E(M', \gamma')$ or $E(M,\gamma)$ is not taut, or supplies
a handle structure ${\cal H}'''$ for a sutured manifold $(M''', \gamma''')$ with properties (i) - (viii) of Theorem 10.1.
The algorithm runs in time that is bounded by a polynomial in $t \log(w(S))$.}

\noindent {\sl Proof.} Each component of $\partial M$ is a torus triangulated using two triangles. This hypothesis
is made to ensure that $\partial S$ is very controlled. It is shown in Lemma 3.5 in [14] that two normal curves
in such a torus are isotopic if and only if they are normally isotopic. Furthermore, in Theorem 3.6 of [14],
a parametrisation is given for all possible normal curves in the torus, and a consequence of this is
that the intersection between any essential normal simple closed curve and any edge in the torus is
a collection of points all with the same sign. We are assuming that the intersection between $\partial S$
and any component of $\partial M$ is a collection of coherently oriented essential curves. Thus,
we deduce that the points of intersection between $\partial S$ and any edge in $\partial M$ all have the
same sign. 

We now possibly modify ${\cal T}$ so that each tetrahedron $\Delta$ of ${\cal T}$ contains at most one
2-simplex lying in $\partial M$. The only other alternative is that $\Delta$ contains two 2-simplices 
which comprise the entirety of a component of $\partial M$. In this case, let $e$ be the edge
of $\Delta$ that does not lie in these two 2-simplices. If $e$ also lies in $\partial M$, then a loop
encircling one of the faces of $\Delta$ not lying in $\partial M$ is essential in $\partial M$ but homotopically
trivial in $M$. Hence, in this case $E(M,\gamma)$ is not taut,
and the algorithm terminates by declaring this. On the other hand, if $e$ does not lie in $\partial M$,
then we may remove this tetrahedron $\Delta$ from ${\cal T}$, 
and the result is a new triangulation ${\cal T}'$ for a 3-manifold $M'$ which is a copy of $M$. 
This new triangulation again satisfies the condition that every boundary component of $M$ 
is a torus that is triangulated using just two triangles. The surface $S \cap M'$ is normal in
${\cal T}'$. We claim that there is a homeomorphism
from $M$ to $M'$ which takes $S$ to $M' \cap S$. This will imply that $M' \cap S$ satisfies
the hypotheses of the theorem.

We can think of $\Delta$ as forming a collar on $\partial M'$. In particular, there is a natural
product structure on $\Delta$. We will show that $S \cap \Delta$ respects this product structure,
which will establish the claim. There are four possible triangle types in $\Delta$, and each of these
respects its product structure. There are three possible square types, two of which respect the
product structure. However, the third square type gives rise to a boundary-compression disc for $S$,
which contradicts the hypothesis that $S$ is taut.

So, we may assume that each tetrahedron $\Delta$ of ${\cal T}$ contains at most one
2-simplex lying in $\partial M$. We now dualise ${\cal T}$ to form a handle structure ${\cal H}$.
Recall that each $i$-handle of ${\cal H}$ arises from a $(3-i)$-simplex of ${\cal T}$ that does
not lie entirely in $\partial M$. Hence, every 0-handle of ${\cal H}$ is subtetrahedral.
In particular, a tetrahedron that intersects $\partial M$ in a triangle is dual to a 0-handle $H_0$
such that $H_0 \cap {\cal F}$ is as shown in Figure 24.

The normal surface $S$ determines a standard surface in ${\cal H}$. We will call this new surface $S'$,
although it is just a copy of $S$. Note that $S'$ need be neither normal nor regulated with respect to ${\cal H}$.
For example, in Figure 24, an elementary disc in a 0-handle $H_0$ is shown which is dual to a normal triangle of $S$.
This disc fails to satisfy Condition $3'$ in the definition of regulated and it fails to satisfy (i) in the definition 
of normality.

\vskip 6pt
\centerline{
\epsfxsize=3.5in
\epsfbox{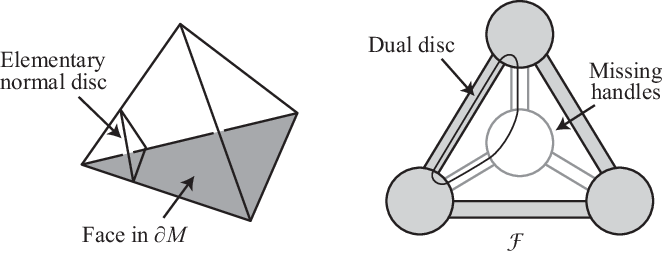}
}
\vskip 6pt
\centerline{Figure 24: Dual handle structure near $\partial M$}

Let $M'$ be the manifold obtained by decomposing along $S'$, and let ${\cal H}'$ be the handle structure that it
inherits from ${\cal H}$. Consider the parallelity handles of ${\cal H}$ for the pair $(M', N(S') \cap M')$. As explained
in Section 9.2, these arise precisely in the space between two parallel normal discs of $S$. Hence,
in each 0-handle $H_0$ of ${\cal H}$, at most $6$ handles of $H_0 \cap {\cal H}$ are not parallelity handles
in this sense. If we let ${\cal B}$ denote the union of these parallelity handles, we can use Theorem 9.3
to compute a handle structure for ${\rm cl}(M' - {\cal B})$ and, for each component $B$ of ${\cal B}$,
the genus and number of boundary components of the base surface for $B$, whether it is a product or
twisted $I$-bundle and the location of $\partial_v B$. 

However, there is another notion of parallelity handle, given in Section 9.1, for the sutured manifold $(M', \gamma')$.
We claim that every parallelity handle $H$ for the pair $(M', N(S') \cap M')$ is also a parallelity handle for $(M', \gamma')$.
In the case where $H$ is disjoint from $\partial M$, this is automatic. In the case where $H$ intersects $\partial M$,
it lies between two elementary normal discs of $S$ that intersect $\partial M$, and we have arranged
that these two discs are therefore coherently oriented. This is because the points of intersection between
$S$ and any edge in $\partial M$ all have the same sign. Hence, $H$ does indeed satisfy (ii) or (iii) in the
definition in Section 9.1. Let ${\cal B}'$ be the union of the parallelity handles for $(M', \gamma')$.
Then, we have shown that ${\cal B} \subseteq {\cal B}'$.

We can readily compute the handles of ${\rm cl}(M' - {\cal B})$ that lie in ${\cal B}'$. Therefore, we can determine
a handle structure for ${\rm cl}(M' - {\cal B}')$ and, for each component $B$ of ${\cal B}'$,
the genus and number of boundary components of the base surface for $B$, whether it is a product or
twisted $I$-bundle and the location of $\partial_v B$. 

This was exactly the information that we required to be able to proceed with the proof of Theorem 10.1.
The remainder of the proof proceeds exactly as in that case. $\square$

\vskip 6pt
\noindent {\caps 10.7. Modifying a triangulation to simplify the boundary tori}
\vskip 6pt

A key hypothesis in Theorem 10.2 was that the triangulation of the 3-manifold restricts to a triangulation of
each boundary torus with only two triangles. In the following result, we show how to arrange that
we have a triangulation of this form.

\noindent {\bf Proposition 10.3.} {\sl Let $M$ be a compact orientable 3-manifold with boundary a union of tori.
Let ${\cal T}$ be a triangulation of $M$ with $t$ tetrahedra. Then there is an algorithm to construct
a triangulation ${\cal T}'$ for $M$ with at most $5t$ tetrahedra, and where every boundary torus is triangulated
using two triangles. Moreover, if $\phi$ is a simplicial 1-cocycle on ${\cal T}$, then the algorithm computes a
simplicial 1-cocycle $\phi$ on ${\cal T}'$ in the same cohomology class as $\phi$ and satisfying $||\phi'||_1 \leq 3^{4t} ||\phi||_1$.
This algorithm runs in time that is bounded by a polynomial function of $t \log ||\phi||_1$.}

\noindent {\sl Proof.} We will modify ${\cal T}$ by attaching tetrahedra to $\partial M$. This will have the
effect of changing the boundary triangulation by a Pachner move. Each time we do this, there is a
simple way to extend $\phi$ to a cocycle $\phi''$ on the larger triangulation. For most ways of attaching a tetrahedron,
$\phi''$ is determined by the fact that has to be a cocycle. But when we attach a tetrahedron along a single triangle,
there is some choice, and so we simply declare that one of the new edges has evaluation zero under $\phi''$.
It is easy to check the sum of the absolute values of the evaluations of the new edges is at most the sum of the values
of $\phi$ on two edges in $\partial M$. Hence, $||\phi''||_1 \leq 3 ||\phi||_1$. We will show that we have to add at most
$4t$ tetrahedra, and hence the final cocycle $\phi'$ satisfies $||\phi'||_1 \leq 3^{4t} ||\phi||_1$. This will also show
that the final triangulation has at most $5t$ tetrahedra.

It is well known that any two triangulations of a surface differ by a sequence of Pachner moves. However, we wish to do
this algorithmically, and with control over the number of moves. Suppose that a triangulation of a component of $\partial M$
has more than one vertex. Our strategy will be to apply at most $4$ Pachner moves to this triangulation to create
a triangulation with fewer vertices. Hence, the total number of Pachner moves that we will apply is at most 4 times
the number of vertices in the boundary. For Euler characteristic reasons, this is half the number of triangles in the
boundary. But each tetrahedron of $M$ contributes at most $2$ faces to the boundary. So, the total number of Pachner
moves will be at most $4t$.

Suppose first that there is an edge in $\partial M$ with the same triangle on both sides. The two edges of the triangle are then
identified. Because $\partial M$ is orientable, the edges are incompatibly oriented as one encircles the triangle. Hence,
their common endpoint ends up as a valence one vertex $v$ in $\partial M$. The remaining edge $e$ of the triangle has
a distinct triangle on the other side of it. Thus we may perform a 2-2 Pachner move centred at $e$. This has the effect of
removing $e$ and increasing the valence of $v$ to $2$. Incident to $v$ there are now two distinct triangles. For each of
these triangles, one of its edges is not incident to $v$. These edges cannot be the same edge of the triangulation,
because this would form a triangulation of the 2-sphere with two triangles. Hence, we may perform a 2-2 Pachner move
along one of these edges. This increases the valence of $v$ to 3. We can then perform
a 3-1 Pachner move at $v$, and this reduces the number of vertices in $\partial M$ by $1$, as required.
So, we may assume that every edge in $\partial M$ is adjacent to two distinct triangles, and hence one may perform
a 2-2 Pachner move along it. 

Using the fact that the torus has zero Euler characteristic, we can find a vertex $v$ of its triangulation with valence at most $6$. 
Our aim is reduce the valence of this vertex to $3$ and then apply a 3-1 Pachner move to remove it.

Suppose first that there is an edge $e$ that starts and ends at $v$. If we perform a 2-2 Pachner move along $e$, then 
the two endpoints of $e$ no longer contribute to the valence. Thus, we reduce the valence of $v$ unless the new edge
also starts and ends at $v$. In this case, the two triangles adjacent to $e$ form a square, and the four sides of the
square must be identified in pairs, in order that the valence of $v$ is at most $6$. Thus, in this case, we deduce
that the triangulation of the torus already has a single vertex, which is contrary to our assumption.

On the other hand, if no edge starts and ends at $v$, then the star of $v$ is an open disc. We may then apply at most 
three 2-2 Pachner moves supported in this disc to reduce the valence of $v$ to $3$. We then can apply a 3-1
Pachner move to remove $v$.

Thus in all cases, after at most $4$ Pachner moves, we can reduce the number of vertices by $1$. $\square$

\vskip 18pt
\centerline{\caps 11. Reduction to the atoroidal and Seifert fibred cases}
\vskip 6pt

At certain points in the previous sections, we made the assumption that the manifold is irreducible and atoroidal.
In this section, we show why it suffices to consider this case and also the case where the manifold is Seifert fibred.

\vskip 6pt
\noindent {\caps 11.1. The canonical tori}
\vskip 6pt

Recall that a surface $S$ properly embedded in a compact orientable 3-manifold $M$ is {\sl essential}
if it is incompressible, boundary-incompressible and no component is parallel to a subsurface of $\partial M$.

A properly embedded torus in $M$ is {\sl canonical} if it is essential and furthermore it can be ambient isotoped off any
other essential torus. 

Canonical tori are also called JSJ tori, since Jaco, Shalen and Johannson developed their theory.
In particular, they proved the following result [15, 17].

\noindent {\bf Theorem 11.1.} {\sl Let $M$ be a compact orientable irreducible 3-manifold with incompressible boundary.
Then up to ambient isotopy, there exists finitely many canonical tori in $M$. The union $T$ of these tori is properly embedded.
Each component of $M - {\rm int}(N(T))$ is either atoroidal or Seifert fibred.}

\vskip 6pt
\noindent {\caps 11.2. The behaviour of Thurston norm upon cutting along spheres and essential tori}
\vskip 6pt

\noindent {\bf Proposition 11.2.} {\sl Let $M$ be a compact orientable 3-manifold with incompressible boundary.
Let $T$ be a union of finitely many disjoint spheres and incompressible tori properly embedded in $M$. 
Let $c$ be a class in $H^1(M)$, and let $i^\ast(c)$ be its image under the homomorphism
$i^\ast \colon H^1(M) \rightarrow H^1(M - {\rm int}(N(T))$ induced by inclusion. Then the Poincar\'e duals of $c$
and $i^\ast(c)$ have equal Thurston norm.}

\noindent {\sl Proof.} Let $S$ be a representative for the Poincar\'e dual of $c$, with minimal Thurston complexity.
If $S \cap T$ contains any simple closed curves that bound discs in $T$, then these may be removed without
increasing the Thurston complexity of $S$. So, we may assume that $S \cap T$ is essential in $T$.
Hence, by the incompressibility of $T$, none of these curves bounds a disc in $M$. Hence, the spheres and discs
of $S$ are precisely the spheres and discs of $S - {\rm int}(N(T))$. Therefore,
$\chi_-(S)  = \chi_-(S- {\rm int}(N(T)))$.
Hence, the Thurston norm of $[S - {\rm int}(N(T))]$ is at most that of $[S]$.
Note that $S - {\rm int}(N(T))$ is a representative of the Poincar\'e dual to $i^\ast(c)$.

Now consider a surface $S'$ representing the Poincar\'e dual to $i^\ast(c)$ with minimal Thurston complexity.
We may assume that $S' \cap \partial N(T)$ is a collection of essential curves. In particular, $S'$ is then disjoint
from the spherical components of $\partial N(T)$. Moreover, we may arrange,
by attaching annuli if necessary, that $S'$ intersects each toral component of $\partial N(T)$ in a collection
of coherently oriented curves. We wish to use $S'$
to build a representative for the Poincar\'e dual to $S$. We first note that if $T_1$ and $T_2$
are the components of $\partial N(T)$ coming from the same component of $T$, then
$S' \cap T_1$ and $S' \cap T_2$ patch together under the gluing map $T_1 \rightarrow T_2$.
This is because $S' \cap T_1$ and $S' \cap T_2$ are determined up to isotopy by the restriction of $c$ to $H^1(T)$.
Hence, when we patch together the tori in $\partial N(T)$ to form $M$, we can get a properly embedded
surface $S$. However, it is not immediately clear that $S$ is Poincar\'e dual to $c$.
Indeed, all we can say is that $c$ and the Poincar\'e dual to $[S]$ have the same image under $i^\ast$.
We have an exact sequence
$$H^1(M, M - {\rm int}(N(T))) \rightarrow H^1(M) \buildrel i^\ast \over \rightarrow H^1(M - {\rm int}(N(T))).$$
Now, $H^1(M, M - {\rm int}(N(T)))$ is isomorphic to $H^1(N(T), \partial N(T))$ by excision, and this is
isomorphic to $H_2(N(T))$ by Poincar\'e duality. Hence, the difference between $[S]$ and the Poincar\'e dual to $c$ 
is given by some copies of toral components of $T$. But we can introduce these components, and form their
double-curve sum with $S$, forming a surface $S''$, say. Then $\chi_-(S'') = \chi_-(S) = \chi_-(S')$.
Hence, the Thurston norm of the Poincar\'e dual to $c$ is at most that of the Poincar\'e dual to $i^\ast(c)$.
$\square$

\vskip 6pt
\noindent {\caps 11.3. The canonical tori in normal form}
\vskip 6pt

The canonical tori may be placed in normal form with respect to any triangulation. In fact, it has long been known
that one may control the weight of this normal surface. For example, an algorithm was given by Jaco and Tollefson
for constructing the canonical tori [16]. An explicit upper
bound on the weight of the canonical tori in normal form was given by Mijatovi\'c (Proposition 2.4 in [24]), as follows.

\noindent {\bf Theorem 11.3.} {\sl Let $M$ be a compact orientable irreducible 3-manifold with incompressible boundary.
Let ${\cal T}$ be a triangulation of $M$ with $t$ tetrahedra. Then the union of the canonical tori may be placed
in normal form with respect to ${\cal T}$, so that it contains at most $2^{80t^2}$ elementary normal discs.}

\vfill\eject
\noindent {\caps 11.4. A triangulation for the exterior of normal tori and spheres}
\vskip 6pt

\noindent {\bf Theorem 11.4.} {\sl Let $M$ be a compact orientable 3-manifold.
Let ${\cal T}$ be a triangulation of $M$ with $t$ tetrahedra. Let $T$ be a union of disjoint properly embedded normal spheres and tori.
Then there is an algorithm to build a triangulation ${\cal T}'$ for $M' = M - {\rm int}(N(T))$,
which runs in time at most a polynomial function of $t \log w(T)$, and where this triangulation has
at most $200t$ tetrahedra.}

\noindent {\sl Proof.} Let ${\cal H}$ be the handle structure of $M$ that is dual to ${\cal T}$.
The normal surface $T$ in ${\cal T}$ determines a standard surface, also called $T$, in ${\cal H}$.
So $M - {\rm int}(N(T))$ inherits a handle structure ${\cal H}'$. Let ${\cal B}$ be the parallelity 
bundle of ${\cal H}'$ for the pair $(M', \partial N(T))$. Let $\partial_v {\cal B}$ be its vertical boundary.

For each 0-handle $H_0$ of
${\cal H}$, all but at most 6 handles of $H_0 \cap {\cal H}'$ are parallelity handles.
Similarly, for each 1-handle $H_1$ of
${\cal H}$, all but at most 4 handles of $H_1 \cap {\cal H}'$ are parallelity handles.
We may triangulate the union of these handles using at most 60 (respectively, 40) tetrahedra, so that the
triangulations match on the intersection between handles and so that the intersection
with $\partial_v{\cal B}$ is simplicial.

It is easy to check that, for each 0-handle $H_0$ of ${\cal H}$, the number of
components of $H_0 \cap \partial_v{\cal B}$ is at most 10. Hence, $\partial_v{\cal B}$ has
at most $10 t$ components.

We now apply Theorem 9.3, which provides an algorithm that determines, for each component $B$
of ${\cal B}$, the genus and number of boundary components of its base surface, whether 
$B$ is a product or twisted I-bundle, and the location of $\partial_v B$ in ${\rm cl}(M' - {\cal B})$.
Using this information, we may build a triangulation for each component of ${\cal B}$
that agrees with the one that we have built for $\partial_v {\cal B}$. Since the genus of the base
surface of each component of ${\cal B}$ is at most one, and the number of boundary components is at most
$10 t$, we may arrange that the number of tetrahedra in this triangulation of ${\cal B}$ is
at most $100 t$. Attaching this to the triangulation of ${\rm cl}(M' - {\cal B})$ already
constructed, we obtain the required triangulation of $M'$. 

As a result of Theorem 9.3, each of these steps takes time that is bounded above by a polynomial function of 
$t \log w(T)$. $\square$

Theorem 11.4 gives an algorithm to construct a triangulation for the
exterior of the canonical tori, as long as that these tori are provided non-deterministically
as a normal surface.

Note however, we are not at this stage claiming to be able to verify that a given collection
of normal tori $T$ is the canonical collection. In particular, we have not checked that $T$ is
incompressible, nor have we checked that each component of the exterior of $T$ is atoroidal or Seifert fibred.

We will use Theorem 11.4 to reduce the main theorem to the case where the
3-manifold is irreducible and atoroidal or Seifert fibred. But to be able to do this, we need to be able
to restrict the given 1-cocycle on $M$ to an explicit 1-cocycle on $M'$.

\noindent {\bf Addendum 11.5.} {\sl Let $M$, ${\cal T}$, $t$, ${\cal T}'$ and $M'$ be
as in Theorem 11.4. Let $\phi$ be a simplicial 1-cocycle on ${\cal T}$. Let $i \colon M' \rightarrow M$
be inclusion. Then we may construct a simplicial 1-cocycle $\phi'$ on ${\cal T}'$ with $||\phi'||_1 \leq 1200 t ||\phi||_1$, where $|| \ ||_1$ denotes the
$\ell^1$ norm of a cocycle and such that the Poincar\'e duals of $i^\ast([\phi])$ and $[\phi']$ have the same Thurston norm.
Moreover, when $T$ consists solely of spheres, then $i^\ast([\phi])$ and $[\phi']$ are cohomologous.
This can be achieved in time at most a polynomial function of $t \log w(T) \log||\phi||_1$.
}

\noindent {\sl Proof.} A cocycle representing $i^\ast([\phi])$ may be formed by first dualising $\phi$ to form a surface
$S$ in $M$, then intersecting $S$ with $M'$ to form a surface $S'$, then dualising this to
give $i^\ast([\phi])$. In order to complete this successfully, we first control the position of $S$ in $M$.
From the cocycle $\phi$, there is a natural way of building $S$, as follows. For each edge $e$ of the
triangulation, one orients $e$ so that $\phi(e)$ is non-negative and one places $|\phi(e)|$ points in the interior of $e$. For each face of the triangulation,
one inserts arcs joining these points, in such a way that arcs start and end on coherently oriented edges. 
This is possible because $\phi$ is a cocycle,
and so the total evaluation of the edges encircling the face is zero. For each tetrahedron,
we have now specified a collection of simple closed curves in its boundary. Each bounds
an elementary normal disc. The union of these discs we take to be $S$.

We perform a normal isotopy of $S$ so that its intersection with $T$ is in general position.
More specifically, we arrange that the intersection between any normal disc of $S$ and any normal disc of $T$
is at most one arc. We may also arrange that $S \cap {\cal B}$ is vertical in ${\cal B}$.
Then $S - {\rm int}(N(T))$ is dual to $i^\ast ([\phi])$.

We now need to specify the triangulation ${\cal T}'$ a little more precisely than in the proof of Theorem 11.4.
The parts of ${\rm cl}(M' - {\cal B})$ were triangulated in a relatively simple way. 
It is clear that we may choose this triangulation so that, for each edge in this part of the triangulation, 
the number of times it intersects $S - {\rm int}(N(T))$
is at most the number of elementary discs of $S$ in the tetrahedron containing it.
We take this algebraic intersection number to be $\phi'$ on that edge.
So, the evaluation of each such edge under $\phi'$ is at most $||\phi||_1$.

The remainder of ${\cal T}'$ was formed by triangulating ${\cal B}$, which is an $I$-bundle
over a surface $F$. Using Theorem 9.3, we know the topological type for each component of $F$.
We pick a triangulation for $F$ with a certain number of vertices in its boundary and none in its interior. 
Then, for each simplex $\sigma$ of this triangulation,
its inverse image in ${\cal B}$ is of the form $\sigma \times I$. We can triangulate this in a simple way.
When $\sigma$ is a 1-simplex, we triangulate $\sigma \times I$ using two 2-simplices.
When $\sigma$ is a 2-simplex, we triangulate $\sigma \times I$ using three 3-simplices.
Note that $\partial_v {\cal B}$ then inherits a triangulation, which agrees with the one
arising from the triangulation on ${\rm cl}(M' - {\cal B})$, as long the number of
vertices in each component of $\partial F$ is chosen correctly. 

Thus, we have triangulated $M'$ and have defined $\phi'$ on ${\rm cl}(M' - {\cal B})$.
We still need to define $\phi'$ on ${\cal B}$. Now $S - {\rm int}(N(T))$ is vertical in ${\cal B}$.
But to determine this precise location of this surface might take too long. So, instead,
we pick any vertical transversely oriented surface in ${\cal B}$ that agrees with $S$ on $\partial_v {\cal B}$.
This can be achieved by choosing arcs in $F$ transverse to the 1-skeleton of $F$ and that join
the points of $S \cap F$ with the correct orientations. These arcs can be chosen so that the number of times they intersect
each edge of $F$ is at most $|S \cap \partial_v{\cal B}|$. They then specify a collection of transversely oriented vertical discs
in ${\cal B}$. We define $\phi'$ on the edges in ${\cal B}$ to be the intersection numbers of these edges with these discs.
Hence, $||\phi'||_1 \leq 1200 t ||\phi||_1$, because there are at most $1200 t$ edges of ${\cal T}'$ and each
edge has evaluation at most $||\phi||_1$ under $\phi'$.

Now, $[\phi']$ and $i^\ast([\phi])$ need not be equal cohomology classes, because their dual surfaces may differ
in ${\cal B}$. But their difference is represented by a union of vertical annuli in ${\cal B}$. Hence,
the duals have the same Thurston norm. Moreover, when $T$ consists solely of spheres, these vertical
annuli in ${\cal B}$ are homologically trivial, and hence $[\phi']$ and $i^\ast([\phi])$ are indeed the
same cohomology classes. $\square$

\vskip 18pt
\centerline {\caps 12. Certification for products and Seifert fibre spaces}
\vskip 6pt

\noindent {\caps 12.1. Products}
\vskip 6pt

Throughout this paper, we have used allowable hierarchies, which terminate in product sutured manifolds.
It will be important that we can certify that a given sutured manifold is indeed a product. The existence
of such a certificate is presented in this subsection.

\noindent {\bf Theorem 12.1.} {\sl Let $(M, \gamma)$ be a product sutured manifold with handle structure
${\cal H}$ of uniform type. Let $h$ be the number of 0-handles of ${\cal H}$. Suppose that $(M, \gamma)$
is a product sutured manifold. Then there is a certificate that proves that $(M, \gamma)$ is indeed a
product. This can be verified in time that is bounded above by a polynomial function of $h$.}

A {\sl complete collection} of product discs for a product sutured manifold $(M, \gamma)$ is a collection
of disjoint product discs $D$, such that decomposing along $D$ gives a collection of taut 3-balls.

\noindent {\bf Theorem 12.2.} {\sl Let $(M, \gamma)$ be a sutured manifold with handle structure
${\cal H}$ of uniform type. Then $(M, \gamma)$ contains a complete collection $D$ of product discs in normal form
with weight satisfying $w(D) \leq c^{h^2}$, where $h$ is the number of 0-handles of ${\cal H}$
and $c$ is a universal computable constant.}

We will prove this result in the next subsection. But first, we show how it can be used to prove Theorem 12.1.

We are given a handle structure ${\cal H}$ for the product sutured manifold $(M, \gamma)$. We need
to certify that it is indeed a product. The certificate that we use will be of the following form:
\item{(i)} a complete collection $D$ of product discs as in Theorem 12.2;
\item{(ii)} a handle structure ${\cal H}'$ for a sutured $(M', \gamma')$ with at most $ch$ handles; 
in fact, this will be the sutured manifold obtained by decomposing $(M, \gamma)$ along $D$;
\item{(iii)} a certificate that $M'$ is a collection of 3-balls, as provided by Schleimer [26] or Ivanov [12].

The algorithm to verify this certificate is as follows:
\item{(1)} Verify that $D$ is a collection of product discs, using the algorithm of Agol-Hass-Thurston to
determine the topological types of the components of $D$, and to verify that each component
intersects $\gamma$ exactly twice.
\item{(2)} Denote the manifold obtained by decomposing $(M, \gamma)$ along $D$ by $(M_D, \gamma_D)$.
This inherits a handle structure ${\cal H}_D$. Denote the parallelity bundle of the pair $(M_D, M_D \cap N(D))$
by ${\cal B}$.
Apply the algorithm given in Theorem 9.3 (adapted in the obvious way to handle structures of uniform type)
to determine a handle structure for ${\rm cl}(M_D - {\cal B})$ and to determine, for each component $B$
of the parallelity bundle ${\cal B}$, the number of boundary components of its base surface,
whether $B$ is a product or twisted $I$-bundle, and the location of $\partial_v B$
in ${\rm cl}(M_D - {\cal B})$. The total number of 0-handles in the handle structure of 
${\rm cl}(M_D - {\cal B})$ is at most $ch$ for some universal computable constant $c$.
\item{(3)} Using the information provided in (2), express ${\cal B}$ as a union of
parallelity 1-handles and 2-handles, and thereby form a handle structure for $(M_D, \gamma_D)$.
Verify that this is equal to the given handle structure ${\cal H}'$. 
\item{(4)} A verification of the certificate that $M'$ is a collection of 3-balls, using the algorithm of Schleimer [26] or Ivanov [12],
together with a verification that each component of $R_\pm(M')$ is a disc.

The time taken to complete this verification is at most a polynomial function of $h \log(w(D))$,
which by Theorem 12.2, at most a polynomial function of $h$.

\vskip 6pt
\noindent {\caps 12.2. Product discs in normal form}
\vskip 6pt

In this subsection, we prove Theorem 12.2. Our method is rather similar to that used by Mijatovi\'c in the proof 
of Theorem 11.3, but somewhat more straightforward.

Let $(M, \gamma)$ be a connected product sutured manifold with a handle structure ${\cal H}$
of uniform type. We may suppose that it is not a 3-ball,
as in this case, we may take $D$ to be empty.

The first step in the proof is to show that there is some essential product disc $P_1$ that is
in normal form with respect to ${\cal H}$ and that is fundamental. This is well known. A proof,
in the related case where $M$ is triangulated, is given in Lemma 4.1.18 of [23]. Hence, by Theorem 8.1,
the weight $w(P_1)$ is at most $(c_1)^h$, where $h$ is the number of 0-handles
of ${\cal H}$ and $c_1$ is a universal computable constant.

We now decompose $(M, \gamma)$ along $P_1$, giving a manifold $(M_2, \gamma_2)$. This inherits
a handle structure ${\cal H}_2$. We denote by ${\cal B}$ the parallelity bundle for the pair
$(M_2, M_2 \cap N(P_1))$ with handle structure ${\cal H}_2$.

Our next goal is to find an essential product disc $P_2$ in $(M_2, \gamma_2)$ that is disjoint from
the copies of $P_1$ in $\partial M_2$. It therefore is an essential product disc in $(M, \gamma)$ also.
Our aim to ensure that $P_2$ is in fact normal in ${\cal H}$, and that there is a
good upper bound on its weight. Simple normal surface theory applied to the handle structure
${\cal H}_2$ will not suffice. This is because the bound on the weight of a 
fundamental normal surface in ${\cal H}_2$ given by Theorem 8.1 is of the
form $(c_1)^{h_2}$, where $h_2$ is the number of 0-handles of ${\cal H}_2$,
but $h_2$ may not be bounded by a linear function of $h$.

So, we develop a modified normal surface theory in ${\cal H}_2$. We consider normal surfaces
$S$ properly embedded in ${\cal H}$ that are disjoint from $P_1$. For any such surface $S$, 
its intersection with ${\cal B}$
is horizontal and so, for each component $B$ of ${\cal B}$, $B \cap S$ intersects every fibre of $B$
the same number of times. We call this the number of {\sl sheets} of $B \cap S$.

One can encode such a surface by means of a vector, which counts the number
of elementary discs of each type in each 0-handle in ${\rm cl}(M_2 - {\cal B})$ and the number of sheets of
$B \cap S$ for each component $B$ of ${\cal B}$. The equations
come in two types. Equations of the first type arise from 1-handles in ${\rm cl}(M_2 - {\cal B})$,
and they simply specify that the elementary discs in the adjacent 0-handles can be
attached to suitable elementary discs in the 1-handles. Equations of the second type
are associated with components $B$ of ${\cal B}$. They assert that, for each 0-handle $H_0$ in ${\rm cl}(M_2 - {\cal B})$ and each component
of $\partial_v B \cap H_0$, the number of arcs of intersection between
$S$ and this component is equal to the number of sheets of $B \cap S$. 

Because of this algebraic structure, one may speak of such a surface being
fundamental. We may find an essential product disc $P_2$
which is normal in ${\cal H}$ and fundamental in the above sense, using the proof of Lemma 4.1.18 of [23].
Then Theorem 8.1 gives an upper bound $(c_1)^h$ on the size of each
co-ordinate in the vector for $P_2$. However, care must be taken when
translating this to a bound on the weight of $P_2$ as a normal surface in ${\cal H}$,
because the part of $P_2$ running through ${\cal B}$ also contributes to this weight.
For each component $B$ of ${\cal B}$, the number of sheets of $B \cap P_2$ is bounded above by 
$(c_1)^h$. Each sheet of $B \cap P_2$
gives rise to a number of elementary discs of $P_2$, which is at most the sum of
the co-ordinates of $P_1$. So, $w(P_2) \leq (c_1)^h(c_1)^h$.

We now repeat by cutting $M_2$ along $P_2$ to give another 3-manifold $M_3$.
We can view this as the result of cutting $M$ along $P_1 \cup P_2$. The above
argument gives another essential product disc $P_3$, and the number of elementary
discs of $P_3$ is at most $(c_1)^{3h}$.

We repeat this process until we have no more essential product discs,
in other words, when we have decomposed $(M, \gamma)$ to a collection of balls.
The number of times that we iterated this procedure was at most $c_2 h$, for some universal
computable constant $c_2$, since this is an upper bound for the number
of disjoint non-parallel properly embedded normal surfaces in ${\cal H}$.
So, we end with a complete collection $D$ of essential normal product discs in 
$(M, \gamma)$ with total weight at most $((c_1)^h)^{c_2h+1}$, as required. $\square$

\vskip 6pt
\noindent {\caps 12.3. Seifert fibre spaces}
\vskip 6pt

As explained in the previous section, we will cut our given compact orientable irreducible 3-manifold along its canonical tori.
Each component of the resulting 3-manifold $M'$ will be atoroidal or Seifert fibred. The reason for doing this is
that many of the arguments given in Section 3 and 4 required the manifold to be atoroidal. However, that
still leaves the case of Seifert fibred manifolds. In this subsection, we will explain how our main theorem
is proved in this case.

\noindent {\bf Theorem 12.3.} {\sl Let $M$ be a connected Seifert fibre space other than a solid torus, with a triangulation ${\cal T}$ having $t$ tetrahedra.
Let $\phi$ be a simplicial 1-cocycle on $M$, and let $m$ be the Thurston norm of the Poincar\'e dual of $[\phi]$. 
Then there is a certificate that certifies that the Thurston norm of this class is indeed $m$ and that $M$ is irreducible and
has incompressible boundary. The algorithm to
verify this certificate can be completed in time that is bounded above by a polynomial function of $t \log ||\phi||_1 \log (m+2)$.}

Note that we do not require the Seifert fibre space structure to be provided to us in any way. All we need for
the existence of the certificate is that there is some Seifert fibration on $M$.

The first stage in the algorithm to apply Proposition 10.3 to obtain a triangulation (which we will also call ${\cal T}$)
for $M$ in which every boundary torus has precisely two triangles.

Any class in $H_2(M, \partial M)$ is represented by a compact oriented incompressible surface. There is a well known
classification of compact orientable incompressible surfaces in a Seifert fibre space $M$. Any such surface is
isotopic to a surface that is {\sl horizontal} or {\sl vertical}. By definition, a surface is {\sl horizontal} if it is everywhere
transverse to the fibres. A surface is {\sl vertical} if it is a union of fibres. The proof divides into these two cases.

\noindent {\sl Case 1.} The Poincar\'e dual of $\phi$ is represented by a horizontal surface $S$.

Then $M - {\rm int}(N(S))$ is a product $I$-bundle over $S$. This implies
that $S$ is in fact a fibre in a fibration of $M$ over the circle. Hence, it has minimal Thurston complexity
in its homology class.

Using Theorem 8.3, there is a taut normal surface $S$ representing this class such that
$w(S) \leq k^{t^2} ||\phi||_1$, where $k$ is a universal computable constant. The vector $(S)$ will form part of our certificate. 
Using Theorem 9.4, we can verify that $\partial S$ contains no component that bounds a disc
in $\partial M$, in time at most a polynomial function of $t \log w(S)$. Let $(M' ,\gamma')$
be the product sutured manifold that is obtained by decomposing along $S$. 
Theorem 10.2 is applied, and it supplies a positive handle structure ${\cal H}''$
for a sutured manifold $(M'', \gamma'')$ that is obtained from $(M' ,\gamma')$ by
decomposing along some non-trivial annuli disjoint from $\gamma'$. According to Theorem 10.2,
$(M'',\gamma'') = E(M'', \gamma'')$ is taut. The time taken to complete the algorithm in Theorem 10.2
is at most a polynomial function of $t \log w(S)$. Since each annulus is vertical in 
the product structure or boundary-parallel, $(M'', \gamma'')$ is
again a product sutured manifold. Theorem 10.2 provides a handle structure
for $(M'', \gamma'')$ of uniform type, with at most $ct$ 0-handles, for a universal
computable constant $c$.

Using Theorem 12.1, there is a certificate which establishes that $(M'', \gamma'')$
is a product sutured manifold. The time that it takes to verify this certificate is bounded above by 
a polynomial function of the number of 0-handles, and so by a polynomial
function of $t$. Assuming that this certificate is correctly verified, then by Theorem 10.2, $(M' ,\gamma')$ is also taut.
Therefore, $S$ is taut and $(M, \gamma)$ is taut.
So, by Lemma 2.2, we know that $S$ has minimal Thurston complexity in
its class in $H_2(M, \partial M)$. So, we have established that the Thurston norm
of this class is indeed $m$.

\noindent {\sl Case 2.} The Poincar\'e dual of $\phi$ is represented by a vertical surface.

This is a union of annuli and tori. In particular, the Thurston norm of this class is zero.
By Theorem 8.3, there is a taut normal surface $S$ representing this class such that
$w(S) \leq k^{t^2} ||\phi||_1$. This surface also has zero Thurston complexity. Since
the Thurston complexity of $S$ may be computed, using the Agol-Hass-Thurston algorithm,
in time bounded above by $t \log(w(S))$, we may readily verify that it is zero.
Thus, the part of Theorem 12.3 concerned with Thurston norm is completed.

However, we will still need to verify that $M$ is irreducible and has incompressible boundary,
or, in other words, that  $(M,\gamma)$ is taut, where $R_-(M) = \partial M$ and $R_+(M) = \emptyset$. In fact, we will verify that
the decomposition $(M, \gamma) \buildrel S \over \longrightarrow (M', \gamma')$ is taut.
We apply Theorem 10.2 to find a handle structure ${\cal H}''$ for a sutured manifold 
$(M'', \gamma'')$ that is obtained from $(M' ,\gamma')$ by
decomposing along some annuli disjoint from $\gamma'$.
This handle structure is of uniform type and has at most $ct$ 0-handles.
According to Theorem 10.2, this is taut if and only if $(M' ,\gamma')$ is.
Note that $M''$ is Seifert fibred. Since it has non-empty boundary, there is
some properly embedded, orientable surface that is horizontal in $M''$.
We may find such a surface $S_2$ that is fundamental. In particular, $w(S_2)$
is at most $k^{ct}$, by Theorem 8.2. We may now revert to Case 1, which provides a certificate to establish
that $(M'', \gamma'')$ is taut. This is verifiable in time at most a polynomial
function of $ch \log w(S_2)$, and this is at most a polynomial function of $t$.
$\square$

\vskip 18pt
\centerline {\caps 13. The certificate and its verification}
\vskip 6pt

In this section, we complete the proof of Theorem 1.5, and of Theorems 1.1 and 1.3, in the case of compact
orientable irreducible 3-manifolds with (possibly empty) toroidal boundary. We describe
the NP algorithm for determining the Thurston norm of a homology class.
We also give the proof of Theorem 1.6 by showing how to certify that a compact orientable irreducible
3-manifold with toroidal boundary and positive first Betti number is indeed irreducible.

\vfill\eject
\noindent {\caps 13.1. A description of the certificate}
\vskip 6pt

In Theorem 1.5, we are given, as an input, a triangulation ${\cal T}$ for the 3-manifold $M$, a simplicial 1-cocycle $\phi$
and an integer $m$. Our aim is determine whether the Thurston norm of the Poincar\'e dual of $[\phi]$ is $m$.
In Theorem 1.6, we are given a triangulation ${\cal T}$ for a 3-manifold $M$, and our aim is to certify that
$M$ is irreducible. In both cases, let $t$ be the number of tetrahedra in ${\cal T}$. 

In the situation of Theorem 1.6, we are assuming that $b_1(M) > 0$, and so there is a simplicial cocycle
representing a non-trivial class in $H^1(M)$. We may find such a cocycle $\phi$ such that $||\phi||_1$ is at most
at an exponential function of $t$, using Lemma 8.5. Hence, the Thurston norm $m$ of its Poincar\'e dual is at most an exponential
function of $t$. In Theorem 1.6, $\phi$ and $m$ will form part of our certificate. In both theorems, we will show
how to certify that $M$ is irreducible and that the Thurston norm of the Poincar\'e dual of $[\phi]$ is $m$.

The following is our certificate. This consists of various pieces of data, which are provided non-deterministically:
\item{(i)} a solution $(T)$, which may be zero, to the normal surface equations in ${\cal T}$ satisfying the quadrilateral constraints,
and with weight at most $2^{2+80t^2}$;
in fact, $T$ will be the canonical tori for $M$;
\item{(ii)} a triangulation ${\cal T}'$ for a 3-manifold $M'$ with at most $1000t$ tetrahedra, and such that every boundary component is
triangulated using just two triangles; in fact, this 3-manifold will be $M - {\rm int}(N(T))$; 
\item{(iii)} a simplicial 1-cocycle $\phi'$ on ${\cal T}'$ such that $||\phi'||_1 \leq 3^{800t} 1200 t ||\phi||_1$; it will in fact have
the same Thurston norm as the dual of $[\phi]$;
\item{(iv)} a decomposition of $M'$ into two subsets $M'_1$ and $M'_2$, which are unions of components of $M'$; these will in fact be the
atoroidal and Seifert fibred components, respectively; if any component of $M'$ is both atoroidal and Seifert fibred, it is
placed in $M_2'$ rather than $M'_1$;
\item{(v)} two integers $m_1$ and $m_2$ that sum to $m$; these will be the Thurston norms of classes in $H_2(M'_1, \partial M'_1)$
and $H_2(M'_2, \partial M_2')$;
\item{(vi)} a certificate, verifiable by the algorithm in Theorem 12.3, that the Thurston norm of the Poincar\'e dual of
$\phi'|M_2'$ is $m_2$ and that $M'_2$ is irreducible and has incompressible boundary;
\item{(vii)} a solution $S_1$ for the normal surface equations in the triangulation of $M_1'$,
satisfying the compatibility constraints, such that each co-ordinate is at most $c^{t^2} ||\phi||_1$;
here $c$ is some universal computable constant;
\item{(viii)} handle structures ${\cal H}_2, \dots, {\cal H}_{n+1}$ for sutured manifolds $(M_2, \gamma_2),
\dots, (M_{n+1}, \gamma_{n+1})$; these must be positive and of uniform type and the number of handles
in each ${\cal H}_i$ is at most $c t$; here $n$ is at most $ct$;
\item{(ix)} for each $i$ between $2$ and $n$, a solution $(S_i)$ for the boundary-regulated matching equations
for ${\cal H}_i$, satisfying the compatibility constraints, such that each co-ordinate of each solution is at most $c^t$;
\item{(x)} for each $i$ between $1$ and $n$, a choice of transverse orientations on at most $ct$ elementary discs of $S_i$;
\item{(xi)} a certificate that $(M_{n+1}, \gamma_{n+1})$ is a product sutured manifold, as provided by Theorem 12.1;
\item{(xii)} certificates that certify that certain 3-manifolds with uniform handle structures are 3-balls, as provided by the work
of Schleimer [26] or Ivanov [12]; the number of such 3-manifolds is bounded by a polynomial function of $t$;
the number of handles in each manifold is bounded by a linear function of $t$.

\vfill\eject
\noindent {\caps 13.2. The algorithm to verify this certificate}
\vskip 6pt

This is as follows:
\item{(1)} Verification that $T$ is a union of tori, using the Agol-Hass-Thurston algorithm [2].
\item{(2)} Verification that ${\cal T}'$ is the triangulation of $M' = M - {\rm int}(N(T))$ provided by
Theorem 11.4 and Proposition 10.3.
\item{(3)} Verification that $\phi'$ is the cocycle provided by Addendum 11.5 and Proposition 10.3,
where $[\phi']$ and $i^\ast([\phi])$ have duals with the same Thurston norm;
\item{(4)} Verification that $m_1 + m_2 = m$.
\item{(5)} Verification of the certificate in (vi) using Theorem 12.3.
\item{(6)} Verification that no component of $\partial S_1$ bounds a disc in $\partial M'$, using
Theorem 9.4.
\item{(7)} The determination of the components of $S_1$ and their Euler characteristic, using the Agol-Hass-Thurston algorithm,
and the verification that each component of $S_1$ has non-positive Euler characteristic.
\item{(8)} Verification that $\chi_-(S_1) = m_1$.
\item{(9)} For each $i$ between $1$ and $n$, the verification that each component of $S_i$ is normally parallel to a unique component that contains an elementary
disc with a pre-assigned transverse orientation given in (x), again using the Agol-Hass-Thurston algorithm.
\item{(10)} The verification that $S_i$ has a transverse orientation compatible with the transverse orientations given in (x) and (for $i > 1$)
the transverse orientations on the discs intersecting $\partial M_i$, using Theorem 9.6 or Theorem 9.7. 
\item{(11)} Using Theorem 9.7, the computation of the signed intersection number between each oriented edge of ${\cal T}'$
and $S_1$.
\item{(12)} Verification that $S_1$ is Poincar\'e dual to $\phi'|M_1'$,
by checking that, for each simplicial loop $\ell$ in some generating set for $H_1(M'_1)$,
$\phi(\ell) = \ell . S_1$.
\item{(13)} For each $i$ between $2$ and $n$, verification that no component of $S_i$ is a planar surface that is disjoint
from $\gamma_i$ and that has all but one of its boundary curves trivial, using Corollary 9.5.
\item{(14)} For each $i$ between $2$ and $n$, verification that each decomposing surface $S_i$ is allowable,
by checking that it is disjoint from the u-sutures of $M_i$ and that any trivialising planar surface of a trivial
curve of $\partial S_i$ is correctly oriented, using Corollary 9.5.
\item{(15)} For each $i$ between $1$ and $n$, the application of the algorithm in Theorem 10.1 or 10.2 to ${\cal H}_i$
and $S_i$. If the algorithm declares `not taut', then the certificate is invalid. On the other hand, if the algorithm
produces a handle structure ${\cal H}''$, then the next stage is a verification that this is combinatorially
equal to ${\cal H}_{i+1}$. In order to complete the algorithm in Theorem 10.1 or 10.2, certain 3-manifolds need to be
certified as 3-balls. These are the manifolds supplied in (xii) and their certificates are verified using
the algorithm of Schleimer [26] or Ivanov [12].
\item{(16)} A verification of the certificate in (xi) using Theorem 12.1, together with a verification
that no component of $R_\pm(M_{n+1})$ is planar and has each boundary component a u-suture.

\vskip 6pt
\noindent {\caps 13.3. Proof of Theorems 1.5 and 1.6.}
\vskip 6pt

In order to prove Theorem 1.6 and Theorem 1.5 in the case of compact orientable irreducible 3-manifolds with (possibly empty) toroidal boundary, 
we need to establish two claims:
\item{(A)} If there is a certificate that is verified by the above procedure, then the Thurston norm of the class in $H_2(M, \partial M)$ that is dual to $[\phi]$ is $m$,
and $M$ is irreducible.
\item{(B)} If Thurston norm of the class in $H_2(M, \partial M)$ that is dual to $[\phi]$ is $m$ and $M$ is irreducible, 
then there is a certificate that is verified by the above procedure
in time at most a polynomial function of $t$, $\log(m+2)$ and $\log||\phi||_1$.

\noindent We will establish these in the next two subsections.

\vskip 6pt
\noindent {\caps 13.4. A proof that the certificate does certify Thurston norm and irreducibility}
\vskip 6pt

We start with (A). So, suppose that there is a certificate as above. The algorithm verifies in (1) of Section 13.2 that $T$ is a union of tori.
It also checks in (2) that ${\cal T}'$ is a triangulation of $M' = M - {\rm int}(N(T))$. It also checks in (3) that
the Poincar\'e duals of $i^\ast([\phi])$ and $[\phi']$ have the same Thurston norm. We will show below that
the latter is $m$, and hence so is the former. By Proposition 11.2, this is equal to the Thurston norm
of the Poincar\'e dual to $\phi$, as long as $T \cup \partial M$ is incompressible. So, assuming that the certificate is
verified, all that we need to establish is that $T \cup \partial M$ is incompressible, $M$ is irreducible and the Thurston norm of the Poincar\'e
dual to $[\phi']$ is $m$.

Now, $m = m_1 + m_2$, by (4). By (5), the Poincar\'e dual of $\phi'|M_2'$ has Thurston norm $m_2$.
So, we must check that the  Poincar\'e dual of $\phi'|M_1'$ has Thurston norm $m_1$.
The algorithm checks in (11) and (12) that $S_1$ is dual to $\phi'_1|M_1'$.
The algorithm also checks that $\chi_-(S_1) = m_1$ in (7) and (8). So, we must establish that $S_1$ has minimal Thurston complexity in its class in $H_2(M', \partial M')$. 
The algorithm checks in (6) that no component of $\partial S_1$ bounds a disc in $\partial M'$.
Hence, by Lemma 2.1, $S_1$ has minimal Thurston complexity in its class in $H_2(M', \partial M')$ if and only if it has minimal Thurston complexity
in its class in $H_2(M', N(\partial S_1))$. So, we check that the latter holds. In fact, we check that $S_1$ is taut. 

Since no component of $\partial S_1$ bounds a disc in $\partial M'$ and $\partial M'$ has no u-sutures, the canonical extension
$E(S_1)$ is equal to $S_1$.  We will show that, assuming that the certificate has been correctly verified, there
is a hierarchy of the form
$$(M_1', \emptyset) = E(M_1, \gamma_1) \buildrel E(S_1) \over \longrightarrow E(M_2', \gamma_2') \buildrel A_1 \over \longrightarrow E(M_2, \gamma_2)
\buildrel E(S_2) \over \longrightarrow E(M_3', \gamma_3') \buildrel A_2 \over \longrightarrow \dots \buildrel A_n \over \longrightarrow E(M_{n+1}, \gamma_{n+1}). \eqno{(\ast)}$$
No component of $E(S_i)$ will be a disc disjoint from the sutures. Also, no boundary curve of $E(S_i)$ bounds
a disc in $\partial E(M_i)$ disjoint from the sutures, by Lemma 3.2. 
Each surface $A_i$ is a union of product discs and non-trivial annuli in $M_i$ disjoint from the sutures $\gamma_i'$.
The final manifold $E(M_{n+1}, \gamma_{n+1})$ is a product
sutured manifold, with no $S^2 \times I$ components, by (16).
Hence, by Theorem 2.3, each manifold in this sequence is taut
and each decomposing surface is taut. In particular, $E(S_1) = S_1$ is taut, and $T$ and $\partial M$ are incompressible.
Also, as $M'$ is irreducible and $T$ is incompressible, then $M$ is irreducible.

For $i > 1$, the manifolds $(M_i, \gamma_i)$ are provided by giving their handle structures ${\cal H}_i$ in (viii) of the certificate. 
For $i > 1$, each surface $S_i$ is a boundary-regulated surface
provided by (ix) in the certificate. It has a well-defined orientation, because (9) and (10) in the verification
complete successfully. The manifold $(M'_{i+1}, \gamma'_{i+1})$ is obtained by decomposing $(M_i, \gamma_i)$
along $S_i$. By (14), this decomposition is allowable. Hence, it has a canonical extension as in $(\ast)$.
By (13), no component of $E(S_i)$ is a disc disjoint from the sutures.
The annuli $A_i$ are provided by Theorem 10.1, together with the successful verification of (15).

Hence, we have shown that if the certificate is verified successfully, then the Thurston norm of the dual of $[\phi]$
is $m$, and $M$ is irreducible, as required.

\vskip 6pt
\noindent {\caps 13.5. The existence of the certificate}
\vskip 6pt

We now check the other direction (B). We start with the 3-manifold $M$ with the given triangulation ${\cal T}$, the cocycle $\phi$ and the integer $m$. 
We suppose that the Thurston norm of the dual of $[\phi]$ is $m$, and we have to show that there is a certificate as in Section 13.1.

By Theorem 11.3, the canonical tori $T$ for $M$ may be realised as a normal surface with at most $2^{80t^2}$ elementary discs.
Hence, it has weight at most $2^{2+80t^2}$. This forms part (i) in our certificate. Theorem 11.4 and Addendum 11.5 then provide
a triangulation for $M' = M - {\rm int}(N(T))$ with at most $200t $ tetrahedra. By Theorem 10.5, we may find such a triangulation
with at most $1000t$ tetrahedra, and where each boundary torus is triangulated using just two triangles. This is (ii) in the certificate.
Addendum 11.5 and Theorem 10.5 provide a 1-coycle $\phi'$
on ${\cal T}'$ with $||\phi'||_1 \leq 3^{800t} 1200 t ||\phi||_1$, and
such that the duals of $i^\ast[\phi]$ and $[\phi']$ have the same Thurston norm. This is (iii) in the certificate.
By construction, steps (1), (2) and (3) in the verification may therefore be completed. 

We let $M'_2$ be the union of the Seifert fibred components of $M'$, and let $M'_1$ be the union
of the remaining components. Therefore $M'_1$ is atoroidal. This decomposition into $M'_1$ and $M'_2$ is
(iv) of the certificate. We denote the inclusion maps of $M'_1$ and $M'_2$ into $M$ by $i_1$ and $i_2$.
We denote the Thurston norms of the Poincar\'e duals of $i^\ast_1([\phi])$ and $i^\ast_2([\phi])$ by
$m_1$ and $m_2$. These form (v) in the certificate. Then $m_1 + m_2 = m$, by Proposition 11.2. 
So step (4) verifies successfully.

Note that Theorem 12.3 provides the certificate used in (vi) and this also is verifiable.
This is step (5). 

So, we now focus on the atoroidal components $M'_1$ of $M'$. Let $\phi'_1$ be the restriction of $\phi'$ to $M'_1$.
By Theorem 8.3, there is a compact oriented lw-taut normal surface $S_1$ such that
$[S_1,\partial S_1] $ is Poincar\'e dual to $[\phi'_1]$, satisfying $w(S_1) \leq k^{(1000 t)^2} ||\phi'_1||_1$, where $k$ is
a universal computable constant. This is at most $c^{t^2}||\phi||_1$ for some computable constant $c$.
Because $S_1$ is lw-taut, whenever two components of $S_1$ are normally parallel, they are compatibly oriented, as otherwise the union of these
two components is null-homologous. The vector $(S_1)$ representing $S_1$ as a normal surface is
(vii) in our certificate. Since $S_1$ is lw-taut, no component of $\partial S_1$ bounds a disc in $\partial M'$, and so step (6)
in the verification may be completed. As $S_1$ is lw-taut, $\chi_-(S_1)$ is equal to the Thurston
norm of $[S_1]$. So, $\chi_-(S_1) = m_1$, and step (8) in the verification completes successfully. Step (7) trivially
completes successfully, because $T$ and $\partial M$ are incompressible, and therefore
no component of $S_1$ can be a disc. Also, no component of $S_1$ is a sphere, because this would have to be
null-homologous, by the irreducibility of $M_1'$.

For each collection of normally parallel components of $S_1$, one elementary disc
from one of these components is chosen. The transverse orientations on these discs form part of the certificate.
Note that there are at most $8000t$ of these discs, because this is an upper bound for the number of disjoint
normal surfaces that are not normally parallel. The reason for including these transverse orientations
as part of the certificate is to ensure that there is no ambiguity about the transverse orientation on $S_1$.
This collection of transverse orientations is (x) in our certificate, in the case $i = 1$. By the
way that these elementary discs are chosen, step (9) in the verification completes successfully.
Also, because they arose from the orientation on $S_1$, so does step (10) in the case $i=1$.

As $[S_1,\partial S_1] $ is Poincar\'e dual to $[\phi'_1]$, steps (11) and (12) also complete successfully.

We now need to establish the existence of the hierarchy as in $(\ast)$. We will establish the existence of
the manifolds $(M_i', \gamma_i')$ and $(M_i, \gamma_i)$ one at a time. We will also prove inductively
that they have the following properties. Our first property is that $(M_i, \gamma_i)$ and $(M'_i, \gamma_i')$ 
admit some allowable hierarchy. Our second property is that $(M_i, \gamma_i)$ has a handle
structure ${\cal H}_i$ of uniform type. Our third property is that $E(M_i, \gamma_i)$ and
$E(M'_i, \gamma_i')$ are atoroidal and no component is a Seifert fibre space other than
a solid torus or a copy of $T^2 \times I$. Our fourth property is that, for each component
$X$ of $(M'_i, \gamma'_i)$, $\partial X$ is not a single torus with no sutures, unless
$X$ is a component of $M_1'$.
The induction starts with $(M_2', \gamma_2')$. Note that Lemma 2.2 gives that $(M_2', \gamma_2')$
is taut. It has no u-sutures, and so by Theorem 2.4, it does admit an allowable hierarchy.
By Lemma 3.11, $E(M_2', \gamma_2')$ is atoroidal and no component is a Seifert fibre space other than
a solid torus or a copy of $T^2 \times I$. For each component $Y$ of $M'_1$ such that
$Y \cap S_1$ is non-empty, then $Y \cap M_2'$ cannot contain a component with boundary
a single torus with no sutures. This is because such a component would have arisen from
a toral component of $S_1$ or a union of annular components, which would have been
homologically trivial in $M'_1$. This contradicts the fact that $S_1$ is lw-taut.

We now consider the inductive step. So, consider the decomposition
$$(M_{i}, \gamma_{i}) \buildrel S_i \over \longrightarrow (M_{i+1}', \gamma_{i+1}').$$
We apply Theorem 10.1 or 10.2 to obtain a decomposition 
$$E(M_{i+1}', \gamma_{i+1}') \buildrel A_{i+1} \over \longrightarrow E(M_{i+1}, \gamma_{i+1})$$
where $A_{i+1}$ is a collection of product discs and oriented annuli disjoint from $\gamma_{i+1}'$ and with non-trivial boundary curves.
By Theorem 10.1 or 10.2, this decomposition is taut, and in particular, no component of $(M_{i+1}, \gamma_{i+1})$ is
a solid torus with no sutures.
Inductively, we are assuming that $(M_{i+1}', \gamma_{i+1}')$ admits an allowable hierarchy,
and hence, by Theorem 10.1 or 10.2, so does $(M_{i+1}, \gamma_{i+1})$.
By Theorem 10.1 or 10.2, $E(M_{i+1}, \gamma_{i+1})$ is atoroidal and no component
is a Seifert fibre space other than a solid torus or a copy of $T^2 \times I$. By Lemma 3.10
and the tautness of $(M_{i+1}, \gamma_{i+1})$, no component of $(M_{i+1}, \gamma_{i+1})$
has boundary a single torus with no sutures unless it is a component of $M'_1$.
Theorem 10.1 or 10.2 provides a handle structure ${\cal H}_{i+1}$ for $(M_{i+1}, \gamma_{i+1})$
of uniform type. Each 0-handle of ${\cal H}_{i+1}$ lies within a 0-handle of ${\cal H}_i$. Moreover, for each 0-handle $H$ of ${\cal H}_i$, 
the 0-handles of ${\cal H}_{i+1}$ lying
within it have, in total, no greater complexity. Indeed, if the complexity is unchanged, then 
${\cal H}_{i+1} \cap H$ contains a single 0-handle isotopic to $H$. Finally, ${\cal H}_{i+1}$ has no parallelity handles other than 2-handles.
We now apply Theorem 7.9 to $(M_{i+1}, \gamma_{i+1})$, which is possible because
$(M_{i+1}, \gamma_{i+1})$ admits an allowable hierarchy. Theorem 7.9 gives
a regulated surface $S_{i+1}$ which is fundamental as a boundary-regulated surface and which also extends to
an allowable hierarchy. Decomposition along $S_{i+1}$ reduces the complexity
of the handle structure by Theorem 7.9. The four properties of $(M_{i+1}, \gamma_{i+1})$
required by the induction are, as above, easily seen to hold.

By construction, steps (13) and (14) of the verification process complete successfully, because the surfaces $S_i$
are part of an allowable hierarchy. Step (15) also
completes successfully, because the handle structure ${\cal H}_{i+1}$ is constructed
using Theorem 10.1.

Repeat this process until we end with a product sutured manifold $(M_{n+1}, \gamma_{n+1})$,
no component of which is pre-spherical. Note that the process is guaranteed to terminate because
the complexity of the handle structure ${\cal H}_{i+1}$ is strictly less than that of ${\cal H}_i$,
for each $i$. Thus, if we apply step (16) of the verification procedure to the final manifold, this also 
completes successfully.

In this way, we obtain the certificate, and we have shown that the verification completes successfully
We need to check that the verification can be completed in polynomial time. In order to do this, we will check that
the length $n$ for the hierarchy is at most $k t$, where $k$ is a universal computable constant, and we will check that the algorithms applied
to the $i$th stage of the hierarchy can be completed in time that is at most a polynomial function of $t$.

We first show that the length of the hierarchy is at most $kt$. Consider a 0-handle $H$ of ${\cal H}_2$, and the 0-handles of
${\cal H}_i$ that lie inside $H$. As $i$ increases, the complexity of these handles does not go up. Moreover, if it stays constant,
the handle remains unchanged up to isotopy. By Theorem 6.1, the number of different possible configurations for ${\cal H}_i \cap H$ is bounded
by some universal computable constant $k$. Hence, the complexity can decrease at most $k$ times. Now, at each stage of the
hierarchy, the complexity must strictly decrease, for some 0-handle $H$ by Theorem 7.9. So, the number of steps in the hierarchy
is at most $k$ times the number of 0-handles of ${\cal H}_2$. But the number of 0-handles of ${\cal H}_2$ is at most $6000t$, since $1000t$
is an upper bound for the number of tetrahedra of ${\cal T}'$, and each tetrahedron of ${\cal T}'$ contains at most six 0-handles
of ${\cal H}_2$ that are not parallelity handles.

We now show that, at each stage of the hierarchy, the algorithms can be completed in time that is at most a polynomial
function of $t$. The weight of each surface $S_i$ is at most $c^{t^2}$. The logarithm of this is $t^2 \log(c)$. Whenever the 
algorithms of Section 9 or Theorem 10.1 are applied, it is this logarithm that appears in the argument. Hence, the running time is at most a
polynomial function of $t$. The final step in the algorithm is the verification of the certificate that $(M_{n+1}, \gamma_{n+1})$ is a product
sutured manifold with no pre-spherical components.
The handle structure on $(M_{n+1}, \gamma_{n+1})$ is of uniform type and has at most $ct$ handles. Hence, the time required to complete this
verification is also at most a polynomial function of $t$.

This completes the proof of Theorem 1.6. It also completes the proof of Theorem 1.5 in the case of compact orientable
irreducible 3-manifolds with (possibly empty) toroidal boundary.

\vskip 18pt
\centerline{\caps 14. The case of closed irreducible 3-manifolds}
\vskip 6pt

In this section, we show that, in the proof of Theorem 1.5, it suffices to focus on the case of closed orientable irreducible 3-manifolds. In conjuction with the argument in the previous section, this will complete the proof of Theorem 1.5. We define the decision problem
{\caps Thurston norm in closed orientable irreducible 3-manifolds} to be the restriction of {\caps Thurston norm of a homology class} to closed orientable irreducible 3-manifolds.

\noindent {\bf Theorem 14.1.} {\sl Suppose that {\caps Thurston norm in closed orientable irreducible 3-manifolds} lies in NP. Then so does {\caps Thurston norm of a homology class}.}

\noindent {\bf Lemma 14.2.} {\sl Let $M$ be a compact orientable 3-manifold and let $z$ be a class in $H_2(M, \partial M)$. Let $\tilde M$ be the double of $M$, and let $\tilde z \in H_2(\tilde M)$ be the doubled class. Then the Thurston norm of $\tilde z$ is twice that of $z$.}

\noindent {\sl Proof.} Let $S$ be a norm-minimising representative for $z$. Its double $\tilde S$ represents $\tilde z$. Each component of $S$ is doubled to form a surface with twice the Thurston complexity. Hence, $\chi_-(\tilde S)=2 \chi_-(S)$. So, the Thurston norm of $\tilde z$ is at most twice that of $z$.

We now establish the opposite inequality. It suffices to consider the case where $M$ has incompressible boundary. For if $\partial M$ is compressible, then we can boundary-reduce $M$, to form a 3-manifold $N$ with incompressible boundary. Thus, $M$ is equal to $N$ with 1-handles attached. Any class $z$ in $H_2(M, \partial M)$ is Poincar\'e dual to a class in $H^1(M)$, which restricts to a class in $H^1(N)$, and this is dual to a class $y$ in $H_2(N, \partial N)$. The classes $z$ and $y$ have equal Thurston norms. Moreover, if $\tilde M$ is the double of $M$, and $\tilde N$ is the double of $N$, then the doubled classes $\tilde z$ and $\tilde y$ have equal Thurston norms.

So, we may assume that $\partial M$ is incompressible. Let $\tilde S$ be a norm-minimising representative for $\tilde z$. We may assume that $\tilde S$ is incompressible. We may also assume that the curves $\tilde S \cap \partial M$ are essential in $\tilde S$. Let $S_1$ be the restriction of $\tilde S$ to one half of $\tilde M$; denote this half by $M$. Let $S_2$ be the image in $M$ of ${\rm cl}(\tilde S - S_1)$ under the involution of $\tilde M$. Thus, $S_1$ and $S_2$ are properly embedded surfaces in $M$ with equal boundary curves. Since the restriction of $\tilde S$ to each component of $\tilde M - \partial M$ is incompressible, $S_1$ and $S_2$ are both incompressible surfaces in $M$. We may therefore compress $S_2$, keeping its boundary fixed, so that no curve of $S_1 \cap S_2$ bounds a disc in either $S_1$ or $S_2$. 
%After these compressions, each curve of $S_1 \cap S_2$ is $\pi_1$-injective in $S_1$ and $S_2$, and hence $\pi_1$-injective in $M$ and $\tilde M$. 
We may perform the corresponding compressions to $\tilde S$, which do not change its Thurston complexity. We denote this new embedded surface also by $\tilde S$. Let $\tilde S'$ be its image under the involution of $\tilde M$. This also represents $\tilde z$. Resolve the intersections between $\tilde S$ and $\tilde S'$ to form an embedded oriented surface $\tilde S''$. Since we resolve only along essential curves, $\chi_-(\tilde S'') = \chi_-(\tilde S) + \chi_-(\tilde S') = 2\chi_-(\tilde S)$. Note that $\tilde S''$ is an embedded surface that is invariant under the covering involution and that represents $2 \tilde z$. Hence, its restriction to each copy of $M$ is a surface $S$ representing $2z$. So, $2 \chi_-(S) = \chi_-(\tilde S'') = 2 \chi_-(\tilde S)$. So the Thurston norm of $2z$ is at most that of $\tilde z$. But Thurston norm is linear along rays in $H_2(M, \partial M)$ and so the Thurston norm of $2z$ is twice that of $z$. Therefore, the Thurston norm  of $\tilde z$ is at least twice that of $z$. $\square$

\noindent {\sl Proof of Theorem 14.1.}
The input to {\caps Thurston norm of a homology class} is a triangulation ${\cal T}$ of a compact orientable 3-manifold $M$ with $t$ tetrahedra, a simplicial 1-cocycle $\phi$ representing an element $[\phi]$ in $H^1(M)$ and a non-negative integer $n$. The problem asks whether the Poincar\'e dual of $[\phi]$ has Thurston norm $n$. The certificate that we will use to solve this problem is as follows:
\item{(i)} a triangulation $\tilde {\cal T}$ of a 3-manifold $\tilde M$, which will in fact be the double of $M$ with its natural triangulation;
\item{(ii)} a 1-cocycle $\tilde \phi$ on $\tilde {\cal T}$, which is in fact the double of $\phi$;
\item{(iii)} a vector $(S)$ of a normal surface $S$ in $\tilde {\cal T}$ with weight at most $2^{740t^2}$; this will in fact be a collection of disjoint embedded 2-spheres;
\item{(iv)} a triangulation ${\cal T}'$ of a 3-manifold $M'$ with at most $400t$ tetrahedra; this will be $\tilde M - {\rm int}(N(S))$;
\item{(v)} a 1-cocycle $\phi'$ on $M'$ satisfying $||\phi'||_1 \leq 2400t ||\tilde \phi||_1$; this will be provided by Addendum 11.5;
\item{(vi)} a triangulation ${\cal T}''$ of a 3-manifold $M''$, obtained from ${\cal T}'$ by attaching a 3-ball, triangulated as a cone, to each component of $\partial M'$;
\item{(vii)} a 1-cocycle $\phi''$ on $M''$ that extends $\phi'$ and satisfies $||\phi''||_1 \leq (1+ 1600t) ||\phi'||_1$;
\item{(viii)} the certificate for {\caps Irreducibility of a compact orientable 3-manifold with toroidal \break boundary and} $b_1 > 0$ giving that $M''$ is irreducible;
\item{(ix)} the certificate for {\caps Thurston norm in closed orientable irreducible 3-manifolds} that the Thurston norm in $M''$ of the dual of $[\phi'']$ is $2n$.

The algorithm to verify the certificate is as follows:
\item{(1)} Verification that $\tilde {\cal T}$ is indeed the double of ${\cal T}$ and that $\tilde \phi$ is the double of $\phi$.
\item{(2)} Verification that $S$ is a union of spheres, using the Agol-Hass-Thurston algorithm [2].
\item{(3)} Verification that ${\cal T}'$ is the triangulation of $\tilde M - {\rm int}(N(S))$ provided by Theorem 11.4.
\item{(4)} Verification that $\phi'$ is the cocycle provided by Addendum 11.5.
\item{(5)} Verification that $\phi''$ extends $\phi'$ and satisfies $||\phi''||_1 \leq (1 + 1600t) ||\phi'||_1$;
\item{(6)} Verification of the certificates in (8) and (9) above.

We now show that the above certificate exists and that its correct verification is equivalent to the Thurston norm of the dual of $[\phi]$ being $n$. We double ${\cal T}$ to form $\tilde {\cal T}$, with $2t$ tetrahedra, and we double $\phi$ to give $\tilde \phi$.  By Lemma 14.2, the Thurston norm of the dual of $\tilde \phi$ is double that of $\phi$. It was shown by King (Lemma 4 in [18]) that there is a maximal collection disjoint essential normal spheres $S$ in $\tilde {\cal T}$ with weight at most $2^{740t^2}$. Applying Theorem 11.4 and Addendum 11.5, we obtain the triangulation ${\cal T}'$ for $\tilde M - {\rm int}(N(S))$, and a 1-cocycle $\phi'$ satisfying  $||\phi'||_1 \leq 2400t ||\tilde \phi||_1$. By Addendum 11.5, the Thurston norms of the duals of $[\phi']$ and $[\tilde \phi]$ are equal. We can extend $\phi'$ to a 1-cocycle $\phi''$ on $M''$ as follows. Dual to $\phi'$ is a surface properly embedded in $M'$. Since each component of $\partial M'$ is a 2-sphere, we can cap off each component of this surface with a disc lying in the attached 3-ball. Each such disc intersects each edge in the interior of the ball at most once. Hence, the total number of intersection points between the interior of these discs and the edges of ${\cal T}''$ is at most $||\phi'||_1$ times the number of vertices of ${\cal T}'$, which is at most $1600t ||\phi'||_1$. Therefore, $||\phi''||_1 \leq (1+1600t ) ||\phi'||_1$. Since $S$ was maximal, the manifold $M''$ obtained by attaching a 3-ball to each component of $\tilde M - {\rm int}(N(S))$ is irreducible. Hence, {\caps Thurston norm in closed orientable irreducible 3-manifolds} can be used to certify the Thurston norm of the dual of $[\phi'']$ and hence the Thurston norm of the dual of $[\phi].$ 
$\square$

\vskip 18pt
\centerline{\caps 15. A certificate for incompressibility}
\vskip 6pt

\noindent {\bf Theorem 1.7.} {\sl {\caps Incompressible boundary} is in NP $\cap$ co-NP.}

\noindent {\sl Proof.} It was shown by Ivanov (Theorem 4 in [12]) that {\caps Incompressible boundary} is in co-NP. So we focus on showing that it is in NP. Let ${\cal T}$ be a triangulation
of a compact orientable 3-manifold $M$ with $t$ tetrahedra. Our certificate for establishing the incompressibility of $\partial M$ is as follows:
\item{(i)} a vector $(S)$ of a normal surface $S$ in ${\cal T}$ with weight at most $2^{185t^2}$; this will in fact be a collection of disjoint embedded 2-spheres;
\item{(ii)} a triangulation ${\cal T}'$ of a 3-manifold $M'$ with at most $1000t$ tetrahedra; this will be $M - {\rm int}(N(S))$ with a 3-ball attached to each component of $\partial N(S)$;
\item{(iii)} a certificate, provided by Theorem 1.6, that the double of each component of $M'$ with non-empty boundary is irreducible.

This certificate is verified in a similar way as in the proof of Theorem 14.1. One checks that $S$ is a collection of 2-spheres using the algorithm of Agol-Hass-Thurston. One checks that ${\cal T}'$ is the triangulation provided by Theorem 11.4, but with a coned 3-ball attached to each component of $\partial N(S)$. One then verifies the certificate provided by Theorem 1.6.

We must show that a verifiable certificate exists if and only if $\partial M$ is incompressible. Suppose first that $\partial M$ is compressible. Then there is a compression disc that avoids the 2-spheres $S$. This therefore forms a compression disc for some component of $\partial M'$.  The double of this disc then forms a reducing sphere for the double of this component of $M'$. Hence, the certificate in (iii) cannot exist.

Suppose now that $\partial M$ is incompressible. By King's result (Lemma 4 in [18]), there is a maximal collection disjoint essential normal spheres $S$ in ${\cal T}$ with weight at most $2^{185t^2}$. (King assumes in his lemma that the 3-manifold is closed, but this is not actually required for his argument.)  Cutting along $S$ and then attaching 3-balls gives the irreducible 3-manifold $M'$. Theorem 11.4 gives the triangulation ${\cal T}'$, which is constructible in polynomial time. The double of each component of $M'$ then must be irreducible. This is because a reducing sphere could be chosen to be disjoint from the incompressible surface $\partial M$, and hence would form a reducing sphere for $M'$, contradicting the irreducibility of this manifold. Thus the above certificate verifies correctly in polynomial time as a function of the number of tetrahedra in ${\cal T}$. $\square$

This has Theorem 1.8 as a rapid consequence.

\noindent {\bf Theorem 1.8.} {\sl {\caps Knottedness in 3-manifolds} is in NP $\cap$ co-NP.}

\noindent {\sl Proof.} Let $\tilde M$ be the double of $M - {\rm int}(N(K))$ along $\partial M$.
Then $K$ is unknotted if and only if $K$ is homologically trivial and $\partial \tilde M$ is compressible.
Whether or not $K$ is homologically trivial can can be easily determined in polynomial time
using linear algebra. Thus, the certificate for knottedness follows from the fact that 
{\caps Incompressible boundary} is in NP, and the certificate for unknottedness uses that
{\caps Incompressible boundary} is in co-NP. $\square$

\vskip 18pt
\centerline{\caps 16. The genus of knots in 3-manifolds}
\vskip 6pt

In this section, we give a proof of Theorem 1.4.

\noindent {\bf Theorem 1.4.} {\sl If {\caps Knot genus in compact orientable 3-manifolds} in is NP, then NP $=$ co-NP.}

 This is a well known consequence of work of Agol, Hass and Thurston [2], but we give a proof here because it does
 not appear to have been put into print before, and because of its great relevance to the theme of this paper.
 
 Recall that co-NP consists of the class of decision problems for which a negative answer can be certified in polynomial
 time. Since in practice, there is no reason to believe that certifying a negative answer is as easy or as hard as certifying
 a positive answer to a problem, then it is widely believed that NP $\not=$ co-NP. However, if a single problem in
 co-NP were NP-complete, then NP and co-NP would be equal. It is this observation that is the basis for the
 proof of Theorem 1.4.

\noindent {\sl Proof.} The main theorem of Agol, Hass and Thurston in [2] is that the problem of deciding whether
a simplicial knot in a triangulated 3-manifold has genus {\sl at most} $g$ is NP-complete. (Somewhat confusingly,
they called this problem `{\caps 3-manifold knot genus}' but we do not do so here.) Therefore, if it is also in co-NP,
then NP $=$ co-NP.  We will show that if {\caps Knot genus in compact orientable 3-manifolds} in is NP, then there is a method of
certifying in polynomial time that the genus of a knot $K$ in a compact orientable 3-manifold $M$ is more than $g$.
So, consider a simplicial knot $K$ in a triangulated 3-manifold $M$. We may easily check, in polynomial time,
whether $K$ is homologically trivial. If it is not, then $K$ does not bound a Seifert surface, and the genus
of $K$ is then, by convention, infinite. On the other hand if $K$ is homologically trivial, then it is easy
to prove, using normal surface theory, that the genus $g(K)$ of $K$ is at most an exponential function of
the number of tetrahedra $t$ in the triangulation of $M$. Recall that we are given an integer $g$, and are trying
to certify that the genus of $K$ is greater than $g$. But assuming that {\caps Knot genus in compact orientable 3-manifolds} is in NP,
there is a certificate that confirms the genus of $K$, which is verifiable in time that is bounded above by a polynomial
function of $t$ and the number of digits of $g(K)$ in binary. Since $g(K)$ is at most an exponential function of $t$,
its number of digits is bounded by a linear function of $t$. Hence, the time required to verify the certificate is at most a polynomial 
function of $t$. Hence, we have indeed shown that the problem considered by Agol, Hass and Thurston is
in co-NP, as required. $\square$

\vskip 18pt
\centerline{\caps References}
\vskip 6pt

\item{1.} {\caps I. Agol,} {\sl Knot genus is NP}, Conference presentation (2002).
\item{2.} {\caps I. Agol, J. Hass, W. Thurston}, {\sl The computational complexity of knot genus and spanning area},
Trans. Amer. Math. Soc. 358 (2006) 3821--3850.
\item{3.} {\caps D. Gabai,} {\sl Foliations and the topology of $3$-manifolds.}
J. Differential Geom. 18 (1983) 445--503.
\item{4.} {\caps D. Gabai,} {\sl Foliations and genera of links.} Topology 23 (1984), no. 4, 381--394.
\item{5.} {\caps D. Gabai,} {\sl Genera of the arborescent links.} Mem. Amer. Math. Soc. 59 (1986), no. 339, i--viii and 1--98.
\item{6.} {\caps D. Gabai,} {\sl Essential laminations and Kneser normal form. }
J. Differential Geom. 53 (1999), no. 3, 517--574. 
\item{7.} {\caps O. Goldreich}, {\sl Computational complexity. A conceptual perspective}, Cambridge University Press (2008).
\item{8.} {\caps W. Haken,} {\sl Theorie der Normalfl\"achen.} Acta Math. 105 (1961) 245--375.
\item{9.} {\caps J. Hass, J. Lagarias,} {\sl The number of Reidemeister moves needed for unknotting.}
J. Amer. Math. Soc. 14 (2001), no. 2, 399--428 
\item{10.} {\caps J. Hass, J. Lagarias, N. Pippenger,} {\sl The computational complexity of knot and 
link problems.} J. ACM 46 (1999), no. 2, 185--211.
\item{11.} {\caps J. Hass, J. Snoeyink, W. Thurston,} {\sl The size of spanning disks for polygonal curves.}
Discrete Comput. Geom. 29 (2003) 1--17.
\item{12.} {\caps S. Ivanov}, {\sl
The computational complexity of basic decision problems in 3-dimensional topology.}
Geom. Dedicata 131 (2008), 1--26. 
\item{13.} {\caps W. Jaco, U. Oertel,} {\sl An algorithm to decide if a 3-manifold is a Haken manifold.}
Topology 23 (1984), no. 2, 195--209. 
\item{14.} {\caps W. Jaco, E. Sedgwick}, {\sl Decision problems in the space of Dehn fillings.} Topology 42 (2003) 845--906.
\item{15.} {\caps W. Jaco, P. Shalen,} {\sl Seifert fibered spaces in irreducible, sufficiently-large 3-manifolds.}
Bull. Amer. Math. Soc. 82 (1976) 765--767.
\item{16.} {\caps W. Jaco, J. Tollefson,} {\sl Algorithms for the complete decomposition of a closed $3$-manifold.}
Illinois J. Math. 39 (1995), no. 3, 358--406.
\item{17.} {\caps K. Johannson,} {\sl Homotopy equivalences of 3-manifolds with boundaries.}
Lecture Notes in Mathematics, 761. Springer, Berlin, 1979.
\item{18.} {\caps S. King}, {\sl The size of triangulations supporting a given link.} Geom. Topol. 5 (2001) 369--398.
\item{19.} {\caps G. Kuperberg}, {\sl Knottedness is in NP, modulo GRH}, Adv. Math. 256 (2014), 493--506.
\item{20.} {\caps M. Lackenby}, {\sl Dehn surgery on knots in 3-manifolds}, J. Amer. Math. Soc. 10 (1997) 835--864.
\item{21.} {\caps M. Lackenby}, {\sl Exceptional surgery curves in triangulated 3-manifolds,}
Pacific J. Math. 210 (2003), 101--163.
\item{22.} {\caps M. Lackenby}, {\sl The crossing number of composite knots,} J. Topology 2 (2009) 747--768.
%\item{23.} {\caps B. Martelli, C. Petronio,} {\sl Three-manifolds having complexity at most 9.} 
Experiment. Math. 10 (2001) 207--236.
\item{23.} {\caps S. Matveev}, {\sl Algorithmic topology and classification of 3-manifold},
Algorithms and Computation in Mathematics, 9. Springer, Berlin, 2007.
\item{24.} {\caps A. Mijatovi\'c}, {\sl Triangulations of fibre-free Haken 3-manifolds}, Pacific J. Math. 219 (2005) 139--186.
\item{25.} {\caps M. Scharlemann,} {\sl Sutured manifolds and generalized Thurston norms.}
J. Differential Geom. 29 (1989), no. 3, 557--614.
\item{26.} {\caps S. Schleimer}, {\sl Sphere recognition lies in NP}, 
Low-Dimensional and Symplectic Topology, Proceedings of Symposia in Pure Mathematics, 82 (2011), 183--213.
\item{27.} {\caps W. Thurston}, {\sl A norm for the homology of 3-manifolds.} Mem. Amer. Math. Soc. 59 (1986), no. 339, i--vi and 99--130.
\item{28.} {\caps J. Tollefson, N. Wang}, {\sl Taut normal surfaces}. Topology 35 (1996), no. 1, 55--75. 
\item{29.} {\caps F. Waldhausen}, {\sl On irreducible 3-manifolds which are sufficiently large.} Ann. of Math. (2) 87 1968 56--88.

\vskip 12pt
\+ Mathematical Institute, University of Oxford, \cr
\+ Radcliffe Observatory Quarter, Woodstock Road, Oxford OX2 6GG, United Kingdom. \cr

\end